\documentclass[11pt,oneside]{amsart}
\usepackage[a4paper,margin=1in]{geometry}
\usepackage[T1]{fontenc}
\usepackage{lmodern,microtype,amsmath,amssymb,amsthm,amscd,mathtools,mathrsfs}
\usepackage{enumitem,booktabs,array,xcolor}
\usepackage[all]{xy}
\usepackage[colorlinks=true,linkcolor=blue!45!black,citecolor=blue!45!black,urlcolor=blue!45!black]{hyperref}
\hypersetup{pdftitle={Comparing p-adic representations and Langlands parameters by reduction to depth zero},pdfauthor={Manish Mishra}}
\numberwithin{equation}{section}
\newtheorem{theorem}{Theorem}[section]
\newtheorem*{introapplication}{Application}
\newtheorem{proposition}[theorem]{Proposition}
\newtheorem{lemma}[theorem]{Lemma}
\newtheorem{corollary}[theorem]{Corollary}
\theoremstyle{definition}
\newtheorem{definition}[theorem]{Definition}
\newtheorem{hypothesis}[theorem]{Hypothesis}
\newtheorem{example}[theorem]{Example}

\theoremstyle{remark}
\newtheorem{remark}[theorem]{Remark}
\DeclareMathOperator{\Irr}{Irr}
\DeclareMathOperator{\Hom}{Hom}
\DeclareMathOperator{\End}{End}
\DeclareMathOperator{\Ad}{Ad}
\DeclareMathOperator{\Stab}{Stab}
\DeclareMathOperator{\GL}{GL}
\DeclareMathOperator{\SL}{SL}

\DeclareMathOperator{\IC}{IC}
\DeclareMathOperator{\Ind}{Ind}
\DeclareMathOperator{\Sc}{Sc}
\DeclareMathOperator{\uSc}{uSc}
\DeclareMathOperator{\pr}{pr}
\DeclareMathOperator{\Spec}{Spec}
\newcommand{\C}{\mathbb C}
\newcommand{\Z}{\mathbb Z}
\newcommand{\WF}{W_F}
\newcommand{\IF}{I_F}
\newcommand{\PF}{P_F}
\newcommand{\hG}{\widehat G}
\newcommand{\hM}{\widehat M}
\newcommand{\LG}[1]{{}^{L}#1}
\newcommand{\Xnr}{X_{\mathrm{nr}}}
\newcommand{\Phie}{\Phi_{\mathrm e}}
\newcommand{\Phic}{\Phi_{\mathrm{e,cusp}}}
\newcommand{\svee}{\mathfrak s^\vee}
\newcommand{\fs}{\mathfrak s}
\newcommand{\cA}{\mathcal A}
\newcommand{\cF}{\mathcal F}
\newcommand{\cE}{\mathcal E}
\newcommand{\cM}{\mathcal M}
\newcommand{\rig}{\mathfrak g}
\newcommand{\extquot}{\mathbin{/\mkern-6mu/}}
\newcommand{\Fr}{\mathrm{Fr}}
\newcommand{\hH}{\widehat H}
\newcommand{\cO}{\mathcal O}
\newcommand{\ssH}{\mathfrak s_0^\vee}
\newcommand{\calD}{\mathcal D}
\newcommand{\red}{R_{\calD}^{+}}
\newcommand{\inc}{J_{\calD}^{+}}
\newcommand{\rigH}{\mathfrak h}
\newcommand{\scrE}{\mathscr E}
\newcommand{\pZ}{p_Z}

\newcommand{\cc}{\operatorname{cc}}
\newcommand{\cusp}{\operatorname{cusp}}
\newcommand{\aff}{\mathrm{aff}}
\newcommand{\src}{\mathrm{src}}

\newcommand{\rb}{\mathrm{b}}
\newcommand{\R}{\mathbb R}

\newcommand{\Rep}{\mathrm{Rep}}
\newcommand{\cind}{\operatorname{c-ind}}
\newcommand{\nr}{\mathrm{nr}}

\newcommand{\Rad}{\operatorname{Rad}}
\newcommand{\Mod}{\mathrm{Mod}}
\newcommand{\op}{\mathrm{op}}
\newcommand{\Jac}{\mathrm{Jac}}
\newcommand{\ft}{\mathfrak{t}}
\newcommand{\ff}{\mathfrak{f}}

\newcommand{\rL}{\mathrm{L}}

\newcommand{\BBB}{\mathcal{B}}
\newcommand{\mH}{\mathcal H}

\newcommand{\cP}{\mathcal P}

\newcommand{\cl}{\operatorname{Cliff}}
\newcommand{\typ}{\operatorname{Type}}
\newcommand{\Res}{\operatorname{Res}}
\newcommand{\cN}{\mathcal N}
\newcommand{\cW}{\mathcal W}

\title[Comparing representations and Langlands parameters]{Comparing $p$-adic representations and Langlands parameters by reduction to depth zero}
\author{Manish Mishra}
\email{manish@iiserpune.ac.in}
\address{Department of Mathematics\\
Indian Institute for Science Education and Research\\
Dr.\ Homi Bhabha Road, Pashan\\
Pune 411 008, India}
\thanks{The author was partially supported by SERB Core Research Grant
(CRG/2022/000415).}
\date{28 November 2026}

\begin{document}
\begin{abstract}
We study how to attach an enhanced Langlands parameter to each
irreducible representation in a tame Bernstein block. We pass to
depth zero, where representations can be studied using finite reductive
groups over the residue field. Suppose the reduced block has
non-singular supercuspidal support in the sense of Solleveld and Xu.
Fix reductions on the two sides that use the same characters at each
stage. We lift the Solleveld--Xu correspondence
to a bijection between the irreducible representations in the original
block and a block of rigidly enhanced Langlands parameters. This applies in positive depth, in particular to regular Bernstein
blocks equipped with such data.

The general reduction uses the twisted extended quotients proposed by
Aubert, Baum, Plymen and Solleveld. In each quotient, an orbit of
cuspidal data can carry several labels, one for each irreducible
projective representation of its stabilizer with the prescribed
multiplier. We show that reduction to depth zero preserves the
stabilizer algebras on each side and the identifications of labels
at conjugate points. When the reduced groups, supporting Levi
subgroups and rigid inner twists are matched, a comparison of the
stabilizer algebras at depth zero determines a unique compatible
bijection between the original quotients. We give a finite trace
criterion for deciding when a prescribed matching of some labels
determines the rest.
\end{abstract}
\maketitle
\setcounter{tocdepth}{1}
\tableofcontents
\clearpage

\section{Introduction}\label{sec:introduction}
Let $G$ be a connected reductive group over a nonarchimedean local
field $F$, and let $\fs=[L,\pi]_G$ be an inertial class. A description
of the irreducible representations in this block starts with the
unramified orbit of the supercuspidal representation $\pi$ of $L(F)$.
The inertial Weyl group acts on this orbit. A point with nontrivial
stabilizer can carry several labels, given by the simple modules of a
twisted group algebra of the stabilizer. An enhanced cuspidal Langlands
parameter for a Levi subgroup gives corresponding data on the Galois
side. To compare the resulting parametrizations one must understand the
algebras responsible for these labels and the maps relating them over
conjugate points.

Solleveld and Xu construct the enhanced correspondence for depth-zero
blocks with non-singular supercuspidal support. At the end of their
introduction they suggest using the reduction of Adler, Fintzen, Mishra
and Ohara to extend their result to positive depth \cite{SXDepthZero}.
Reducing the enhanced parameters as well gives the following application.

\begin{introapplication}[Theorem~\ref{thm:regular-block}]
Let $\fs=[L,\pi]_G$ be a tame Bernstein block, with fixed matching
reductions. The two reductions use the same tower of groups and the
same characters at every stage, including the final character of
$G(F)$. Write
$\fs_{-1}=[L^0,\pi_{-1}]_{G^0}$ for the reduced block. Suppose that
$\pi_{-1}$ is an $F$-non-singular depth-zero supercuspidal representation
of $L^0(F)$. Then the Solleveld--Xu correspondence on $\fs_{-1}$
lifts uniquely through these reductions to a bijection
\[
 \mathcal L_\Sigma^{\rm ns}:\Irr(G(F))^{\fs}
       \xrightarrow{\ \sim\ }\Phie^Z(\rig)^{\svee}.
\]
To obtain the parameter of an irreducible representation in $\fs$,
reduce the representation to depth zero, apply the Solleveld--Xu
correspondence, and lift the resulting enhanced parameter. The target
block on the chosen rigid inner form $\rig$ is indexed by $\svee$,
the inertial class of the lifted parameter of $\pi_{-1}$.
\end{introapplication}

For the trivial tower with trivial character and identity embedding,
the bijection is the Solleveld--Xu correspondence itself. In positive
depth it assigns a parameter to each irreducible constituent of a
reducible parabolic induction in the block. For a regular Bernstein block in the sense of
Adler--Mishra \cite{AdlerMishraRegular}, the reduced supercuspidal
support is non-singular in the regular presentation. The application
therefore covers such a block whenever that presentation admits the
matching character data. Aubert--Xu had previously compared the irreducible
classes of a regular block and its depth-zero counterpart
\cite[Theorem~2.2.14]{AubertXuBlocks}.

The proposal to construct the enhanced local Langlands correspondence
by comparing two twisted extended quotients is due to
Aubert--Baum--Plymen--Solleveld (ABPS). In
\cite[\S3, diagram~(39)]{ABPSConjectures} they place their conjectural
parametrizations of Bernstein blocks (Conjecture~2) and of enhanced
parameter blocks (Conjecture~4) on the two sides of a square. The
correspondence for supercuspidal representations of Levi subgroups is
to induce its bottom arrow. The representation-side conjecture goes
back to Aubert--Baum--Plymen \cite{ABPHecke}, who also related the
extended quotient of the second kind to Kazhdan--Lusztig--Reeder
parameters in the split principal series \cite{ABPPrincipal}.
An earlier formulation of ABPS gives a triangle with a finite-to-one
map to Langlands parameters and a bijective version with enhancements
for split principal-series blocks under their residue-characteristic
hypothesis \cite[Conjecture~1.1 and Theorem~1.3]{ABPSGeometric}.

Solleveld constructed the representation-side parametrization using
endomorphism algebras of progenerators
\cite[Theorems~9.7 and~9.9]{SolEndo}. Aubert--Moussaoui--Solleveld
constructed the parameter-side parametrization by generalized Springer
theory \cite[Theorem~9.3]{AMS}. For rigid enhancements, Dillery--Schwein
give the corresponding parametrization and formulate the comparison
of the two quotients \cite[Theorems~4.8.7 and~5.2.5 and
Conjecture~5.3.10]{DS}. Aubert--Xu obtain a block bijection by comparing
the two extended quotients under equivariance and cocycle hypotheses
\cite[Theorem~3.1.32]{AubertXuBlocks}. They use this comparison for
the intermediate series in their explicit correspondence for split
$G_2(F)$ in residue characteristic greater than three
\cite[\S10.2]{AubertXuGtwo}. For the specified tame type and parameter
data in this paper, our main result reduces the comparison of stabilizer
algebras with prescribed projective labels to depth zero.

\subsection*{The two reductions}
On the representation side fix a tame type datum $\Sigma$ of
Adler--Fintzen--Mishra--Ohara (AFMO) for $\fs$, with its
Fintzen--Kaletha--Spice twist. Write $G^0$ for its initial twisted Levi
and $\fs_{-1}=[L^0,\pi_{-1}]_{G^0}$ for the depth-zero block. The
AFMO Hecke-algebra isomorphism gives a block equivalence
$\mathcal F_\Sigma$ \cite[Theorem~4.3.11]{AFMO-reduce}.
Put $P=\Irr(G(F))^{\fs}$, $P_0=\Irr(G^0(F))^{\fs_{-1}}$ and
$D_{\rm p}=\Irr(\mathcal F_\Sigma)$. We work over $\C$. Our tame
parametrization uses tame splitting, $p\ne2$, the full compact
depth-zero Levi stabilizer, and the depth-zero hypotheses of
Theorem~\ref{thm:fiber-irr} for $G^0$, including a $0$-generic
admissible embedding and $p\nmid|W_{\rm abs}(G^0)|$.

On the Galois side fix a rigid inner twist $\rig=(G^*,G,\xi,z)$, a
finite central subgroup $Z$ through which the rigidification is
represented, and a wild reduction datum $\calD$ as in
Definition~\ref{sq:datum}. Its ingredients include a tame
$L$-embedding $j:{}^LH^*\to{}^LG^*$, a full central Weil cocycle
$\varkappa\in Z^1(W_F,Z(\widehat H))$, and a common dual cover. Here
$\widehat H=Z_{\widehat G}(\varkappa(P_F))$ is a connected Levi. A
tame cuspidal enhanced seed for a supporting Levi in $H^*$ determines
the reduced and ambient inertial classes. Denote the two sets of
relevant enhanced parameters by $Y_0$ and $Y$, respectively. On the
slice with fixed wild restriction the map $D_\vee=R^+_{\calD}:Y\to Y_0$
divides the Weil homomorphism by $\varkappa$ in $^LH^*$ and keeps the
unipotent element and its enhancement in a common centralizer.
Dat--Fintzen construct the underlying wild reduction and discuss its
relation to the local Langlands correspondence \cite[\S\S2.2.6 and~4.3]{DF}.
Here we identify its effect on the rigid Springer multiplicity spaces.
For marked quotient-character towers, with $p\ne2$ and
$p\nmid |W_{\rm abs}(G)|$, Theorem~\ref{sq:Sigma-matching}
identifies $H$ and its supporting Levi with the prescribed $G^0$ and
$L^0$, including their relative rigid realization. For a general
specified wild datum its square holds under the hypotheses of
Theorem~\ref{sq:theorem}.

Write $\mathscr E,\mathscr E_0$ for the representation-side extended
quotients and $\mathscr E^\vee,\mathscr E_0^\vee$ for the Galois-side
quotients. At the upper level their descriptions are
\begin{equation}\label{eq:intro-target}
 \mathscr E=(T^{\fs}\extquot W^{\fs})_{\kappa^{\rm p}},\qquad
 \mathscr E^\vee=(T^{\svee}\extquot W^{\svee})_{\kappa^\vee}.
\end{equation}
Here each $T$ is the orbit of the chosen cuspidal datum under connected
unramified twists, and $W$ is its finite inertial normalizer group.
Theorems~\ref{thm:tame-global-extended-quotient}
and~\ref{sq:theorem} give the following two commutative squares of
bijections, each for its own data:
\begin{equation}\label{eq:two-squares}
\begin{array}{c@{\qquad\qquad}c}
\begin{CD}
 P @>{b}>> \mathscr E\\
 @V{D_{\rm p}}VV @VV{d_{\rm p}}V\\
 P_0 @>{b_0}>> \mathscr E_0
\end{CD}
&
\begin{CD}
 Y @>{b^\vee}>> \mathscr E^\vee\\
 @V{D_\vee}VV @VV{d_\vee}V\\
 Y_0 @>{b_0^\vee}>> \mathscr E_0^\vee .
\end{CD}
\end{array}
\end{equation}
On the left, we construct $b_0$ from the depth-zero type algebra and
define $b$ using the AFMO equivalence. The quotient map $d_{\rm p}$ is
realized by an isomorphism of the normalizer algebras of the given
types. On the right, $D_\vee$ is defined on enhanced parameters before
applying the Springer parametrizations $b^\vee$ and $b_0^\vee$.

\subsection*{Comparison of the blocks}
Let $(T,W)$ and $(T_0,W_0)$ be the representation-side orbit torsors
and their inertial Weyl groups. Add a superscript $\vee$ for the
corresponding parameter data. For $t\in T$ let $A^{\rm p}_t$ be
the stabilizer algebra at the endpoint where the reflection parameters
equal one. It is the centralizer of the neutral matrix algebra in the
semisimple orbit corner. For $s\in T^\vee$ let $A^\vee_s$ be the
endomorphism algebra of the Springer induction object. Use a
subscript $0$ for the two reduced families. All simple modules are
left modules. By a comparison of these families we mean an
equivariant isomorphism of orbit torsors and inertial groups together
with algebra isomorphisms respecting the stabilizer gradings and
commuting with orbit transport on simple classes
(Definition~\ref{def:twisted-extended-quotient} and
Hypothesis~\ref{hyp:cocycle}).
\begin{theorem}[Transfer of a depth-zero comparison]\label{thm:joint-transfer}
Fix the data and choices of Theorems~\ref{thm:tame-global-extended-quotient}
and~\ref{sq:theorem}. Suppose their reduced groups and supporting Levis
are matched with the specified rigid realization. Supply a comparison
of the two reduced stabilizer-algebra families over
\[
 f_0:T_0\xrightarrow{\sim}T_0^\vee,\qquad
 \iota_0:W_0\xrightarrow{\sim}W_0^\vee ,
\]
and let $Q_0:\mathscr E_0\xrightarrow{\sim}\mathscr E_0^\vee$ be its
map on extended quotients. Then the two reductions transport the given
family to a comparison of the original stabilizer algebras. Its map
$Q:\mathscr E\to\mathscr E^\vee$ is the unique bijection such that
$d_\vee Q=Q_0d_{\rm p}$. The resulting bijection of the original blocks is
\begin{equation}\label{eq:joint-transfer}
 \mathcal L_\Sigma=(b^\vee)^{-1}Qb
     =D_\vee^{-1}\mathcal L_0D_{\rm p},\qquad
 \mathcal L_0=(b_0^\vee)^{-1}Q_0b_0 .
\end{equation}
Conversely, transporting along the inverse reductions recovers the
reduced comparison for the fixed vertical algebra maps.
\end{theorem}

The proof in \S\ref{sec:transfer-proof} gives the transported algebra
map at each point. This is the step that identifies the projective
labels at a point where the stabilizer jumps. Over the free locus the
orbit and group identification alone suffices. At a point with
stabilizer $S_2$, even two trivial multipliers admit a change of
coordinates by the sign character, interchanging the two simple
labels. In particular a depth-zero cuspidal correspondence
equivariant for unramified twists and normalizers supplies the orbit
data in the theorem. The comparison of the algebras over those orbits
is a further input.

For a pinned comparison we must fix two character ambiguities. The
first occurs in removing the Hecke parameters, on the stabilizer
$Q_V$ of a neutral simple module. The second occurs in comparing the
resulting endpoint label with its Springer label, on the stabilizer
$(W_0)_{t_0}$ of an unramified twist. At either step, a character of
the stabilizer changes a graded algebra comparison while preserving
its multiplier class. Proposition~\ref{prop:pinned-trace} determines
exactly when independently prescribed pairs of simple modules fix the
map on all labels. It suffices to detect the degrees on which some
simple character has nonzero trace. Theorem~\ref{thm:pinned-removal}
applies this criterion to the representation-side multiplicity spaces.
Theorem~\ref{thm:pinned-cross} applies it to the endpoint and Springer
algebras in the reduced block.

In the application to non-singular reduced support, Solleveld--Xu
supply $\mathcal L_0$ on the entire reduced block
\cite[(10.12), Theorem~10.4]{SXDepthZero}. The proof of
Theorem~\ref{thm:regular-block} identifies its target with $Y_0$ by
comparing the toral Springer induction series in the two support
conventions. The two reductions then give
$\mathcal L_\Sigma^{\rm ns}=D_\vee^{-1}\mathcal L_0D_{\rm p}$ directly.
Identifying this bijection with a prescribed comparison of graded
stabilizer algebras is the further question considered above.

\begin{remark}[Relation to the local Langlands correspondence]
\label{rem:jd-llc}
The enhanced Jordan decomposition for Bruhat--Tits quotients
\cite{JDII} describes how the finite cuspidal labels change with the
pinning. In forthcoming work we plan to use these labels to compare
depth-zero supercuspidal representations with enhanced Langlands
parameters. To apply the present results to tame blocks, one must
realize the prescribed tests on the two stabilizers above and
verify the multiplier and covariance hypotheses. Once these hold,
\eqref{eq:joint-transfer} transports the resulting depth-zero
comparison to the tame block.
\end{remark}

There is also a distinction between the orbit coordinate and ordinary
cuspidal support. On representations the orbit coordinate records the
central character after removing the affine Hecke parameters. On
enhanced parameters it records unnormalized cuspidal support. The
normalized map changes the Weil part by a cocharacter determined by
the Springer summand. Proposition~\ref{prop:support-comparison} gives
the exact compatibility condition for these corrections, and
Corollary~\ref{cor:joint-support} explains how it behaves under
reduction. These properties matter already for the two labels above
the scalar orbit in split $\GL_2$ (Example~\ref{ex:GL2}).

\subsection*{Outline of the proofs}
We first describe twisted extended quotients using the conjugation
maps on the spectra of the stabilizer algebras. On the representation
side, the Levi type algebra is a twisted lattice algebra. Its
commutator pairing recovers the Bernstein lattice, and its Azumaya
fibres give the matrix factor in the orbit corners of the full
normalizer algebra (Theorem~\ref{thm:source-support-realization}).
Ohara's comparison with an affine factor of a unipotent type
\cite[Theorem~4.4.1]{Ohara25} supplies the connected geometry used to
remove the reflection parameters. In a local graded chart the
semisimple coordinate changes from $\sigma$ to $\sigma-rh$, where $h$
comes from an $\mathrm{SL}_2$-triple for the nilpotent label. We use the
graded Hecke and affine reduction results of
\cite{AMS-GHA,Sol-red,Sol-survey} and keep the residual action and
multiplication defects of the original type. The full normalizer
comparison of Proposition~\ref{prop:tame-normalizer-isom} then gives
$d_{\rm p}$, with its affine orbit action and corner maps.

On the parameter side the proof is shorter. A conjugator between two
parameters in the fixed wild slice belongs to $\widehat H$. In this
group the central factors cancel, so the Weil centralizers agree
before and after reduction, as do their full inverse images in the
common cover. Generalized Springer induction therefore takes place
in the same group with the same cuspidal local system. Its
endomorphisms, their products and the multiplicity spaces identify
under reduction (Proposition~\ref{sq:algebras}). For a marked tower
the successive generic characters identify the full wild centralizer
and determine the smaller group with its prescribed enhanced seed.
The last part combines the two algebra comparisons and works out a
split $\GL_4$ family with nontrivial wild restriction and jumping
stabilizers (Example~\ref{ex:GL4-joint}).

\section{Twisted extended quotients and Bernstein orbits}\label{sec:common}
Both parametrizations use a finite group acting on an unramified orbit.
We describe their common quotient in terms of the algebras that supply
the projective labels.

\begin{definition}[Twisted extended quotient]
\label{def:twisted-extended-quotient}\label{def:TEQ}
Let $W$ act on $T$, and write $W_t$ for the stabilizer of $t$. At each
$t$ let $A_t$ be an algebra presented as a twisted group algebra
$\C_{\kappa_t}[W_t]$. Suppose we are given transports
$\Irr(A_t)\to\Irr(A_{wt})$ for $w\in W$, compatible with composition,
such that a stabilizer element acts trivially on $\Irr(A_t)$. The
\emph{twisted extended quotient} and its projection are
\begin{equation}\label{eq:TEQ}
 (T\extquot W)_\kappa=
 \left(\coprod_{t\in T}\{t\}\times\Irr(A_t)\right)\!/W,
 \qquad \pr[t,V]=[t]\in T/W .
\end{equation}
After choosing one point in each orbit this is the disjoint union of
the spectra $\Irr(A_t)$ over the chosen points. We retain the transports
when the points vary. For an isomorphism $a:A_t\to A'_{t'}$, transport
of a left module sends $V$ to $V\circ a^{-1}$.
\end{definition}
\begin{lemma}[Orbit representatives suffice for existence]\label{lem:orbit-transport}
Let $W$ act on $T$ and $W'$ on $T'$, each with a system of twisted
stabilizer algebras and coherent transports as above. Let
$\iota:W\xrightarrow{\sim}W'$ be an isomorphism and
$f:T\xrightarrow{\sim}T'$ an $\iota$-equivariant bijection.
Suppose that for one representative $t$ of each $W$-orbit,
\[
 [\kappa_t]=(\iota|_{W_t})^*[\kappa'_{f(t)}]
 \quad\text{in }H^2(W_t,\C^\times).
\]
Then the extended quotients admit a bijection over the map
$T/W\xrightarrow{\sim}T'/W'$ induced by $f$. The bijection depends
on the chosen algebra identifications at the representatives.
\end{lemma}
\begin{proof}
Equivariance identifies $W_t$ with $W'_{f(t)}$ and carries the chosen
orbit representatives to representatives in $T'$. At each such $t$,
choose a cochain realizing the displayed equality of classes.
It gives an isomorphism of the two twisted group algebras and hence
a bijection of their simple spectra.
Every quotient label has a representative over that $t$, unique up
to a stabilizer action, which is trivial on the simple spectrum.
The bijections therefore define the desired map orbit by orbit.
Alternatively, transport the chosen maps to all other points.
Two choices of transporter differ by a stabilizer element and
give the same map on simple classes. Changing the initial algebra
identifications can change these bijections.
\end{proof}

The remainder of this section identifies the orbit and the finite group
for a representation block. The Galois-side orbit is described in
\S\ref{sec:orbits}.

Let $F$ be a non-archimedean local field with residue field $\ff$
of characteristic $p$, and let $G$ be a connected reductive $F$-group.
For a connected reductive group $H$, write $W_{\mathrm{abs}}(H)$
for its absolute Weyl group.
Fix an $F$-Levi subgroup $L$ and a supercuspidal representation $\pi$ of $L(F)$.
Write $\fs=[L,\pi]_G$ for the inertial class and $\Rep(G(F))^\fs$ for its Bernstein block.

\subsection{Unramified characters and the Bernstein torus}
Let
\[
X_{\nr}(L) := \Hom_{\mathrm{cont}}(L(F)/{}^{\circ}{L(F)},\C^\times)
\]
be the group of unramified characters, where ${}^{\circ}{L(F)}$ is the subgroup generated by compact open subgroups.
Let
\[
X_{\nr}(L,\pi) := \{\chi\in X_{\nr}(L): \pi\otimes \chi \simeq \pi\}.
\]
Then the \emph{Bernstein torus} attached to $\fs$ is
\[
X_{\nr}(L)/X_{\nr}(L,\pi),
\]
a complex torus. We write $T^\fs=\mathcal O^\fs=\{\pi\otimes\chi:
\chi\in X_{\nr}(L)\}$ for the orbit torsor. The chosen $\pi$ identifies
this torsor with the quotient torus, and we use this identification
whenever we write lattice coordinates on $T^\fs$.

\subsection{The finite inertial normalizer group}
Define the inertial normalizer
\[
N^{\fs} := \{g\in G(F): gLg^{-1}=L\ \text{and}\ \exists\chi\in X_{\nr}(L)\ \text{with}\ \pi^g \simeq \pi\otimes \chi\}.
\]
Then
\[
W^{\fs} := N^{\fs}/L(F)
\]
is finite and acts algebraically on $T^{\fs}$ (equivalently on $\mathcal{O}^{\fs}$).

\subsection{The Bernstein center}
Let $\ft=[L,\pi]_L$. We write
\[
\Rep^{\fs}(G(F)):=\Rep(G(F))^{\fs},
\qquad
\Rep^{\ft}(L(F)):=\Rep(L(F))^{\ft},
\]
where $\Rep(L(F))^{\ft}$ denotes the Bernstein block of $L(F)$ attached to $\ft$.
We also write
\[
\Irr^{\ft}(L(F))
:=\Irr\bigl(\Rep^{\ft}(L(F))\bigr)
=\{\pi\otimes\chi:\chi\in X_{\nr}(L)\}
=\mathcal O^{\fs}.
\]
The center $\mathfrak{Z}$ (resp.\ $\mathfrak{Z}^{\fs}$, resp.\ $\mathfrak{Z}^{\ft}$) of the category $\Rep(G(F))$ (resp.\ $\Rep^{\fs}(G(F))$,
resp.\ $\Rep^{\ft}(L(F))$) is called its \emph{Bernstein center}. Recall here that the \emph{center} of an abelian category is the endomorphism ring of the identity functor.

\begin{theorem}[Bernstein]\label{thm:Bernstein}
The center of the specified block is
\[
\mathfrak Z^\fs\simeq
\mathcal O\bigl(\Irr^{\ft}(L(F))/W^\fs\bigr).
\]
\end{theorem}

\clearpage
\part{Tame Bernstein blocks}

\section{Depth-zero Hecke algebras}\label{sec:hecke-input}
\subsection{Setup: a depth-zero cuspidal datum in a Levi}\label{sec:notations}

Fix an $F$-Levi subgroup $L\subset G$.
Let $x$ be a point in the Bruhat--Tits building of $L$ whose image $[x]$ in the reduced building is a vertex.
Let $\rL^{\circ}(\ff) := L(F)_{x,0}/L(F)_{x,0+}$ be the reductive quotient.
Put $K_L=L(F)_x$ and $K_L^+=L(F)_{x,0+}$. Choose
\[
\underline\tau_{\rb}\in\Irr(K_L/K_L^+)
\]
whose restriction to $L(F)_{x,0}/K_L^+$ has cuspidal irreducible
constituents, and let $\tau_{\rb}$ be its inflation to $K_L$. If starting
from one connected cuspidal constituent, form its inertia group in
$K_L/K_L^+$ and use Clifford theory, with its projective multiplicity
space, to obtain $\underline\tau_{\rb}$. No extension of an arbitrary
connected constituent to the full quotient is assumed.

Set $H=\Stab_{L(F)_{[x]}}(\tau_{\rb})$. The full enlarged stabilizer $K_L$
is normal in $L(F)_{[x]}$, so this is an inertia subgroup containing $K_L$.
Choose an irreducible $\tau_H$ lying over $\tau_{\rb}$ such that
\[
\pi=\cind_H^{L(F)}\tau_H
\]
is irreducible. This is the usual compact-induction description of a
depth-zero supercuspidal representation. The restriction
$\tau_H|_{K_L}$ is $\tau_{\rb}$-isotypic and may have multiplicity larger
than one.

\subsection{Presentation of a depth-zero Hecke algebra}

Choose an AFMO depth-zero $\fs$-type $(K_G,\rho_G)$ covering the
Levi type $(K_L,\tau_{\rb})$, with a $0$-generic admissible building
embedding as in \cite[\S5]{AFMO-structure}. Write
\[
P_G=\cind_{K_G}^{G(F)}\rho_G,\qquad
\mathcal H^{\mathrm R}_G=\End_{G(F)}(P_G),\qquad
\mH_G=(\mathcal H^{\mathrm R}_G)^{\op}.
\]
Throughout, $\Mod$ denotes left modules. The type equivalence sends $V$
to $\Hom_{G(F)}(P_G,V)$, with the action by precomposition.
We also denote this left-module algebra by
$\mathcal H(G(F),(K_G,\rho_G))$.

\begin{theorem}[AFMO structure]\label{thm:AFMO-structure}
In the stated AFMO setup, the source/right-module algebra has a presentation
\[
\mathcal H^{\mathrm R}_G\cong
\C_{\mu^{\mathrm R}}[\Omega(\tau_{\rb})]
\ltimes\mathcal H_{\mathrm{aff}}(W_{\mathrm{aff}},\mathbf q).
\]
The affine factor is the Coxeter Hecke algebra of $W_{\mathrm{aff}}$.
Here $\Omega(\tau_{\rb})$ is the length-zero subgroup; the full support
indexing group is $\mathcal W=\Omega(\tau_{\rb})\ltimes W_{\mathrm{aff}}$.
The group $\Omega(\tau_{\rb})$ need not be finite. The full Levi
translation subgroup need not lie in it.
\end{theorem}
\begin{proof}
This is \cite[Theorem~5.3.6]{AFMO-structure}, with the support notation of
its Definition~3.4.15 and Proposition~5.3.2.
\end{proof}

\begin{lemma}[Opposite convention for scalar multipliers]
\label{lem:opposite-cocycle}
If $u_gu_h=\mu(g,h)u_{gh}$, put $v_g=u_g^{-1}$ in the opposite algebra,
where the inverse is taken in the original twisted group algebra. Then
\[
v_g\circ v_h=u_h^{-1}u_g^{-1}
=\mu(g,h)^{-1}u_{gh}^{-1}.
\]
Thus $\C_\mu[Q]^{\op}\simeq\C_{\mu^{-1}}[Q]$. With this indexing,
$v_g$ has original support $g^{-1}$.
\end{lemma}
\begin{proof}
Invert $u_gu_h=\mu(g,h)u_{gh}$ in the original algebra. The displayed
identity is the multiplication rule in the opposite algebra.
\end{proof}

\begin{lemma}[Affine reflection parameters]\label{lem:q-parameters-powers}
For the depth-zero AFMO datum with full compact Levi stabilizer and
$0$-generic admissible embedding, let $s$ be a simple reflection of the
\emph{relevant} affine Coxeter group. In the positive AFMO normalization,
\[
 (T_s+1)(T_s-q_s)=0,\qquad
 q_s=q^{m_s},\qquad m_s\in\Z_{\geq1},
\]
where $q$ is the residue-field cardinality.
\end{lemma}
\begin{proof}
Ohara's affine comparison identifies these quadratic parameters with
unipotent affine Hecke parameters
\cite[\S4.2 and Theorem 4.4.1]{Ohara25}. The unipotent parameter
classification gives integral powers of the residue-field cardinality
in this Iwahori--Matsumoto convention; see
\cite[Table 1 and its preceding classification references]{Sol2025}.
Passing to an unramified extension of degree $f$ replaces an exponent
$m$ by $fm$ and preserves integrality.

The relevance condition excludes parameter $1$: in
\cite[Proposition 3.8.19]{AFMO-structure} the coefficient of the
reflection operator is nonzero, and the normalization in Proposition
3.8.23 and Choice 3.10.3 gives coefficient $q_s-1$.
In characteristic zero, Proposition 3.9.1 identifies $q_s$ with a
positive ratio of dimensions. Choose the larger member of the pair
$\{q_s,q_s^{-1}\}$, as in the positive normalization. Thus $q_s>1$,
and the integral exponent is strictly positive.
\end{proof}

This assertion concerns source Coxeter parameters. It does not exclude
zero parameters in a local graded reduction or the auxiliary geometric
label $Q_\beta^*=1$ in the second row of Proposition~\ref{prop:local-coefficients}.
In that proposition the source parallel-wall label satisfies
$Q_\alpha^*>1$, so
$a(Q_\alpha^*)/a(Q)>0$, where $a(Q)=Q^{1/2}-Q^{-1/2}$.

The root-indexed reducibility parameters $q_\alpha,q_{\alpha*}$,
which may involve half-powers, are different from the Coxeter parameters
$q_s$ used here.

\paragraph{Scope of the hypotheses.}
The main results are stated for the tame setting used in the AFMO
reduction. Ohara's depth-zero Theorem~4.4.1 itself requires neither
tame splitting nor a Weyl-order restriction: it uses the depth-zero
data in its \S4.2, including the full-stabilizer option. Tame splitting
and $p\ne2$ enter the AFMO positive-depth construction, while
$p\nmid|W_{\mathrm{abs}}(G)|$ is used in its exhaustion of Bernstein
blocks; see \cite[\S\S4.1, 4.5]{AFMO-reduce} and
\cite[\S\S4.5--4.6]{Ohara25}.

\section{The depth-zero parametrization}\label{sec:rep-block}

\subsection{The endpoint and its coordinates}
Let $\pi$ be depth-zero and let \(\ft = [L,\pi]_L\) and \(\fs=[L,\pi]_G\) as before. Put \(A:=\mathfrak{Z}^\ft\) and \(Z:=\mathfrak{Z}^\fs\). Then $Z=A^{W^\fs}$.
\begin{definition}\label{def:cusp-map}
Let \(\mathfrak s=[L,\pi]_G\) be a Bernstein block. Denote by
\[
\operatorname{cusp} :
\Irr\bigl(\Rep(G)^{\mathfrak s}\bigr)\longrightarrow
\Irr^{\ft}(L(F))/W^{\mathfrak s}
\]
the cuspidal-support map.

For \(t\in \Irr^{\ft}(L(F))\), we write
\[
\Irr\bigl(\Rep(G)^{\mathfrak s}\bigr)_{[t]}
\]
for the fibre of this map above the class
\([t]\in\Irr^{\ft}(L(F))/W^{\mathfrak s}\); that is,
\[
\Irr\bigl(\Rep(G)^{\mathfrak s}\bigr)_{[t]}
=
\left\{
\pi'\in \Irr\bigl(\Rep(G)^{\mathfrak s}\bigr)
\;\middle|\;
\operatorname{cusp}(\pi')=[t]
\right\}.
\]
\end{definition}

Put $T^\fs=\Irr^\ft(L(F))$. A choice of $\pi$ identifies this
unramified orbit with a complex torus; its $W^\fs$-action includes
the translation part determined by ${}^n\pi$. Let $B_1$ be the fixed left AFMO presentation of the type algebra and
$B_0$ its source endpoint. Write $D$ for the actual Levi type algebra,
embedded in $B_0$ by Proposition~\ref{prop:Levi-qzero-embedding}. Thus
\[
 Z(D)=\cO(T^\fs),\qquad
 B_0=\bigoplus_{w\in W^\fs}Du_w.
\]
The cocycle $\kappa_t$ is the scalar multiplier obtained from the
$t$-corner in Proposition~\ref{prop:central-fiber-semisimple}, after passing to the
semisimple quotient of the central fiber $B_{0,[t]}$.
Denote the resulting parametrization by
\[
 \cl_0:\Irr(B_0)\longrightarrow(T^\fs\extquot W^\fs)_\kappa.
\]
The type equivalence followed by the chosen presentation gives
\[
 \typ_\fs:\Irr(\Rep(G(F))^\fs)\xrightarrow{\sim}\Irr(B_1).
\]

\paragraph{Neutral invariants.}
The finite crossed presentation of Proposition~\ref{prop:local-finite} has the form
\begin{equation}\label{eq:statement-finite-presentation}
 A_\epsilon=\C[\Lambda_c]\otimes
                  \mathcal H(W_{\aff},q^\epsilon),\qquad
 B_\epsilon=\bigoplus_{g\in Q}A_\epsilon u_g,
 \quad Q=\Omega/\Lambda_c.
\end{equation}
Here $W_{\aff}=X_{\aff}\rtimes W_f$,
$X=X_{\aff}\oplus\Lambda_c$, and $T_X=\Hom(X,\C^\times)$.
For a simple $A_\epsilon$-module $V$, let
\[
 p_\epsilon(V)=([u_V]_{W_f},\tau_V),\qquad
 \tau_V\in\Hom(\Lambda_c,\C^\times),
\]
where $[u_V]_{W_f}$ is the $W_f$-orbit of the unitary parts of its Bernstein weights in
$T_X$ and $\tau_V$ is the full character of $\C[\Lambda_c]$.
For a simple $B_\epsilon$-module $M$, define the finite set
\begin{equation}\label{eq:neutral-invariant}
 \mathcal P_\epsilon(M)
 =\{p_\epsilon(V):V\text{ is a simple constituent of }
                                  \Res_{A_\epsilon}M\}.
\end{equation}
This definition does not require a single constituent to be
distinguished, or an action of $Q$ on the displayed pairs to be
specified independently of the modules.

\begin{theorem}[Depth-zero parametrization]
\label{thm:fiber-irr}
Fix the depth-zero AFMO datum of \S\ref{sec:notations}, with full compact Levi stabilizer and
$0$-generic admissible embedding. Assume that $G$ splits over a
tamely ramified extension and $p\nmid|W_{\mathrm{abs}}(G)|$.
Fix the source presentation, compatible local reductions and sign
conventions, and the global projective-label choices in Propositions~\ref{prop:equiv-full-q}--\ref{prop:assembly-global}.

The corner cocycles $\kappa_t$ of $B_0$, with their homogeneous
corner transports, define a $W^\fs$-equivariant twisted extended
quotient. For this fixed source algebra their classes are independent
of the choices of matrix coordinates and nonzero intertwiners.
There is a bijection
\begin{equation}\label{eq:depth-parametrization}
 \mathfrak b_\fs=\cl_0\circ\Theta_B\circ\typ_\fs:
 \Irr(\Rep(G(F))^\fs)\xrightarrow{\sim}
                      (T^\fs\extquot W^\fs)_\kappa
\end{equation}
with the following properties.
\begin{enumerate}[label=\textup{(\roman*)},leftmargin=*]
\item The neutral constituents of $\typ_\fs(\Pi)$ and of
$\Theta_B\typ_\fs(\Pi)$ correspond under $\Theta_A$, with the
same multiplicities. Each paired constituent has the same unitary
$W_f$-orbit and independent toral character. In particular,
\eqref{eq:neutral-invariant-preservation} holds for the associated
$B_1$-module.
\item In every specified connected geometric chart the neutral map
is the enhanced-parameter shift
\begin{equation}\label{eq:geometric-shift}
 (y,\sigma,\rho)\longmapsto(y,\sigma-rh,\rho_0),
 \qquad [\sigma-rh,y]=0,
\end{equation}
factor by factor when several geometric parameters occur. Here
$[\sigma,y]=2ry$, and $h$ comes from a compatible algebraic
$\mathrm{SL}_2$, with $[h,y]=2y$; the enhancement $\rho_0$
is transported by the component-group identification. The fixed
root normalizations, sign maps, and the endpoint
Iwahori--Matsumoto conventions are included in the chart maps.
\item Residual multiplicity spaces are transported by the cochains
of Lemma~\ref{lem:neutral-invariants}. The effective stabilizer classes,
including central multiplication defects, agree. The resulting
maps on simple classes agree under the specified chart and
full-corner identifications.
\end{enumerate}
If $\pr$ denotes the projection of the twisted extended quotient to
$T^\fs/W^\fs$, then
\begin{equation}\label{eq:endpoint-character}
 \pr\mathfrak b_\fs(\Pi)
 =\cc_{B_0}\bigl(\Theta_B\typ_\fs(\Pi)\bigr).
\end{equation}
This is the central character after parameter removal; it need not
be the ordinary cuspidal support of $\Pi$.
\end{theorem}
\begin{remark}[The support change and finite sheets]
In the geometric coordinates of \eqref{eq:geometric-shift}, the
semisimple label changes from $\sigma$ to $\sigma-rh$; transport
through the fixed reduction maps gives the corresponding local
Bernstein characters. The finite covers, sheet data, and residual
projective labels remain part of the construction. A neutral
character on $T_X$ alone need not determine a point of $T^\fs$.
Thus the shift is not asserted to be a single global translation
on $T^\fs$, nor is the endpoint central character asserted to be
independent of all projective-label choices.
\end{remark}

\section{The depth-zero source algebra}\label{sec:proof-of-theorem}
We keep the depth-zero datum from \S\ref{sec:notations}.  Put
\[
\Gamma_H:=H/K_L.
\]

\paragraph{The Bernstein torus and the Clifford lattice.}
\begin{lemma}[Depth-zero intertwining]\label{lem:depth-zero-intertwining}
The intertwining of $\tau_{\rb}$ in $L(F)$ is
\[
I_{L(F)}(\tau_{\rb})=H.
\]
\end{lemma}
\begin{proof}
Put $K_L^\circ=L(F)_{x,0}$ and $K_L^+=L(F)_{x,0+}$. If $g$ does not
fix $[x]$, the usual apartment description of parahoric intersections
shows that the image of
$K_L^\circ\cap{}^gK_L^+$ in $K_L^\circ/K_L^+$ contains the unipotent
radical of a proper parabolic subgroup. Every connected constituent of
$\tau_{\rb}$ is cuspidal, so it has no invariants under that radical.
On the other hand, ${}^g\tau_{\rb}$ is trivial on ${}^gK_L^+$.
A map in
$\Hom_{K_L\cap{}^gK_L}({}^g\tau_{\rb},\tau_{\rb})$
would therefore have zero image. This excludes $g$ outside
$L(F)_{[x]}$.

Every element of $L(F)_{[x]}$ normalizes $K_L$, since it acts by translation
on the vectorial fibre of the enlarged building. For such $g$,
intertwining is exactly ${}^g\tau_{\rb}\simeq\tau_{\rb}$, which defines
$H$. Conversely every element of $H$ intertwines $\tau_{\rb}$.
\end{proof}

\begin{lemma}\label{lem:reduced-enlarged-stabilizer}
One has
\[
L(F)_{[x]}\cap{}^\circ L(F)=L(F)_x=K_L,
\]
and hence $H\cap{}^\circ L(F)=K_L$.
\end{lemma}
\begin{proof}
An element of $L(F)_{[x]}$ fixes the reduced point $[x]$ and acts on the
vectorial fibre of the enlarged building by its valuation vector $\nu_L(l)$,
characterized by
\[
\langle\chi,\nu_L(l)\rangle=-\operatorname{ord}_F(\chi(l)).
\]
Membership in ${}^\circ L(F)$ is equivalent to $\nu_L(l)=0$. Thus an element
of $L(F)_{[x]}\cap{}^\circ L(F)$ fixes the enlarged point $x$, giving the
claim.
\end{proof}

Define
\[
X_H^{\nr}(\tau_H)
=
\{\chi\in\Hom(\Gamma_H,\C^\times):\tau_H\otimes\chi\simeq\tau_H\}.
\]

\begin{lemma}\label{lem:Ts-as-HomLambda}
Restriction of ordinary unramified characters induces an isomorphism
\[
T^\fs\xrightarrow{\sim}
\Hom(\Gamma_H,\C^\times)/X_H^{\nr}(\tau_H).
\]
Consequently the character lattice
\[
\Lambda_\fs:=X^*(T^\fs)
\]
is
\[
\Lambda_\fs
=
\{\lambda\in\Gamma_H:\chi(\lambda)=1\ \text{for every }\chi\in X_H^{\nr}(\tau_H)\}.
\]
In particular, canonically
\[
T^\fs\cong\Hom(\Lambda_\fs,\C^\times).
\]
\end{lemma}
\begin{proof}
Put $Q=L(F)/{}^\circ L(F)$. By Lemma~\ref{lem:reduced-enlarged-stabilizer},
the image of $H$ in $Q$ is $\Gamma_H$, a free abelian subgroup.
It has finite index: $H$ contains $Z(L)(F)K_L$ because central elements
stabilize $\tau_{\rb}$, and the image of $Z(L)(F)$ has finite index in
the valuation lattice $Q$. Since
$\C^\times$ is divisible, restriction gives a surjection
\[
r:X_{\nr}(L)=\Hom(Q,\C^\times)\twoheadrightarrow\Hom(\Gamma_H,\C^\times).
\]
We claim
\[
r^{-1}(X_H^{\nr}(\tau_H))=X_{\nr}(L,\pi).
\]
One inclusion follows from $\pi\simeq\cind_H^{L(F)}\tau_H$. Conversely, if
$\pi\otimes\nu\simeq\pi$, Mackey theory gives a double coset $HgH$
contributing a nonzero intertwining space between $\tau_H$ and
$\tau_H\otimes\nu|_H$. Since $\nu$ is trivial on compact subgroups,
restriction to $K_L\cap{}^gK_L$ shows that $g$ intertwines $\tau_{\rb}$.
Lemma~\ref{lem:depth-zero-intertwining} gives $g\in H$, and hence
$\tau_H\otimes\nu|_H\simeq\tau_H$. Taking the quotient and then character
lattices proves the assertion.
\end{proof}

\begin{proposition}[The actual Levi support algebra]\label{prop:Levi-support-algebra}
Let $P_L=\cind_{K_L}^{L(F)}\tau_{\rb}$. There is a support-graded
identification
\[
D:=\End_{L(F)}(P_L)^{\op}\simeq\C_{\mu_H}[\Gamma_H]
\]
for a normalized scalar cocycle $\mu_H$. We index a left homogeneous
operator by $\gamma=hK_L$ so that, on $\Hom_{K_L}(\tau_{\rb},V)$, it has
the form $f\mapsto V(h)fA_h$, with a fixed intertwiner
$A_h\tau_{\rb}(k)=\tau_{\rb}(h^{-1}kh)A_h$.
\end{proposition}
\begin{proof}
Lemma~\ref{lem:depth-zero-intertwining} shows that the Hecke support is
$H$. Since $H$ normalizes $K_L$, each relevant double coset is $hK_L$ and
its intertwining line has dimension one by Schur's lemma. Its nonzero
operators are invertible, their products lie in the corresponding product
line, and associativity gives a scalar $2$-cocycle after choosing bases.
The degree-zero line is normalized to the identity.

In the compact-induction model generated by $g\otimes v$, the
endomorphism $T_h(g\otimes v)=gh\otimes A_hv$ is well defined by the
intertwining relation. Precomposition by $T_h$ acts as $V(h)fA_h$.
Composition in the opposite endomorphism algebra multiplies the degrees
in the order $h_1h_2$, establishing the asserted left convention.
\end{proof}

Define
\[
b_{\mu_H}(\lambda,\lambda')=
\frac{\mu_H(\lambda,\lambda')}{\mu_H(\lambda',\lambda)},
\qquad
R_H:=\Rad(b_{\mu_H}).
\]

\begin{proposition}[Bernstein lattice inside the AFMO Levi support group]
\label{prop:AFMO-Bernstein-lattice}
One has
\[
\Lambda_\fs=R_H.
\]
After a cochain normalization on $R_H$,
\[
Z\bigl(\C_{\mu_H}[\Gamma_H]\bigr)\simeq\C[\Lambda_\fs].
\]
\end{proposition}
\begin{proof}
Let $E=\Hom_{L(F)}(P_L,\pi)$, with its left $D$-action.
Mackey theory and Lemma~\ref{lem:depth-zero-intertwining} identify its
underlying vector space with the multiplicity space in
$\tau_H|_{K_L}\simeq\tau_{\rb}\otimes E$. It is a finite-dimensional
irreducible $D$-module; thus its homogeneous operators have multiplier
$\mu_H$ in the convention of Proposition~\ref{prop:Levi-support-algebra}.
Put $\widehat{\Gamma_H}=\Hom(\Gamma_H,\C^\times)$ and
\[
S(E)=\{\chi\in\widehat{\Gamma_H}:E\otimes\chi\simeq E\}.
\]
Twisting $\tau_H$ shows $S(E)=X_H^{\nr}(\tau_H)$. If $r\in R_H$, the twisted
basis element $u_r$ is central and acts on $E$ by a scalar, so
$S(E)\subset R_H^\perp$. Conversely, if $d=\dim E$, determinants of the
projective commutation relations give $d\Gamma_H\subset R_H$. Hence
$\Gamma_H/R_H$ is finite, and the central reduction at a fixed $R_H$-character
is a twisted group algebra of this quotient with nondegenerate commutator
pairing. It is therefore a full matrix algebra. Thus a character trivial on
$R_H$ does not change $E$, and $S(E)=R_H^\perp$. Lemma~\ref{lem:Ts-as-HomLambda}
now gives $\Lambda_\fs=(R_H^\perp)^\perp=R_H$. The center statement follows
because a twisted basis element is central exactly on the radical; the
restriction of the multiplier to the free abelian group $R_H$ is symmetric and
hence cohomologically trivial.
\end{proof}

\paragraph{The genuine affine Bernstein action.}
Choose $\pi$ as base point of the torsor $\Irr^\ft(L(F))$ and write
\[
\jmath_\pi:\Irr^\ft(L(F))\xrightarrow{\sim}T^\fs,
\qquad
\pi\otimes\nu\longmapsto[\nu].
\]
For $w\in W^\fs$, choose $n\in N^\fs$ representing $w$ and
$\chi_n\in X_{\nr}(L)$ with ${}^n\pi\simeq\pi\otimes\chi_n$. Define
\[
w_{\mathrm{lin}}[\nu]=[{}^n\nu],
\qquad
\eta_w=[\chi_n]\in T^\fs.
\]
These are independent of the choices and satisfy
\[
\eta_{w_1w_2}=\eta_{w_1}\,w_{1,\mathrm{lin}}(\eta_{w_2}).
\]
Thus
\[
w\star t:=\eta_w\,w_{\mathrm{lin}}(t)
\]
is an affine action.

\begin{proposition}\label{prop:center}
The map $\jmath_\pi$ is equivariant for the Bernstein action on
$\Irr^\ft(L(F))$ and the affine action $\star$ on $T^\fs$. Consequently
\[
\mathfrak Z^\fs\simeq\mathcal O(T^\fs)^{W^\fs,\star}
=\C[\Lambda_\fs]^{W^\fs,\mathrm{aff}}.
\]
\end{proposition}
\begin{proof}
For a representative $n$ of $w$,
\[
{}^n(\pi\otimes\nu)\simeq\pi\otimes\chi_n\,{}^n\nu.
\]
Applying $\jmath_\pi$ gives
\[
\jmath_\pi({}^n(\pi\otimes\nu))
=\eta_w\,w_{\mathrm{lin}}([\nu]).
\]
The center formula is Bernstein's description of the center as regular
functions on the inertial quotient.
\end{proof}

\subsection{Normalization of the Levi center}\label{subsec:Levi-center-normalization}
By Proposition~\ref{prop:AFMO-Bernstein-lattice}, choose multiplicative
central basis elements $z_\lambda$ of $D$, for $\lambda\in\Lambda_\fs$.
Let $\xi_\pi(\lambda)$ be their scalar on
$E_\pi=\Hom_{L(F)}(P_L,\pi)$ and replace $z_\lambda$ by
$\xi_\pi(\lambda)^{-1}z_\lambda$. This preserves multiplicativity and
makes the central point of $E_\pi$ equal to $1$.

\begin{proposition}[Unramified normalization]\label{prop:Levi-center-normalized}
In these coordinates, the central point of
$E_{\pi\otimes\nu}$ is $\nu|_{\Lambda_\fs}$. Consequently
\[
Z(D)=A_\fs:=\C[\Lambda_\fs]\simeq\mathcal O(T^\fs)
\]
is identified with the actual Bernstein torus based at $\pi$.
\end{proposition}
\begin{proof}
Since $\nu$ is trivial on $K_L$, the spaces
$\Hom_{K_L}(\tau_{\rb},\pi)$ and
$\Hom_{K_L}(\tau_{\rb},\pi\otimes\nu)$ are naturally the same vector
space. Proposition~\ref{prop:Levi-support-algebra} shows that a degree
$\gamma=hK_L$ operator on the latter is $\nu(h)$ times its operator on
the former. Thus normalized $z_\lambda$ acts on $E_{\pi\otimes\nu}$
by $\nu(\lambda)$. Lemma~\ref{lem:Ts-as-HomLambda} identifies the resulting
characters with $T^\fs$; the finite self-twist quotient is included in that
identification. Any other multiplicative homogeneous central basis differs
by a character of $\Lambda_\fs$, and that character cancels on dividing
by its scalar on $E_\pi$. Thus the normalized identification is independent
of this basis choice.
\end{proof}

\begin{proposition}\label{prop:Levi-Azumaya}
The algebra $D$ is Azumaya over $A_\fs$.  For every $t\in T^\fs$,
\[
D_t:=D\otimes_{A_\fs}\C_t\simeq M_d(\C),
\qquad d^2=[\Gamma_H:\Lambda_\fs].
\]
\end{proposition}
\begin{proof}
The quotient $\Gamma_H/\Lambda_\fs$ is finite. Choose representatives
$\gamma$ for its cosets. The corresponding homogeneous elements give
$D=\bigoplus_\gamma A_\fs u_\gamma$, so $D$ is finite free of rank
$[\Gamma_H:\Lambda_\fs]$ over its center. The central fibre is a twisted
group algebra of this finite abelian quotient, and the induced commutator
pairing is nondegenerate by the definition of $\Lambda_\fs$.  Hence the fibre
is finite-dimensional central simple and therefore a full matrix algebra.  Its
dimension gives the displayed value of $d^2$. The finite free algebra
therefore has central simple fibres at all closed points. Equivalently,
the map $D\otimes_{A_\fs}D^{\op}\to\End_{A_\fs}(D)$ is an isomorphism
on all such fibres and hence an isomorphism, by finite freeness and
Nakayama's lemma. This proves the Azumaya assertion.
\end{proof}

\subsection{The source support extension and its Bernstein action}
\label{subsec:support-realization}
Let $B_0$ be the algebra obtained by putting the affine reflection
parameters equal to $1$ in the fixed left AFMO presentation of $\mH_G$.
It has a twisted support presentation $B_0\simeq\C_{\bar\mu}[\mathcal W]$.
We identify its support extension and its Levi subalgebra with the
specified Bernstein data. This construction uses the source AFMO
presentation and the compact inducing datum.

\subsubsection{The source normalizer and the inertial quotient}
Keep the notation of \S\ref{sec:notations}, and abbreviate $K_L$ to $K$. Thus $[x]$ is a reduced vertex for $L$,
\[
 K=K_L=L(F)_x,\qquad K^+=L(F)_{x,0+},\qquad \tau=\tau_{\mathrm b},
\]
and the connected constituents of $\tau$ are cuspidal. Set
\[
 H=\operatorname{Stab}_{L(F)_{[x]}}(\tau),\qquad
 \Gamma_H=H/K,
 \qquad \pi=\cind_{H}^{L(F)}\tau_H,
\]
with $\tau_H|_K$ isotypic of type $\tau$ and $\pi$ irreducible.
The source support normalizer is
\begin{equation}\label{eq:source-normalizer}
 \mathcal N
 =\{n\in N_G(L)(F):[nx]_L=[x]_L,\ {}^n\tau\simeq\tau\}.
\end{equation}
Here $[nx]_L=[x]_L$ means $nx\in x+X_*(A_L)\otimes\mathbb R$, where
$A_L$ is the maximal split subtorus of $Z(L)$. The full point stabilizer
$K$ is normalized by every element preserving this reduced point:
translation in the central directions does not change the stabilizer.
In particular $K\triangleleft\mathcal N$ and $\mathcal N\cap L(F)=H$.
For the depth-zero AFMO family, the actual support indexing group is
\begin{equation}\label{eq:source-support-group}
 \mathcal W=\mathcal N/K.
\end{equation}
This is the normalizer furnished by \cite[Definition~3.4.15 and
Proposition~5.3.2]{AFMO-structure}, with the full compact stabilizer convention fixed above. It is the entire support group, not just its length-zero subgroup.

\begin{lemma}[Cuspidality detects the reduced vertex]
\label{lem:cuspidal-vertex}
Let $x,y\in\mathcal B(L,F)$ have reduced vertex images, and suppose
$[x]\ne[y]$. Let $\tau_x$ be a representation of $L(F)_x$, trivial on
$L(F)_{x,0+}$, with cuspidal connected constituents. If a representation
$\sigma_y$ of a subgroup $J_y\supseteq L(F)_{y,0+}$ is trivial on
$L(F)_{y,0+}$, then
\[
 \Hom_{L(F)_x\cap J_y}(\tau_x,\sigma_y)=0.
\]
\end{lemma}
\begin{proof}
Choose an apartment of $L$ containing $x$ and $y$. The image of
\[
 L(F)_{x,0}\cap L(F)_{y,0+}
 \quad\hbox{in}\quad L(F)_{x,0}/L(F)_{x,0+}
\]
contains the unipotent radical of a proper rational parabolic of the
connected reductive quotient at $x$. Explicitly, take the parabolic
defined by the direction from $x$ to $y$. Its unipotent root groups are
the reductions of affine root groups whose roots vanish at $x$ and
have positive value at $y$. Those groups belong to the displayed
intersection. At a reduced vertex the gradients of the walls through
the vertex span the noncentral apartment directions. Since the
displacement is noncentral, the resulting parabolic is proper. This also proves the parahoric-intersection assertion used in
Lemma~\ref{lem:depth-zero-intertwining}.

On this intersection $\sigma_y$ is trivial. An intertwining map must
therefore factor through the coinvariants of $\tau_x$ for its image
in the finite reductive quotient. These vanish by cuspidality of every
connected constituent. Over $\C$, invariants and coinvariants for the
finite unipotent subgroup have the same vanishing. Hence the map is zero.
\end{proof}

\begin{proposition}[The actual inertial quotient]
\label{prop:inertial-support-quotient}
Writing $N^{\fs}$ for the inertial normalizer of $[L,\pi]$, one has
\[
 \mathcal N\subset N^{\fs},\qquad
 N^{\fs}=L(F)\mathcal N.
\]
Consequently the natural map gives an exact sequence
\begin{equation}\label{eq:B-group-seq}
 1\longrightarrow\Gamma_H\longrightarrow\mathcal W
 \longrightarrow W^{\fs}\longrightarrow1.
\end{equation}
\end{proposition}
\begin{proof}
First take $n\in\mathcal N$. It normalizes $H$, and
$({}^n\tau_H)|_K$ is isotypic of type $\tau$. The compact Clifford
description identifies its multiplicity space with a finite-dimensional
simple module for the actual Levi algebra
\[
 D=\End_{L(F)}(\cind_{K}^{L(F)}\tau)^{\mathrm{op}}
   \simeq\C_{\mu_H}[\Gamma_H].
\]
By Propositions~\ref{prop:AFMO-Bernstein-lattice}--\ref{prop:Levi-Azumaya}, $D$ is Azumaya over its normalized center
$A_\fs=\C[\Lambda_\fs]$, and there is exactly one simple module at
each closed point of that center. Unramified twists of $\tau_H$ attain
every such point: restriction of unramified characters surjects onto
$\Hom(\Gamma_H,\C^\times)$, and restriction from that character group
to the finite-index sublattice $\Lambda_\fs$ is also surjective.
Thus the two multiplicity modules become isomorphic after an unramified
twist. The $K$-isotypic Clifford equivalence gives
\[
 {}^n\tau_H\simeq\tau_H\otimes\nu|_H
\]
for some $\nu\in X_{\nr}(L)$. Compact induction gives
${}^n\pi\simeq\pi\otimes\nu$, proving $\mathcal N\subset N^{\fs}$.

Conversely let $n\in N^{\fs}$. Since unramified characters are trivial
on $K$, the representation ${}^n\pi$ contains $\tau$. Apply compact
Mackey theory to
\[
 {}^n\pi\simeq
 \cind_{{}^nH}^{L(F)}({}^n\tau_H).
\]
There is $\ell\in L(F)$ such that, with $g=\ell n$,
\begin{equation}\label{eq:inertial-support-Mackey}
 \Hom_{K\cap{}^gH}(\tau,{}^g\tau_H)\ne0.
\end{equation}
Indeed the compact Mackey restriction is a direct sum over double
cosets, and the finite-dimensional $K$-module $\tau$ detects one
nonzero summand. The group ${}^gH$ contains $L(F)_{gx,0+}$ and its
representation ${}^g\tau_H$ is trivial on that subgroup. Lemma
\ref{lem:cuspidal-vertex} applied to \eqref{eq:inertial-support-Mackey} yields
$[gx]_L=[x]_L$.

It follows that ${}^gK=K$, and hence $K\subset{}^gH$.
The restriction of ${}^g\tau_H$ to $K$ is isotypic of type
${}^g\tau$. Equation~\eqref{eq:inertial-support-Mackey} now gives
${}^g\tau\simeq\tau$. Therefore $g\in\mathcal N$, proving
$N^{\fs}=L(F)\mathcal N$. The kernel of the quotient map
$\mathcal N/K\to N^{\fs}/L(F)$ is
$(\mathcal N\cap L(F))/K=H/K$, as required.
\end{proof}

\subsubsection{From the normalizer cocycle to the AFMO endpoint}
The preceding quotient does not alone identify the Levi algebra inside
$B_0$. The full normalizer-line cocycle must be compared with the
multiplier obtained by setting the affine Hecke parameters equal to one.
The following normalization provides that comparison.

For this calculation only, start in the AFMO right/endomorphism
convention. Choose the intertwiners
\[
 T_n\in\Hom_K({}^n\tau,\tau),\qquad n\in\mathcal N,
\]
as in \cite[Choice~3.10.3]{AFMO-structure}. Extend the notation $\mu^{\mathrm R}$ from
Theorem~\ref{thm:AFMO-structure} to their full scalar cocycle:
\[
 T_nT_m=\mu^{\mathrm R}(\bar n,\bar m)T_{nm}.
\]
The convention $T_{nk}=T_n\tau(k)$ includes the representative defects
in this formula; its quotient indices lie in $\mathcal W$.
Write
\[
 \mathcal W=\Omega\ltimes W_{\aff},\qquad
 \mu^{\mathrm R}_\Omega=\mu^{\mathrm R}|_{\Omega\times\Omega},
\]
and let $p_\Omega:\mathcal W\to\Omega$ be the projection. Define
\[
 \overline\mu^{\mathrm R}=p_\Omega^*\mu^{\mathrm R}_\Omega.
\]
This is the cocycle of the right $q=1$ AFMO support algebra.

The precise AFMO relations used in the next proof are
\begin{align}
 \mu^{\mathrm R}(s,s)&=q_s &&(s\text{ simple affine}),
 \label{eq:normalizer-square}\\
 \mu^{\mathrm R}(v,w)&=1
 &&(v,w\in W_{\aff},\ \ell(vw)=\ell(v)+\ell(w)),
 \label{eq:normalizer-additive}\\
 \mu^{\mathrm R}(\omega,v)&=\mu^{\mathrm R}(v,\omega)=1
 &&(\omega\in\Omega,\ v\in W_{\aff}).
 \label{eq:normalizer-mixed}
\end{align}
Here \eqref{eq:normalizer-square} follows from AFMO Lemma~3.8.10 and the
normalized quadratic relation; \eqref{eq:normalizer-additive}--\eqref{eq:normalizer-mixed}
are AFMO Proposition~3.10.7. The parameters are invariant under the
source length-zero action. These are stronger data than a bare list of
double-coset supports.

To make the normalization in \eqref{eq:normalizer-square} explicit,
let $x$ be the base generic point for the chosen relevant alcove and
simple system, and write $K_x$ for the ambient compact subgroup in
the AFMO family. The compact Levi subgroup in this calculation is
$K=K_L$. If $w=\bar n$, then
\begin{equation}\label{eq:AFMO-function-value}
 \varphi_{x,w}(n)=D_{x,w}^{-1/2}T_n,
 \qquad D_{x,w}=[K_x:K_x\cap K_{wx}].
\end{equation}
This is \cite[Lemma~3.5.26]{AFMO-structure}, obtained from the normalized
averaging operator in its Definition~3.5.9. Thus $T_n$ is generally
different from the value of the Hecke function at $n$.

For a relevant affine simple reflection $s$, choose a lift $n_s$.
Since $n_s^2\in K$, the representative convention gives
\[
 T_{n_s^{-1}}=T_{n_s}\tau(n_s^2)^{-1},\qquad
 T_{n_s}T_{n_s^{-1}}=\mu^{\mathrm R}(s,s)\operatorname{id}.
\]
Put $D_s=D_{x,s}$. With $\operatorname{vol}(K_x)=1$, the double
coset $K_xn_sK_x$ has $D_s$ right $K_x$-cosets. Formula
\eqref{eq:AFMO-function-value} applies at both $n_s$ and $n_s^{-1}$.
The product of the values is scalar, so covariance and convolution give
\[
 (\varphi_{x,s}*\varphi_{x,s})(1)
 =D_sD_s^{-1}T_{n_s}T_{n_s^{-1}}
 =\mu^{\mathrm R}(s,s)\operatorname{id}.
\]
On the other hand, the quadratic relation gives
\[
 \varphi_{x,s}^2=(q_s-1)\varphi_{x,s}+q_s\varphi_{x,1},
 \qquad \varphi_{x,s}(1)=0,\quad
 \varphi_{x,1}(1)=\operatorname{id}.
\]
Evaluation at $1$ proves \eqref{eq:normalizer-square}.
For an Iwahori type with trivial compact representation, $D_s=q_s$
and the standard Hecke function has value $1$ at $n_s$. The chosen
$T_{n_s}$ is therefore $q_s^{1/2}\operatorname{id}$; choosing
$T_{n_s}=1$ instead rescales the Hecke function and changes its
quadratic normalization.

\begin{lemma}[An explicit cochain for the full normalizer cocycle]
\label{lem:normalizer-cocycle-normalization}
For $v=s_1\cdots s_r\in W_{\aff}$ reduced, put
\[
 q_v=\prod_{j=1}^r q_{s_j},\qquad
 d_q(\omega v)=q_v^{1/2}\quad(\omega\in\Omega),
\]
using positive real square roots. Then $d_q$ is well defined, $d_q|_\Omega=1$,
and
\begin{equation}\label{eq:normalizer-cochain}
 \mu^{\mathrm R}=(\delta d_q)\,\overline\mu^{\mathrm R},
 \qquad
 (\delta d_q)(a,b)=d_q(a)d_q(b)d_q(ab)^{-1}.
\end{equation}
\end{lemma}
\begin{proof}
Parameter invariance on conjugate simple reflections makes $q_v$
independent of the reduced expression: every braid move preserves its
parameter product, including the odd braid relations. Let $U_a$ be the
standard basis of $\C_{\mu^{\mathrm R}}[\mathcal W]$. Set
\[
 V_s=q_s^{-1/2}U_s\qquad(s\text{ simple affine}).
\]
Equation~\eqref{eq:normalizer-square} gives $V_s^2=1$. The length-additive
relations and independence of the parameter product give the braid
relations for the $V_s$. They therefore define a homomorphism
$W_{\aff}\to\C_{\mu^{\mathrm R}}[\mathcal W]^\times$, with
\[
 V_v=q_v^{-1/2}U_v.
\]
The mixed relations give
$U_\omega U_vU_\omega^{-1}=U_{\omega v\omega^{-1}}$.
Since $q_{\omega v\omega^{-1}}=q_v$, the same conjugation identity holds
with $V_v$ in place of $U_v$.

Consequently $V_{\omega v}:=U_\omega V_v=d_q(\omega v)^{-1}U_{\omega v}$
satisfies
\[
 V_aV_b=\overline\mu^{\mathrm R}(a,b)V_{ab}.
\]
Comparing this product with $U_aU_b=\mu^{\mathrm R}(a,b)U_{ab}$ proves
\eqref{eq:normalizer-cochain}. The construction is on the entire group
$\mathcal W$; it does not assume that Levi translations lie in $\Omega$.
\end{proof}

\subsubsection{The actual Levi algebra in the left endpoint}
Return to the left-module convention of \S3. In the opposite of a
twisted group algebra, label the inverse of the original degree-$a$
unit by $a$. Its multiplier is the inverse of the original multiplier,
as in Lemma~\ref{lem:opposite-cocycle}. Thus put
\[
 c=(\mu^{\mathrm R})^{-1},\qquad
 \bar\mu=(\overline\mu^{\mathrm R})^{-1}.
\]
The left full normalizer-line algebra and the left AFMO endpoint have
the support presentations
\[
 \mathscr C_\tau=\C_c[\mathcal W],\qquad
 B_0=\C_{\bar\mu}[\mathcal W].
\]
The first is also
$\End_{\mathcal N}(\cind_{K}^{\mathcal N}\tau)^{\mathrm{op}}$.
Here $K\triangleleft\mathcal N$ and all its normalizer intertwining
lines are one-dimensional; convolution uses their composition and the
compact representative defects.

To fix the geometric indexing explicitly, let $A_n=T_n^{-1}$. On the
multiplicity space $E_\sigma=\Hom_K(\tau,\sigma)$ of an $L(F)$-representation,
the degree-$hK$ Levi operator is
\begin{equation}\label{eq:normalizer-left-operator}
 R_h(f)=\sigma(h)fA_h\qquad(h\in H).
\end{equation}
Indeed $A_h\tau(k)=\tau(h^{-1}kh)A_h$, so this operator preserves
$\Hom_K(\tau,\sigma)$. Its products have multiplier $c$ because
$A_mA_n=(\mu^{\mathrm R}(\bar n,\bar m))^{-1}A_{nm}$.
For $k\in K$ one has $R_k=1$. This is the convention in
Proposition~\ref{prop:Levi-support-algebra}. In the original right support indexing, inverse units
have inverse support; \eqref{eq:normalizer-left-operator} specifies the left
indexing used for the map to $W^\fs$.

\begin{proposition}[A multiplication-compatible embedding of the Levi algebra]
\label{prop:Levi-qzero-embedding}
There is a support-graded algebra isomorphism
\begin{equation}\label{eq:normalizer-qzero-isom}
 J:\mathscr C_\tau\xrightarrow{\sim}B_0,
 \qquad E_a\longmapsto d_q(a)^{-1}e_a,
\end{equation}
where $E_a,e_a$ are the respective twisted support bases. The
$\Gamma_H$-subalgebra of $\mathscr C_\tau$ is the actual Levi algebra
$D$, with operators \eqref{eq:normalizer-left-operator}. In particular $J$ gives
an embedding $j_0:D\hookrightarrow B_0$ and
\begin{equation}\label{eq:qzero-multiplier-restriction}
 [\bar\mu|_{\Gamma_H\times\Gamma_H}]=[\mu_H].
\end{equation}
The quotient grading of $B_0$ is a strong grading
\begin{equation}\label{eq:qzero-strong}
 B_0=\bigoplus_{w\in W^\fs}j_0(D)u_w
\end{equation}
with invertible homogeneous elements $u_w$.
\end{proposition}
\begin{proof}
Lemma~\ref{lem:normalizer-cocycle-normalization} gives
$c=(\delta(d_q^{-1}))\bar\mu$, so \eqref{eq:normalizer-qzero-isom}
respects multiplication and has inverse $e_a\mapsto d_q(a)E_a$.
The intertwining of $\tau$ inside $L(F)$ is exactly $H$ by Lemma~\ref{lem:depth-zero-intertwining}. Its Levi support lines are therefore exactly the
$\Gamma_H$-lines of $\mathscr C_\tau$, and their products are the same
operator products \eqref{eq:normalizer-left-operator}. This identifies their
algebra with $D$, not merely with an abstract Laurent algebra.
Restriction of the cochain identity proves
\eqref{eq:qzero-multiplier-restriction}. If different bases were originally
used to define $\mu_H$, their change of basis supplies the corresponding
additional coboundary.

Choose lifts $\widetilde w\in\mathcal W$ of $w\in W^\fs$ and put
$u_w=J(E_{\widetilde w})$. By \eqref{eq:B-group-seq}, every support degree over
$w$ is a $\Gamma_H$-translate of $\widetilde w$. The twisted group
units are invertible, giving \eqref{eq:qzero-strong}. In general
$u_wu_{w'}u_{ww'}^{-1}$ belongs to $j_0(D)^\times$; it need not be
a scalar. This is the form of grading used in
Proposition~\ref{prop:central-fiber-semisimple}.
\end{proof}

\begin{remark}[Why the normalization matters]
The raw full normalizer cocycle has
$\mu^{\mathrm R}(s,s)=q_s$, whereas the $q=1$ affine group basis has
square one at a simple reflection. Identifying these bases without the
cochain would not preserve products. Formula
\eqref{eq:normalizer-qzero-isom} gives the product comparison.
It is an identification at the $q=1$ endpoint; it does not assert that
the positive-parameter standard basis spans a twisted group subalgebra
on every set of translation supports.
\end{remark}

\subsubsection{Conjugation on the normalized Bernstein center}
Normalize the center of $D$ as in Proposition~\ref{prop:Levi-center-normalized}. Thus homogeneous
central elements $z_\lambda$, $\lambda\in\Lambda_\fs$, multiply as the
lattice, act by one on $E_\pi$, and act on $E_{\pi\otimes\nu}$ by
$\nu(\lambda)$. Use $j_0$ to regard this particular center as a subalgebra
of $B_0$.

For $n\in\mathcal N$, let $\beta_n$ be conjugation by $E_{nK}$ on $D$
inside $\mathscr C_\tau$. Scalar normalizations of this unit do not
change $\beta_n$. The corresponding automorphism in $B_0$ is conjugation
by $J(E_{nK})$. Elements of $H$ give inner automorphisms of $D$, so the
action on $Z(D)$ factors through $W^\fs$.

\begin{proposition}[The genuine Bernstein action]
\label{prop:qzero-Bernstein-action}
Under the normalized identification $Z(D)=\mathcal O(T^\fs)$, the action
on its closed points induced by \eqref{eq:qzero-strong} is
\[
 w\star t=\eta_w\,w_{\mathrm{lin}}(t).
\]
Here $n\in\mathcal N$ represents $w$,
${}^n\pi\simeq\pi\otimes\chi_n$, and $\eta_w=[\chi_n]$, as in Proposition~\ref{prop:center}. In homogeneous central coordinates, equivalently,
\begin{equation}\label{eq:qzero-Bernstein-action}
 \beta_n^{-1}(z_\lambda)
 =\eta_w(\lambda)z_{n^{-1}\lambda n}.
\end{equation}
\end{proposition}
\begin{proof}
For any representation $\sigma$ containing $\tau$, let
${}^n\sigma(l)=\sigma(n^{-1}ln)$. The map
\[
 F_n:E_\sigma\longrightarrow E_{{}^n\sigma},\qquad f\longmapsto fA_n,
\]
is an isomorphism of vector spaces. The intertwining relation for $A_n$
proves that its image belongs to $E_{{}^n\sigma}$. Computing with
\eqref{eq:normalizer-left-operator}, including the same representative defects
as in $\mathscr C_\tau$, gives
\begin{equation}\label{eq:normalizer-module-naturality}
 R_{{}^n\sigma}(b)F_n
   =F_nR_\sigma(\beta_n^{-1}(b))\qquad(b\in D).
\end{equation}
For clarity, on a degree-$h$ operator both sides act as
$f\mapsto\sigma(n^{-1}hn)fA_nA_h$; the conjugated operator before
applying $F_n$ has coefficient $A_nA_hA_n^{-1}$.
Thus if $E_\sigma$ has central character $t$, the central character of
$E_{{}^n\sigma}$ is $t\circ\beta_n^{-1}$.

Now put $\sigma=\pi\otimes\nu$. Proposition~\ref{prop:inertial-support-quotient}
gives
\[
 {}^n(\pi\otimes\nu)
 \simeq\pi\otimes\chi_n\,{}^n\nu.
\]
The normalized center sends this representation to
$\eta_w w_{\mathrm{lin}}([\nu])$. Since these unramified twists exhaust
the closed points of $Z(D)$, this proves the claimed action.
Conjugation sends the degree-$\lambda$ central line to the
degree-$n\lambda n^{-1}$ line. Evaluate the inverse-conjugation formula
at the base point $E_\pi$ to determine its scalar, obtaining
\eqref{eq:qzero-Bernstein-action}. The isomorphism $J$ transports all these
identities to $B_0$.
\end{proof}

\subsubsection{Central fibres and the source extended quotient}
We first record the central-fibre statement for a finite crossed extension
of an Azumaya algebra. The strong grading in
Proposition~\ref{prop:Levi-qzero-embedding} and the center action in
Proposition~\ref{prop:qzero-Bernstein-action} give its application to $B_0$.

\begin{proposition}[Azumaya central fibres]\label{prop:central-fiber-semisimple}
Let $A$ be a finitely generated commutative $\C$-algebra, $D$ an Azumaya
$A$-algebra, and
\[
B=\bigoplus_{w\in W}Du_w
\]
a strongly $W$-graded algebra with invertible homogeneous elements, inducing an
action of the finite group $W$ on $A=Z(D)$.  For a closed complex point $t\in\Spec A$, let $W_t$ be
its stabilizer and set
\[
B_{[t]}=B\otimes_{A^W}\C_{[t]}.
\]
If $D_t\simeq M_d(\C)$, then there is a cocycle
$\kappa_t\in Z^2(W_t,\C^\times)$, well defined up to coboundary, such that
\[
B_{[t]}/\Jac(B_{[t]})
\simeq
M_{|W/W_t|d}(\C_{\kappa_t}[W_t]).
\]
Hence
\[
\Irr(B_{[t]})\longleftrightarrow\Irr(\C_{\kappa_t}[W_t]).
\]
\end{proposition}
\begin{proof}
The commutative fibre $A\otimes_{A^W}\C_{[t]}$ is finite-dimensional.
Its nilradical is $W$-stable and nilpotent, so its product with the
specialized strongly graded algebra is a nilpotent ideal. The reduced
commutative quotient is
$\prod_{x\in Wt}\C e_x$.  In the $t$-corner the neutral algebra is
$D_t\simeq M_d(\C)$ and only degrees in $W_t$ survive.  For $w\in W_t$,
conjugation by $u_w$ is an automorphism of $D_t$, hence is inner by
Skolem--Noether; choose $a_w\in D_t^\times$ implementing it and put
$z_w=a_w^{-1}e_tu_w$.  Then $z_w$ centralizes $D_t$ and
\[
z_{w_1}z_{w_2}=\kappa_t(w_1,w_2)z_{w_1w_2}
\]
for a scalar cocycle, changed only by a coboundary if the $a_w$ are changed.
Thus the corner is $M_d(\C_{\kappa_t}[W_t])$.  The idempotent $e_t$ is full
because the homogeneous invertibles move it through the orbit, which gives the
displayed matrix algebra.  The reduced quotient is semisimple, so the
nilpotent ideal removed at the start is the Jacobson radical.
\end{proof}

\begin{lemma}[Intrinsic orbit corners]\label{lem:comparison-representatives}
In the setting of Proposition~\ref{prop:central-fiber-semisimple}, let
\[
\mathscr C_t=e_t(B_{[t]}/\Jac(B_{[t]}))e_t,
\]
where $e_t$ is the primitive idempotent of the reduced commutative orbit
algebra. Homogeneous conjugation gives coherent bijections
$\Irr(\mathscr C_t)\to\Irr(\mathscr C_{u t})$, independent of the lift and
of the element $u$ carrying $t$ to the specified point. In particular,
stabilizer elements act trivially on these simple-module classes.
\end{lemma}
\begin{proof}
A degree-$u$ unit conjugates $e_t$ to $e_{ut}$ and identifies the corners.
Two units of this degree differ by a neutral unit, giving an inner change
on the target corner. If two degrees carry $t$ to the same point, their
ratio stabilizes $t$ and its homogeneous unit restricts to a unit of the
$t$-corner. Its action there is inner. Inner automorphisms act trivially
on simple isomorphism classes. Products of the units give composition
coherence for the same reason.
\end{proof}

Thus the intrinsic target is
\[
\left(\bigsqcup_{t\in T^\fs}\Irr(\mathscr C_t)\right)/W^\fs.
\]
Choices of matrix decompositions and Skolem--Noether operators identify it
with the twisted extended quotient of
Definition~\ref{def:twisted-extended-quotient}. The induced cocycle classes
are equivariant, but the coordinate transports must accompany them.
For example, the automorphism $u_s\mapsto-u_s$ of $\C[C_2]$ fixes its
cocycle and exchanges its two simple characters. Consequently a
cohomology class alone does not define a choice-free transport of labels.

\begin{lemma}[Graded isomorphisms and effective corner cocycles]
\label{lem:graded-corner-functoriality}
For $i=0,1$, let $D_i$ be Azumaya over $\cO(T_i)$ and let
\[
 C_i=\bigoplus_{w\in W_i}D_i u_{i,w}
\]
be strongly graded algebras with invertible homogeneous units and
finite groups $W_i$. Suppose an algebra isomorphism
$F:C_0\to C_1$ carries $D_0$ to $D_1$, carries degree $w$ to
degree $\jmath(w)$ for an isomorphism $\jmath:W_0\to W_1$,
and restricts on centers to $f^*:\cO(T_0)\to\cO(T_1)$,
where $f:T_1\to T_0$ is equivariant for $\jmath^{-1}$.
Then the intrinsic orbit corners correspond, with coherent
transports, and
\begin{equation}\label{eq:graded-corner-class}
 [\kappa_{1,t}]
 =\bigl(\jmath^{-1}|_{W_{1,t}}\bigr)^*
                       [\kappa_{0,f(t)}].
\end{equation}
Compatible choices of homogeneous units and matrix identifications
make the cocycles equal under this pullback.
\end{lemma}
\begin{proof}
Here a central fiber means
$C_{i,[x]}=C_i\otimes_{\cO(T_i)^{W_i}}\C_{[x]}$.
The isomorphism carries the central fiber over $[f(t)]$ to the
fiber over $[t]$, preserves the Jacobson radical, and carries the
reduced-orbit idempotent at $f(t)$ to the one at $t$. In the first
corner choose units $a_w$ in the neutral matrix algebra implementing
conjugation by $u_{0,w}$. Then
$z_w=a_w^{-1}e_{f(t)}u_{0,w}$ centralizes that matrix algebra, and
$z_wz_v=\kappa_{0,f(t)}(w,v)z_{wv}$.
The units $F(a_w)$ and $F(u_{0,w})$ give exactly the same products
in the target corner, proving \eqref{eq:graded-corner-class}.
Changing their scalar normalizations gives coboundaries.
Finally, $F$ intertwines conjugation by all homogeneous units,
so it intertwines the corner transports as well as their
multipliers. This is stronger than a match of supports alone.
\end{proof}

\begin{theorem}[The source endpoint and its extended quotient]
\label{thm:source-support-realization}
In the stated depth-zero AFMO setup, the sequence
\eqref{eq:B-group-seq}, the embedding $j_0:D\hookrightarrow B_0$, the multiplier
comparison \eqref{eq:qzero-multiplier-restriction}, and the center action
\eqref{eq:qzero-Bernstein-action} give the source support realization. Moreover
\[
 Z(B_0)=\mathcal O(T^\fs)^{W^\fs,\star}.
\]
For every closed orbit $[t]$, Proposition~\ref{prop:central-fiber-semisimple} gives
\[
 B_{0,[t]}/\Jac(B_{0,[t]})
 \simeq M_{|W^\fs/W_t^\fs|d}(\C_{\kappa_t}[W_t^\fs]),
 \qquad d^2=[\Gamma_H:\Lambda_\fs].
\]
The orbit transports are the intrinsic corner transports of Lemma~\ref{lem:comparison-representatives}. Consequently the source endpoint has a twisted
extended-quotient parametrization
\[
 \Irr(B_0)\simeq(T^\fs//W^\fs)_\kappa
\]
after the projective-label coordinates of
Definition~\ref{def:twisted-extended-quotient} are fixed.
\end{theorem}
\begin{proof}
The preceding propositions give the support realization. The action of
$W^\fs$ on $T^\fs$ is faithful: an inertial normalizer acting trivially
on the real unramified directions centralizes $A_L$, whose centralizer
in $G$ is $L$. Passing to a finite self-twist quotient does not change
this linear faithfulness.

Write a central element of $B_0$ using \eqref{eq:qzero-strong}. Commutation
with $A_\fs$ forces every coefficient in nonidentity degree to vanish,
by faithfulness and the torsion-freeness of $D$ over the integral algebra
$A_\fs$. The remaining neutral coefficient is central in $D$, hence in
$A_\fs$, and commutation with the $u_w$ makes it invariant. This proves
the center formula. The central-fibre and orbit-corner statements follow
from the cited abstract results, now with all their application data
identified. Finiteness over the invariant center ensures that all simple
$B_0$-modules occur in these closed central fibres.
\end{proof}

Theorem~\ref{thm:source-support-realization} identifies the $q=1$ endpoint.
The comparison of its simple modules with those at positive parameters
is constructed in the next section.

\section{Geometric parameter removal}\label{subsec:deformation}

The construction in this section follows the local-to-global
parameter-removal strategy used by Solleveld in
\cite[\S8.2, Theorem~8.6 and Corollary~8.7]{SolEndo}: pass to graded
Hecke algebras on strata with fixed unitary part, compare their simple
modules with those at zero parameters, and assemble the comparisons
equivariantly. Our direct inputs are the affine-to-graded reductions
of \cite{Sol-red}, the distinguished finite Weyl types and their
Langlands-classification extension in
\cite[Theorem~6.2, Lemma~6.6 and Proposition~6.10]{Sol-survey}, and
the geometric classification and enhanced-parameter shift of
\cite{AMS-GHA}.

We implement this strategy for the full AFMO source extension. Its
length-zero group can be infinite; absorbing central translations
produces central-unit multiplication defects, and residual elements
can move the special vertex used in a Bernstein presentation. The
source lattice, its finite covers, and these multiplication defects
must therefore be retained through the local reductions. The
additional comparisons below use a finite source lattice enlargement
and an equivariant full corner to identify the effective stabilizer
cocycles, including the mixed cocycle of the enlarged action. Global
assembly then uses one coherent choice of projective-label transport
on each neutral orbit.

\subsection{Affine comparison and the source parameter family}
\label{subsec:affine-input}
Fix the depth-zero AFMO datum of \S\ref{sec:notations}. Thus
$K_L=L(F)_x$ is the full compact Levi stabilizer, the image of $x$ in
the reduced building of $L$ is a vertex, and $\tau_{\rb}$ is an
irreducible depth-zero label whose connected constituents are cuspidal.
The admissible embedding into the building of $G$ is $0$-generic at $x$.
Throughout this section assume that $G$ splits over a tamely ramified
extension and that $p\nmid |W_{\mathrm{abs}}(G)|$.

\paragraph{The source presentation.}
Use the left-module convention and fix a compatible AFMO presentation
\cite[\S3.10 and Theorem~5.3.6]{AFMO-structure}:
\begin{equation}\label{eq:source-presentation}
 \mH_G\simeq A_{\rho,1}\rtimes_{\alpha,\mu}\Omega,
 \qquad A_{\rho,1}=\mH(W_{\rho,\aff},\mathbf q).
\end{equation}
Here $\Omega$ is the source length-zero group, possibly infinite,
and $\mu\in Z^2(\Omega,\C^\times)$ is its source multiplier in this
convention. The multiplication rules are
\[
 e_\omega a=\alpha_\omega(a)e_\omega,
 \qquad e_\omega e_\eta=\mu(\omega,\eta)e_{\omega\eta}.
\]
The relevant affine simple reflections have $q_s>1$. Their normalized
AFMO generators satisfy
$\alpha_\omega(T_s)=T_{\omega s\omega^{-1}}$, and this permutation
preserves $q_s$ \cite[Lemma~3.10.6]{AFMO-structure}.
Passing from a right-module presentation means taking the opposite of
the whole algebra, including its cross relations. The affine
anti-involution identifies the affine opposite; on the group part,
Lemma~\ref{lem:opposite-cocycle} gives the inverse multiplier with
inverse support indexing. All subsequent products use the left convention.

\begin{proposition}[Unipotent realization of the affine factor]
\label{prop:ohara-affine}
There is an auxiliary connected reductive group $G_c$ splitting over
an unramified extension of $F$, with a parahoric unipotent datum,
and a based affine Coxeter isomorphism
\[
 f_c:(W_{\rho,\aff},S_\rho)\xrightarrow{\sim}(W_{c,\aff},S_c),
 \qquad q_{c,f_c(s)}=q_s.
\]
Consequently $T_w\mapsto T_{f_c(w)}$ is an isomorphism of affine
Hecke algebras. The source datum in this assertion uses the full
compact stabilizer $K_L$.
\end{proposition}
\begin{proof}
Apply \cite[Theorem~4.4.1]{Ohara25}, using the full-stabilizer option
in its \S4.2 and the connected-parahoric option on the unipotent side
in its \S4.3. The equality of relevant reflection arrangements and
parameters identifies the based Coxeter systems.
\end{proof}

Fix such a realization and based identification, denoted collectively
by $c$. Only the affine factor is used from this datum. The source
Bernstein lattice remains
$\Lambda_\fs=\Rad(b_{\mu_H})\subset\Gamma_H$; it is not replaced
by a lattice of $G_c$. Likewise the source group $\Omega$, its action,
and its multiplier remain part of the construction.

\begin{example}[A finite-index change can change stabilizers]
\label{ex:lattice-stabilizers}
In the cocharacter space of a split maximal torus of $\mathrm{GL}_2$,
consider the two lattices
\[
 X_{\rm s}=\Z e_1\oplus\Z e_2,
 \qquad X_{\rm a}=\Z(e_1-e_2)\oplus\Z(e_1+e_2).
\]
The second has index two in the first. The transposition $e_1\leftrightarrow
e_2$ preserves both. Restriction gives an equivariant degree-two isogeny
\[
 r:\Hom(X_{\rm s},\C^\times)\longrightarrow\Hom(X_{\rm a},\C^\times),
 \qquad (x_1,x_2)\longmapsto(a,b)=(x_1/x_2,x_1x_2).
\]
The action swaps $x_1,x_2$ on the source and sends $(a,b)$ to
$(a^{-1},b)$ on the target. Nevertheless,
\[
 \Stab_{S_2}(1,-1)=1,
 \qquad \Stab_{S_2}\bigl(r(1,-1)\bigr)
 =\Stab_{S_2}(-1,-1)=S_2.
\]
Indeed, $(1,-1)$ and $(-1,1)$ differ by the nontrivial kernel element
$(-1,-1)$. Thus an equivariant isogeny need not identify stabilizers or
their extended-quotient fibres. A rational apartment identification
does not supply an identification of the source Bernstein torus with
an auxiliary torus.
\end{example}

\begin{proposition}[Source transport through the parameter family]
\label{prop:source-family}
For $0\leq\epsilon\leq1$, put
\[
 A_{\rho,\epsilon}=\mH(W_{\rho,\aff},\mathbf q^\epsilon),
 \qquad A_{c,\epsilon}=\mH(W_{c,\aff},\mathbf q_c^\epsilon).
\]
The based source action defines an action $\alpha^\epsilon$ on
$A_{\rho,\epsilon}$. There are based isomorphisms
\[
 \Psi_{c,\epsilon}:A_{\rho,\epsilon}\xrightarrow{\sim}A_{c,\epsilon},
 \qquad T_w\longmapsto T_{f_c(w)}.
\]
Set $\alpha^{c,\epsilon}_\omega
=\Psi_{c,\epsilon}\alpha^\epsilon_\omega\Psi_{c,\epsilon}^{-1}$ and
\begin{equation}\label{eq:source-retained}
 B^{\src}_{c,\epsilon}
 =A_{c,\epsilon}\rtimes_{\alpha^{c,\epsilon},\mu}\Omega.
\end{equation}
Then
\begin{equation}\label{eq:transport-family}
 \widehat\Psi_{c,\epsilon}(a e_\omega)
 =\Psi_{c,\epsilon}(a)e_\omega
\end{equation}
is an algebra isomorphism from
$A_{\rho,\epsilon}\rtimes_{\alpha^\epsilon,\mu}\Omega$ to
$B^{\src}_{c,\epsilon}$. In particular,
$\mH_G\simeq B^{\src}_{c,1}$ and $B_0\simeq B^{\src}_{c,0}$.
For any fixed collection of choices $c,d,e$, the transitions
\[
 \Theta_{d,c,\epsilon}=\Psi_{d,\epsilon}\Psi_{c,\epsilon}^{-1}
\]
extend by the identity on the source group basis and satisfy
\[
 \Theta_{e,d,\epsilon}\Theta_{d,c,\epsilon}
 =\Theta_{e,c,\epsilon},\qquad
 \Theta_{c,c,\epsilon}=\mathrm{id}.
\]
\end{proposition}
\begin{proof}
The source action preserves the braid relations and the relations
$(T_s+1)(T_s-q_s^\epsilon)=0$, because it permutes the simple
reflections and their parameters. The based map $f_c$ preserves the
same presentations, proving the assertions about $\alpha^\epsilon$
and $\Psi_{c,\epsilon}$, including their inverses.
For homogeneous terms the source product is
\[
 (a e_\omega)(b e_\eta)
 =a\alpha^\epsilon_\omega(b)\mu(\omega,\eta)e_{\omega\eta}.
\]
Applying $\Psi_{c,\epsilon}$ to its affine coefficient gives exactly
the product in $B^{\src}_{c,\epsilon}$. This proves
\eqref{eq:transport-family}; the inverse uses $\Psi_{c,\epsilon}^{-1}$.
At $\epsilon=0$ the affine factor is $\C[W_{\rho,\aff}]$, with the
same source group action and multiplier, hence the endpoint is the
specified support algebra $B_0$.

The transition map satisfies
$\Theta_{d,c,\epsilon}\alpha^{c,\epsilon}_\omega
=\alpha^{d,\epsilon}_\omega\Theta_{d,c,\epsilon}$.
It therefore extends to the crossed products while fixing $e_\omega$.
The displayed composition identities follow by cancellation of the
$\Psi$ maps. This proves coherence for the chosen comparisons.
\end{proof}

The family $B_\epsilon$ provides compatible algebra presentations.
The bijection between its endpoints is constructed below using the
local geometric shift and effective Clifford transport; it does not
use a continuity argument in $\epsilon$ or continuous choices of
comparison cochains.

\begin{remark}[Affine normalization under a change of finite pinning]
\label{rem:pinned-affine}
Fix the source datum and index the relevant walls by its relative
geometry. For the finite Jordan construction cited in Remark~\ref{rem:jd-llc},
use its standing inputs \textup{(F), (C), (R)} on the actual finite
stabilizer extensions and rational pinned-component lifts, the connected
Jordan and cuspidal matching theorems in their stated scope, and the
regular-embedding and rational descent data. Fix an admissible rational
pinning $\cP$ and the global normalization
\begin{equation}\label{eq:pinned-global-data}
 \mathcal D_{\rm tot}
  =(\mathcal D_{\rm conn},\mathcal D_{\rm HL}^{\star},
                                      \mathcal D_{\rm ext}^{\rm all}).
\end{equation}
These entries specify the connected normalizations, the coherent
Howlett--Lehrer system and the relative extension classes, respectively.
For a connected constituent of the pinned Jordan label to which
the finite affine comparison applies, the
normalized maps of Proposition~\ref{prop:source-family} are determined
by the source-indexed simple support lines for $0<\epsilon\leq 1$.
Indeed, if $T$ and $aT$ satisfy the Hecke quadratic relations with
parameters $q,q'>1$ in such a line, comparison of the coefficients of
$T$ and $1$ gives $a(q-1)=q'-1$ and $a^2q=q'$. The two possible
values are $a=1$ and $a=-1/q$. The latter gives $q'=1/q<1$.
Taking $q=q_s^\epsilon$ proves uniqueness. At $\epsilon=0$ we use
the specialization $T_w\mapsto T_{f_c(w)}$ of the same basis formula.

The compatibility of the pinned Jordan labels with the ordinary
connected Jordan orbit gives these same transitions
for its permitted changes of unipotent realization. Equation
\eqref{eq:transport-family} extends them using the original source
units, action and multiplier. The passage between the spectra at
$\epsilon=1$ and $0$ uses the connected shift and the residual
Clifford comparison. For a further change of presentation its action
on that shift must also be checked, with the distinguished finite
Weyl type, real Langlands datum, root and sign conventions of
\S\ref{subsec:source-equivariance}. The remaining choice on each
multiplicity space is described in \S\ref{sec:pinned-neutral}.
\end{remark}

\paragraph{Source multiplication lines.}
The following elementary observation records the quotient compatibility
needed for the source multiplier, including the representative defects.

\begin{lemma}[Inflation and multiplication of intertwining lines]
\label{lem:inflation-lines}
Let $K^+\subset K\subset N$, with both $K^+$ and $K$ normal in $N$
and $K/K^+$ finite, and
let $\rho$ be an irreducible complex representation of $K$ inflated
from $\bar\rho$ on $K/K^+$. Assume that $N$ preserves the isomorphism
class of $\rho$. Put $\Delta=N/K$, choose representatives $n_g$ with
$n_1=1$, and write
\[
 n_gn_h=d(g,h)n_{gh},\qquad d(g,h)\in K.
\]
For each $g\in\Delta$, inflation identifies the one-dimensional spaces
\[
 \mathscr L_g=\{T\in\End_{\C}(V_\rho):
 T\rho(k)=\rho(n_gkn_g^{-1})T\quad(k\in K)\}
\]
with the corresponding spaces for $\bar\rho$ and conjugation by
$\bar n_g\in N/K^+$. This identification preserves the corrected
products and their scalar cocycles, using $\overline{d(g,h)}$ in the
quotient. The indexing group remains $\Delta$.
\end{lemma}
\begin{proof}
Both representations act on the same vector space, and the equations
defining the two spaces are identical since $\rho$ is trivial on $K^+$.
Schur's lemma gives dimension one. Choose $T_g\ne0$ in each space,
with $T_1=1$. There are unique scalars $c(g,h)$ such that
\begin{equation}\label{eq:inflation-line-product}
 T_gT_h=c(g,h)\rho(d(g,h))T_{gh}.
\end{equation}
Indeed, the two sides without the scalar implement the same conjugation
of $K$. The group identity
\[
 d(g,h)d(gh,k)
 =n_gd(h,k)n_g^{-1}d(g,hk)
\]
and associativity of operator composition imply
$c(g,h)c(gh,k)=c(h,k)c(g,hk)$.
Moreover,
$\rho(d(g,h))=\bar\rho(\overline{d(g,h)})$,
so \eqref{eq:inflation-line-product} is literally the same equation
after inflation. Replacing $T_g$ by $a(g)T_g$ changes both cocycles
by the same coboundary. Taking opposite conventions on both sides
preserves this compatibility as well.
\end{proof}

Apply the lemma with the full compact Levi stabilizer and its
pro-unipotent radical, and with the source normalizer preserving the
label. It identifies the multiplication lines used to compute the
source multiplier in a fixed convention. It does not replace the
source indexing group by its finite image in an automorphism group.
The central translation absorption in Proposition~\ref{prop:local-finite}
and the local transports in Proposition~\ref{prop:local-local-source}
keep these products. The endpoint is identified in
Theorem~\ref{thm:source-support-realization}.

\subsection{Clifford transport and multiplication defects}
\begin{lemma}[Clifford transport for a finite action]\label{lem:Clifford-transport}
Let a finite group $Q$ act on complex algebras $\mathscr A_r$ and
$\mathscr A_0$, and let $\omega\in Z^2(Q,\C^\times)$ be the external
crossed-product multiplier. Consider their finite-dimensional simple
modules. Suppose $\Theta:\Irr(\mathscr A_r)\to\Irr(\mathscr A_0)$ is an
equivariant bijection. For a simple $V$ with stabilizer $Q_V$, let
$c_V$ be the multiplier of intertwining operators implementing this action.
If
\begin{equation}\label{eq:intrinsic-Clifford-target}
[c_V]=[c_{\Theta(V)}]\quad\text{in }H^2(Q_V,\C^\times)
\end{equation}
for every $V$, then choices of cochains realizing these equalities, one
for each $Q$-orbit, give a bijection
\[
\Irr(\mathscr A_r\rtimes_\omega Q)
\longrightarrow\Irr(\mathscr A_0\rtimes_\omega Q).
\]
\end{lemma}
\begin{proof}
Choose operators $T_g$ on $V$ implementing the stabilizer action, with
$T_gT_h=c_V(g,h)T_{gh}$. On a multiplicity space $E$, the operators
$T_g\otimes S_g$ satisfy the crossed-product relation precisely when
\[
S_gS_h=\omega(g,h)c_V(g,h)^{-1}S_{gh}.
\]
Clifford theory parametrizes simples over the orbit of $V$ by the
irreducible projective representations of $Q_V$ with this multiplier.
Equation~\eqref{eq:intrinsic-Clifford-target}, with a chosen cochain,
identifies these projective representation sets for $V$ and $\Theta(V)$.
Induction from the inertia crossed product and passage over the orbits
prove the assertion.
\end{proof}

For localized crossed products the chosen residual lifts can have
algebra-valued multiplication defects. The following version includes
those defects in the stabilizer invariant.

\begin{lemma}[Clifford transport for coherent weak actions]
\label{lem:weak-Clifford}
Let $Q$ be finite, and for $i\in\{1,0\}$ let
\[
 \mathscr B_i=\bigoplus_{g\in Q}\mathscr A_i u_{i,g}
\]
be crossed products of complex algebras with invertible homogeneous
elements and normalized data
\begin{align*}
 u_{i,g}a u_{i,g}^{-1}&=\beta_{i,g}(a),&
 u_{i,g}u_{i,h}&=a_i(g,h)u_{i,gh},\\
 \beta_{i,g}\beta_{i,h}&=\Ad(a_i(g,h))\beta_{i,gh},&
 a_i(g,h)a_i(gh,k)&=\beta_{i,g}(a_i(h,k))a_i(g,hk).
\end{align*}
Consider finite-dimensional simple modules. Suppose
$\Theta:\Irr(\mathscr A_1)\to\Irr(\mathscr A_0)$ is $Q$-equivariant.
For a simple $V$ with representation $\pi_V$ and stabilizer $Q_V$,
choose operators $T_{i,g}$ implementing $\beta_{i,g}$ on $V$, and define
the scalar multiplier $\gamma_{i,V}$ by
\begin{equation}\label{eq:weak-obstruction}
 T_{i,g}T_{i,h}
 =\gamma_{i,V}(g,h)\,\pi_V(a_i(g,h))T_{i,gh}.
\end{equation}
Here $V$ is a module for $\mathscr A_i$ in the corresponding formula.
Then $\gamma_{i,V}$ is a scalar cocycle, and the multiplicity-space
multiplier for simples of $\mathscr B_i$ above $V$ is
\begin{equation}\label{eq:effective-kappa}
 \kappa^{\mathrm{eff}}_{i,V}=\gamma_{i,V}^{-1}.
\end{equation}
If
\begin{equation}\label{eq:effective-comparison}
 [\kappa^{\mathrm{eff}}_{1,V}]
 =[\kappa^{\mathrm{eff}}_{0,\Theta(V)}]
 \quad\text{in }H^2(Q_V,\C^\times)
\end{equation}
for every $V$, then choices of comparison cochains on orbit
representatives give a bijection
$\Irr(\mathscr B_1)\to\Irr(\mathscr B_0)$.
\end{lemma}

\begin{proof}
The operators on both sides of \eqref{eq:weak-obstruction} implement the
same automorphism of $\mathscr A_i$ on $V$. Schur's lemma gives the
scalar $\gamma_{i,V}(g,h)$. Associativity of the $T$ products and the
displayed defect identity imply its scalar cocycle identity; changing
$T_{i,g}$ changes it by a coboundary. For an inertia module on $V\otimes E$,
let $u_{i,g}$ act as $T_{i,g}\otimes S_g$. The product relation becomes
\[
 S_gS_h=\gamma_{i,V}(g,h)^{-1}S_{gh}.
\]
Finite crossed-product Clifford theory therefore describes simples
above the orbit of $V$ by
$\Irr(\C_{\kappa^{\mathrm{eff}}_{i,V}}[Q_V])$. This applies to
finite-dimensional simples: their restriction to the neutral algebra
is semisimple, since the sum of the conjugates of a simple neutral
submodule is a nonzero submodule of the whole simple crossed-product
module and hence is the whole module.

Choose $b_V:Q_V\to\C^\times$ such that
$\kappa^{\mathrm{eff}}_{1,V}=(\delta b_V)
\kappa^{\mathrm{eff}}_{0,\Theta(V)}$. Rescaling the operators on $E$
by $b_V(g)^{-1}$ identifies the two projective representation sets.
Induction from the inertia crossed products, followed by passage over
the neutral $Q$-orbits, proves the bijection.
\end{proof}

For a strict action with
$a_i(g,h)=\omega_i(g,h)1$, write
$T_{i,g}T_{i,h}=c_{i,V}(g,h)T_{i,gh}$. Then
\begin{equation}\label{eq:strict-effective}
 \kappa^{\mathrm{eff}}_{i,V}=\omega_i c_{i,V}^{-1}.
\end{equation}
When the two $\omega_i$ are identified, condition
\eqref{eq:effective-comparison} is precisely \eqref{eq:intrinsic-Clifford-target}. Without that identification, the
effective condition \eqref{eq:effective-comparison} is the appropriate
one. Matching only the intrinsic $c$ classes is then insufficient.

\paragraph{Residual products after lattice absorption.}
The multiplier $\mu$ in Proposition~\ref{prop:source-family} is fixed
before localization. A residual group can arise only after absorbing
translations. If representatives of this quotient multiply with lattice
defect $\ell(g,h)$, evaluation of that defect at a neutral central character
contributes to the residual product. It is not simply the restriction of
$\mu$ to a subgroup. The character can move under parameter removal.

\begin{example}[A scalar defect created by lattice absorption]
The group algebra $\C[\Z]=\C[u,u^{-1}]$ has trivial original multiplier.
View it as a crossed extension of $\C[2\Z]=\C[z,z^{-1}]$, with $z=u^2$,
by $Q=C_2$. A lift of the nontrivial element has square $z$, so the
specialization $z=t\in\C^\times$ gives a residual multiplier with value
$t$ on that pair. This multiplier varies although the original one was
trivial. Its class is pointwise trivial because $H^2(C_2,\C^\times)=0$,
but choosing its trivialization requires a square root of $t$; there is
no single-valued continuous square root on all of $\C^\times$.
This example does not assert a nonzero obstruction class. It shows why
product identifications and compatible projective-label transports
are needed in a parameter-removal construction.
\end{example}

\subsection{Source local models and connected geometry}\label{subsec:source-local-models}
\subsubsection{Absorbing the central translations}
Write the source presentation from Proposition~\ref{prop:source-family} as
\[
 B_\epsilon=\mH(W_{\rm aff},q^\epsilon)
                    \rtimes_{\alpha,\mu}\Omega,\qquad 0\leq\epsilon\leq1.
\]
Here $\mu$ is the \emph{left} multiplier: the inverse of the normalized
right multiplier of AFMO on $\Omega$. Thus $e_\omega e_{\omega'}=
\mu(\omega,\omega')e_{\omega\omega'}$ and
$e_\omega T_s e_\omega^{-1}=T_{\omega s\omega^{-1}}$.
At $\epsilon=0$, the group cocycle is inflated from $\Omega$.
Use the source groups
\[
 1\longrightarrow\Gamma\longrightarrow\mathcal W
   \longrightarrow W^\fs\longrightarrow1,\qquad
 \Gamma=\Gamma_H,\qquad \Lambda=\Lambda_\fs=\operatorname{Rad}(b_{\mu_H}).
\]
The lattice $\Lambda$ has finite index in $\Gamma$ and carries the
normalized central basis $z_\lambda$ of \S\ref{subsec:Levi-center-normalization}. We use the same
notation for its image under $j_0$ in $B_0$.

Choose a special vertex for $W_{\rm aff}$ and a finite Weyl complement,
so that
\[
 W_{\rm aff}=X_{\rm aff}\rtimes W_f.
\]
Here $X_{\rm aff}$ denotes its translation subgroup, not the whole
translation group of $\mathcal W$. Put
\begin{equation}\label{eq:local-lattices}
 \Gamma_0=\Gamma\cap\Omega,\qquad
 \Lambda_c=\Lambda\cap\Omega,\qquad
 X=X_{\rm aff}\oplus\Lambda_c,\qquad Q=\Omega/\Lambda_c.
\end{equation}

\begin{proposition}[The source translation lattices]\label{prop:local-lattices}
The inclusions $X_{\rm aff}\subset\Lambda$ and $X\subset\Lambda$ hold.
The subgroup $\Lambda_c$ is normal in $\mathcal W$ and centralizes
$W_{\rm aff}$. Both $[\Lambda:X]$ and $|Q|$ are finite. More precisely,
\[
 1\longrightarrow\Gamma_0/\Lambda_c\longrightarrow Q
 \longrightarrow\operatorname{im}(\Omega\longrightarrow W^\fs)
 \longrightarrow1.
\]
These assertions include the empty-root case, when $X_{\rm aff}=0$,
$X=\Lambda$ and $Q=\mathcal W/\Lambda$.
\end{proposition}

\begin{proof}
Let $V_c$ be the common kernel of the gradients of the relevant affine
walls. AFMO's reflection realization acts faithfully on the quotient
apartment by $V_c$, and the affine reflections fix the directions $V_c$
\cite[Axiom~3.7.1, Proposition~3.7.4]{AFMO-structure}. Hence a translation in
$W_{\rm aff}$ has trivial linear part on the full apartment. By Theorem~\ref{thm:source-support-realization},
the linear action of $W^\fs$ is faithful, so this translation belongs to
$\Gamma$.

For $x\in X_{\rm aff}$ and $\gamma\in\Gamma$, the inflated cocycle of
$B_0$ has value $1$ on both $(x,\gamma)$ and $(\gamma,x)$. The underlying
group elements commute. The embedding $j_0$ therefore shows that $x$
lies in the radical of the Levi commutator, namely $\Lambda$.

An element of $\Gamma_0$ is a translation preserving the chosen alcove
in the quotient apartment. That alcove is bounded, so its translation
there is zero. Its vector lies in $V_c$ and it commutes with every
relevant affine reflection on the full apartment. The reflection
realization is faithful on $W_{\rm aff}$; thus $\Gamma_0$ centralizes
$W_{\rm aff}$. The radical $\Lambda$ is $\mathcal W$-stable, and
$\Omega$ normalizes $\Gamma_0$, so $\Lambda_c$ is normal in $\Omega$.
Together with centralization this proves normality in $\mathcal W$.

The group $\Gamma_0/\Lambda_c$ injects into the finite group
$\Gamma/\Lambda$, and $\Omega/\Gamma_0$ injects into $W^\fs$. Thus $Q$
is finite. The intersection $X_{\rm aff}\cap\Lambda_c$ is trivial
because $W_{\rm aff}\cap\Omega=1$. Finally $X$ is normal in
$\mathcal W$, and $\mathcal W/X$ has order $|W_f||Q|$.
Since $X\subset\Lambda\subset\mathcal W$, the claimed index is finite.
\end{proof}

For $c\in\Lambda_c$, write
\[
 z_c=b(c)e_c,\qquad b(c)\in\C^\times.
\]
The normalization by $E_\pi$ fixes $b(c)$, and
$b(c)b(d)\mu(c,d)=b(cd)$. We use multiplicative notation for
$\Lambda_c$ inside $\Omega$, and additive notation when it is viewed
as a character lattice. These same $z_c$ define a copy of
$\C[\Lambda_c]$ in every $B_\epsilon$.

\begin{proposition}[A finite crossed presentation]\label{prop:local-finite}
Set
\begin{equation}\label{eq:local-neutral}
 A_\epsilon=\C[\Lambda_c]\otimes
                        \mH(W_{\rm aff},q^\epsilon)\subset B_\epsilon.
\end{equation}
Choose a normalized section $\sigma:Q\to\Omega$ and put
$u_g=e_{\sigma(g)}$. Then
\[
 B_\epsilon=\bigoplus_{g\in Q}A_\epsilon u_g,
 \qquad u_g a=\beta_g(a)u_g,
 \qquad u_gu_h=a(g,h)u_{gh}.
\]
The action and defects are
\begin{align}
 \beta_g(T_s)&=T_{\sigma(g)s\sigma(g)^{-1}},\label{eq:local-beta}\\
 \beta_g(z_c)&=\chi_{\sigma(g)}(c)z_{\sigma(g)c\sigma(g)^{-1}},\notag\\
 \chi_\omega(c)&=
 \frac{b(c)\mu(\omega,c)}
 {b(\omega c\omega^{-1})\mu(\omega c\omega^{-1},\omega)},\notag\\
 c(g,h)&=\sigma(g)\sigma(h)\sigma(gh)^{-1}\in\Lambda_c,\notag\\
 a(g,h)&=
 \frac{\mu(\sigma(g),\sigma(h))}
 {\mu(c(g,h),\sigma(gh))b(c(g,h))}\,z_{c(g,h)}.
 \label{eq:local-defect}
\end{align}
In particular $a(g,h)\in\C[\Lambda_c]^\times\subset Z(A_\epsilon)^\times$,
$\beta_g\beta_h=\beta_{gh}$, and
\[
 a(g,h)a(gh,k)=\beta_g(a(h,k))a(g,hk).
\]
The formulas are independent of $\epsilon$ in the standard source
presentations.
\end{proposition}

\begin{proof}
The preceding proposition and the source action show that the $z_c$
commute with the affine factor. Disjoint support lines give the tensor
embedding in \eqref{eq:local-neutral} and the direct-sum decomposition.
The identity $e_\omega e_c=\mu(\omega,c)e_{\omega c}$ gives the
conjugation formula after comparing with $e_{\omega c\omega^{-1}}e_\omega$.
Similarly, compare $e_{\sigma(g)}e_{\sigma(h)}$ with
$e_{c(g,h)}e_{\sigma(gh)}$ to obtain \eqref{eq:local-defect}.
These products involve only length-zero elements. The defects are
central in $A_\epsilon$, so conjugation gives a genuine $Q$-action;
associativity gives the displayed identity.
\end{proof}

This $Q$ is explicit and finite even when $\Omega$ is infinite. Its
possibly ineffective subgroups must be retained. The extension is a
crossed product with \emph{central-unit} defects, not necessarily a
crossed product with a fixed scalar cocycle. On a neutral simple module
with central character $z$, the external scalar factor is $a(g,h)(z)$
for $g,h$ in its stabilizer.

\subsubsection{The full neutral Bernstein presentation}
Choose the reduced root system $R\subset X_{\rm aff}\otimes\R$ whose
root lattice is $X_{\rm aff}$ and whose Weyl group is the chosen $W_f$.
Thus the affine group in this Bernstein convention is
$\mathbb ZR\rtimes W_f$; relative to the apartment root system this
uses the dual convention. Extend every coroot by zero on $\Lambda_c$.
The full root datum and torus are
\begin{equation}\label{eq:local-rootdatum}
 \mathcal R_X=(X,R,Y,R^\vee,\Delta),\qquad
 Y=\Hom(X,\mathbb Z),\qquad T_X=\Hom(X,\C^\times).
\end{equation}
This is the root-lattice semisimple datum together with the actual
absorbed source toral lattice.

Here are explicit parameter conventions. Let $Q_\alpha$ and
$Q_\alpha^*$ be the positive Bernstein labels for the wall through the
chosen special vertex and the other parallel reflection class. They
are obtained from the affine simple labels by the root-parameter
extension; they are \emph{not} values of the length-multiplicative
function on arbitrary nonsimple reflections. Set $Q_\alpha^*=Q_\alpha$
when $\alpha^\vee\notin2Y$. These labels are $W_f$-invariant on $R$.

Write $N_{s,\epsilon}=q_s^{-\epsilon/2}T_s$. For dominant
$x\in X_{\rm aff}$ set $\vartheta_{x,\epsilon}=N_{t_x,\epsilon}$,
extend multiplicatively to $X_{\rm aff}$, and put
$\vartheta_{x+c,\epsilon}=\vartheta_{x,\epsilon}z_c$.
The Bernstein presentation \cite[Theorem~1.2.1]{Sol-red} gives
\[
 A_\epsilon=\bigoplus_{w\in W_f}\C[X]N_{w,\epsilon},\qquad
 Z(A_\epsilon)=\C[X]^{W_f}.
\]
The finite generators satisfy their braid relations and
$(N_{s,\epsilon}-q_s^{\epsilon/2})
 (N_{s,\epsilon}+q_s^{-\epsilon/2})=0$.
For $\alpha\in\Delta$, $s=s_\alpha$ and $f\in\C[X]$, the cross relation is
\begin{equation}\label{eq:local-bernstein}
 \begin{aligned}
 &fN_{s,\epsilon}-N_{s,\epsilon}s(f)\\
 &\quad=\begin{cases}
 (Q_\alpha^{\epsilon/2}-Q_\alpha^{-\epsilon/2})
 \dfrac{f-s(f)}{1-\vartheta_{-\alpha,\epsilon}},
       &\alpha^\vee\notin2Y,\\[6pt]
 \left(Q_\alpha^{\epsilon/2}-Q_\alpha^{-\epsilon/2}
 +((Q_\alpha^*)^{\epsilon/2}-(Q_\alpha^*)^{-\epsilon/2})
                 \vartheta_{-\alpha,\epsilon}\right)
 \dfrac{f-s(f)}{1-\vartheta_{-2\alpha,\epsilon}},
       &\alpha^\vee\in2Y.
 \end{cases}
\end{aligned}
\end{equation}
The quotients in this formula combine to elements of $\C[X]$.
In particular \eqref{eq:local-rootdatum} and
\eqref{eq:local-bernstein} specify the whole neutral affine algebra,
including its toral directions and both possible rank-one labels.

\paragraph{Comparison with the normalized source torus.}
At $\epsilon=0$ put $E_{x+c}=e_xz_c$ for
$x\in X_{\rm aff}$ and $c\in\Lambda_c$. These elements are a
multiplicative central family in $j_0(D)$. Let
$\xi_\pi\in\Hom(X,\C^\times)$ be their scalar character on $E_\pi$;
then $\xi_\pi|_{\Lambda_c}=1$ and
\begin{equation}\label{eq:local-origin}
 z_x=\xi_\pi(x)^{-1}E_x\quad(x\in X).
\end{equation}
Consequently the map from the actual Bernstein torsor in its fixed
normalization to the standard neutral torus is
\begin{equation}\label{eq:local-isogeny}
 r_\pi:T^\fs\longrightarrow T_X,\qquad
 r_\pi(t)(x)=\xi_\pi(x)t(x),\qquad
 K_X=\Hom(\Lambda/X,\C^\times).
\end{equation}
It is the restriction isogeny followed by translation by $\xi_\pi$;
its fibres are $K_X$-torsors. The self-twist quotient is already present
in $\Lambda=X^*(T^\fs)$. Neither $X$ nor an auxiliary realization
lattice can be substituted for $\Lambda$ without this correction.

Equation~\eqref{eq:local-isogeny} is an endpoint coordinate comparison.
It does not assert that the homogeneous Levi elements $z_x$ for
$x\in X_{\rm aff}$ are positive-parameter Bernstein generators. The
positive-parameter coordinates are the specified
$\vartheta_{x,\epsilon}$.

On the neutral center the action \eqref{eq:local-beta} is independent of
$\epsilon$ under the Bernstein identifications. Indeed, if $D_\omega$
is the linear part of $\omega$, it sends a Bernstein presentation for
a chamber $P$ to that for $D_\omega P$ and sends $x$ to $D_\omega x$.
The $W_f$-invariant Bernstein elements are independent of that chamber.
Combine this action on the semisimple invariant functions with the
explicit formula for $z_c$. In particular, set
\[
 Z=\C[X]^{W_f},\qquad S=Z^Q\subset Z(B_\epsilon).
\]
Every $B_\epsilon$ is finite over $S$. No assertion that $S$ is its full
center is needed: the finite crossed extension retains the further
central and projective data. At the endpoint its full center is the
one established in Theorem~\ref{thm:source-support-realization}.

\paragraph{Index-two check.}
For the untwisted $\mathrm{GL}_2$ Iwahori example,
\[
 \Lambda=\mathbb Z^2,\quad
 X_{\rm aff}=\mathbb Z(e_1-e_2),\quad
 \Lambda_c=\mathbb Z(e_1+e_2),\quad Q=C_2.
\]
With base character $\xi_\pi=1$, the map is
$(t_1,t_2)\mapsto(t_1/t_2,t_1t_2)$, of degree two.
A standard generator $\omega$ of $\Omega\simeq\mathbb Z$ satisfies
$\omega^2=t_{e_1+e_2}$. Thus its residual lift satisfies
$u^2=z_{e_1+e_2}$, and conjugates the two affine simple generators.
Discarding this residual extension would discard the missing finite
information. In the torus case $W_{\rm aff}=1$, the same construction
retains $\C[\Lambda]$ and the finite twisted extension by
$\Gamma/\Lambda$, which accounts for the Levi Azumaya fibres.

\subsubsection{Local algebraic graded reductions at both endpoints}
Let
\[
 \mathfrak t=\Hom(X,\C),\quad \mathfrak a=\Hom(X,\R),\quad
 T_{X,\mathrm{un}}=\Hom(X,S^1).
\]
For $u\in T_{X,\mathrm{un}}$ define
\begin{align}
 R_u&=\{\alpha\in R:s_\alpha u=u\},\label{eq:local-unitary-roots}\\
 W_u&=\operatorname{Stab}_{W_f}(u),\qquad
 \Delta_u=\text{the basis of $R_u$ determined by $R^+$},\notag\\
 \Xi_u&=\{w\in W_u:w(\Delta_u)=\Delta_u\}.
 \notag
\end{align}
Here $\Xi_u$ is the finite stabilizer complement; it is distinct
from the source translation lattice $\Gamma$. Then $W_u=W(R_u)\rtimes\Xi_u$.
The reflection test is lattice sensitive:
\[
 s_\alpha u=u\quad\Longleftrightarrow\quad
 \begin{cases}
 \alpha(u)=1,&\alpha^\vee\notin2Y,\\
 \alpha(u)\in\{1,-1\},&\alpha^\vee\in2Y.
 \end{cases}
\]
In particular the condition is not universally $\alpha(u)=1$.
Define
\begin{equation}\label{eq:local-graded-parameters}
 k_u(\alpha)=\tfrac12\bigl(\log Q_\alpha+
                 \alpha(u)\log Q_\alpha^*\bigr),\qquad
 k_{\epsilon,u}=\epsilon k_u.
\end{equation}
All logarithms of positive real parameters here are real logarithms.
Write
\[
 G_{\epsilon,u}=
 \mathbb H(\mathfrak t,R_u,\epsilon k_u)
                         \rtimes\Xi_u.
\]
Our graded convention is
$f s_\alpha-s_\alpha s_\alpha(f)
 =\epsilon k_u(\alpha)(f-s_\alpha(f))/\alpha$ for
$f\in\mathcal O(\mathfrak t)$.
The full vector space $\mathfrak t$, including the directions dual to
$\Lambda_c$ and all other directions on which $R_u$ vanishes, is retained.

\begin{proposition}[Neutral algebraic local models]\label{prop:local-reduction}
For each $u\in T_{X,\mathrm{un}}$ and $0\leq\epsilon\leq1$, there is
an equivalence
\begin{equation}\label{eq:local-unitary-equivalence}
 \operatorname{Mod}_{f,W_fu\exp(\mathfrak a)}(A_\epsilon)
 \simeq \operatorname{Mod}_{f,\mathfrak a}(G_{\epsilon,u}).
\end{equation}
Subscripts specify all generalized weights for the respective
commutative polynomial or Laurent subalgebras. More explicitly, for
$t=u\exp(\nu)$, $\nu\in\mathfrak a$, put $z=W_ft\in\operatorname{Spec}Z$.
There is a choice of algebra isomorphism
\begin{equation}\label{eq:local-completion}
 \Psi_{\epsilon,z}:\widehat A_{\epsilon,z}
 \xrightarrow{\sim}
 M_{n_u}\left(\widehat G_{\epsilon,u,W_u\nu}\right),
 \qquad n_u=[W_f:W_u].
\end{equation}
Completions are taken at the indicated central characters. On their
finite-dimensional continuous modules these give the corresponding
central-character equivalences. The models cover all simple modules
of $A_\epsilon$ at either endpoint.
\end{proposition}

\begin{proof}
Apply the two reduction theorems to the \emph{actual} datum
\eqref{eq:local-rootdatum}, with no external diagram group.
The first reduction retains the stabilizer complement
$\Xi_u$ and a matrix factor of size $[W_f:W_u]$.
The second uses $\exp_u(\lambda)=u\exp\lambda$ and sends normalized
rational intertwiners to graded intertwiners. Their analytic restrictions
and completions give \eqref{eq:local-completion}; the parameter formula
is \eqref{eq:local-graded-parameters}
\cite[Theorems~2.1.2, 2.1.4 and Corollary~2.1.5]{Sol-red}.
Suitable logarithmic neighborhoods exist around every real $\nu$:
shrink their imaginary directions to avoid nonzero multiples of $\pi i$
in the root evaluations and their real parameter translates.
The assertions include $q^0=1$. Finally every character of $X$ has a
unique polar decomposition $u\exp\nu$, proving coverage.
\end{proof}

The parameter $k_u$ may vanish, and a difference of labels at
$\alpha(u)=-1$ may be negative. A sign normalization therefore
includes the corresponding transport of the source action and
coordinates. The complement $\Xi_u$ can also be nontrivial even
when the neutral affine algebra has no external finite group.

\paragraph{Zero parameters belong to a finite residual factor.}
There is a useful further algebraic refinement before asking for
connected geometric data. Define, using the parameters at $\epsilon=1$,
\[
 R_u^\times=\{\alpha\in R_u:k_u(\alpha)\ne0\},\qquad
 F_u=\{w\in W_u:w((R_u^\times)^+)=(R_u^\times)^+\},
\]
where $(R_u^\times)^+=R_u^\times\cap R^+$.
These groups and root systems are held fixed as $\epsilon$ varies.

\begin{proposition}[Separating the identically zero parameters]
\label{prop:local-zero}
The set $R_u^\times$ is a $W_u$-stable root subsystem, and
\[
 W_u=W(R_u^\times)\rtimes F_u,\qquad \Xi_u\subset F_u.
\]
For all $0\leq\epsilon\leq1$, the polynomial subalgebra and the
subgroups of $W_u$ give an isomorphism
\begin{equation}\label{eq:local-nonzero}
 G_{\epsilon,u}\simeq
 \mathbb H(\mathfrak t,R_u^\times,
                       \epsilon k_u|_{R_u^\times})\rtimes F_u.
\end{equation}
The action of $F_u$ on this connected factor is the ordinary
linear action on $\mathfrak t$ and conjugation on $W(R_u^\times)$.
\end{proposition}

\begin{proof}
The function $k_u$ is $W_u$-invariant. Hence $R_u^\times$ is stable
under its root reflections, contains both signs of each root, and
inherits the remaining reduced root-system axioms from $R_u$. Its
Weyl group is normal in $W_u$. The simple transitivity of this Weyl
group on its positive systems gives the asserted semidirect product.
The group $\Xi_u$ preserves $R_u^+$ and $k_u$, hence lies in $F_u$.

In the polynomial commutator with a Weyl element, the terms in the
graded cross relations are indexed by its inversion roots, each with
the corresponding parameter as a factor. This follows by induction
on a reduced expression. An element preserving $(R_u^\times)^+$ has
only zero-parameter inversions in its $W(R_u)$ component; its
$\Xi_u$ component already acts linearly. Thus elements of $F_u$
commute with polynomials by their ordinary linear action. For a simple
reflection of $R_u^\times$, the only nonzero-parameter inversion root
is its simple root. Its cross relation is precisely that for the
restricted parameter function. These observations define the
homomorphism. The PBW bases and the displayed group decomposition
make it bijective, including at $\epsilon=0$.
\end{proof}

There is also an explicit sign normalization. Put
$k_u^+(\alpha)=|k_u(\alpha)|$ on $R_u^\times$.
The assignments
$\varepsilon_u(s_\alpha)=\operatorname{sgn}(k_u(\alpha))$
for its simple reflections extend to a character of $W(R_u^\times)$:
odd braid relations join conjugate reflections, whose signs agree.
This character is $F_u$-invariant. Therefore the map fixing polynomials
and $F_u$, and sending $w\in W(R_u^\times)$ to
$\varepsilon_u(w)w$, gives
\begin{equation}\label{eq:local-positive}
 G_{\epsilon,u}\simeq
 \mathbb H(\mathfrak t,R_u^\times,\epsilon k_u^+)\rtimes F_u
                   \quad(0\leq\epsilon\leq1).
\end{equation}
The cross relation verifies this isomorphism directly.
The sign character is fixed throughout the family, including the
zero endpoint, where it can change Weyl-representation labels.
It must be included in the coordinate and action transports if this
positive presentation is used. No change of the torus origin is involved.

For example, take $X=\mathbb Z\alpha$ in rank one with equal affine
labels $Q_\alpha=Q_\alpha^*=q>1$ and $\alpha(u)=-1$.
Then $\alpha^\vee\in2Y$, $R_u=\{\pm\alpha\}$, but $k_u=0$.
The local algebra is $\mathcal O(\mathfrak t)\rtimes C_2$.
Here $R_u^\times$ is empty: the connected factor is toral and the
reflection is retained in $F_u=C_2$.

\subsubsection{Retaining the complete residual action locally}
A local neutral model by itself is insufficient for $B_\epsilon$.
The following construction gives its actual completed algebra and
specifies how to carry every source product through the reductions.
It avoids discarding matrix factors before their actions have been
transported.

\begin{proposition}[The full source completion]\label{prop:local-local-source}
Choose $z\in\operatorname{Spec}Z=T_X/W_f$, put
$\mathcal O=Qz$, $Q_z=\operatorname{Stab}_Q(z)$, and let $\bar z$ be
its image in $\operatorname{Spec}S$. There are algebra isomorphisms
\begin{align}
 \widehat B_{\epsilon,\bar z}
 &\simeq
 \left(\prod_{z'\in\mathcal O}\widehat A_{\epsilon,z'}\right)
                           \rtimes_{\beta,a}Q,
                         \label{eq:local-orbit-completion}\\
 &\simeq M_{|\mathcal O|}\left(
    M_{n_u}(\widehat G_{\epsilon,u,W_u\nu})
                          \rtimes_{\widetilde\beta,\widetilde a}Q_z
                          \right).
                         \label{eq:local-full-local}
\end{align}
In the second line choose $t=u\exp\nu$ over $z$ and the isomorphism
$\Psi=\Psi_{\epsilon,z}$ from \eqref{eq:local-completion}. The residual
data in that formula are, explicitly,
\begin{equation}\label{eq:local-local-data}
 \widetilde\beta_g=\Psi\beta_g\Psi^{-1},\qquad
 \widetilde a(g,h)=\Psi(a(g,h))\quad(g,h\in Q_z),
\end{equation}
where \eqref{eq:local-beta}--\eqref{eq:local-defect} specify the arguments
on the source generators. The multiplication in the local model is
\[
 (mU_g)(nU_h)=m\,\widetilde\beta_g(n)\,
                  \widetilde a(g,h)U_{gh}.
\]
These are finite residual extensions of the complete neutral matrix
algebras. They are equivalences of the full source categories, including
all projective actions on the retained factors.
\end{proposition}

\begin{proof}
The algebra $Z$ is finite over $S=Z^Q$, and the points over $\bar z$
form the orbit $\mathcal O$. Its completion is the product of its
completions at those points. Since $A_\epsilon$ is finite over $Z$ and
$Q$ is finite, the crossed presentation commutes with completion.
This proves \eqref{eq:local-orbit-completion}.

Let $p_z$ be the central idempotent of the completed neutral algebra
corresponding to $z$. The units $u_g$ conjugate these idempotents
transitively. Thus $p_z$ is full in the completed crossed product, and
its corner is
\[
 p_z\widehat B_{\epsilon,\bar z}p_z
       =\widehat A_{\epsilon,z}\rtimes_{\beta,a}Q_z.
\]
Choosing coset representatives and their inverse units gives matrix
units and the outer matrix factor $|\mathcal O|$. Applying $\Psi$ to
this corner gives \eqref{eq:local-full-local} and
\eqref{eq:local-local-data}. Because this is an algebra transport, it
preserves all products and actions, including those on the inner matrix
factor. The usual full-idempotent bimodules supply the asserted Morita
equivalences. Every simple $B_\epsilon$-module has a character of $S$
since $B_\epsilon$ is finite over $S$, so these completions cover its
simple spectrum.
\end{proof}

Composing $\Psi$ with \eqref{eq:local-positive} gives the same statement
with the inner algebra written using the positive-parameter connected factor
and $F_u$. The source data are still transported by
\eqref{eq:local-local-data}, including the fixed sign character at zero.

The groups $Q_z$, $W^\fs$ and the stabilizer $Q_V$ of a neutral
simple have distinct roles. The model retains
$M_{n_u}(\widehat G)$ together with its transported action and
defects; any Morita reduction transports these data through its
bimodules. The Levi degree $[\Gamma:\Lambda]^{1/2}$ need not equal
$n_u$, since the finite crossed data also carry multiplicities.

On a neutral simple $V$ of central character $z$, choose implementing
operators with intrinsic multiplier $c_V$. The explicit external
factor from \eqref{eq:local-defect} gives
\begin{equation}\label{eq:local-effective}
 \kappa^{\mathrm{eff}}_V(g,h)=a(g,h)(z)c_V(g,h)^{-1},
                \qquad g,h\in Q_V.
\end{equation}
For the chosen section this is the strict-action form of Lemma~\ref{lem:weak-Clifford}. Transport by \eqref{eq:local-local-data} preserves its class
at each fixed $\epsilon$. It does not compare that class with the
class of a different simple module after parameter removal.

The source scalars $\chi_\omega(c)$ need not have modulus one under
an arbitrary fixed normalization of $E_\pi$. Accordingly a $Q$-orbit
of neutral central characters is treated point by point with its own
polar decomposition. No preservation of the chosen compact torus by
the full residual affine action is assumed.

\subsubsection{Auxiliary Bernstein coordinates and root normalization}
The connected geometric recognition uses a comparison of Bernstein
cross relations. The source root system $R$, whose root lattice is
$X_{\rm aff}$, is kept throughout. Write $R_{\rm g}$ for the roots
in a geometric affine presentation. These roots can differ from the
source roots in the same rational vector space. In particular, the
root lattice of $R_{\rm g}$ need not be $X_{\rm aff}$.

Use the depth-zero hypotheses of \S\ref{subsec:affine-input}, the
chosen affine comparison of Proposition~\ref{prop:ohara-affine}, and
the source normalization $q_s>1$ on relevant affine simple reflections.
We use an auxiliary full \emph{unipotent} algebra to recognize the
connected factor. The residual group and multiplier of the original
source algebra are those in \eqref{eq:source-presentation}.

\begin{proposition}[An auxiliary Bernstein embedding]
\label{prop:local-auxiliary}
There is a choice of the finite Weyl complement $W_f$ and a
$W_f$-equivariant finite-index lattice inclusion
\[
 \iota:X=X_{\rm aff}\oplus\Lambda_c
       \hookrightarrow L=\Bigl(\bigoplus_i L_i\Bigr)\oplus\Lambda_c
\]
with the following properties. There are geometric affine algebras
$\mathcal H_i^{\rm g}(v_i^\epsilon)$ on the lattices $L_i$, where
$v_i=q_{F_i}^{1/2}$. Put
\[
 \mathcal A^{\rm g}_\epsilon
 =\Bigl(\bigotimes_i\mathcal H_i^{\rm g}(v_i^\epsilon)\Bigr)
       \otimes\C[\Lambda_c].
\]
Its full character lattice is $L$. At $\epsilon=1$ there is an
injective algebra map
\begin{equation}\label{eq:local-aux-embedding}
 I_1:A_1\hookrightarrow\mathcal A^{\rm g}_1,
 \qquad \vartheta_{x,1}\longmapsto\theta_{\iota(x),1},
 \qquad N_{s,1}\longmapsto N^{\rm g}_{s,1}
             \quad(s\text{ simple in }W_f).
\end{equation}
Here the generators on both sides have the same quadratic parameter.
The roots $\alpha_s\in R$ and $\beta_s\in R_{\rm g}$ for such an $s$
are positive rational multiples of one another; they are not asserted
to be equal. The image has PBW decomposition
$\bigoplus_{w\in W_f}\C[\iota(X)]N^{\rm g}_{w,1}$.
\end{proposition}
\begin{proof}
Proposition~\ref{prop:ohara-affine} supplies a parahoric unipotent
datum for an auxiliary connected group $G_c$, together with a based,
parameter-preserving identification of its affine factor with the source
factor. Passing to the adjoint group preserves this factor and its labels
\cite[Lemma~3.5]{Sol-unip}. Work factor by factor with absolutely simple
adjoint groups, allowing restriction of scalars. The extra toral
factor will be added at the end.

For an auxiliary parahoric unipotent type $(P_{\mathfrak f},\sigma)$,
use the groups in \cite[\S3.1, equations~(27)--(34)]{Sol-unip}:
\[
 W_{\rm aff}(J,\sigma)=X(J)\rtimes W^\circ(J,\sigma),\qquad
 W(J,\sigma)/\Omega_{\mathfrak f,\rm tor}
                =L_i\rtimes W^\circ(J,\sigma).
\]
Choose the finite Weyl subgroup as in that presentation and pull it
back to the source. Thus $X(J)=X_{{\rm aff},i}$ and the finite group
is the corresponding factor of $W_f$. The lattice $L_i$ is the
\emph{translation lattice of this auxiliary extended group}.
In particular it is not defined to be the root lattice of the
normalized geometric root system.

Choose an extension of $\sigma$ and a character of
$\Omega_{\mathfrak f,\rm tor}$. The corresponding algebra contains
$\mathcal H(W_{\rm aff},q)$ injectively, by
\cite[equations~(39)--(42)]{Sol-unip}. Its Bernstein translation
monomials restrict to the source monomials: for a dominant
$x\in X(J)$ they are the same normalized translation element, and
multiplicativity gives the assertion for every $x$. Its finite
Bernstein generators are the same affine support generators for
$W^\circ(J,\sigma)$.

Now use the arithmetic/geometric \emph{Bernstein-presentation}
comparison in \cite[Lemma~4.5, equations~(63)--(66) and its concluding
paragraph]{Sol-unip}, based on \cite[Theorem~10.11]{Lus-unip}.
Use the Laurent coordinates and basepoints of that presentation
comparison. It identifies the Laurent monomials through a lattice
isomorphism and sends each finite generator $T_s$ to $T_s^{\rm g}$.
The common positive quadratic parameter then gives the same assertion
for $N_s$. This is the generator-level input in
\eqref{eq:local-aux-embedding}; an arbitrary equivalence of module
categories or an arbitrary choice of LLC origin would not suffice.
We do not identify the affine simple generators associated with the
\emph{normalized geometric roots} with all the source affine simple
generators. The changes between those presentations are precisely
what the next proposition checks.

In the common rational character space the $(-1)$-eigenspace of a
reflection is one-dimensional. The two positive roots for that
reflection therefore lie on the same rational ray, after choosing
corresponding finite chambers. The PBW bases give the stated image.
Finally take products and adjoin the torus with character lattice
$\Lambda_c$, on which the embedding is the identity.
\end{proof}

The finite lattice information in this construction is explicit.
Place the chosen finite vertex at the origin in the auxiliary
apartment. Write an auxiliary length-zero element as
$\omega=t_{v_\omega}w_\omega$ with $w_\omega\in W_f$, as in
\cite[Proposition~3.1 and equation~(33)]{Sol-unip}. Then
\begin{equation}\label{eq:local-geo-lattice}
 L_i=X_{{\rm aff},i}+\sum_{\omega\in
       \Omega_{\mathfrak f}/\Omega_{\mathfrak f,\rm tor}}
                   \mathbb Z v_\omega
       \ \subset X_{{\rm aff},i}\otimes\mathbb Q.
\end{equation}
The quotient is finite because the auxiliary ambient group is
semisimple adjoint. Equivalently,
$L_i/X_{{\rm aff},i}\simeq
\Omega_{\mathfrak f}/\Omega_{\mathfrak f,\rm tor}$.
Formula~\eqref{eq:local-geo-lattice} specifies the inclusion as well as
its index. Transport it through the fixed geometric character-lattice
identification, and suppress $\iota$ when writing elements of $X$ in $L$.
No source or geometric root lattice is substituted for $L_i$.

Consequently the finite torus diagram is
\begin{equation}\label{eq:local-geo-covers}
 T^\fs\xrightarrow{\ r_\pi\ }T_X
           \xleftarrow{\ p\ }T_L,
 \qquad K_{\rm g}=\ker(p)=\Hom(L/X,\C^\times).
\end{equation}
The underlying isogeny of $r_\pi$ has kernel
$K_X=\Hom(\Lambda/X,\C^\times)$; $r_\pi$ itself includes the origin
translation in \eqref{eq:local-isogeny}. These two finite groups are
retained separately. The auxiliary intersection lattice is $X\cap L=X$.

\subsubsection{The exact comparison of cross relations}
Put $a(Q)=Q^{1/2}-Q^{-1/2}$. For every source root set
\begin{equation}\label{eq:local-Bernstein-coefficient}
 B^{\rm s}_{\alpha,\epsilon}
 =\frac{a(Q_\alpha^\epsilon)
           +a((Q_\alpha^*)^\epsilon)\theta_{-\alpha}}
         {1-\theta_{-2\alpha}}\in\C(T_L).
\end{equation}
This is also the coefficient in the primitive-coroot case: our
convention $Q_\alpha^*=Q_\alpha$ cancels $1+\theta_{-\alpha}$.
Define $B^{\rm g}_{\beta,\epsilon}$ by the same formula for the
geometric Bernstein labels, setting $Q_\beta^*=Q_\beta$ when its
coroot is primitive. The quadratic label for corresponding finite
simple generators is a common number $Q>1$.

\begin{proposition}[Bernstein comparison with different roots]
\label{prop:local-coefficients}
For every corresponding positive root ray, the two possibilities are
\[
\begin{array}{c|c|c}
\text{roots in }L\otimes\mathbb Q
   &\text{source labels}&\text{geometric labels}\\ \hline
\alpha=\beta &(Q,Q^*)&(Q,Q^*)\\
\alpha=2\beta &(Q,Q)&(Q,1).
\end{array}
\]
In the second row the geometric coroot is divisible by $2$.
These comparisons imply, for every $0\leq\epsilon\leq1$,
\begin{equation}\label{eq:local-coefficient-equality}
 B^{\rm s}_{\alpha,\epsilon}
           =B^{\rm g}_{\beta,\epsilon}\quad\text{in }\C(T_L).
\end{equation}
The generator maps in \eqref{eq:local-aux-embedding} therefore extend
to injective algebra maps $I_\epsilon:A_\epsilon\hookrightarrow
\mathcal A^{\rm g}_\epsilon$ throughout the family.
\end{proposition}
\begin{proof}
Initially let $s$ be simple in $W_f$. Apply $I_1$ to the source cross
relation for an $f\in\C[X]$ with $f-s(f)\ne0$ and compare with the
geometric cross relation. Cancellation in $\C(T_L)$ gives
$B^{\rm s}_{\alpha,1}=B^{\rm g}_{\beta,1}$. This step uses only the
embedding already proved, not an equality of roots or of the second
labels.

Choose a primitive lattice vector $\chi$ on the common positive ray,
and write $\alpha=m\chi$, $\beta=n\chi$ with positive integers $m,n$.
Set $z=\theta_{-\chi}$ and
\[
 c=\frac{a(Q_\alpha^*)}{a(Q)},\qquad
 d=\frac{a(Q_\beta^*)}{a(Q)}.
\]
All source parallel-wall labels are $>1$, so $c>0$. Geometric labels
are positive, and $d$ is a real number. The identity becomes
\begin{equation}\label{eq:local-ray-identity}
 \frac{1+cz^m}{1-z^{2m}}
             =\frac{1+dz^n}{1-z^{2n}}.
\end{equation}
If $d\ne0$, the first nonconstant terms in the formal power series
at $z=0$ have degrees $m$ and $n$. Hence $m=n$, and then $c=d$.
Since $a$ is injective on the positive real numbers, the second labels
agree. If $d=0$, the first nonconstant term on the right has degree
$2n$ and coefficient $1$. Hence $m=2n$ and $c=1$; therefore
$Q_\alpha^*=Q$ and $Q_\beta^*=1$. The latter is incompatible with the
primitive-coroot convention because $Q>1$.
This proves the table for simple reflections. Weyl covariance of the
root-parameter functions and of the torus inclusion proves it on
all corresponding root rays.

In the first row the coefficient equality persists on raising every
label to $\epsilon$. In the second row it is the identity
\[
 \frac{a(Q^\epsilon)(1+\theta_{-2\beta})}
      {1-\theta_{-4\beta}}
       =\frac{a(Q^\epsilon)}{1-\theta_{-2\beta}}.
\]
Thus the assertion for the parameter family follows from explicit
cross relations, rather than from an inference about a deformation
from one numerical specialization. The finite braid and quadratic
relations agree as well. The Bernstein presentations consequently
define $I_\epsilon$, and the PBW bases prove injectivity, including
at $\epsilon=0$.
\end{proof}

\begin{proposition}[The nonzero local data are preserved as ratios]
\label{prop:local-lattice-recognition}
Let $u\in T_{X,\rm un}$ and choose a unitary lift
$\widetilde u\in T_{L,\rm un}$. The nonzero local reflection sets agree:
\begin{equation}\label{eq:local-reflection-match}
 \{s_\alpha:\alpha\in R_u,\ k_u(\alpha)\ne0\}
 =\{s_\beta:\beta\in(R_{\rm g})_{\widetilde u},\
                                      k^{\rm g}_{\widetilde u}(\beta)\ne0\}.
\end{equation}
For the corresponding positive roots the precise parameter identity is
\begin{equation}\label{eq:local-ratio-match}
 \frac{k_u(\alpha)}{\alpha}
      =\frac{k^{\rm g}_{\widetilde u}(\beta)}{\beta}
          \quad\text{as rational functions on }\mathfrak t.
\end{equation}
The full torus stabilizers need not agree. Every source reflection
lost on lifting has zero graded parameter and remains in the source
residual factor $F_u$.
\end{proposition}
\begin{proof}
A unitary lift exists for a finite-index lattice inclusion. Restriction
is Weyl equivariant, so any reflection fixing $\widetilde u$ fixes $u$.
Consider the two rows of Proposition~\ref{prop:local-coefficients}.

If $\alpha=\beta$, their values at the two points agree. A source
reflection can be lost only when $\alpha(u)=-1$ and its coroot becomes
primitive on $L$. In that case the geometric primitive-coroot convention
and the first row force $Q_\alpha^*=Q_\alpha$, so $k_u(\alpha)=0$.
For every retained reflection the label formulas give the same $k$,
including when the root value is $1$ and coroot divisibility changes.

If $\alpha=2\beta$, the labels are $(Q,Q)$ and $(Q,1)$, and
$\alpha(u)=\beta(\widetilde u)^2$. When the source root value is $1$,
the geometric value is $1$ or $-1$. Its coroot is divisible, so the
geometric reflection fixes $\widetilde u$. The two parameters are
\[
 k_u(\alpha)=\log Q,\qquad
 k^{\rm g}_{\widetilde u}(\beta)=\tfrac12\log Q,
\]
which give \eqref{eq:local-ratio-match}. A source reflection at root
value $-1$, when allowed by its coroot divisibility, has $k_u=0$;
its geometric root value is $i$ or $-i$, so that reflection is lost.
These cases prove the assertions in both directions. Equivalently,
\eqref{eq:local-ratio-match} is the common limit of the coefficients
in \eqref{eq:local-coefficient-equality} at
$\widetilde u\exp(\epsilon\nu)$ as $\epsilon\downarrow0$, for generic
real $\nu$. The case analysis also treats the vanishing limits.
\end{proof}

\paragraph{A rank-one root normalization.}
Let $X=\mathbb Z\alpha\subset L=\mathbb Z\omega$ with
$\alpha=2\omega$ and source labels $(q,q)$, $q>1$.
The normalized geometric root is $\beta=\omega$ with labels $(q,1)$.
For $z=\theta_\omega$ the two coefficients are
\[
 \frac{a(q^\epsilon)(1+z^{-2})}{1-z^{-4}}
                =\frac{a(q^\epsilon)}{1-z^{-2}}.
\]
The kernel of $T_L\to T_X$ is $\mu_2$. At $\alpha(u)=1$ its two
lifts have $\omega(\widetilde u)=\pm1$; the local parameters are
$\log q$ on $\alpha$ and $\tfrac12\log q$ on $\omega$.
Their ratios agree. At $\alpha(u)=-1$ the lifts have values $\pm i$,
the source graded parameter is zero, and the lost reflection remains
in $F_u$. Thus both the additional root change and the finite torus
kernel are visible, including on the exceptional unitary stratum.

\subsubsection{Connected geometric recognition}
\begin{proposition}[Connected geometric datum on each source stratum]
\label{prop:local-geometric}
For every $u$, the connected factor
\[
 C_{\epsilon,u}=\mathbb H(\mathfrak t,R_u^\times,\epsilon k_u^+)
\]
has a geometric realization by connected reductive groups, Levi
subgroups and cuspidal nilpotent local systems. The realizations can
be chosen uniformly for $0\leq\epsilon\leq1$, with the real polynomial
space and all toral coordinates retained. 
\end{proposition}
\begin{proof}
Fix a unitary lift as in Proposition~\ref{prop:local-lattice-recognition}.
For each auxiliary factor use the notation
$(G\vartheta,G_J,\mathcal C,\mathcal F)$ of
\cite[\S\S7--9]{Lus-unip}, with $K=I\setminus J$. The notation $J$
here is an affine-node set, not a lattice. Represent the unitary Weyl
orbit by $p'(d)$, with $d$ on the corresponding face $C_S$,
$\varnothing\ne S\subset K$. One may take the auxiliary logarithmic
complement in \cite[\S9.1]{Lus-unip} to contain $i\mathbb R$, so these
unitary points lie in its $C_S^\diamond$.
The connected cuspidal datum is
\begin{equation}\label{eq:local-cuspidal-datum}
 G_u=Z_G(p(d))=G_{I\setminus S},\qquad
 L_u=G_J,\qquad (\mathcal C,\mathcal F).
\end{equation}
Connectedness is \cite[Lemma~6.20]{Lus-unip}; the Levi assertion is
\cite[Lemma~7.6]{Lus-unip}. The original cuspidal local system on $G_J$
is retained.

Record the Weyl element used to choose the unitary orbit representative.
If the positive basis of the local root system in the source reduction
and the face basis in Lusztig's description differ, also choose the
unique local Weyl element sending the first positive system to the
second. Apply these Weyl maps to polynomial coordinates and to
reflection generators together. They are real linear and fix every
independent toral direction. They are included in $J_{\epsilon,u}$
and in the subsequent transports of the source actions. The following
identities are written after these specified coordinate changes; no
claim that a change of positive chamber fixes individual polynomial
coordinates is needed.

In the normalized geometric affine presentation, write the local
simple roots as $\beta_k=\widehat\gamma_k$.
Lusztig's calculation in \cite[\S9.3]{Lus-unip} gives
\[
 \mu(\beta_k)=d_k\bar z_k c_k,\qquad
 \gamma_k^*=(d_k\bar z_k)^{-1}\beta_k,\qquad
 \mu'(\gamma_k^*)=c_k>0.
\]
With $r_i=\log v_i=\tfrac12\log q_{F_i}$, the normalized geometric
local parameter is $k^{\rm g}_{\widetilde u}=r_i\mu$.
Here one uses $Q_\beta=v_i^{2\lambda(\beta)}$ and
$Q_\beta^*=v_i^{2\lambda^*(\beta)}$ in the divisible-coroot case,
with the primitive convention otherwise. These parameters are indexed by the geometric root $\beta$.

Combine this calculation with \eqref{eq:local-ratio-match}. If $\alpha_k$
is the corresponding source local simple root, the complete comparison is
\begin{equation}\label{eq:local-root-rescaling}
 \frac{k_u(\alpha_k)}{\alpha_k}
   =\frac{r_i\mu(\beta_k)}{\beta_k}
   =\frac{r_i c_k}{\gamma_k^*}.
\end{equation}
The first equality incorporates the preliminary Bernstein root change
proved above; the second is the geometric graded root rescaling.
They identify the graded algebras by the identity on polynomials and
corresponding reflection generators. Equality of individual roots
or of individual $k$ values is unnecessary and is not used.

The auxiliary local Weyl group is $W_{K\setminus S}$ and its
parameters are positive, by \cite[\S9.3]{Lus-unip}.
Equation~\eqref{eq:local-reflection-match} identifies its reflections
with exactly the nonzero source reflections. Their positive systems
are compared through the common reflection arrangement; hence their
simple reflections correspond. The geometric algebra in
\eqref{eq:local-cuspidal-datum} consequently realizes the entire
connected source factor. The fixed sign map
\eqref{eq:local-positive} is included if that presentation is used.
Every cross-relation identity persists on replacing $r_i$ by
$\epsilon r_i$. At zero the same reflection subgroup and the same
polynomial space are retained, so the identification includes that
endpoint as well.

Take products over the auxiliary factors, and multiply both each
geometric group product and its Levi by the torus with character
lattice $\Lambda_c$, with its trivial cuspidal datum. This supplies
$\operatorname{Sym}(\Hom(\Lambda_c,\C)^*)$. The other central
directions already present in $G_u$ are retained. The lattice
inclusion and the recorded Weyl maps identify the real Lie spaces.
The root rescalings themselves leave polynomial coordinates fixed. No descent of a cuspidal local
system through an unspecified central isogeny is required.
The groups $F_u$, $Q$ and their source defects remain in the local
source extension, as before.
\end{proof}

\subsubsection{The connected shift and its inverse}
Write $\mathcal H_{\epsilon,u}^{\rm geo}$ for the connected geometric
algebra just constructed, including its toral factor, and fix the
isomorphisms
\[
 J_{\epsilon,u}:C_{\epsilon,u}\xrightarrow{\sim}
                      \mathcal H_{\epsilon,u}^{\rm geo}.
\]
Here a product with several parameters means a tensor product of
geometric factors, specialized at $(\epsilon r_i)_i$; it does not
assert that all residue-field parameters are equal.

\begin{proposition}[Connected parameter removal on real strata]
\label{prop:local-connected-shift}
There is a fixed-choice bijection
\[
 \Theta_u^\circ:\Irr_{\rm real}(C_{1,u})
                   \xrightarrow{\sim}\Irr_{\rm real}(C_{0,u}).
\]
It uses the generalized-Springer support condition and the
distinguished zero-parameter simple specified in the proof. It preserves the independent toral
coordinate. It need not preserve the remaining central character.
\end{proposition}
\begin{proof}
Transport first through $J_{1,u}$. For one geometric factor a simple
is labelled by an admissible enhanced triple $(y,\sigma,\rho)$ with
$[\sigma,y]=2r y$. Choose a compatible $\mathfrak{sl}_2$ triple with
semisimple element $h$, normalized by $[h,y]=2y$, and put
\[
 \sigma_0=\sigma-rh,\qquad [\sigma_0,y]=0.
\]
There is a natural component-group identification
\[
 \pi_0 Z_{G_u}(\sigma,y)\simeq
                   \pi_0 Z_{G_u}(\sigma_0,y)
\]
and hence a transported enhancement $\rho_0$.
The precise independence on enhanced conjugacy classes is
\cite[Lemma~3.6]{AMS-GHA}. Admissibility means that the generalized-Springer
cuspidal support of $(y,\rho_0)$ in $Z_{G_u}(\sigma_0)$ is conjugate
to $(L_u,\mathcal C,\mathcal F)$. This generalized-Springer support condition is part of the
admissibility requirement.

The classification in \cite[Theorem~3.20(c),(d) and
Proposition~3.8]{AMS-GHA}, with the erratum understood, now identifies the
positive simple with the distinguished zero-parameter simple labelled
by $(y,\sigma_0,\rho_0)$. At zero it is the simple supplied by the
appropriate generalized-Springer representation and polynomial
character, as specified in that proposition. It is not the entire
standard homology module.

Conversely, start with that admissible enhanced zero label. Choose an
algebraic $\mathrm{SL}_2$ homomorphism for $y$ inside
$Z_{G_u}(\sigma_0)$ and restore $\sigma=\sigma_0+rh$.
The same component-group identifications restore $\rho$.
Lemma~3.6 and Theorem~3.20 show that these constructions are inverse
on simple classes. They also give independence of the chosen compatible
triple on those classes.

The real-stratum assertion is needed here because the preceding affine
reduction fixes a unitary part. For real $r$, the geometric
classification has real central character if and only if $\sigma_0$
is conjugate into the real polynomial space; see
\cite[Theorem~5.7(c)]{Sol-survey}. Its use of the
Iwahori--Matsumoto involution does not change this criterion, since
negation preserves the real space. Consequently both directions of
the shift stay in the required real stratum. An algebraic
$\mathrm{SL}_2$ homomorphism has trivial image in the independent
central torus, so its coordinate is unchanged. Apply this argument
factor by factor, with the respective $r_i$, and transport back
through $J_{0,u}$. If a sign normalization was made, its zero-endpoint
map is included in this last transport, rather than resetting Weyl
labels at zero.
\end{proof}

\subsubsection{Exact interface with the source residual action}
The geometric identification is an algebra isomorphism, so it also
specifies how the already constructed source local model is to be
written geometrically. Transport the actual $F_u$-action through
$J_{\epsilon,u}$ and set
\[
 D_{\epsilon,u}=\mathcal H_{\epsilon,u}^{\rm geo}
                    \rtimes_{\eta_{\epsilon,u}}F_u,
 \qquad \eta_{\epsilon,u}(f)=J_{\epsilon,u}fJ_{\epsilon,u}^{-1}.
\]
Together with the fixed zero-separation and sign maps this gives an
isomorphism $\mathcal J_{\epsilon,u}:G_{\epsilon,u}\simeq D_{\epsilon,u}$.
Compose its entrywise matrix map with $\Psi_{\epsilon,z}$ in
\eqref{eq:local-local-data}, including the corresponding central
completion, and denote the composite by $\Phi_{\epsilon,z}$. Then
\begin{equation}\label{eq:local-geometric-source-action}
 \beta^{\rm geo}_{g,\epsilon,z}
   =\Phi_{\epsilon,z}\beta_g\Phi_{\epsilon,z}^{-1},\qquad
 a^{\rm geo}_{\epsilon,z}(g,h)=\Phi_{\epsilon,z}(a(g,h)).
\end{equation}
These formulas, in the matrix algebra over the completed
$D_{\epsilon,u}$, retain exactly the original source products and
matrix multiplicities in \eqref{eq:local-full-local}. Toral functions
and their defect evaluations are transported through the specified
coordinates. The automorphisms are defined by these algebra
transports; the finite-group comparison below applies to the original
neutral affine simples through an equivariant full corner.

\subsection{Residual equivariance and Clifford classes}\label{subsec:source-equivariance}
\subsubsection{An equivariant form of the connected shift}
Write $W_u^\times=W(R_u^\times)$. The action of $F_u$ on
$C_{\epsilon,u}$ preserves the real polynomial space, the positive
root system and the positive parameter function. In choosing the
geometric models in \S\ref{subsec:source-local-models}, use the same model and root normalization on
components permuted by $F_u$. This is possible because these components
have isomorphic based root data with the same parameters. Root
normalizations within a component are constant on its diagram-automorphism
orbits. After conjugating by the fixed chamber changes in \S\ref{subsec:source-local-models}, the
transported action on the normalized geometric root datum is therefore
a based linear action. This assertion concerns the graded presentation;
it does not require lifting the original source group to an automorphism
group of the chosen cuspidal local system.

We use the following representation-theoretic substitute for such a
geometric lift. For a positive geometric graded Hecke algebra $H$,
let $\zeta_H'$ denote the real-central-character bijection of
\cite[Theorem~6.2, Lemma~6.6 and Proposition~6.10]{Sol-survey}.
On tempered representations it picks the distinguished finite Weyl
type. That type occurs with multiplicity one, and the restriction
matrix is triangular with diagonal entries~$1$. For all real central
characters it is obtained from this map by the Langlands classification.
The geometric parameter functions here meet the restrictions in those
results, including the type $F_4$ restriction; we apply them after the
root normalizations of \S\ref{subsec:source-local-models}.

\paragraph{The involution convention matters.}
The map in \S\ref{subsec:source-local-models} uses the untilded AMS convention. The survey uses
\[
 \widetilde M_{y,\lambda,r,\rho}
       =\operatorname{IM}_r^*M_{y,rh-\lambda,r,\rho},\qquad
 \operatorname{IM}(f)(\nu)=f(-\nu),\quad
 \operatorname{IM}(w)=\operatorname{sgn}(w)w.
\]
Here $\operatorname{IM}$ fixes the parameter. Thus the connected map
here is, in the geometric coordinates,
\begin{equation}\label{eq:equiv-im}
 \Theta^\circ=\operatorname{IM}_0^*\,\zeta_H'\,
                                      \operatorname{IM}_1^*.
\end{equation}
In particular $M_{y,\sigma_0+rh,r,\rho}$ maps to
$M_{y,\sigma_0,0,\rho_0}$. Formula~(5.14), Lemma~6.3 and the
geometric compatibility discussed in \S6.2 of \cite{Sol-survey}
give this identification: in the Langlands datum, the centralizer of
the real semisimple label supplies the Levi, the tempered Weyl type
supplies the generalized-Springer label, and the real central
translation is retained. Applying the same construction inside that
Levi gives the distinguished zero simple of \cite[Proposition~3.8]{AMS-GHA}.
Both endpoint involutions are necessary. They commute with based
linear actions, as do the fixed sign maps of \S\ref{subsec:source-local-models}.

\begin{proposition}[Equivariance and intrinsic cocycles for $F_u$]
\label{prop:equiv-fu}
The connected bijection $\Theta_u^\circ$ can be chosen as above so that
it is $F_u$-equivariant. For every real-central-character connected
simple $V$, put $V_0=\Theta_u^\circ(V)$ and
$F_{u,V}=\operatorname{Stab}_{F_u}[V]$. Then
\[
 F_{u,V}=F_{u,V_0},\qquad
 [c_{1,V}]=[c_{0,V_0}]\quad\hbox{in }
                    H^2(F_{u,V},\C^\times).
\]
The convention here is $T_gT_h=c_V(g,h)T_{gh}$ for a strict action.
The same assertions hold for any finite group of based linear
automorphisms preserving this positive graded datum and its toral
decomposition, including an ineffective action.
\end{proposition}
\begin{proof}
First use the survey convention and take $V$ tempered with real
central character and zero independent central translation. Its
distinguished $W_u^\times$-type $\tau=\zeta_H'(V)$ is characterized
by the triangular restriction bijection. Uniqueness makes this
assignment equivariant, hence identifies its stabilizers. On the
stabilizer, every intertwiner $T_g$ of $V$ preserves the unique
$\tau$-isotypic subspace, which is one copy of $\tau$. Restriction
therefore gives intertwiners of $\tau$ with exactly the same products:
\[
 (T_g|_\tau)(T_h|_\tau)
                 =c_V(g,h)(T_{gh}|_\tau).
\]
This proves equality of cocycle representatives for these choices,
without inferring it from a bijection of labels. It is the
multiplicity-one argument in \cite[Lemma~6.6]{Sol-survey}.

For a general real simple, use its unique Langlands datum
$(P,\tau,\lambda)$, with $\tau$ tempered of real central character
for the semisimple Levi factor and $\lambda$ real and strictly
positive outside $P$. An element stabilizing the simple stabilizes
this datum, because the action preserves the chosen positive chamber.
Intertwiners of $\tau$ induce intertwiners of the standard module
with the same scalar products. They descend to its unique irreducible
quotient, still with those products. At zero, replace $\tau$ by its
distinguished $W_P$-type and retain $\lambda$. The resulting induced
zero module is irreducible: the Weyl stabilizer of $\lambda$ is
$W_P$, so this is precisely Mackey induction. The same induced
intertwiners again have the same products. These constructions are
the ones in \cite[Corollary~6.8 and Proposition~6.10]{Sol-survey}.
The independent real toral coordinate is retained throughout; a
stabilizer must fix it. This proves the assertion for all real simples.

Finally apply \eqref{eq:equiv-im} at both endpoints and transport through
the identifications of \S\ref{subsec:source-local-models}. Pullback by an equivariant algebra isomorphism
does not change an intertwiner product. All root, chamber and sign
normalizations are used at both endpoints. This proves the assertion
in the convention of $\Theta_u^\circ$.
\end{proof}

\begin{proposition}[Parameter removal for the full neutral local model]
\label{prop:equiv-local-full}
There is a fixed-choice bijection
\begin{equation}\label{eq:equiv-local-full}
 \Theta_u:\Irr_{\rm real}(G_{1,u})
                    \xrightarrow{\sim}\Irr_{\rm real}(G_{0,u}).
\end{equation}
It includes $F_u$, all zero-parameter reflections and all toral
directions. It preserves the toral coordinate and has an inverse
given by the reverse connected map and the inverse comparison
cochains.
\end{proposition}
\begin{proof}
Use $G_{\epsilon,u}\simeq C_{\epsilon,u}\rtimes F_u$ with the
fixed zero-endpoint sign identification. Choose one $V$ in every
$F_u$-orbit of connected simples. Its effective Clifford multiplier
is $\kappa_{\epsilon,V}=c_{\epsilon,V}^{-1}$, since this extension
has no external defect. Proposition~\ref{prop:equiv-fu} identifies
the stabilizers and their effective classes. Choose a normalized
cochain $b_V$ with
\[
 \kappa_{1,V}=\delta b_V\,\kappa_{0,\Theta_u^\circ V},\qquad
 \delta b_V(g,h)=b_V(g)b_V(h)b_V(gh)^{-1}.
\]
For a projective label with operators $R_g$, replace them by
$b_V(g)^{-1}R_g$, then induce the resulting Clifford module from
the common stabilizer. This is the strict case of Lemma~\ref{lem:weak-Clifford}.
Use conjugation transport from the chosen representatives for every
other description. The inverse uses $b_V$, rather than $b_V^{-1}$.
For the ordinary based extension one may fix the distinguished-type
normalizations in the preceding proof, obtaining equal representatives
and $b_V=1$. Arbitrarily renormalized intertwiners are handled by
the displayed formula.
\end{proof}

\paragraph{A compatible choice for nested finite groups.}
When another based finite action is present, fix the preceding
normalizations by the distinguished finite type, rather than choosing
unrelated cochains for each subgroup. More generally, let $F$ be a
normal subgroup of a finite group $E$ acting on the based positive
graded datum. Apply \cite[Theorem~6.2 and Proposition~6.10]{Sol-survey}
to $H\rtimes F$ and $H\rtimes E$. The distinguished
$W\rtimes F$-type occurs with multiplicity one in a real tempered
$H\rtimes F$-simple. Intertwiners for its $E/F$ stabilizer restrict
to that same subspace. If representatives multiply with a defect
$f(g,h)\in F$, restriction retains its operator as well as the
intertwiners. Hence it identifies the \emph{effective} cocycle of
this weak quotient action. The Langlands argument in the proof above
extends the assertion to all real simples. Consequently this choice
of parameter removal commutes with Clifford induction in stages.
The assertion also applies to a fixed scalar twist: first use one
cochain at both endpoints to replace its cocycle by a cohomologous
root-of-unity-valued representative. Realize that representative in
a finite central extension and restrict to the prescribed central
character. Multiplicity one preserves that character; then undo the
same cochain change at both endpoints. This supplies
the compatible finite-group choices used below; it does not apply to
an arbitrary action on a completed matrix algebra.

\subsubsection{Preservation of the external source defects}
Put $T_c=\Hom(\Lambda_c,\C^\times)$. If a neutral simple $V$ has
toral character $\tau\in T_c$, formula~\eqref{eq:local-defect} gives
\begin{equation}\label{eq:equiv-external}
 a(g,h)(V)=d(g,h)\tau(c(g,h)),\qquad
 d(g,h)=\frac{\mu(\sigma(g),\sigma(h))}
              {\mu(c(g,h),\sigma(gh))b(c(g,h))}.
\end{equation}
Every connected shift above preserves $\tau$: polar reduction writes
it as $u_c\exp(\nu_c)$ and both factors are retained. Therefore
\begin{equation}\label{eq:equiv-external-equal}
 a(g,h)(V)=a(g,h)(V_0)
\end{equation}
for corresponding neutral simples, even when their semisimple central
characters differ. This is equality of values, not just of classes.
No unitarity assumption on the scalars $\chi_\omega(c)$ is needed.

When a $Q$-equivariant neutral bijection is available, its common
stabilizer $Q_V$ thus has
\begin{equation}\label{eq:equiv-reduced-obstruction}
 \kappa^{\rm eff}_{\epsilon,V}(g,h)
  =d(g,h)\tau(c(g,h))c_{\epsilon,V}(g,h)^{-1}.
\end{equation}
The remaining effective-class comparison for the original strict
action on $A_\epsilon$ is exactly the comparison of its intrinsic
$c_V$ classes. This simplification is specific to the absorbed
\emph{central} translation lattice. It does not justify deleting
defect units in an unrelated weak presentation. Changing to matrix
corners must transport the same products through the Morita bimodules.

\subsubsection{The vertex-preserving source action}
Let $o$ be the chosen special vertex in the semisimple quotient
apartment. Define
\[
 Q^{\rm vert}=\{g\in Q:\sigma(g)o=o\}.
\]
This is well-defined because $\Lambda_c$ acts trivially on that
apartment. It is a subgroup; no normality in $Q$ is asserted.
It includes any ineffective subgroup. Its elements preserve both the
alcove and its chosen vertex, so act by automorphisms of the finite
based Bernstein datum. On the independent torus they still act by
the actual formula~\eqref{eq:local-beta}.

\begin{proposition}[The vertex-preserving source extension]
\label{prop:equiv-vertex}
Use the compatible finite-type normalizations above. The neutral
parameter-removal map assembled in \S\ref{subsec:source-assembly} is $Q^{\rm vert}$-equivariant
and preserves the intrinsic and effective stabilizer classes for that
action. There is consequently a fixed-choice bijection
\[
 \Irr\bigl(A_1\rtimes_{\beta,a}Q^{\rm vert}\bigr)
       \simeq
 \Irr\bigl(A_0\rtimes_{\beta,a}Q^{\rm vert}\bigr).
\]

\end{proposition}
\begin{proof}
First fix a toral character $\tau$ and use Clifford induction over
its finite $Q^{\rm vert}$-orbit. Its stabilizer acts on the
semisimple neutral affine algebra by based root-datum automorphisms;
the external cocycle is the fixed scalar specialization
\eqref{eq:equiv-external}. Work in this stabilizer and fix a unitary
semisimple character $u_s$ up to $W_f$.

The first reduction with diagram automorphisms retains
\[
 E=\operatorname{Stab}_{W_f\rtimes(Q^{\rm vert})_\tau}(u_s).
\]
After separating the nonzero roots, let
$E^+=\{e\in E:e((R_u^\times)^+)=(R_u^\times)^+\}$. Then
\[
 E=W_u^\times\rtimes E^+,
 \qquad
 1\longrightarrow F_u\longrightarrow E^+
 \longrightarrow (Q^{\rm vert})_{\tau,W_fu_s}\longrightarrow1.
\]
These identities follow by moving each positive chamber with a unique
element of $W_u^\times$. The reductions of
\cite[Theorems~2.1.2, 2.1.4]{Sol-red} give the local crossed graded
algebra with this full group, with matrix factors handled by their
reduction equivalence. The scalar source cocycle is inflated along
$E^+\to(Q^{\rm vert})_{\tau,W_fu_s}$. A finite central extension
realizes it, after the common cochain normalization just specified,
as a central-character summand, so the same reductions
apply with that extension as an ineffective diagram group.

The action of $E^+$ is based after zero separation. Choose the geometric root normalizations of \S\ref{subsec:source-local-models} on $E^+$-orbits of components. The
nested-group multiplicity-one argument above applies to
$F_u\triangleleft E^+$ and identifies the intrinsic quotient
classes with their defects. Equivalently, apply the distinguished-type
map in the same scalar-twisted $E^+$ extension at both endpoints.
It restricts to the neutral map $\Theta_A$ and preserves
the projective labels in Clifford induction. Equation~
\eqref{eq:equiv-external-equal} identifies the external factors.
The reduction equivalences transport these products back to the
original source algebra. Orbit induction over $u_s$ and $\tau$
completes the proof. Notice that $\tau$ was fixed as a complex torus
point throughout; its orbit need not preserve a chosen compact torus.
\end{proof}

\paragraph{The full source comparison.}
For the full source group the required identities are
\begin{equation}\label{eq:equiv-source-target}
 \Theta_A(gV)=g\Theta_A(V),\qquad
 [c_{1,V}]=[c_{0,\Theta_A(V)}]
       \quad\hbox{in }H^2(Q_V,\mathbb C^\times).
\end{equation}
The following construction proves these identities for a coherent choice of the
neutral map, while retaining the actual source action and products.

\subsubsection{Enlarging the source lattice to include the vertex movers}
The finite cover used here is constructed from the source action.
It is independent of the auxiliary geometric cover in \S\ref{subsec:source-local-models}. Write
$X_s=X_{\rm aff}$ and let $V_s=X_s\otimes_{\mathbb Z}\mathbb R$.
The action of $Q$ on the semisimple quotient apartment is by affine
automorphisms normalizing $W_{\rm aff}$ and preserving its base alcove.
The action can have a kernel, which is retained throughout.

\paragraph{A compatible choice of the component vertices.}
Choose the special vertices used in the Bernstein presentations as
follows. In each $Q$-orbit of irreducible affine components, choose
one component and a special vertex compatible with its auxiliary
Bernstein presentation in Proposition~\ref{prop:local-auxiliary}.
Choose transporters to the other components and transport that vertex. All the constructions below use
it, including the finite positive systems. Transport the local graded recognitions in \S\ref{subsec:source-local-models} and connected shifts by these same
parameter-preserving component isomorphisms. Thus the geometric
input on each representative component supplies it on the whole
orbit. These component presentations are fixed throughout
the local reductions and the construction of the bijection.

Let $R$ and $\Delta$ be the resulting source finite root system and
base, with $\mathbb ZR=X_s$, and put
\[
 P=\{v\in V_s:\langle v,\alpha^\vee\rangle\in\mathbb Z
                                      \text{ for every }\alpha\in R\}.
\]
This is the weight lattice in the source Bernstein convention.
The alcove stabilizer in $P\rtimes W_f$ is denoted $\Omega_P$;
it is canonically isomorphic to the finite abelian group $P/X_s$.

\begin{proposition}[A parameter-compatible source enlargement]
\label{prop:equiv-enlargement}
There are a homomorphism
$\delta:Q\to\operatorname{Aut}(R,\Delta)$ and elements
$\ell_g\in\Omega_P$ such that the affine action of $g$ is
$\ell_g\delta_g$ and
\begin{equation}\label{eq:equiv-lengthzero-cocycle}
 \ell_1=1,\qquad \ell_g\delta_g(\ell_h)=\ell_{gh}.
\end{equation}
Both $\delta_g$ and $\ell_g$ preserve the source affine parameters.
Let $S$ be the subgroup of $\Omega_P$ generated by the $\ell_g$,
and let $L^\sharp$ be its inverse image in $P$ under
$P\to P/X_s\simeq\Omega_P$. Then
\[
 X_s\subset L^\sharp\subset P,\qquad L^\sharp/X_s\simeq S,
 \qquad \delta_g(L^\sharp)=L^\sharp.
\]
For every $\epsilon$, the algebra
\begin{equation}\label{eq:equiv-enlarged-algebra}
 H_\epsilon^\sharp
   =\mathcal H(W_{\rm aff},q^\epsilon)\rtimes S
\end{equation}
is the affine Hecke algebra for the datum
$(L^\sharp,R,(L^\sharp)^*,R^\vee,\Delta)$.
Its length-zero units $n_s$ satisfy $n_sn_t=n_{st}$.
Its inclusion of the original affine algebra preserves the finite
Hecke generators and the Bernstein functions $\vartheta_x$ for
$x\in X_s$.
\end{proposition}
\begin{proof}
The derivative $D_g$ preserves the translation lattice $X_s$ and
normalizes the finite reflection group. It therefore preserves $R$:
the roots are the primitive lattice vectors on its reflection-normal
lines. Write $D_g=w_g\delta_g$, where $w_g\in W_f$ and
$\delta_g$ preserves $\Delta$. If $v_g=g(0)$, conjugation of a
finite reflection shows that
$(1-s_\alpha)v_g\in X_s$ for every $\alpha\in R$.
Since $\alpha$ is primitive in $X_s$, this gives
$\langle v_g,\alpha^\vee\rangle\in\mathbb Z$, so $v_g\in P$.
Thus $\ell_g=t_{v_g}w_g\in P\rtimes W_f$. It preserves the
alcove because $g$ and $\delta_g$ do. Uniqueness in the semidirect
product $\Omega_P\rtimes\operatorname{Aut}(R,\Delta)$ proves
\eqref{eq:equiv-lengthzero-cocycle} and that $\delta$ is a homomorphism.

It remains to check the parameters separately; preservation by their
product is not enough. On one irreducible component, a nontrivial
finite based root automorphism occurs only in the simply laced
types $A_n$ ($n\geq2$), $D_n$ and $E_6$. All simple affine
reflections in each of these components are conjugate in its affine
Weyl group, so their parameters are equal. In the other irreducible
types the finite based automorphism is trivial. It follows that the
finite based part of a parameter-preserving alcove automorphism of
one component preserves its parameters. For components permuted by
$Q$, apply this observation after the chosen vertex transporters.
Those transporters identify the labelled finite and affine
presentations. Consequently every $\delta_g$ preserves the complete
affine parameter function. So does $\ell_g=g\delta_g^{-1}$.

Equation~\eqref{eq:equiv-lengthzero-cocycle} gives
$\delta_g(\ell_h)=\ell_g^{-1}\ell_{gh}$, hence $S$ is
$\delta(Q)$-stable. The Weyl group acts trivially on $P/X_s$, so
$L^\sharp$ is also $W_f$-stable. The resulting root datum is
integral. Every element of $S$ preserves the parameters, so its
Iwahori--Matsumoto crossed presentation is the affine Hecke algebra
asserted in \eqref{eq:equiv-enlarged-algebra}
\cite[Theorems~1.8 and 1.10]{Sol-survey}. Dominant translations
$x\in X_s$ have the same normalized support elements in both
presentations; their Bernstein functions therefore agree. These
facts hold for $q^\epsilon$, including $\epsilon=0$.
\end{proof}

\paragraph{Component vertices and parameter preservation.}
Two components can have affine parameter lists which agree only
after an alcove rotation. For unrelated initial special vertices,
their finite based exchange need not preserve those lists.
Transporting the component vertices first removes exactly this
problem. Nontrivial vertex movers within a component remain, and
are the elements $\ell_g$ treated below.

\subsubsection{An equivariant full corner and the mixed cocycle}
Put
\[
 K^\sharp=\operatorname{Hom}(S,\mathbb C^\times),\qquad
 \widetilde A_\epsilon=\mathbb C[\Lambda_c]\otimes H_\epsilon^\sharp,
 \qquad \mathscr D_\epsilon=\widetilde A_\epsilon\rtimes K^\sharp.
\]
The action of $k\in K^\sharp$ fixes the original affine algebra and
$\mathbb C[\Lambda_c]$, and sends $n_s$ to $k(s)n_s$.
Equivalently it multiplies $\vartheta_x$, $x\in L^\sharp$, by the
character of $L^\sharp/X_s$ determined by $k$. In particular it is
trivial on all roots. This group is distinct from both $K_X$ and
$K_{\rm g}$.

Set
\[
 e=\frac1{|K^\sharp|}\sum_{k\in K^\sharp}[k].
\]
Extend the source action to $\widetilde A_\epsilon$ by
\begin{equation}\label{eq:equiv-lifted-action}
 \widehat\beta_g=\operatorname{Ad}(n_{\ell_g})\delta_g
       \quad\hbox{on }H_\epsilon^\sharp,\qquad
 \widehat\beta_g(z_c)=\beta_g(z_c).
\end{equation}
Here and below $\delta_g$ on the toral factor means this specified
source action when the two factors are considered together. On
$K^\sharp$ write $gk=k\circ\delta_g^{-1}$.

\begin{proposition}[Equivariant Morita descent]\label{prop:equiv-corner}
The action \eqref{eq:equiv-lifted-action}, together with
$\widehat\beta_g([k])=[gk]$, is a strict $Q$-action on
$\mathscr D_\epsilon$. It fixes $e$ and restricts on the full corner
to the original source action:
\begin{equation}\label{eq:equiv-full-corner}
 e\mathscr D_\epsilon e\simeq A_\epsilon,\qquad
 e(\mathscr D_\epsilon\rtimes_{\widehat\beta,a}Q)e
                           \simeq B_\epsilon.
\end{equation}
In particular $M\mapsto eM$ is an equivariant Morita equivalence.
It identifies stabilizers and preserves intrinsic intertwining
cocycles, with equality of representatives when the intertwiners
are restricted to $eM$.
\end{proposition}
\begin{proof}
The $S$-grading of $H_\epsilon^\sharp$ has the invertible
homogeneous elements $n_s$. Thus
\[
 e\widetilde A_\epsilon e=A_\epsilon e,\qquad
 \sum_{s\in S}n_s e n_s^{-1}=1.
\]
More explicitly the elements $n_s e n_t^{-1}$ are matrix units:
the middle factor $e n_t^{-1}n_{s'}e$ is zero unless $t=s'$.
They identify $\mathscr D_\epsilon$ with $M_{|S|}(A_\epsilon)$
and prove that $e$ is full.

Equation~\eqref{eq:equiv-lengthzero-cocycle} and the exact products
of the $n_s$ show that $\widehat\beta$ is strict. Scalar factors in
$k(n_{\ell_g})$ cancel in inner conjugation, so
$\widehat\beta_g k\widehat\beta_g^{-1}=gk$ on
$\widetilde A_\epsilon$. This proves the assertion on the crossed
product. The average $e$ is fixed. On the original affine generators,
conjugation by $n_{\ell_g}$ followed by $\delta_g$ is precisely
the original affine-wall permutation; on $z_c$ it is
\eqref{eq:local-beta}. The defects $a(g,h)$ remain in
$\mathbb C[\Lambda_c]^\times$. Hence taking the corner gives
both isomorphisms in \eqref{eq:equiv-full-corner}.

If $T_g$ intertwines a simple $M$ with its $g$-translate, then
$T_ge=eT_g$, since $\widehat\beta_g(e)=e$. Its restriction is a
nonzero intertwiner on the simple $eM$, and
$(T_gT_h)|_{eM}=c_M(g,h)T_{gh}|_{eM}$.
The inverse Morita equivalence identifies the stabilizers, proving
the cocycle assertion. This retains the projective action of the
matrix factor explicitly.
\end{proof}

For the intrinsic comparison, temporarily form the untwisted
extension $\mathscr D_\epsilon\rtimes_{\widehat\beta}Q$ and write
$U_g$ for its group units. This does not change or renormalize the
source defects; it isolates the intrinsic action whose classes are
needed in \eqref{eq:equiv-reduced-obstruction}. Put
\[
 V_g=n_{\ell_g}^{-1}U_g,\qquad
 \eta_g=\operatorname{Ad}(V_g)|_{\mathscr D_\epsilon}.
\]
An exact computation gives
\begin{align}
 V_gV_h&=V_{gh},&
 \eta_g|_{H_\epsilon^\sharp}&=\delta_g,&
 \eta_g([k])&=(gk)(\ell_g)[gk].
 \label{eq:equiv-untwist}
\end{align}
Thus $\eta$ is another strict $Q$-action. It differs from
$\widehat\beta$ by the displayed inner units, but includes a
nontrivial character factor on $K^\sharp$.

Let $\mathcal E=K^\sharp\rtimes_\delta Q$, with multiplication
$(k,g)(l,h)=(k(gl),gh)$. Its action on $\widetilde A_\epsilon$ is
the Bernstein action $k\delta_g$, including the specified toral
action. With group units $[k]V_g$ the multiplier is
\begin{equation}\label{eq:equiv-mixed}
 \xi((k,g),(l,h))=(gl)(\ell_g).
\end{equation}
It is root-of-unity-valued, independent of $\epsilon$, and trivial
on $K^\sharp\times K^\sharp$ and on $Q\times Q$. It need not
be cohomologically trivial on $\mathcal E$. Its cocycle identity is
exactly \eqref{eq:equiv-lengthzero-cocycle}. We have
\begin{equation}\label{eq:equiv-mixed-presentation}
 \mathscr D_\epsilon\rtimes_{\widehat\beta}Q
      \simeq\widetilde A_\epsilon\rtimes_\xi\mathcal E.
\end{equation}
If the actual source products $U_gU_h=a(g,h)U_{gh}$ are used
instead, the multiplier on the right is
$\xi((k,g),(l,h))a(g,h)$. In particular the mixed factor is not
discarded when the external source defects are restored.

\subsubsection{Finite character translations in the local reductions}
The action in \eqref{eq:equiv-mixed-presentation} preserves the
Bernstein presentation, although it can translate its torus by a
finite central character. We spell out why the reduction argument
and the nested-group cocycle comparison extend to these actions.

\begin{proposition}[Reduction with finite central-character translations]
\label{prop:equiv-character-reduction}
Let a finite group $\mathcal E'$ act on an affine Hecke algebra
with positive parameters by based root-datum automorphisms and
finite torus characters trivial on the roots. Ineffective elements
are allowed. At a unitary torus point $\widetilde u$, put
\[
 E_{\widetilde u}=(W_f\rtimes\mathcal E')_{\widetilde u},
 \qquad
 R_{\widetilde u}=\{\alpha:s_\alpha\widetilde u=\widetilde u\}.
\]
The first and second reductions give the real-stratum local model
with connected graded root system $R_{\widetilde u}$ and the
full stabilizer complement in $E_{\widetilde u}$. After zero
separation the model has the form
\begin{equation}\label{eq:equiv-character-local}
 C_{\epsilon,\widetilde u}\rtimes
       E_{\widetilde u}^{+},\qquad
 E_{\widetilde u}^{+}
 =\{e'\in E_{\widetilde u}:
       D_{e'}((R_{\widetilde u}^{\times})^+)
                         =(R_{\widetilde u}^{\times})^+\}.
\end{equation}
Here $D_{e'}$ is the derivative, and all finite elements remain
labelled even when their derivatives coincide. The reductions
intertwine actions and products, including products of
intertwiners and the matrix-corner Morita equivalences.
For a scalar twist inflated from $\mathcal E'$, its restriction
to the displayed complement is retained. The same assertions
hold simultaneously for a normal subgroup of $\mathcal E'$.
\end{proposition}
\begin{proof}
The Bernstein formulas show that a character trivial on $R$
fixes each normalized rational intertwiner $\iota_w^\circ$.
A based automorphism sends $\iota_w^\circ$ to
$\iota_{\delta(w)}^\circ$. Consequently the units
$\iota_w^\circ[\gamma]$ multiply according to
$W_f\rtimes\mathcal E'$ in the meromorphic localization.
In the twisted case they have exactly the scalar multiplier
inflated from $\mathcal E'$. This follows directly from the
intertwiner formula in \cite[Proposition~1.12]{Sol-survey}.

Choose the first-reduction neighborhoods around the finite orbit
of $\widetilde u$ and its real translates. The characteristic
idempotents and the intertwiner maps between their summands
transform under these same group units. The analytic invertibility
checks in \cite[Theorem~2.1.2]{Sol-red} involve only the root
functions and their parameters. A character trivial on $R$
does not change any of them. Hence the same corner construction
works, retaining $E_{\widetilde u}$. Its reflection subgroup
$W(R_{\widetilde u})$ is normal; taking positive chambers gives
its stabilizer complement.

Every $e'\in E_{\widetilde u}$ satisfies
\[
 e'(\widetilde u\exp\nu)=
                      \widetilde u\exp(D_{e'}\nu).
\]
The finite translation has thus disappeared from the local
logarithmic action, not from the group. Under the meromorphic
second-reduction map, functions pull back through
$\exp_{\widetilde u}$ and normalized affine intertwiners go to
the corresponding graded intertwiners. These assignments preserve
the displayed products. The root-wise analytic argument and
parameter formula in \cite[Theorem~2.1.4]{Sol-red} are unchanged,
so this map restricts to the required analytic and completed
algebras. Its real-stratum form follows as in that reference's
Corollary~2.1.5.

Finally apply the zero-separation proof from \S\ref{subsec:source-local-models}. The nonzero
reflection subgroup is normal, and the chamber stabilizer in
\eqref{eq:equiv-character-local} contains the zero reflections and
all other residual elements. The fixed sign character is invariant.
Every construction used the same functions, intertwiners and
idempotents; it can therefore be made for a normal subgroup
inside the full construction. Equivalently one can first take
that subgroup's orbit corner and then the larger orbit corner.
The intervening Morita bimodules have their actual multiplication
maps, so induction in stages retains the weak quotient defects.
For a scalar twist, use a finite central extension after one
common cochain normalization, retaining its prescribed central
character at every step.
\end{proof}

\begin{proposition}[Parameter removal before Morita descent]
\label{prop:equiv-upstairs}
For the action $\eta$ in \eqref{eq:equiv-untwist}, there is a
$Q$-equivariant fixed-choice bijection
\[
 \Theta_{\mathscr D}:\Irr(\mathscr D_1)
                           \xrightarrow{\sim}\Irr(\mathscr D_0)
\]
which preserves the independent toral character and the intrinsic
stabilizer classes for $\eta$. Its local maps use the connected
shift and the sign and involution conventions of \S\S\ref{subsec:source-local-models}--\ref{subsec:source-equivariance}.
\end{proposition}
\begin{proof}
First fix a toral character $\tau$ and its stabilizer $Q_\tau$;
the other toral fibres are handled by orbit transport. This is
necessary because the source action on $T_c$ need not preserve
its compact torus. All the finite translations on the semisimple
cover are unitary.

Choose a finite cyclic group $Z_m$ containing the values of
$\xi$, and form the extension of
$K^\sharp\rtimes Q_\tau$ defined by \eqref{eq:equiv-mixed}.
Its inverse image of $K^\sharp$ is $Z_m\times K^\sharp$,
and the $Q_\tau$ units still multiply strictly. The central
character $Z_m\hookrightarrow\mathbb C^\times$ recovers
$H_\epsilon^\sharp\rtimes K^\sharp$ with precisely the quotient
action $\eta$. In particular the character factor in
$\eta_g([k])$ is part of this extension.

Apply Proposition~\ref{prop:equiv-character-reduction} to this
finite extension and to its normal subgroup
$Z_m\times K^\sharp$. On a chosen unitary orbit the local
complements form a normal-subgroup pair. They act by based
linear automorphisms on the same connected nonzero graded
factor. The geometric input for that factor is supplied by
Proposition~\ref{prop:local-geometric} and the lattice comparison below. Choose an invariant
complement to the root span by averaging an inner product
over the finite stabilizer, and choose the geometric root
normalizations on its orbits of components. As in the
$F_u$ argument, these choices preserve the full based action.
Thus the nested-group
multiplicity-one argument following
Proposition~\ref{prop:equiv-local-full} applies. On real tempered
modules, restriction of intertwiners to the distinguished
finite type preserves their scalar products and every weak
quotient defect. On all real simples, induction from the
Langlands datum and passage to its unique simple quotient
give the same comparison. This is
\cite[Theorem~6.2, Lemma~6.6 and Proposition~6.10]{Sol-survey}
with the prescribed central-character summand retained.
The two endpoint involutions in \eqref{eq:equiv-im} commute
with these based actions.

More explicitly, denote the local complements in this finite
central extension by $E_0^+\triangleleft E^+$. Choose lifts
$r_g$ of their quotient. Their defects are the actual group
elements $d(g,h)=r_gr_hr_{gh}^{-1}\in E_0^+$, and the local
units satisfy
\[
 U_{r_g}U_{r_h}=U_{d(g,h)}U_{r_{gh}}
\]
at both endpoints. In a stabilizer representation the weak
intertwiner equation is therefore
\[
 T_gT_h=\gamma(g,h)\pi(U_{d(g,h)})T_{gh}.
\]
The distinguished finite type retains the operator
$\pi(U_{d(g,h)})$, including its central character, so
restriction preserves $\gamma$, not merely the permutation
of labels. The quotient and its defects are independent of
the moving real central character.

The distinguished type has the same $Z_m$ character as the
module, so the construction stays in the required summand.
The subgroup and quotient comparisons therefore give an
equivariant bijection on the $\mathscr D_\epsilon$ simples,
not just a bijection of the full crossed-product spectra.
Their effective quotient classes agree. Globally this
quotient action is the strict action $\eta$, with units
$V_gV_h=V_{gh}$; its effective class is $c_M^{-1}$.
Consequently the intrinsic $c_M$ classes agree as asserted.
The simultaneous reduction equivalences transport this
statement back from local corners without discarding their
inner defects. Indeed, with their multiplication maps retained,
these equivalences preserve composition of the intertwining
isomorphisms. They thus identify the central extensions of a
simple's stabilizer by $\mathbb C^\times$ which define the
effective Clifford classes.

Use one set of unitary orbit representatives and the inverse
construction at zero to assemble the bijection. Choose the
finite-type normalizations under the full stabilizer before
passing to its subgroup. Transport them to the remaining
$Q_\tau$-orbit representatives and then to the other toral
fibres. The orbit and cochain argument of \S\ref{subsec:source-assembly} gives consistent
maps and inverses. The toral character is unchanged.
\end{proof}

\subsubsection{Comparison with the original neutral shift}
The source lattice enlargement retains the roots $R$. The
following local comparison identifies its connected geometric
factor and finite stabilizers with those of the original neutral
model.

\begin{proposition}[The covering corner has the original local model]
\label{prop:equiv-local-descent}
Let $\widetilde u\in\operatorname{Hom}(L^\sharp,S^1)$ restrict
to $u_s\in\operatorname{Hom}(X_s,S^1)$. Under the differential
identification of their real polynomial spaces, the nonzero
graded roots and parameters at the two points are identical.
Moreover projection onto $W_f$ gives an isomorphism
\begin{equation}\label{eq:equiv-stabilizer-cover}
 (W_f\rtimes K^\sharp)_{\widetilde u}
                         \xrightarrow{\sim}(W_f)_{u_s}.
\end{equation}
After zero separation its complement is the original $F_u$.
The Morita equivalence $M\mapsto eM$ identifies these local
models, with their root intertwiners and residual products.
Consequently the descended map from $\Theta_{\mathscr D}$
is a coherent choice of the neutral map $\Theta_A$ in
Proposition~\ref{prop:assembly-neutral}.
\end{proposition}
\begin{proof}
The torus cover has kernel $K^\sharp$, acting freely on its
fibres. For $w\in(W_f)_{u_s}$ there is a unique $k\in K^\sharp$
such that $(w,k)\widetilde u=\widetilde u$. This proves
\eqref{eq:equiv-stabilizer-cover}, including its multiplication.

The common Bernstein embedding was proved in
Proposition~\ref{prop:equiv-enlargement}. A reflection that fixes
$\widetilde u$ fixes $u_s$. The only reflection that can be
lost on lifting has $\alpha(u_s)=-1$, with its coroot
divisible by $2$ on $X_s$ but primitive on $L^\sharp$.
In that case the enlarged Bernstein relation has equal
labels $Q_\alpha^*=Q_\alpha$. Comparison with the original
relation, or its common restriction to
$\mathbb C[X_s]$, forces the same equality in the source.
Thus $k_u(\alpha)=0$. Retained reflections have identical
root values, labels and graded parameters. This is precisely
the same-root case of Propositions~\ref{prop:local-coefficients} and~\ref{prop:local-lattice-recognition}. It proves equality of the nonzero connected
data and supplies their geometric recognition from \S\ref{subsec:source-local-models}.
The independent toral directions are unchanged.

For completeness, the algebra comparison is not merely a
comparison of stabilizer cardinalities. The covering algebra
with $K^\sharp$ is the finite Galois descent presentation of
the original affine algebra; its averaging corner is
\eqref{eq:equiv-full-corner}. In the first-reduction orbit
corners, take one sheet over $u_s$. A lift of a residual
$w$ is represented by the actual unit
$\iota_w^\circ[k]$ with the unique $k$ in
\eqref{eq:equiv-stabilizer-cover}. Multiplication is the
group multiplication because $K^\sharp$ fixes all
$\iota_w^\circ$. Descending through the averaging corner
gives the original $\iota_w^\circ$. Invariant Bernstein
functions descend to the same functions on $T_{X_s}$.
These assertions can also be checked on the matrix units
$n_s e n_t^{-1}$ and the sheet idempotents: multiplication
followed by the corner compression is their canonical Morita
map, not a newly chosen projective action.

The second reduction takes these same functions and
intertwiners to the original graded ones. Hence the residual
complement in \eqref{eq:equiv-stabilizer-cover} is identified
with $F_u$, and lost reflections remain in that complement
with zero parameter. On the connected nonzero factor use
the same root normalizations, fixed sign map and two endpoint
involutions as in \S\ref{subsec:source-equivariance}. Extend the involutions by fixing the
chosen finite complement at both endpoints. The connected
shift, its finite-type extension and its inverse are
therefore exactly the local maps used in \S\ref{subsec:source-assembly}, for these
compatible choices.
\end{proof}

\begin{proposition}[The full source action and its intrinsic classes]
\label{prop:equiv-full-q}
Under the depth-zero hypotheses in \S\ref{subsec:affine-input},
there is a coherent choice of the neutral parameter-removal
bijection such that, for every $g\in Q$ and every neutral simple $V$,
\begin{equation}\label{eq:equiv-full-q}
 \Theta_A(gV)=g\Theta_A(V),\qquad
 [c_{1,V}]=[c_{0,\Theta_A(V)}]
          \quad\hbox{in }H^2(Q_V,\mathbb C^\times).
\end{equation}
This includes the elements moving the special vertex and all
ineffective subgroups. The effective source classes also agree,
with the external factors in \eqref{eq:equiv-external} retained.
\end{proposition}
\begin{proof}
Apply Proposition~\ref{prop:equiv-upstairs} and descend through
Proposition~\ref{prop:equiv-corner}. We verify the intervening
inner change of action explicitly. On a simple
$\mathscr D_\epsilon$-module $M$, choose intertwiners $T_g^\eta$
for $\eta$ on its stabilizer and put
\[
 T_g^{\widehat\beta}=\pi_M(n_{\ell_g})T_g^\eta.
\]
Since $\eta_g(n_{\ell_h})=\delta_g(n_{\ell_h})$, equation
\eqref{eq:equiv-lengthzero-cocycle} gives
\[
 n_{\ell_g}\eta_g(n_{\ell_h})=n_{\ell_{gh}},\qquad
 T_g^{\widehat\beta}T_h^{\widehat\beta}
                      =c_M^\eta(g,h)T_{gh}^{\widehat\beta}.
\]
The intrinsic cocycle representatives for these choices are
therefore equal. The stabilizers are equal because the two
actions differ by inner automorphisms. Restrict
$T_g^{\widehat\beta}$ to $eM$; the full-corner proposition
gives the original source intertwiners and exactly the same
scalar products. This proves both assertions in
\eqref{eq:equiv-full-q}. Proposition~\ref{prop:equiv-local-descent}
identifies the descended map with the required choice of
$\Theta_A$. Any cochain comparing the upstairs endpoint
products consequently compares the downstairs products as
well.

Finally $\Theta_A$ preserves $\tau$. Thus
\eqref{eq:equiv-external-equal} identifies the external defect
values exactly, and \eqref{eq:equiv-reduced-obstruction}
proves equality of the effective classes. No assertion that
the mixed cocycle or the source multiplier is trivial has
been used.
\end{proof}

\paragraph{Rank-one check of the mixed factor.}
For the equal-parameter $A_1$ affine algebra with its wall-swapping
$Q=C_2$ action, take
$L^\sharp=\mathbb Z(\alpha/2)$, $S=C_2$ and $K^\sharp=C_2$.
Let $\ell$ and $k$ be the nontrivial elements of $S$ and
$K^\sharp$, so $k(\ell)=-1$. Then $\delta=1$ and the new
residual unit satisfies
\[
 V[k]=-[k]V,\qquad V^2=1.
\]
Thus the mixed cocycle on $K^\sharp\times Q$ is nontrivial:
its group units anticommute. The original unit is $U=n_\ell V$;
it commutes with $e=(1+[k])/2$ because both $n_\ell$ and $V$
anticommute with $[k]$. For the $\mathrm{GL}_2$ source example,
restoring the external defect changes the square to
$V^2=U^2=z_{e_1+e_2}$, with the anticommutation unchanged.
This value is preserved by the shift. At $\alpha(u_s)=-1$
the lifts have values $\pm i$ on $\alpha/2$; the missing
reflection has zero parameter and is restored by its combined
$(w,k)$ stabilizer. Both the matrix corner and the exceptional
unitary stratum are therefore visible in the construction.

\subsection{Global assembly on simple classes}\label{subsec:source-assembly}
\subsubsection{Polar strata and compatible logarithms}
For $t\in T_X$ write uniquely
\[
 t=u\exp\nu,\qquad
 u(x)=\frac{t(x)}{|t(x)|},\quad
 \nu(x)=\log|t(x)|,
 \quad u\in T_{X,\rm un},\quad\nu\in\mathfrak a.
\]
Choose a set $\mathcal U$ of representatives of
$T_{X,\rm un}/W_f$. Each simple $A_\epsilon$-module has exactly one
unitary $W_f$-orbit in its Bernstein weights. Thus the real-stratum
equivalences~\eqref{eq:local-unitary-equivalence} give
\begin{equation}\label{eq:assembly-strata}
 \Irr(A_\epsilon)
   =\coprod_{u\in\mathcal U}
       \Irr_{W_fu\exp\mathfrak a}(A_\epsilon)
   \simeq\coprod_{u\in\mathcal U}\Irr_{\rm real}(G_{\epsilon,u}).
\end{equation}

Here is a compatible choice of the analytic charts behind these
equivalences. For fixed $u$, choose a sufficiently small $W_u$-invariant
convex ball $B_u\subset\mathfrak a$ around zero. Require that
\begin{itemize}
\item $(B_u-B_u)$ contain no nonzero point of
$2\pi\Hom(X,\mathbb Z)$;
\item all root imaginary parts in $B_u$ be sufficiently small to
satisfy the logarithmic restrictions in the second reduction;
\item the finitely many distinct unitary translates by $W_f/W_u$
have disjoint imaginary neighborhoods, and the root phases excluded
in the first reduction stay excluded.
\end{itemize}
Such a ball exists: there are only finitely many root phases and
unitary translates to separate. On the tube
$\mathfrak a+iB_u$ the map $\lambda\mapsto u\exp\lambda$ is
injective. Use its single logarithm on every smaller neighborhood of
a real point. Normalized rational intertwiners and their analytic
extensions in the reduction theorems define one set of algebra maps;
restrict those maps when a neighborhood is shrunk. Consequently two
charts describing the same real-stratum module use the same logarithm
and the same reduction equivalence on their intersection. For a
different representative $wu$, use the transported charts and
equivalences. Matrix corners use the corresponding matrix-unit Morita
bimodules; they do not introduce independent projective labels.

These analytic neighborhoods implement a map on \emph{simple classes
stratified by the exact unitary orbit}. Points of a neighboring unitary
stratum are processed with their own polar part. No assertion of a
holomorphic or continuous parameter-removal map across all those
strata is being made.

\begin{proposition}[Global neutral parameter removal]
\label{prop:assembly-neutral}
Fix the reductions and the compatible finite-type choices in \S\ref{subsec:source-equivariance}.
The maps $\Theta_u$ assemble to a bijection
\begin{equation}\label{eq:assembly-neutral}
 \Theta_A:\Irr(A_1)\xrightarrow{\sim}\Irr(A_0).
\end{equation}
It preserves the unitary $W_f$-orbit and the independent toral
character. Its inverse is assembled from $\Theta_u^{-1}$. The local
descriptions agree on their common simple classes when transported
through the specified reduction equivalences.
\end{proposition}
\begin{proof}
For each summand in \eqref{eq:assembly-strata}, compose the equivalence
at $1$, the bijection~\eqref{eq:equiv-local-full}, and the inverse
equivalence at $0$. The intermediate map and its inverse preserve
real central characters and the independent toral coordinate. They
therefore stay in the same summand. Disjointness proves injectivity;
the inverse construction proves surjectivity. The choices of charts
just described give the same equivalence on a common stratum module.
If an orbit representative or a matrix corner is changed, transport
the map through that change rather than reconstructing its labels
independently. The resulting local maps and inverses agree by
composition. This proves the asserted compatibility. It makes no
claim that a fixed original semisimple central character is preserved.
\end{proof}

\paragraph{Coordinates on the source torsor.}
The endpoint comparison is the already specified
$r_\pi:T^\fs\to T_X$, with fibres the $K_X$-torsors in
\eqref{eq:local-isogeny}. Charts on $T^\fs$ are inverse images of
the above charts, with a local lift chosen on each component.
On an overlap two lifts differ by the recorded element of $K_X$.
For the auxiliary geometric cover the same instruction uses
$K_{\rm g}$, separately. The source-action cover
$T_{L^\sharp\oplus\Lambda_c}\to T_X$ uses $K^\sharp$ and the
full-corner descent of Proposition~\ref{prop:equiv-corner}.
Its corner is the original $A_\epsilon$, and its extended corner is
the original $B_\epsilon$; in particular at zero it retains the origin normalization of \S\ref{subsec:Levi-center-normalization} and the same $r_\pi$.
These lift changes are coordinate transitions;
they are not identifications that discard simple modules. The
endpoint Clifford construction and Theorem~\ref{thm:source-support-realization} retain their
finite and projective information. The affine comparison is not
being used to replace $X$, $\Lambda$ or the auxiliary lattice by one
another.

\subsubsection{One global choice of projective-label transport}
The following proposition gives the global assembly, with its source
comparison supplied by Proposition~\ref{prop:equiv-full-q}. It explains exactly
which coherence is needed, and why arbitrary pointwise choices on
overlapping charts should not be made.

\begin{proposition}[Global Clifford assembly with central defects]
\label{prop:assembly-global}
Use the neutral bijection $\Theta_A$ and source comparison of
Proposition~\ref{prop:equiv-full-q}. There is a fixed-choice bijection
\begin{equation}\label{eq:assembly-global}
 \Theta_B:\Irr(B_1)\xrightarrow{\sim}\Irr(B_0)
\end{equation}
whose local Clifford descriptions and inverses agree on common
simple classes. All comparison cochains can be chosen on one set
of global neutral orbit representatives; no continuous choice of a
cochain on the torus is required.
\end{proposition}
\begin{proof}
Choose a set $\mathcal V$ of representatives of
$\Irr(A_1)/Q$ and choose actual modules and normalized intertwiners
there. Use $\Theta_A(\mathcal V)$ at zero; equivariance makes it a
set of representatives of $\Irr(A_0)/Q$ with the same stabilizers.
For $V\in\mathcal V$, the external factors are identical by
\eqref{eq:equiv-external-equal}. The intrinsic comparison in
\eqref{eq:equiv-source-target} supplies a normalized $b_V$ with
\[
 \kappa^{\rm eff}_{1,V}
              =\delta b_V\,\kappa^{\rm eff}_{0,\Theta_A(V)}.
\]
If $R_gR_h=\kappa^{\rm eff}_{1,V}(g,h)R_{gh}$, put
$R_g^0=b_V(g)^{-1}R_g$. Define
\begin{equation}\label{eq:assembly-induced}
 \Theta_B\!\left(\operatorname{Ind}_{A_1\rtimes Q_V}^{B_1}
                                  (V\otimes R)\right)
   =\operatorname{Ind}_{A_0\rtimes Q_V}^{B_0}
                                  (\Theta_A(V)\otimes R^0).
\end{equation}
The notation on both sides includes the restrictions of the actual
defects $a(g,h)$. With $u_g$ acting as $T_g\otimes R_g$ their
products are exactly $a(g,h)$, by
\eqref{eq:equiv-reduced-obstruction}. Clifford theory classifies all
simples of a finite crossed extension by these induced modules,
uniquely up to orbit and projective-label equivalence. It applies
here because the algebras are finite over their indicated affine
centers, so the relevant simple modules are finite dimensional.
This proves that \eqref{eq:assembly-induced} is a well-defined bijection.
Its inverse uses $\Theta_A^{-1}$ and the inverse cochain change.

For any other representative of the orbit of $V$, choose a transporter
from $V$, and transport both endpoint modules and their label
coordinates through the actual units $u_g$. A different transporter
differs by an element of $Q_V$; its action is already implemented in
the stabilizer module. Induction therefore gives an isomorphic
module. Associativity of the source crossed product is exactly the
coherence relation for these transports.

On a logarithmic chart, completion, or matrix corner, define the
local map by this same global representative and $b_V$, followed by
the stipulated reduction equivalence. Two charts that describe the
same simple consequently yield the same induced endpoint module.
The argument applies to the inverse as well. Coverage follows from
\eqref{eq:assembly-strata} and Clifford induction; these same two
classifications account for every duplicate description. There is
no extra choice of a square root or a cochain independently on each
chart. A different global choice can change label coordinates, as
described by the cochain rescaling above, while the displayed choices
define the stated bijection.
\end{proof}

\begin{lemma}[Transfer of neutral invariants and projective labels]
\label{lem:neutral-invariants}
Use the equivariant bijection $\Theta_A$ of Propositions~\ref{prop:equiv-full-q}--\ref{prop:assembly-neutral} and the effective cocycle comparison there. The Clifford
bijection $\Theta_B$ of Proposition~\ref{prop:assembly-global} satisfies
\begin{equation}\label{eq:neutral-restriction}
 [\Res_{A_0}\Theta_B(M)]
 =\sum_V m_V[\Theta_A(V)]
 \quad\text{if}\quad
 [\Res_{A_1}M]=\sum_Vm_V[V].
\end{equation}
Both restrictions are semisimple. In particular,
\begin{equation}\label{eq:neutral-invariant-preservation}
 \mathcal P_0(\Theta_B(M))=\mathcal P_1(M).
\end{equation}
If $Q_V$ is the stabilizer of $V$ and
$\kappa^{\mathrm{eff}}_{1,V}=(\delta b_V)
\kappa^{\mathrm{eff}}_{0,\Theta_A(V)}$, the multiplicity operators
are transported by $R_g\mapsto b_V(g)^{-1}R_g$.
\end{lemma}
\begin{proof}
Every simple under consideration is finite dimensional. Choose a
simple neutral submodule $V\subset M$. The sum of its translates by
the invertible homogeneous units is a nonzero $B_1$-submodule of
$M$, hence equals $M$. Each translate is simple over $A_1$, so this
sum is semisimple and its constituents form one $Q$-orbit.

Let $Q_V=\Stab_Q[V]$. Clifford theory writes
\[
 M\simeq\Ind_{B_{1,V}}^{B_1}(V\otimes E),\qquad
 B_{1,V}=\bigoplus_{g\in Q_V}A_1u_g,
\]
where $E$ is a simple projective $Q_V$-module with the effective
multiplier. Its restriction to $A_1$ contains each distinct $gV$
with multiplicity $\dim E$. Equivariance of $\Theta_A$ identifies
the stabilizers and the constituent orbits. Rescaling $R_g$ by
$b_V(g)^{-1}$ changes its multiplier to
$\kappa^{\mathrm{eff}}_{0,\Theta_A(V)}$, without changing
$\dim E$. The same induction formula at zero proves
\eqref{eq:neutral-restriction}. Proposition~\ref{prop:assembly-neutral} gives
$p_0(\Theta_A(V))=p_1(V)$ for every constituent, which proves
\eqref{eq:neutral-invariant-preservation}.

The comparison uses the effective multiplier, including the evaluated
central defects $a(g,h)$ of \eqref{eq:local-defect}. These evaluations
agree at the endpoints because $a(g,h)$ lies in
$\C[\Lambda_c]^\times$ and $\Theta_A$ preserves $\tau_V$.
Thus no external defect is discarded in the rescaling.
\end{proof}

\subsection{The depth-zero parametrization}
\begin{proposition}[Parameter removal]\label{prop:bridge-stratified}
Under the depth-zero hypotheses in \S\ref{subsec:affine-input}, fix the
source datum and the coherent comparison choices of
\S\S\ref{subsec:source-local-models}--\ref{subsec:source-assembly}.
There is a bijection
\[
 \Irr(\mH_G)\xrightarrow{\sim}\Irr(B_0).
\]
It need not preserve the same individual Bernstein point.
\end{proposition}
\begin{proof}
Proposition~\ref{prop:source-family} identifies the source presentation
with $B_1$. Proposition~\ref{prop:equiv-full-q} supplies equivariance
and comparison of the effective stabilizer cocycles for the neutral
parameter-removal map. Proposition~\ref{prop:assembly-global} therefore
gives $\Irr(B_1)\simeq\Irr(B_0)$ with compatible global label
transports. Their composition is the required bijection.
\end{proof}

\begin{proof}[Proof of Theorem~\ref{thm:fiber-irr}]
Theorem~\ref{thm:source-support-realization} identifies the source endpoint, its center, and
its orbit corners, giving $\cl_0$. The type equivalence gives
$\typ_\fs$. Proposition~\ref{prop:equiv-full-q} supplies equivariance for the
full source group $Q$, including elements moving the chosen special
vertex, and compares the intrinsic and effective stabilizer classes.
Propositions~\ref{prop:assembly-neutral}--\ref{prop:assembly-global} therefore supply the bijection
$\Theta_B$. Their composition proves \eqref{eq:depth-parametrization}.

Lemma~\ref{lem:neutral-invariants} proves (i), including multiplicities.
For (ii), the local connected construction is precisely the shift
in Proposition~\ref{prop:local-connected-shift}, using
\cite[Lemma 3.6, Proposition 3.8, Theorem 3.20 and erratum]{AMS-GHA}.
In the conventions of \cite{Sol-survey}, the same map is
$\operatorname{IM}_0^*\zeta'_H\operatorname{IM}_1^*$, with the
sign and root maps already fixed in \S\ref{subsec:source-equivariance}. The local
real-central-character assertion ensures that both endpoints lie
in the same unitary stratum. An algebraic $\mathrm{SL}_2$ has
trivial image in an independent central torus, so the independent
toral character is unchanged. The finite residual extensions use
the same effective Clifford transport as in
Lemma~\ref{lem:neutral-invariants}.

For (iii), choose $b_V$ once on each global neutral $Q$-orbit.
Transport its intertwiners and labels by the actual homogeneous
units. Changing a transporter by a stabilizer element gives an
inner action on the corresponding inertia module, hence the same
induced simple. Associativity includes the defects $a(g,h)$.
The chart and corner compatibilities are those of Propositions~\ref{prop:equiv-full-q}--\ref{prop:assembly-global}; no independent cochains are chosen on
overlapping charts. This proves the asserted coherence.
Finally $\cl_0$ is constructed in the central fiber at its first
coordinate, so \eqref{eq:endpoint-character} follows from
Theorem~\ref{thm:source-support-realization}.
\end{proof}

\subsection{A finite trace criterion for projective labels}
\label{sec:pinned-traces}

In Proposition~\ref{prop:assembly-global} a choice of cochain
identifies two twisted stabilizer algebras. Given independently
prescribed pairs of simple modules, we determine when such a comparison
realizes the prescribed pairs and when they determine its action on
all simple classes.

Let $K$ be a finite group, $C_i=\C_{\kappa_i}[K]$ for $i=0,1$, and
$e_g^i e_h^i=\kappa_i(g,h)e_{gh}^i$ in normalized homogeneous bases.
A comparison over the identity of $K$ has the form
\begin{equation}\label{eq:pinned-cochain}
 a_b(e_g^1)=b(g)e_g^0,\qquad
 \kappa_1(g,h)=\frac{b(g)b(h)}{b(gh)}\,\kappa_0(g,h).
\end{equation}
These equations, with $b(1)=1$, are equivalent to multiplicativity.
They admit a solution exactly when $[\kappa_1]=[\kappa_0]$.
For left modules $(a_b)_* E=E\circ a_b^{-1}$. Replacing $b$ by
$\chi b$ for a character $\chi$ of $K$ tensors the image by
$\chi^{-1}$, and all other solutions arise this way.

For $C=\C_\kappa[K]$, write $\bar g$ for the image of $g$ in
$K_{\rm ab}$ and put
\begin{equation}\label{eq:pinned-visible}
 B(C)=\langle\bar g:\operatorname{tr}_E(e_g)\ne 0
                              \text{ for some }E\in\Irr(C)\rangle
                           \ \subseteq K_{\rm ab}.
\end{equation}
Rescaling a homogeneous basis preserves this subgroup.
\begin{lemma}\label{lem:pinned-invisible}
Tensoring by $\chi$ fixes every simple class of $C$ precisely when
$\chi|_{B(C)}=1$. If $K$ is abelian, $B(C)$ is the radical of the
commutator pairing $b_\kappa(g,h)=\kappa(g,h)/\kappa(h,g)$.
\end{lemma}
\begin{proof}
The twisted group algebra is semisimple over $\C$ and its simple
characters distinguish its simple modules. Tensoring by $\chi$
multiplies the trace at $e_g$ by $\chi(g)$, giving the first assertion.
For abelian $K$, if $g$ lies outside the radical, conjugation by some
$e_h$ multiplies $e_g$ by a nontrivial scalar. Its trace vanishes on
every module. Inside the radical, $e_g$ is central and invertible.
Its trace on each simple is its nonzero scalar times the dimension.
\end{proof}

Prescribe a finite collection $(E_j,F_j)\in\Irr(C_1)\times\Irr(C_0)$.
Its partners are part of the input and are specified independently of
the comparison being sought. Choose any initial graded comparison $a$
and put
\[
 m_j(g)=\operatorname{tr}_{a_*E_j}(e_g^0),\qquad
 n_j(g)=\operatorname{tr}_{F_j}(e_g^0),\qquad
 B_{\rm test}=\langle\bar g:n_j(g)\ne0\text{ for some }j\rangle .
\]

\begin{proposition}[Calibration by finite traces]\label{prop:pinned-trace}
Assume $[\kappa_1]=[\kappa_0]$. A graded comparison carries each
$E_j$ to $F_j$ if and only if the following two conditions hold.
\begin{enumerate}[label=\textup{(\alph*)},leftmargin=*]
\item For all $j,g$, one has $m_j(g)=0$ if and only if $n_j(g)=0$.
\item The values $m_j(g)/n_j(g)$ wherever $n_j(g)\ne0$ define a single
character of $B_{\rm test}$. In particular they agree whenever $\bar g$
is the same and satisfy the relations in $B_{\rm test}$.
\end{enumerate}
When these conditions hold, the comparison cochains form a torsor
under $\Hom(K_{\rm ab}/B_{\rm test},\C^\times)$. They all induce
the same bijection on simple classes if and only if
\begin{equation}\label{eq:pinned-detecting}
                       B_{\rm test}=B(C_0).
\end{equation}
In this case the bijection is independent of the initial comparison
and of the chosen homogeneous coordinates. We call the tests
\emph{detecting} when \eqref{eq:pinned-detecting} holds.
\end{proposition}
\begin{proof}
A character correction $\chi$ changes $m_j(g)$ to
$\chi(g)^{-1}m_j(g)$. Thus the two conditions are necessary. Every
character of a subgroup of a finite abelian group extends to the
whole group, by adjoining one generator at a time and taking a root
in $\C^\times$. Extend the character in (b) to $K_{\rm ab}$. After
correction all traces agree on the homogeneous basis, so the modules
are isomorphic. The extensions form the stated torsor.

Two solutions differ by a character trivial on $B_{\rm test}$.
If \eqref{eq:pinned-detecting} holds, Lemma~\ref{lem:pinned-invisible}
makes their actions on every simple class equal. If
$B_{\rm test}\subsetneq B(C_0)$, characters of the finite quotient
$K_{\rm ab}/B_{\rm test}$ separate its points. One is nontrivial on
$B(C_0)$ and changes a simple class while preserving all tests.
The property of carrying the prescribed pairs to each other is
intrinsic, giving the asserted independence of coordinates.
\end{proof}

\begin{lemma}[Naturality of the calibrated simple classes]
\label{lem:pinned-naturality}
Let $r_i:C_i\to C'_i$ be graded isomorphisms over the same group
isomorphism $K\to K'$. Suppose $(r_1,r_0)$ transports a collection
of compatible detecting tests to the prescribed collection for the
primed pair. If $\phi,\phi'$ are the calibrated simple-class maps, then
\[
                 \phi'(r_1)_*=(r_0)_*\phi .
\]
This also holds for coherent transports on simple classes admitting
such graded representatives up to inner automorphisms.
\end{lemma}
\begin{proof}
For any calibrated comparison $a$, the conjugate $r_0 a r_1^{-1}$ is
graded and satisfies the primed tests. Apply uniqueness in
Proposition~\ref{prop:pinned-trace}. Inner automorphisms fix simple
isomorphism classes, so the same argument applies to compositions
and to a return to the original presentation.
\end{proof}

\subsection{Pinned parameter removal at depth zero}\label{sec:pinned-neutral}
Fix the finite data \eqref{eq:pinned-global-data}. For an admissible
pinned presentation $d$ let $A_{\epsilon,d}$ be the neutral algebra
and $B_{\epsilon,d}=\bigoplus_{g\in Q_d}A_{\epsilon,d}u_g$ the finite
crossed presentation of Proposition~\ref{prop:local-finite}, for
$\epsilon=0,1$. Write $a_d(g,h)=u_gu_hu_{gh}^{-1}$ for the full central
unit of \eqref{eq:local-defect}, including its scalar factor. For a
neutral simple $V$ with stabilizer $K=Q_{d,V}$ and intertwiners
$T_gT_h=c_{\epsilon,V}(g,h)T_{gh}$, its multiplicity algebra is
\begin{equation}\label{eq:pinned-effective}
 C_{\epsilon,d,V}=\C_{\kappa^{\rm eff}_{\epsilon,V}}[K],
 \qquad \kappa^{\rm eff}_{\epsilon,V}(g,h)
                 =a_d(g,h)(V)c_{\epsilon,V}(g,h)^{-1}.
\end{equation}
Here $a_d(g,h)(V)$ is evaluation at the full absorbed toral character
of $V$. The finite group $B_s$ stabilizing a rational Lusztig series,
the group $Q_{d,V}$ in this formula, and the endpoint stabilizer
$W_t^\fs$ in Theorem~\ref{thm:source-support-realization} enter at
different stages.

We recall the change-of-pinning formula for the normalized cuspidal
correspondence cited in Remark~\ref{rem:jd-llc}. Fix a finite source irreducible
$\rho$ over the cuspidal orbit and let $\cP_0,\cP_1$ be admissible
pinnings. The two enhanced targets have a common fibre product
$\widehat H_1\times_H\widehat H_0$, with projections $p_i$. On this
group the finite theorem gives
\begin{equation}\label{eq:pinned-jd-change}
 p_1^*\mathcal J_{\cP_1}(\rho)
  \simeq \vartheta_{10}\otimes p_0^*\mathcal J_{\cP_0}(\rho).
\end{equation}
The characters for successive changes multiply after pullback to
the triple fibre product. Both projections and their prescribed
central sectors enter this identity. To use it on $Q_{d,V}$ one must
pass through the actual restrictions and multiplicity spaces with
their normalizer operators. A comparison on the endpoint groups
$W_t^\fs$ requires its own realization, as in
\S\ref{sec:pinned-mixed}.

Consider a specified collection of presentations and changes between
them. We use the following hypotheses.
\begin{enumerate}[label=\textup{(N\arabic*)},leftmargin=*]
\item\label{pin:N1}
For each $d$ take the $Q_d$-equivariant neutral map $\Theta_{A,d}$
of Propositions~\ref{prop:equiv-full-q} and~\ref{prop:assembly-neutral}.
It preserves the full absorbed toral character and identifies the
classes in \eqref{eq:pinned-effective} for $V$ and $V_0=\Theta_{A,d}(V)$.
For each stipulated change $h:d\to d'$ there are algebra transports
$F_{\epsilon,h}:B_{\epsilon,d}\to B_{\epsilon,d'}$ preserving the
neutral subalgebras and their crossed products over the same
isomorphism $Q_d\to Q_{d'}$. Their actions on simple classes are
coherent and their neutral actions commute with $\Theta_A$.
\item\label{pin:N2}
For each neutral orbit representative $V$, the pinned finite Jordan
normalization and the actual block normalizer operators independently
prescribe test pairs for $C_{1,d,V}$ and $C_{0,d,V_0}$.
They satisfy (a), (b) and \eqref{eq:pinned-detecting} of
Proposition~\ref{prop:pinned-trace}.
\item\label{pin:N3}
The test pairs transport to the prescribed pairs under changes of
orbit representative and under every stipulated $h$, including
returns. Use the actual induced graded transports of the multiplicity
algebras, with their representative defects and the action induced
by \eqref{eq:pinned-jd-change} for changes of pinning. These
transports are coherent on simple classes.
\end{enumerate}
The fixed-presentation assertions in \ref{pin:N1} follow from the
cited propositions. For additional changes one must check their
action on the neutral shift, finite sheets and full corners. The
tests in \ref{pin:N2} are finite for each stabilizer. Their independent
realization on the entire family, and covariance in \ref{pin:N3},
are additional hypotheses.

\begin{theorem}[Pinned parameter removal]\label{thm:pinned-removal}
Assume \ref{pin:N1}--\ref{pin:N3}. Among the Clifford bijections
covering $\Theta_{A,d}$ and induced by graded comparisons on
multiplicity spaces, exactly one map on simple classes satisfies the
prescribed tests. It obeys
\begin{equation}\label{eq:pinned-neutral-natural}
 \Theta^{\rm pin}_{B,d}:\Irr(B_{1,d})\xrightarrow{\sim}\Irr(B_{0,d}),
 \qquad
 \Theta^{\rm pin}_{B,d'}(F_{1,h})_*=(F_{0,h})_*\Theta^{\rm pin}_{B,d}.
\end{equation}
If the changes intertwine the type equivalences and transport the
embedded Levi algebras, then
\begin{equation}\label{eq:pinned-depth-parametrization}
 \mathfrak b^{\rm pin}_{\fs,d}
     =\cl_{0,d}\circ\Theta^{\rm pin}_{B,d}\circ\typ_{\fs,d}
\end{equation}
is natural for the induced endpoint corner transports and satisfies
the conclusions of Theorem~\ref{thm:fiber-irr}.
\end{theorem}
\begin{proof}
Let $V_0=\Theta_{A,d}(V)$. Condition \ref{pin:N1} identifies the
stabilizers and gives equality of the effective multiplier classes.
By \ref{pin:N2} and Proposition~\ref{prop:pinned-trace} there is a
unique calibrated map $\phi_V$ on the corresponding simple
multiplicity modules. In \eqref{eq:assembly-induced} set
\[
 \Theta^{\rm pin}_{B,d}\!\left(\Ind_{B_{1,d,V}}^{B_{1,d}}(V\otimes E)\right)
   =\Ind_{B_{0,d,V_0}}^{B_{0,d}}(V_0\otimes\phi_V(E)),\qquad
 B_{\epsilon,d,V}=\bigoplus_{g\in Q_{d,V}} A_{\epsilon,d}u_g .
\]
On the input $u_g$ acts by $T_g\otimes R_g$ and its product with
$u_h$ includes exactly $a_d(g,h)(V)$ by \eqref{eq:pinned-effective}.
Condition \ref{pin:N3} and Lemma~\ref{lem:pinned-naturality} give the
same induced class for all choices of orbit representative. Clifford
classification in Proposition~\ref{prop:assembly-global} proves
bijectivity and uniqueness. The same lemma identifies the labels
after $h$, and $F_{\epsilon,h}$ transports induction from the inertia
subalgebra. This proves \eqref{eq:pinned-neutral-natural}, including
for returns to $d$.

The multiplicity of each distinct neutral constituent is $\dim E$,
which graded transport preserves. Its toral character, unitary Weyl
orbit and connected shift are those of $\Theta_{A,d}$, as in
Lemma~\ref{lem:neutral-invariants}. Finally a transport of the
embedded Levi algebra carries its central fibres and reduced orbit
idempotents, the neutral matrix algebra in each corner and its
centralizer to their counterparts. The intrinsic $\cl_{0,d}$ thus
intertwines the induced corner transports. Composing with the type
equivalences proves the last assertion.
\end{proof}

\begin{corollary}\label{cor:pinned-automatic}
On a transport-stable union of neutral orbits satisfying \ref{pin:N1},
suppose $B(C_{0,d,V_0})=1$ at every orbit. The simple-class map of
Theorem~\ref{thm:pinned-removal} is then unique and natural without
test pairs. This applies when the stabilizer is trivial or its
twisted group algebra is simple.
\end{corollary}
\begin{proof}
The empty collection is detecting and is preserved by every transport.
If the algebra is simple it has a single simple class, fixed by all
characters. Lemma~\ref{lem:pinned-invisible} and separation of
points by the characters of $K_{\rm ab}$ give $B(C)=1$.
\end{proof}
For a return fixing the neutral class and its stabilizer degrees,
the discrepancy between two graded comparisons is a character of
$Q_{d,V}$. By Lemma~\ref{lem:pinned-invisible} it acts trivially on
all labels exactly when it annihilates $B(C_{0,d,V_0})$. This is the
extra covariance condition beyond choosing a cochain on each neutral
orbit.

\section{The tame reduction square}\label{sec:tame}

For any connected reductive $F$-group $H$, write $\BBB(H,F)$ for its enlarged Bruhat--Tits building.

Assume $G$ splits over a tamely ramified extension and that the residue
characteristic is different from $2$. Throughout this section the type
$(K_G,\rho_G)$ and its depth-zero counterpart $(K^0,\rho^0)$ are constructed
as in \cite[\S4.1]{AFMO-reduce}, including the Fintzen--Kaletha--Spice
quadratic twist. In particular, the Levi factor uses the restriction of
the \emph{ambient} twist; it need not be the twist constructed intrinsically
for the Levi, as explained in \cite[Remark~4.1.5]{AFMO-reduce}.
The datum $\pi_{-1}$ below is the depth-zero inducing representation
corresponding to this chosen convention.

Consider a Kim--Yu datum (see \cite{KY17}), with the twist convention just specified,
\[
\Sigma=((\vec{G},L^{0}),(y,\{\iota\}),\vec{r},(K_{L^{0}},\rho_{L^{0}}),\vec{\phi}),
\]
attached to $\fs=[L,\pi]_{G}$, where $\pi$ is a tame
supercuspidal of $L(F)$. So $\vec{G}=(G^{0},\ldots,G^{d})$ is a
tamely ramified twisted Levi sequence in $G$, $L^{0}$ is a Levi
subgroup of $G^{0}$, $\vec{r}=(r_{0},\ldots,r_{d})$ is a sequence
of real numbers, $y$ is a point in $\BBB(L^{0},F)$, $\{\iota\}$
is a commutative diagram
\[
\xymatrix{
\BBB(G^{0},F) \ar@{^{(}->}[r]
& \BBB(G^{1},F) \ar@{^{(}->}[r]
& \cdots \ar[r]
& \BBB(G^{d},F) \\
\BBB(L^{0},F) \ar@{^{(}->}[r] \ar@{^{(}->}[u]
& \BBB(L^{1},F) \ar@{^{(}->}[r] \ar@{^{(}->}[u]
& \cdots \ar[r]
& \BBB(L^{d},F) \ar@{^{(}->}[u]
}
\]
where $L^{i}$ is the centralizer in $G^{i}$ of the split part of
the connected center of $L^{0}$. Write $\vec L=(L^{0},\ldots,L^{d})$.
Let $\pi_{-1}$ be the depth-zero inducing representation of $L^0(F)$
matched to $\pi$ by \eqref{eq:tame-seed}, and put
\[
\ft_{-1}=[L^0,\pi_{-1}]_{L^0},\qquad
\fs_{-1}=[L^0,\pi_{-1}]_{G^0}.
\]
Write $(K_{L^0},\rho_{L^0})$ for its compact Levi type and
$(K_L,\rho_L)$ for the positive-depth Levi type obtained by restriction
of the ambient twisted construction. These symbols in this section refer
to the tame construction; the depth-zero compact pair of
\S\ref{sec:notations} is recovered on the $G^0$-side.

\subsection{Type data and Levi algebras}
\label{subsec:tame-type-data}

Work over $\C$ with the fixed datum above, under the hypotheses
of \cite[\S4.1]{AFMO-reduce}. In particular, $G$ splits tamely and $p\ne2$.
Write $(K_0,\rho_0)=(K_{L^0},\rho_{L^0})$ and use the full
compact depth-zero Levi stabilizer $K_0=L^0(F)_y$, with $\rho_0$ trivial on
$L^0(F)_{y,0+}$ and with cuspidal connected constituents.
The depth-zero pair is a type for a single supercuspidal block.
For the positive-depth pair, the finite-type input of
\cite[\S4.1 and the opening of \S4.3]{AFMO-reduce}, together with
Proposition~\ref{prop:tame-single-block} below, gives the single-block
assertion before any use of inertial-normalizer exhaustion.
Write $\ft=[L,\pi]_L$,
$\fs=[L,\pi]_G$ and let $\ft_{-1}$ be the supercuspidal block
of $L^0(F)$ determined by $(K_0,\rho_0)$.
Its associated block in $G^0(F)$ is denoted $\fs_{-1}$.
The matched representative $\pi_{-1}$ is fixed below; these
inertial classes do not depend on that base-point choice. Set
\[
 K=K_L,\qquad \rho=\rho_L,\qquad
 (K_x,\rho_x)=(K_G,\rho_G),\qquad K^0_x=K^0 .
\]
We use the AFMO family for all generic $x$ in
$y+X_*(A_L)\otimes_{\mathbb Z}\R$, and initially fix one such $x$.
Write $K_+$ for the group denoted $K_{M,0+}$ in
\cite[(4.1.3)]{AFMO-reduce}; this notation is confined to this section.
The construction gives
\begin{equation}\label{eq:tame-compact-data}
 K=K_0K_+,\quad K\cap G^0(F)=K_0,\quad
 \rho=\operatorname{inf}(\rho_0)\otimes\eta_\Sigma,\quad
 \eta_\Sigma=\kappa_L^{\mathrm{nt}}\otimes
                   \epsilon_y^{\vec G}|_K .
\end{equation}
Here $\kappa_L^{\mathrm{nt}}$ is the untwisted Heisenberg
representation and $\epsilon_y^{\vec G}$ is the ambient
Fintzen--Kaletha--Spice character. Its restriction to $K$ is part of
the fixed construction.

Put
\begin{equation}\label{eq:tame-N}
 \mathcal N=N(\rho_0)_{[y]_{L^0}}
       \subset N_{G^0}(L^0)(F),\qquad H_0=\mathcal N\cap L^0(F).
\end{equation}
The results used below are the following concrete AFMO inputs.
\begin{enumerate}[label=(\roman*)]
\item The pairs $(K_x,\rho_x)$ are quasi-$G$-covers of $(K,\rho)$
and satisfy the family axioms, including
$K_{nx}=nK_xn^{-1}$ for $n\in\mathcal N$.
\item The support identity is
\[
 I_{G(F)}(\rho_x)=K_x\mathcal N K_x .
\]
Both assertions follow from \cite[Lemma~4.3.2,
Proposition~4.3.5]{AFMO-reduce}. Their depth-zero analogues hold for
$(G^0,L^0,K_0,\rho_0)$.
\item The group $\mathcal N$ normalizes $K$, $K_+$ and $\rho$,
and $\eta_\Sigma$ has a genuine extension
\[
 \widetilde\eta_\Sigma:\mathcal N K\longrightarrow
                      \operatorname{GL}(V_{\eta_\Sigma})
\]
given in \cite[Definition~4.3.3, Proposition~4.3.4]{AFMO-reduce}.
Fix this extension once and for all.
\end{enumerate}
Here $I_J(\sigma)$ denotes the set of elements of $J$ intertwining
the indicated compact representation, not its stabilizer alone.
These inputs refer to the actual ambient-twisted type.

\begin{proposition}[An absorption lemma for normalizer double cosets]
\label{prop:tame-absorption}
For $n\in\mathcal N$ one has
\begin{equation}\label{eq:tame-absorption}
 K_xnK_x\cap N_G(L)(F)=Kn .
\end{equation}
Consequently
\begin{equation}\label{eq:tame-intertwining}
 I_{N_G(L)(F)}(\rho)=K\mathcal N .
\end{equation}
\end{proposition}
\begin{proof}
We spell out the strengthening of the proof of
\cite[Proposition~3.4.12]{AFMO-structure} that is needed here.
For $g=k_1nk_2$ in the left-hand side of
\eqref{eq:tame-absorption},
\[
 gn^{-1}\in K_xK_{nx}.
\]
The last assertion of AFMO Axiom~3.4.1(5) supplies opposite
unipotent radicals $U,\bar U$ of parabolics with Levi $L$ such that
\[
 K_x\cap U\subset K_{nx}\cap U,\qquad
 K_{nx}\cap\bar U\subset K_x\cap\bar U .
\]
In this proof these symbols denote groups of $F$-points.
The decompositions in Axiom~3.4.1(4a), with common Levi
factor $K$, then give
\[
 K_xK_{nx}
 =(K_x\cap\bar U)(K_{nx}\cap U)K.
\]
It follows that $gn^{-1}=vuk$ with
$v\in\bar U$, $u\in U$ and $k\in K$.
The element $v u=gn^{-1}k^{-1}$ normalizes $L$.
The elementary big-cell identity
\[
 N_G(L)(F)\cap\bar U(F)U(F)=\{1\}
\]
is \cite[Lemma~3.4.13]{AFMO-structure}. Hence $g=kn$.
The converse inclusion in \eqref{eq:tame-absorption} is immediate.
Notice that $g$ was only assumed to normalize $L$; it was not
assumed to preserve $[y]_L$ or to lie in $\mathcal N$.

For $g\in N_G(L)(F)$ the pair
$({}^gK_x,{}^g\rho_x)$ is a quasi-$G$-cover of
$({}^gK,{}^g\rho)$. The quasi-cover intertwining identity
\cite[Lemma~3.3.5]{AFMO-structure} therefore gives
\[
 \Hom_{K_x\cap{}^gK_x}(\rho_x,{}^g\rho_x)
   =\Hom_{K\cap{}^gK}(\rho,{}^g\rho).
\]
If the latter is nonzero, the support identity places $g$ in
$K_xnK_x$ for some $n\in\mathcal N$.
Equation~\eqref{eq:tame-absorption} gives $g\in K\mathcal N$.
Every element of $K\mathcal N$ normalizes the isomorphism class of
$(K,\rho)$, proving the reverse inclusion.
\end{proof}

\begin{lemma}[Levi intertwining support]\label{lem:tame-Levi-support}
For the fixed compact Levi data,
\begin{equation}\label{eq:tame-Levi-support}
 I_{L(F)}(\rho)=J:=KH_0,\qquad
 I_{L^0(F)}(\rho_0)=H_0,\qquad J/K\simeq H_0/K_0.
\end{equation}
\end{lemma}
\begin{proof}
The relation $L\cap G^0=L^0$ gives $\mathcal N\cap L(F)=H_0$.
Intersect \eqref{eq:tame-intertwining} with $L(F)$ to obtain the
first equality. The same absorption argument for the depth-zero
family gives $K_0H_0=H_0$. Finally,
$K\cap H_0=K_0$ by \eqref{eq:tame-compact-data}, and $H_0$
normalizes $K$, proving the quotient isomorphism. These statements
use the quasi-cover and intertwining inputs, without a single-block
assumption for $(K,\rho)$.
\end{proof}

The support calculation also provides an explicit comparison of the
two Levi Hecke algebras. Initially use the convolution convention of AFMO and set
\[
 D_0^{\mathrm R}=\mathcal H(L^0(F),\rho_0),\qquad
 D^{\mathrm R}=\mathcal H(L(F),\rho).
\]
Normalize Haar measures by $\operatorname{vol}(K_0)=
\operatorname{vol}(K)=1$. Their opposites $D_0,D$ are the
left-module algebras used throughout.

\begin{proposition}[An equivariant Levi algebra isomorphism]
\label{prop:tame-Levi-isom}
There is a support-index-preserving isomorphism
\[
 \mathcal I_L:D_0^{\mathrm R}\xrightarrow{\sim}D^{\mathrm R}
\]
defined, on the line supported at $K_0h$ with $h\in H_0$, by
\begin{equation}\label{eq:tame-Levi-formula}
 \mathcal I_L(F_0)(h)=F_0(h)\otimes\widetilde\eta_\Sigma(h),
 \qquad \operatorname{supp}\mathcal I_L(F_0)=Kh .
\end{equation}
It intertwines the conjugation actions of $\mathcal N$ on these
algebras. It also commutes with unramified twisting:
twisting on the $L$-side by $\nu$ corresponds to twisting on the
$L^0$-side by $\nu|_{L^0(F)}$.
The same conclusions hold for the opposite algebras.
\end{proposition}
\begin{proof}
By \eqref{eq:tame-Levi-support} the two algebras have one-dimensional
support lines indexed by the same group $\Gamma=H_0/K_0\simeq J/K$.
Indeed their support groups normalize their respective compact
irreducible types, so Schur's lemma applies. For $h\in H_0$,
$F_0(h)\otimes\widetilde\eta_\Sigma(h)$ is a nonzero intertwiner for
$\rho$ whenever $F_0(h)$ is nonzero.
Its covariance under a change of representative by $K_0$ follows
from \eqref{eq:tame-compact-data} and
$\widetilde\eta_\Sigma|_K=\eta_\Sigma$.
This proves that \eqref{eq:tame-Levi-formula} is well defined and is
a linear isomorphism on every support line.

Convolution of two such lines is a single support line, because
$K_0\triangleleft H_0$ and $K\triangleleft J$.
With the chosen Haar measures there is no double-coset index factor.
If representatives satisfy $h_\gamma h_\delta=
k_{\gamma,\delta}h_{\gamma\delta}$, with
$k_{\gamma,\delta}\in K_0$, the factor contributed on the positive
side by $\eta_\Sigma$ is exactly
\[
 \widetilde\eta_\Sigma(h_\gamma)\widetilde\eta_\Sigma(h_\delta)
   =\eta_\Sigma(k_{\gamma,\delta})
                     \widetilde\eta_\Sigma(h_{\gamma\delta}).
\]
Thus the scalar multiplication defects of the two support-line
bases agree. This proves multiplicativity, including representative
defects, rather than just equality of their supports.

For $n\in\mathcal N$ choose a type intertwiner $S_n^0$ implementing
conjugation by $n$ on $\rho_0$, and choose
\[
 S_n=S_n^0\otimes\widetilde\eta_\Sigma(n)
\]
on $\rho$. The corresponding algebra automorphisms are
$F(h)\mapsto S_nF(n^{-1}hn)S_n^{-1}$, and similarly at depth zero.
The identity
\[
 \widetilde\eta_\Sigma(n)\widetilde\eta_\Sigma(n^{-1}hn)
                     \widetilde\eta_\Sigma(n)^{-1}
       =\widetilde\eta_\Sigma(h)
\]
proves that $\mathcal I_L$ commutes with these automorphisms.
Changing $S_n^0$ by a scalar does not change them.
The products $S_nS_m$ have exactly the same scalar defects as
$S_n^0S_m^0$, since $\widetilde\eta_\Sigma$ is a genuine representation.

An unramified character is trivial on $K$ and $K_0$.
The two twisting automorphisms multiply the line labelled by $h$
by the same scalar $\nu(h)$. Hence they commute with
$\mathcal I_L$. Taking opposites preserves all of these identities;
the inverse-multiplier convention is applied on both sides.
\end{proof}

\begin{proposition}[The supercuspidal Levi type]
\label{prop:tame-single-block}
The pair $(K,\rho)$ is a type for the single supercuspidal Bernstein
block $[L,\pi]_L$.
\end{proposition}
\begin{proof}
The AFMO construction makes $(K,\rho)$ a type for a finite set
$\mathfrak S$ of Bernstein blocks. By
Proposition~\ref{prop:Levi-Azumaya}, the depth-zero algebra $D_0$
is Azumaya over the coordinate ring of its unramified orbit. This
Laurent polynomial algebra has no nontrivial idempotents.
Proposition~\ref{prop:tame-Levi-isom} identifies $D_0$ with $D$,
so $D$ has no nontrivial central idempotents either.

The type equivalence with the finite product of blocks indexed by
$\mathfrak S$ would supply a central idempotent for every nontrivial
decomposition of $\mathfrak S$. Hence $\mathfrak S$ is a singleton.
The fixed supercuspidal representation $\pi$ contains $\rho$, so
this block is $[L,\pi]_L$.
\end{proof}

\subsection{Inertial normalizers}\label{subsec:tame-normalizers}

\begin{proposition}[Types detect inertial normalizers]
\label{prop:tame-types}
Let $(K,\rho)$ be a type for a supercuspidal block
$\ft=[L,\pi]_L$. If $n\in N_G(L)(F)$ satisfies
$ {}^n\pi\simeq\pi\otimes\nu$ for an unramified character $\nu$,
then some element of $L(F)n$ intertwines $\rho$.
Conversely, an element of $N_G(L)(F)$ normalizing the pair
$(K,\rho)$ belongs to $N^\fs$, where $\fs=[L,\pi]_G$.
\end{proposition}
\begin{proof}
Let
\[
 P=\operatorname{c\!-\!Ind}_{K}^{L(F)}\rho,\qquad
 P_n=\operatorname{c\!-\!Ind}_{{}^nK}^{L(F)}{}^n\rho.
\]
There are nonzero maps $P\to\pi$ and $P_n\to\pi$.
For the second map, first use Frobenius reciprocity for
$ {}^n\pi$, and then the fact that $\nu$ is trivial on the compact
group $ {}^nK$. Both maps are surjective since $\pi$ is irreducible.
The object $P$ is projective: over $\C$, the functor
$\Hom_K(\rho,-)$ is exact by compact averaging.
Lift $P\to\pi$ through $P_n\to\pi$ to obtain a nonzero map
$P\to P_n$. Compact Mackey decomposition and Frobenius reciprocity
then give $\ell\in L(F)$ with
\[
 \Hom_{K\cap{}^{\ell n}K}(\rho,{}^{\ell n}\rho)\ne0.
\]
The finite-dimensionality of $\rho$ ensures that a nonzero map into
the Mackey direct sum has a nonzero component.

Conversely, if $n$ normalizes the pair, then $ {}^n\pi$ contains
$\rho$. The single-block type property implies that $ {}^n\pi$
lies in $\operatorname{Rep}(L)^\ft$. Its irreducibles are precisely
the unramified twists of $\pi$, which proves the claim.
\end{proof}

\begin{proposition}[Exhaustion of the inertial normalizers]
\label{prop:tame-exhaustion}
With the fixed AFMO types,
\begin{equation}\label{eq:tame-exhaustion}
 N^\fs=L(F)\mathcal N,\qquad
 N^{\fs_{-1}}=L^0(F)\mathcal N.
\end{equation}
In particular
\begin{equation}\label{eq:tame-W}
 \mathcal N/H_0\xrightarrow{\sim}W^\fs,\qquad
 \mathcal N/H_0\xrightarrow{\sim}W^{\fs_{-1}} .
\end{equation}
\end{proposition}
\begin{proof}
The group $\mathcal N$ normalizes $(K,\rho)$, so
Proposition~\ref{prop:tame-types} gives $\mathcal N\subset N^\fs$.
If $n\in N^\fs$, the same proposition and
\eqref{eq:tame-intertwining} give $\ell n\in K\mathcal N$ for
some $\ell\in L(F)$. Thus $n\in L(F)\mathcal N$.
Repeat the argument for the depth-zero quasi-cover family to obtain
the second equality in \eqref{eq:tame-exhaustion}.

The identity $\mathcal N\cap L(F)=H_0$ from
Lemma~\ref{lem:tame-Levi-support} identifies both kernels in
\eqref{eq:tame-W}.
\end{proof}

In particular,
\[
 N^{\fs_{-1}}\subset N^\fs,\qquad
 L(F)\cap N^{\fs_{-1}}=L^0(F),
\]
under the standing AFMO hypotheses and the single-block type
convention.

\subsection{The unramified orbit map}
\label{subsec:tame-orbit}

\paragraph{Fixing the depth-zero base representation.}
The compact pair $(K_0,\rho_0)$ specifies an inertial block, not a
unique member of its unramified orbit. The datum $\pi_{-1}$ must
therefore include a compatible choice of the representation on the
compact-mod-center inducing group. Here is an explicit prescription
in the ambient convention.

For the fixed $\pi$ containing $\rho$, let
$E_\pi=\Hom_K(\rho,\pi)$, viewed as a simple left $D$-module.
Pull it back along the opposite of $\mathcal I_L$, and let
$\pi_{-1}$ be the depth-zero representation corresponding to that
simple $D_0$-module. This is equivalently the following compact
induction prescription. For the corresponding irreducible
$H_0$-representation $\tau_0$, whose restriction to $K_0$ is
$\rho_0$-isotypic, use
\begin{equation}\label{eq:tame-seed}
 \begin{split}
 \pi_{-1}&\simeq\operatorname{c\!-\!Ind}_{H_0}^{L^0(F)}\tau_0,\\
 \pi&\simeq\operatorname{c\!-\!Ind}_{J}^{L(F)}
       \bigl(\operatorname{inf}\tau_0
                      \otimes\widetilde\eta_\Sigma|_J\bigr).
 \end{split}
\end{equation}
The inflation is along
$J=H_0K_+\to H_0/(H_0\cap K_+)$.
Its kernel acts trivially on $\tau_0$, since
$H_0\cap K_+\subset K_0\cap\ker\rho_0$.
The equivalence of these prescriptions follows from compact Clifford
theory and induction in stages: compact induction from $K$ to $J$
has the same support-line endomorphism algebra as compact induction
from $K$ to $L(F)$, by \eqref{eq:tame-Levi-support}.
The same statement holds for $K_0,H_0,L^0(F)$.
The corresponding compact projective generators identify their
categories of representations generated by the indicated compact
type with modules over that algebra; induction in stages corresponds
to the identity functor on these modules.
Their simple modules are therefore precisely the multiplicity
modules of the inducing representations in \eqref{eq:tame-seed}.
Formula~\eqref{eq:tame-Levi-formula} tensors these representations by
the fixed $\widetilde\eta_\Sigma$.

We use this matched base representation throughout. The convention
depends on the fixed extension $\widetilde\eta_\Sigma$ and retains the
restriction of the ambient twist on the Levi.

\begin{theorem}[Orbit and symmetry comparison]
\label{prop:tame-orbit}\label{thm:isom-center-supercuspidal-blocks}
\label{thm:supercuspidal-block-equivalence}\label{thm:weyl-isom}
For the matched inducing data \eqref{eq:tame-seed}, there is an
isomorphism of algebraic unramified orbits
\begin{equation}\label{eq:tame-orbit}
 f_\Sigma:\Irr^\ft(L(F))\xrightarrow{\sim}
                \Irr^{\ft_{-1}}(L^0(F)),\qquad
 f_\Sigma(\pi\otimes\nu)=\pi_{-1}\otimes\nu|_{L^0(F)} .
\end{equation}
It is $\mathcal N$-equivariant. Consequently it is equivariant for
the isomorphism
\begin{equation}\label{eq:tame-iota}
 \iota_\Sigma:W^\fs\xrightarrow{\sim}W^{\fs_{-1}}
\end{equation}
obtained by identifying both sides with $\mathcal N/H_0$.
\end{theorem}
\begin{proof}
The type equivalences and $\mathcal I_L$ give an equivalence of the
two supercuspidal blocks. It sends $\pi$ to $\pi_{-1}$ by
construction. Its compatibility with twisting proves the formula
in \eqref{eq:tame-orbit}. The irreducibles of a supercuspidal block
form exactly one unramified orbit, so this formula is a bijection
on the displayed sets. In particular
\[
 \pi\otimes\nu\simeq\pi
 \quad\Longleftrightarrow\quad
 \pi_{-1}\otimes\nu|_{L^0(F)}\simeq\pi_{-1}.
\]
This proves well-definedness and injectivity on the self-twist
quotients, and the equivalence gives surjectivity.
The induced isomorphism of centers is an algebraic isomorphism of
these orbits, since their Bernstein centers are their coordinate
rings. 

Conjugation by $n\in\mathcal N$ on representations corresponds under
the type equivalences to the algebra automorphisms in
Proposition~\ref{prop:tame-Levi-isom}. Their intertwining identity gives
\[
 f_\Sigma({}^n\sigma)={}^nf_\Sigma(\sigma)
                 \qquad(\sigma\in\Irr^\ft(L(F))).
\]
The subgroup $H_0$ acts by inner conjugation on both groups' simple
representations, so it acts trivially on their isomorphism classes.
Now \eqref{eq:tame-W} proves equivariance for \eqref{eq:tame-iota}.

For clarity, if $ {}^n\pi\simeq\pi\otimes\nu_n$, the same identity
gives
\[
 {}^n\pi_{-1}\simeq\pi_{-1}\otimes\nu_n|_{L^0(F)} .
\]
Thus both the linear part of the unramified-character action and
the translation part of the affine action agree. Equivariance is
not limited to the action on the tangent space of the torus.
\end{proof}

\subsection{Normalizer algebras and endpoint cocycles}
\label{subsec:tame-endpoint-cocycles}
Retain the notation above, and put
\[
 \widehat\cN=\cN K,\qquad
 \cW=\cN/K_0\simeq\widehat\cN/K,\qquad
 \Gamma=H_0/K_0\simeq KH_0/K.
\]
The compact subgroups are normal in these normalizers. Define the
raw normalizer algebras
\[
 C_0^{\mathrm R}=\End_{\cN}
                  (\cind_{K_0}^{\cN}\rho_0),\qquad
 C_\Sigma^{\mathrm R}=\End_{\widehat\cN}
                  (\cind_K^{\widehat\cN}\rho),
\]
and write $C_0,C_\Sigma$ for their opposites. The subalgebras in
degrees $\Gamma$ are $D_0,D$.

\begin{proposition}[Normalizer comparison, including representative defects]
\label{prop:tame-normalizer-isom}
There is a $\cW$-graded isomorphism
\[
 \mathcal I_N:C_0\xrightarrow{\sim}C_\Sigma
\]
whose restriction is $\mathcal I_L$. Its induced map on the Levi centers is $f_\Sigma^*$, where
\[
 f_\Sigma(\pi\otimes\nu)
     =\pi_{-1}\otimes\nu|_{L^0(F)}
\]
is the matched orbit map of Theorem~\ref{thm:weyl-isom}.
Consequently the endpoint cocycle family intrinsically constructed
from this fixed tame type and its full normalizer is the pullback
of the depth-zero family, including its corner transports.
\end{proposition}
\begin{proof}
Choose representatives $n_a\in\cN$ for $a\in\cW$, with
$n_1=1$, and set
$n_an_b=d(a,b)n_{ab}$, where $d(a,b)\in K_0$.
Choose nonzero type intertwiners $S_a^0$, normalized at $1$, with
\begin{equation}\label{eq:tame-normalizer-products}
 S_a^0S_b^0=c(a,b)\rho_0(d(a,b))S_{ab}^0.
\end{equation}
Set $S_a=S_a^0\otimes\widetilde\eta_\Sigma(n_a)$.
The construction of the genuine extension makes $S_a$ an
intertwiner of the positive-depth type. Since it is genuine,
\begin{align*}
 S_aS_b
 &=c(a,b)\rho_0(d(a,b))S_{ab}^0\otimes
            \eta_\Sigma(d(a,b))\widetilde\eta_\Sigma(n_{ab})\\
 &=c(a,b)\rho(d(a,b))S_{ab}.
\end{align*}
Each support line is one dimensional by Schur's lemma. The two
normalizer convolution algebras therefore have identical
multiplication constants with these bases; no ambient double-coset
index occurs because $K_0\triangleleft\cN$ and
$K\triangleleft\widehat\cN$. Normalizing their compact Haar
volumes to $1$ gives the isomorphism in the right convention, and
taking both opposites gives $\mathcal I_N$.
The representative rule in \eqref{eq:tame-normalizer-products} ensures
well-definedness under $n_a\mapsto n_ak_0$ as well.

On the $\Gamma$-lines this is exactly the tensor map defining
$\mathcal I_L$. Its effect on simple Levi modules is $f_\Sigma$,
so its center map is $f_\Sigma^*$, with the centers normalized by
the matched base representations. The quotient grading is
identified with $W^\fs$ and $W^{\fs_{-1}}$ by Proposition~\ref{prop:tame-exhaustion}. Lemma~\ref{lem:graded-corner-functoriality} now applies to the raw
normalizer algebras.

To express the assertion using the source endpoints, let
\[
 J_0:C_0\xrightarrow{\sim}B_{0,-1},\qquad
 J_\Sigma:C_\Sigma\xrightarrow{\sim}B_{0,\Sigma}
\]
be the normalization maps of Lemma~\ref{lem:normalizer-cocycle-normalization} and Proposition~\ref{prop:Levi-qzero-embedding}, applied to the two AFMO presentations. The same argument
applies at positive depth because the AFMO relations used in
that lemma are axiomatic and hold for the datum $\Sigma$.
Then
\begin{equation}\label{eq:tame-endpoint-isom}
 F_{\Sigma,0}=J_\Sigma\mathcal I_NJ_0^{-1}:
             B_{0,-1}\xrightarrow{\sim}B_{0,\Sigma}
\end{equation}
preserves the quotient grading and the actual embedded Levi
algebras. Applying Lemma~\ref{lem:graded-corner-functoriality} once more proves
the endpoint assertion and coherence of its transports.
\end{proof}

\begin{remark}[Tensor bases and normalized AFMO bases]
The tensor intertwiners in \eqref{eq:tame-normalizer-products} are used
to compare the raw multiplication lines. They are not asserted to
be the normalized AFMO generators on every affine support.
\cite[Proposition~4.4.4]{AFMO-reduce} prescribes this tensor normalization on
length-zero supports; the remaining normalized generators are
fixed by \cite[Choice~3.10.3]{AFMO-structure}. Formula
\eqref{eq:tame-endpoint-isom} includes the necessary basis
changes. It does not identify raw tensor bases with normalized
ambient Hecke-function values.
\end{remark}

\subsection{Transfer to the tame block}
\begin{theorem}[AFMO reduction, {\cite[Theorem~4.3.11]{AFMO-reduce}}]
\label{thm:AFMO-depth-zero-counterpart-block}\label{thm:AFMO-depthzero}
Under the standing assumptions of \cite[\S4.1]{AFMO-reduce}, in
particular tame splitting and residue characteristic different from $2$,
let $(K_G,\rho_G)$ be its Fintzen--Kaletha--Spice twisted type and
$(K^0,\rho^0)$ its depth-zero counterpart. Then there is
a support-preserving algebra isomorphism
\[
\mathcal H\bigl(G(F),(K_G,\rho_G)\bigr)
\;\xrightarrow{\ \sim\ }\;
\mathcal H\bigl(G^{0}(F),(K^{0},\rho^{0})\bigr).
\]
If the coefficient field admits a nontrivial involution, the isomorphism may be
chosen to preserve the corresponding anti-involutions.
\end{theorem}
\begin{corollary}\label{cor:isom-B-center}
Under the standing AFMO assumptions, there is an algebra isomorphism
\[
 \mathfrak Z^\fs\simeq\mathfrak Z^{\fs_{-1}}.
\]
\end{corollary}
\begin{proof}
Theorem~\ref{thm:AFMO-depthzero} identifies the centers of the two
type algebras. The type equivalences identify these centers with the
centers of the corresponding Bernstein blocks.
\end{proof}

Fix the matched datum $\Sigma$ and write
$\mathcal F_\Sigma:\Rep(G(F))^\fs\to
\Rep(G^0(F))^{\fs_{-1}}$ for the chosen AFMO block equivalence.
The orbit and symmetry maps are $f_\Sigma$ above and
$\iota_\Sigma:W^\fs\to W^{\fs_{-1}}$, from Theorem~\ref{thm:weyl-isom}.
Write $T^{\fs_{-1}}=\Irr^{\ft_{-1}}(L^0(F))$, and let
$\kappa_{-1}$ be the depth-zero endpoint family from
Theorem~\ref{thm:source-support-realization}, with its intrinsic
corner transports. Choose corner coordinates compatibly with
\eqref{eq:tame-endpoint-isom}, and set
\begin{equation}\label{eq:tame-endpoint-cocycle}
 \kappa_{\Sigma,t}(w,v)
 =\kappa_{-1,f_\Sigma(t)}
                  (\iota_\Sigma(w),\iota_\Sigma(v)).
\end{equation}
Let
\[
 \mathfrak q_\Sigma:
 (T^{\fs_{-1}}\extquot W^{\fs_{-1}})_{\kappa_{-1}}
 \longrightarrow (T^\fs\extquot W^\fs)_{\kappa_\Sigma}
\]
send $[t_0,E]$ to
$[f_\Sigma^{-1}(t_0),\iota_\Sigma^*E]$, with the corresponding
corner transports.

The two parametrizations and the vertical maps form the square
\begin{equation}\label{eq:tame-comparison-diagram}
\xymatrix@C=40pt{
\Irr(\Rep(G(F))^\fs) \ar[r]^-{\mathfrak b_{\fs,\Sigma}} \ar[d]_{\Irr(\mathcal F_\Sigma)}
& (\Irr^{\ft}(L(F))//W^\fs)_{\kappa_\Sigma}
\ar[d]^{(f_\Sigma,\iota_\Sigma)_*}_{\sim}\\
\Irr(\Rep(G^0(F))^{\fs_{-1}}) \ar[r]^-{\mathfrak b_{\fs_{-1}}}
& (\Irr^{\ft_{-1}}(L^0(F))//W^{\fs_{-1}})_{\kappa_{-1}}.
}
\end{equation}

\begin{theorem}[Tame parametrization]\label{thm:tame-global-extended-quotient}
Use the ambient-twisted AFMO datum $\Sigma$, over $\C$, with
tame splitting, $p\ne2$, and full compact depth-zero Levi
stabilizer. Assume its depth-zero counterpart satisfies the
hypotheses of Theorem~\ref{thm:fiber-irr}, with the specified
$0$-generic admissible embedding and
$p\nmid|W_{\mathrm{abs}}(G^0)|$.
Fix the AFMO block equivalence and the comparison and
projective-label choices in that theorem.

The family \eqref{eq:tame-endpoint-cocycle} is the intrinsic endpoint
corner family for this fixed tame type, up to compatible
coboundary changes. Its coherent transports are those of the
actual tame normalizer algebra. The map
\begin{equation}\label{eq:tame-parametrization}
 \mathfrak b_{\fs,\Sigma}
 =\mathfrak q_\Sigma\circ\mathfrak b_{\fs_{-1}}
                         \circ\Irr(\mathcal F_\Sigma)
\end{equation}
is a bijection
\[
 \Irr(\Rep(G(F))^\fs)\xrightarrow{\sim}
                          (T^\fs\extquot W^\fs)_{\kappa_\Sigma}.
\]
After applying $\mathcal F_\Sigma$, its neutral constituents,
unitary orbits, independent toral characters, local geometric
shift, and effective projective labels satisfy properties
\textup{(i)--(iii)} of Theorem~\ref{thm:fiber-irr}. Its endpoint central character
is obtained from the depth-zero endpoint character by the induced
map $f_\Sigma^{-1}$ on orbit spaces. Thus the square
\eqref{eq:tame-comparison-diagram} commutes, and its right vertical
map is compatible with the projections of the twisted extended
quotients to their orbit spaces.
\end{theorem}
\begin{proof}
Proposition~\ref{prop:tame-single-block} justifies the single-block Levi
input before inertial-normalizer exhaustion. Propositions~\ref{prop:tame-exhaustion} and~\ref{prop:tame-Levi-isom}, together with Theorem~\ref{thm:weyl-isom}, then give the actual orbit map and the
full symmetry isomorphism, with their affine translation parts.
The AFMO algebra isomorphism
\cite[Theorem 4.3.11]{AFMO-reduce} gives
$\mathcal F_\Sigma$ for the specified single-block types.

Proposition~\ref{prop:tame-normalizer-isom} identifies the endpoint
corners and their multiplication, not just their indexing sets.
Thus \eqref{eq:tame-endpoint-cocycle} and $\mathfrak q_\Sigma$
represent the intrinsic tame endpoint data in compatible
coordinates. Apply Theorem~\ref{thm:fiber-irr} to
$\mathcal F_\Sigma(\Pi)$ and then $\mathfrak q_\Sigma$.
Every map is a bijection, proving \eqref{eq:tame-parametrization}.
The neutral and projective-label properties follow from the
depth-zero theorem. For the last assertion, $f_\Sigma^*$ is
the actual Levi center map of
\eqref{eq:tame-endpoint-isom}, so the first coordinate of
$\mathfrak q_\Sigma([t_0,E])$ is precisely
$f_\Sigma^{-1}(t_0)$. The same statement descends to orbit
spaces by equivariance.
\end{proof}

\begin{remark}[Cocycles and comparison choices]
The construction uses three multipliers on generally different groups.
The source multiplier $\mu$ is defined on $\Omega$ in the left AFMO
presentation. The effective multiplier
$\kappa^{\mathrm{eff}}_{\epsilon,V}$ is defined on the stabilizer
$Q_V$ of a neutral constituent; it includes the evaluated central
defects and the inverse intrinsic obstruction. The endpoint
multiplier $\kappa_t$ is defined on $W_t^\fs$ and comes from the
centralizer of the neutral matrix algebra in the orbit corner.
Passing from the neutral description to the endpoint corner retains
the finite lattice and sheet information through the full algebra.

For the fixed type and source algebra, changing a nonzero intertwiner
changes its scalar multiplier by a coboundary, and changing matrix
coordinates identifies the corresponding projective labels. The
coherent corner transports are part of the target. Two comparison
cochains with the same coboundary differ by a character, which tensors
the residual projective label by that character. The parametrizations
use the fixed type datum and comparison choices; no independence from
arbitrary changes of type datum is asserted.

The endpoint cocycles are not identified here with those in
endomorphism algebras of individual parabolically induced
representations. Such a comparison involves analytic $R$-groups,
which can be smaller than the full Weyl stabilizer; see
\cite[Theorem~4.2 and Example~4.3]{Sol-survey}. The structural
properties of Theorems~\ref{thm:fiber-irr}
and~\ref{thm:tame-global-extended-quotient} concern the stated
endpoint comparison. Section~\ref{sec:pinned-neutral} specifies which
changes of presentation the pinned refinement accommodates, and
Proposition~\ref{prop:support-comparison} gives the criterion for
preserving the original cuspidal-support fibres.
\end{remark}

\begin{corollary}[Transport of pinned parameter removal]\label{cor:pinned-tame}
In Theorem~\ref{thm:tame-global-extended-quotient}, suppose that the
reduced block satisfies Theorem~\ref{thm:pinned-removal}, with the
specified type and Levi transports. Its pinned parametrization
determines the unique upper parametrization making
\eqref{eq:tame-comparison-diagram} commute. Explicitly,
\[
 \mathfrak b^{\rm pin}_{\fs,\Sigma,d}
  =\mathfrak q_\Sigma\,\mathfrak b^{\rm pin}_{\fs_{-1},d}
                           \,\Irr(\mathcal F_\Sigma).
\]
It is covariant under changes of the pinned datum for which the AFMO
reduction and the actual normalizer corner maps intertwine the given
transports.
\end{corollary}
\begin{proof}
The right vertical map of \eqref{eq:tame-comparison-diagram} is
$\mathfrak q_\Sigma^{-1}$ by
Proposition~\ref{prop:tame-normalizer-isom}. Inverting it gives the
formula and uniqueness, with the projective label transported through
the actual corner. Covariance follows by composing the three
intertwining identities.
\end{proof}

\clearpage
\part{Rigidly enhanced Langlands parameters}

\section{Rigid enhancements, support and inertial orbits}\label{sec:rigid}

Let $G^*$ be a quasi-split connected reductive group over a
nonarchimedean local field $F$, and let
$\rig=(G^*,G,\xi,z)$ be a fixed rigid inner twist. Fix a finite central
$F$-subgroup scheme $Z\subset G^*$ through which the rigidification
is represented. Write $\hG=\widehat{G^*}$ and
$\hG^+=\widehat{G^*/Z}$. For a rigid Levi datum
$\mathfrak m$ of $\rig$, let $\cM\subset\LG{G^*}$ be the corresponding
Levi $L$-subgroup, with connected part $\hM$.

\subsection{The common cover and the parameter model}

Put $p_Z:\hG^+\to\hG$. For any subgroup $J\subset\hG$, write
$J^+=p_Z^{-1}(J)$. In particular $\hM^+=p_Z^{-1}(\hM)$ is
the dual group of $M^*/Z$ for a compatible Levi $M^*$ of $G^*$.
This is a full inverse image. It is not obtained by dividing out
the center before taking a simply connected cover.

We use the commuting-unipotent model of a parameter:
\[
 (\psi,u),\qquad
 \psi:\WF\longrightarrow\LG{G^*}\text{ admissible},\quad
 u\in Z_{\hG}(\psi)\text{ unipotent}.
\]
It corresponds to an $\SL_2$-parameter
$\varphi:\WF\times\SL_2(\C)\to\LG{G^*}$ by taking
$\psi=\varphi|_{\WF}$ and the usual upper-triangular unipotent
element in its $\SL_2$ factor. Conversely an $\SL_2$ enrichment
exists, and the resulting conjugacy class is independent of its choice.
The Weil--Deligne homomorphism obtained by inserting the diagonal
$|w|^{1/2}$ is a different model and must not be substituted for
$\psi$ in the centralizers below.

Define
\begin{equation}\label{eq:rigid-groups}
 H_\psi^+=Z_{\hG}(\psi)^+,
 \qquad
 S_{\psi,u}^+=Z_{\hG}(\psi,u)^+,
 \qquad
 A_{\psi,u}^+=\pi_0(S_{\psi,u}^+).
\end{equation}
The unipotent $u$ has a unique unipotent lift $u^+$ to $\hG^+$.
Thus $S_{\psi,u}^+=Z_{H_\psi^+}(u^+)$. Equivalently, the
component group in \eqref{eq:rigid-groups} is canonically identified
with $\pi_0(Z_{\hG}(\varphi)^+)$; see
\cite[Corollary~4.2.3]{DS}. A rigid enhancement is
$\rho\in\Irr(A_{\psi,u}^+)$.

The groups written as inverse images in \eqref{eq:rigid-groups}
do not require a lift of $\psi$ to an $L$-group with connected
dual group $\hG^+$. This observation also applies to twisted-Levi
comparisons in Section~\ref{sec:reduction}.

\subsection{Rigid central characters}

Let
\begin{equation}\label{eq:central-character-group}
 C_G^Z=\pi_0\bigl((Z(\hG)^{\Gamma_F})^+\bigr).
\end{equation}
Rigid duality attaches to $z$ a character
$\zeta_\rig:C_G^Z\to\C^\times$; see
\cite{KalethaRigid,Dillery}. There is a natural map
$a_{\psi,u}:C_G^Z\to A_{\psi,u}^+$ with central image.
The relevance condition is
\begin{equation}\label{eq:relevance}
 \rho(a_{\psi,u}(c))=\zeta_\rig(c)\operatorname{Id}
 \qquad(c\in C_G^Z).
\end{equation}
The map $a_{\psi,u}$ need not be injective. We use
\eqref{eq:relevance}, not a restriction to an implicitly embedded
copy of $C_G^Z$.

For a rigid Levi $\mathfrak m$, define $C_M^Z$ in the same way.
The common cover gives $C_G^Z\to C_M^Z$, and the character of
$\mathfrak m$ pulls back to $\zeta_\rig$. The Levi includes its
rigidification as well as its rational realization in $\rig$.
An abstract inner form of $M^*$ is insufficient for this bookkeeping.

\subsection{Cuspidality}

The enhancement is an $H_\psi^+$-equivariant local system on
the orbit of $u^+$. We use Dillery--Schwein's cuspidality convention:
the parameter is discrete in its $L$-group and the local system is
cuspidal, equivalently its constituents on the identity component
are cuspidal. Quasi-cuspidal support is taken in their theory of
quasi-Levi subgroups. For a disconnected reductive group $H$, these
are subgroups
\[
 K=Z_H(Z(L)^\circ),\qquad L\text{ a Levi of }H^\circ.
\]
Write $\Phic^Z(\mathfrak m)$ for the relevant cuspidal enhanced
parameters, up to $\hM$-conjugacy.

\begin{remark}[The precise Springer convention]\label{rem:DS-AMS}
We use the quasi-cuspidal support and parametrization of \cite{DS}
throughout. Their Remark~3.4.6 does not establish equality with every
label in the original AMS parametrization \cite{AMS}. In their endomorphism
calculation, a Fourier-associated cuspidal local system may occur;
its equality with the original local system is Conjecture~3.4.2
there. Neither that equality nor equality of the two disconnected
support constructions is an input here. Any application importing
an \emph{AMS-normalized} multiplier from another paper must compare
it with the actual algebra defined in Section~\ref{sec:algebras}.
\end{remark}

\subsection{The two cuspidal-support maps}\label{sec:support}

For $(\psi,u,\rho)$ take its quasi-cuspidal support in $H_\psi^+$:
\[
 [K^+,v^+,\epsilon]_{H_\psi^+}.
\]
Let $\cM$ be the centralizer in $\LG{G^*}$ of the image of
$Z(K^+)^\circ$. The construction gives
\begin{equation}\label{eq:support-Levi}
 K^+=H_\psi^+\cap\hM^+,
 \qquad
 \uSc(\psi,u,\rho)=[\cM,\psi,v,\epsilon]_{\hG}.
\end{equation}
The right side is a cuspidal datum, well defined up to conjugacy;
these assertions are \cite[Construction~4.6.2 and Lemma~4.6.4]{DS}.

If $\varphi$ is an $\SL_2$ enrichment, its infinitesimal character
is the conjugacy class of
\begin{equation}\label{eq:infinitesimal}
 \lambda_\varphi(w)=
 \varphi\left(w,
 \begin{pmatrix}|w|^{1/2}&0\\0&|w|^{-1/2}\end{pmatrix}\right).
\end{equation}
In general \eqref{eq:support-Levi} changes this character.
There is an unramified central twist $z_{\psi,u,\rho}$ of
$\cM$ such that
\begin{equation}\label{eq:normalized-support}
 \Sc(\psi,u,\rho)
  =[\cM,z_{\psi,u,\rho}\psi,v,\epsilon]_{\hG}
\end{equation}
has the same infinitesimal character as $(\psi,u)$.
Here multiplication denotes the unramified cocycle twist, rather
than multiplication by a constant on every element of $\WF$.
The correcting element has a construction by an algebraic
cocharacter; the output class, including its enhancement, is
well defined by \cite[Lemma~4.6.6 and Definition~4.6.7]{DS}.
Our convention for $|w|$ fixes the sign in
\eqref{eq:infinitesimal}; it also fixes the corresponding
correcting cocharacter convention.

\begin{proposition}\label{prop:inertial-support}
The maps $\Sc$ and $\uSc$ have the same inertial class. Both preserve
the ambient rigid central character. In particular, they define the
same partition into rigid Galois Bernstein blocks.
\end{proposition}
\begin{proof}
Equation~\eqref{eq:normalized-support} differs from
\eqref{eq:support-Levi} by an unramified twist in the supporting
Levi. This is precisely an allowed operation in inertial equivalence.
For the central character, the center of $H_\psi^+$ acts by the
same character on a simple summand and its quasi-cuspidal support;
this is \cite[Lemma~4.6.8]{DS}. Pull back that equality along
$C_G^Z$ to obtain \eqref{eq:relevance} for the support.
A central twist in the Levi changes neither its enhancement nor
its rigid central-character map. Hence the assertion holds for
both support maps.
\end{proof}

\begin{example}[The correction for $\GL_2$]\label{ex:GL2}
Take the split group $\GL_2$, $Z=1$, and the parameter with
trivial Weil part and standard $\SL_2$ factor. Its unipotent is
regular and its enhancement is trivial. Ordinary Springer support
in $\GL_2(\C)$ is the diagonal torus with trivial local system.
Thus its unnormalized support is the torus parameter $(1,1)$,
whereas its normalized support is
\[
 (|\cdot|^{1/2},|\cdot|^{-1/2})
 \quad\text{up to permutation}.
\]
The parameter with the same trivial Weil part and trivial
$\SL_2$ factor also has unnormalized support $(1,1)$, but its
normalized support is $(1,1)$. The two labels above the orbit of
$(1,1)$ in $(\C^\times)^2/S_2$ therefore need not have the same
normalized support. This proves directly that the quotient
projection cannot generally be replaced by $\Sc$.
\end{example}

\subsection{Unramified orbits and their inertial symmetries}\label{sec:orbits}

For a standard realization of $\cM$, put $D=Z(\hM)^{I_F}$ and
write $D_{\Fr}=D/(\Fr-1)D$. The weakly unramified character group
and the Bernstein unramified character torus are respectively
\begin{equation}\label{eq:Xnr}
 X_{\mathrm{w.nr}}(\cM)=D_{\Fr},\qquad
 \Xnr(\cM)=(D_{\Fr})^\circ .
\end{equation}
Indeed, evaluation at Frobenius identifies $H^1(\WF/I_F,D)$ with
$D_{\Fr}$. Inflation identifies this group with the central cocycle
classes admitting a representative trivial on $I_F$. Bernstein
inertial equivalence uses its identity component
\cite[\S3.3.1]{Haines}. The natural map
$(D^\circ)_{\Fr}\to(D_{\Fr})^\circ$ is an isogeny. Its kernel is
\[
 \bigl(D^\circ\cap(\Fr-1)D\bigr)/(\Fr-1)D^\circ,
\]
which can be nontrivial. For instance, take inertia to act trivially on
\[
 \hM=(\C^\times\times\SL_2(\C)\times\SL_2(\C))/
             \langle(-1,-I,-I)\rangle
\]
and let Frobenius exchange the two $\SL_2$ factors. Writing central
elements as $[(t,I,\epsilon I)]$ identifies $D$ with
$\C^\times\times\mu_2$. Frobenius sends $(t,\epsilon)$ to
$(\epsilon t,\epsilon)$, so the displayed kernel is
$\{(1,1),(-1,1)\}$. For a nonstandard realization we transport
\eqref{eq:Xnr} by conjugacy. A class in $\Xnr(\cM)$ twists $\psi$
and leaves $v$ and $\epsilon$ unchanged. Centrality preserves the
centralizer inside $\hM$, cuspidality, and relevance.

Fix $t_*=[\psi_*,v_*,\epsilon_*]_{\hM}$ and set
\begin{equation}\label{eq:selftwists}
 U=\Xnr(\cM),\quad
 U_{t_*}=\{a\in U:a t_*=t_*\},\quad T=U\cdot t_*.
\end{equation}

\begin{lemma}\label{lem:orbit-torus}
The group $U_{t_*}$ is finite. The orbit $T$ is a torsor under
$U/U_{t_*}$, and a choice of $t_*$ identifies it with that torus.
In particular it need not be a torsor under $U$ itself.
\end{lemma}
\begin{proof}
The self-twist group of the enhanced parameter is a subgroup of
the self-twist group of its underlying parameter. Finiteness of
the latter is \cite[Lemma~4.5.7]{DS}: projection to the abelian
quotient of $\hM$ forces a self-twist into the finite kernel of
the induced map of unramified tori. Since $U$ is abelian, all its
point stabilizers on this orbit equal $U_{t_*}$.
The quotient therefore acts freely and transitively. The quotient
of a complex torus by a finite subgroup is again a complex torus.
\end{proof}

Put $W(\LG{G^*},\cM)=N_{\hG}(\cM)/\hM$. The ordinary data underlying
the extended quotient are
\begin{equation}\label{eq:intro-TW}
 T=\Xnr(\cM)t_*,\qquad W=\{w\in W(\LG{G^*},\cM): wT=T\}.
\end{equation}
The group $W$ is finite. In coordinates
$U/U_{t_*}$ its action can have a translation part. Indeed, if
$w t_*=a_w t_*$, then
\begin{equation}\label{eq:affine-action}
 w(a t_*)=w(a)a_w t_*.
\end{equation}
The element $a_w$ is determined modulo $U_{t_*}$. The relation
$a_{ww'}=w(a_{w'})a_w$ holds in that quotient, as follows by
applying both sides to $t_*$. Thus replacing this action by the
linear action on $U$ loses data in general.

The ambient class $\svee=[\cM,t_*]$ consists of conjugates and
unramified twists of the cuspidal datum. Its component in the
dual Bernstein variety is $T/W$. Define
\begin{equation}\label{eq:block}
 \Phie^Z(\rig)^{\svee}
 =\{x\in\Phie^Z(\rig):[\Sc(x)]_{\mathrm{inert}}=\svee\}.
\end{equation}
By Proposition~\ref{prop:inertial-support}, $\Sc$ can be replaced
by $\uSc$ in this formula.

\begin{remark}[Ambient centralizers can jump]\label{rem:jump}
An unramified twist is central in $\hM$, not in general in $\hG$.
Consequently $Z_{\hG}(a\psi_*)$ can vary with $a$, even though
$Z_{\hM}(a\psi_*)$ is unchanged. For the diagonal torus of split
$\GL_2$, the centralizer of $\chi_1\oplus\chi_2$ is $\GL_2(\C)$
when the two Weil characters agree and is the diagonal torus otherwise.
In the unramified principal block, equality can be tested on Frobenius.
The algebras below must
therefore be formed at each point $t$. A single centralizer at
$t_*$ is not a model for every point of $T$.
\end{remark}

\section{The rigid Galois block and its Springer algebras}\label{sec:algebras}

We recall the rigid generalized Springer parametrization of
\cite[Theorems~4.8.7 and~5.2.5]{DS} in the common cover of
\S\ref{sec:rigid}. We use its actual induction object to identify the
multiplication and the simple-module transports under wild reduction.

For $t=[\psi_t,v_t,\epsilon_t]_{\hM}\in T$, set
\begin{equation}\label{eq:HKt}
 H_t^+=Z_{\hG}(\psi_t)^+,
 \qquad K_t^+=H_t^+\cap\hM^+.
\end{equation}
Let $\cE_t$ be the cuspidal local system on the $K_t^+$-orbit
of $v_t^+$ corresponding to $\epsilon_t$. Work with the
nilpotent model obtained by logarithm of unipotents. Choose a
quasi-parabolic $P_t^+$ with quasi-Levi $K_t^+$ and define the
Springer induction object and its algebra by
\begin{equation}\label{eq:intrinsic-A}
 \cF_t=\Ind_{K_t^+\subset P_t^+}^{H_t^+}
              \IC(\mathcal O_t,\cE_t),
 \qquad \cA_t=\End(\cF_t).
\end{equation}
The coefficients are in characteristic zero. One may use the
sheaf formalism of \cite{DS}, with a fixed identification of its
coefficient field with $\C$. Parabolic, Fourier, and Springer
conventions are fixed compatibly when coordinates are required.

\subsection{Left modules and multiplicity spaces}

We use left modules, with the map from an irreducible summand
$\mathcal L$ of $\cF_t$ given by
\begin{equation}\label{eq:multiplicity}
 \mathcal L\longmapsto\Hom(\mathcal L,\cF_t).
\end{equation}
The algebra acts by postcomposition. Semisimplicity of the
Springer induction object and its decomposition
$\cF_t=\bigoplus_\mathcal L\mathcal L\otimes V_\mathcal L$
give $\cA_t\simeq\bigoplus_\mathcal L\End(V_\mathcal L)$,
so \eqref{eq:multiplicity} is a bijection onto $\Irr(\cA_t)$.
Using $\Hom(\cF_t,\mathcal L)$ instead gives right modules and
requires an opposite-algebra conversion.

\begin{proposition}[The stabilizer identification]\label{prop:stabilizer}
There is an identification
\begin{equation}\label{eq:Springer-W}
 W_t\simeq
 N_{H_t^+}(K_t^+,\mathcal O_t,\cE_t)/K_t^+.
\end{equation}
Under the generalized Springer construction, $\cA_t$ is a twisted
group algebra of this group. The cocycle class is independent of
the nonzero scalar choices of homogeneous intertwining operators.
\end{proposition}
\begin{proof}
A representative of an element of $W_t$ can be adjusted by
$\hM$ so that it fixes $\psi_t$ and preserves the conjugacy
class of $(v_t,\epsilon_t)$ inside $K_t^+$. Lift the adjusted
element through the central isogeny to obtain the right side of
\eqref{eq:Springer-W}. The kernel is $K_t^+$: downstairs an
element with trivial relative Weyl class lies in $\hM$, and its
centralizer condition puts it in $H_t^+\cap\hM^+$ upstairs.
Conversely an element of the normalizer on the right normalizes
$Z(K_t^+)^\circ$, hence its centralizer $\cM$ in the $L$-group,
and fixes $t$. These maps are inverse. This is the construction
of \cite[Lemma~4.7.2]{DS}.

The algebra assertion is \cite[Proposition~3.4.4 and
Lemma~3.4.5]{DS}, applied to \eqref{eq:intrinsic-A}, with the
left-module convention above. Choose homogeneous units $e_w$,
$e_1=1$, in such a presentation. Define
$e_we_{w'}=\kappa_t(w,w')e_{ww'}$. Associativity is precisely the
cocycle identity. Replacing $e_w$ by $b(w)e_w$ changes
$\kappa_t$ by the coboundary $\delta b$.
\end{proof}

\subsection{Transport and the extended quotient}

For $g\in W$, conjugation of the parameter and induction data,
together with the compatible Springer identifications, transports
simple $\cA_t$-modules to simple $\cA_{gt}$-modules. We retain
these maps as part of the quotient data. In twisted-algebra
coordinates they can be chosen as in \cite[Lemma~4.8.6]{DS};
on simple spectra they compose exactly, and elements of $W_t$
act trivially. Changes of object identifications produce inner
automorphisms and do not affect these simple classes.

These data give the extended quotient of Definition~\ref{def:TEQ},
with $A_t=\cA_t$.

The intrinsic description of a fiber is $\Irr(\cA_t)$.
Writing it as $\Irr(\C_{\kappa_t}[W_t])$ uses a presentation of
$\cA_t$; even two cochains with the same coboundary can differ
by a character and hence permute projective labels. The standard
Springer normalization removes part, but not automatically all,
of this ambiguity. The character ambiguity on the disconnected
relative quotient is recorded in \cite[Proposition~4.7.1]{DS}.

\subsection{The block parametrization for a fixed rigid form}\label{sec:main-proof}

\begin{theorem}[Dillery--Schwein]
\label{thm:main}
There are twisted stabilizer algebras
$\cA_t\simeq\C_{\kappa_t}[W_t]$, with coherent transports of their
simple modules, and a bijection
\begin{equation}\label{eq:main}
 \mathfrak b^\vee_{\svee}:
 \Phie^Z(\rig)^{\svee}\xrightarrow{\ \sim\ }
 (T\extquot W)_\kappa.
\end{equation}
The following statements hold.
\begin{enumerate}[label=\textup{(\roman*)},leftmargin=*]
\item The projection of \eqref{eq:main} to $T/W$ is the
\emph{unnormalized} cuspidal-support map $\uSc$. Consequently
\[
 \uSc^{-1}([t])\simeq\Irr(\cA_t).
\]
\item Every parameter in the displayed block is relevant for the
fixed rigid form $\rig$. Its rigid central character is retained.
\item The bijection identifies bounded parameters with quotient
labels whose cuspidal parameter is bounded.
\item The induction endomorphism algebras and their multiplicity modules
give the parametrization before any bases of intertwining lines are
chosen. Twisted-group-algebra coordinates use the specified Springer
and projective-label conventions.
\end{enumerate}
\end{theorem}

\begin{proof}[Proof of Theorem~\ref{thm:main}]
First work without fixing a rigid form. Suppose
$x=[\psi,u,\rho]$ has unnormalized support in $\svee$. Choose
a representative $t\in T$ of that support and conjugate the
support datum to $(\cM,\psi_t,v_t,\epsilon_t)$. The remaining
unipotent pair is a simple constituent in the induction series
of \eqref{eq:intrinsic-A}. Send it to its multiplicity module
\eqref{eq:multiplicity}, and take its $W$-orbit. Different
representatives of the support are related by $W$; the transports
of Section~\ref{sec:algebras} identify the resulting labels.

Conversely, given $[t,V]$, decompose $\cF_t$ and take the unique
simple constituent corresponding to $V$. It gives a pair
$(u^+,\rho)$ in $H_t^+$. Its image $u$ in $\hG$, together with
$\psi_t$, defines an enhanced parameter with support $[t]$.
An $\SL_2$ enrichment recovers the usual parameter model if
desired. The two procedures are inverse, by the Springer
correspondence and the multiplicity-space decomposition. This is
\cite[Theorem~4.8.7]{DS} in the present notation. Its first
coordinate is $\uSc$ by construction.

Now take $t_*$ relevant for the rigid Levi $\mathfrak m$.
Unramified twisting preserves its rigid character. Its character
on $C_M^Z$ pulls back to $\zeta_\rig$. The central-character
assertion of Proposition~\ref{prop:inertial-support} shows that
every simple constituent produced above satisfies
\eqref{eq:relevance}. Conversely every relevant parameter in
the block was already included in the construction. This proves
the fixed-rigid-form statement. Existence of the corresponding
rational rigid Levi for an arbitrary relevant parameter is
\cite[Proposition~5.2.4]{DS}; the global restriction is their
Theorem~5.2.5.

The inverse construction leaves the Weil homomorphism $\psi_t$
unchanged. Boundedness, understood for that Weil part, is
therefore preserved. Finally, \eqref{eq:intrinsic-A} and
\eqref{eq:multiplicity} give the assertion about intrinsic algebras;
Proposition~\ref{prop:stabilizer} and the transport discussion give
the stated dependence on coordinates.
\end{proof}

\begin{corollary}[Rigid Bernstein decomposition]\label{cor:all-blocks}
Choose one representative of each relevant inertial class for
$\rig$. Then
\begin{equation}\label{eq:all-blocks}
 \Phie^Z(\rig)\simeq
 \bigsqcup_{\svee}(T_{\svee}\extquot W_{\svee})_{\kappa_{\svee}}.
\end{equation}
One may first group the classes by the rational Levi subgroup of
$G$ carrying their cuspidal support.
\end{corollary}
\begin{proof}
Cuspidal support exists uniquely up to dual conjugacy. Relevance
realizes its supporting Levi in $G$, with the compatible rigid
structure, by \cite[Proposition~5.2.4]{DS}. Its inertial orbit
is then unique up to the corresponding relative Weyl action.
Theorem~\ref{thm:main} applies to each disjoint class. This is
also the block refinement of \cite[Corollary~5.2.6]{DS}.
\end{proof}

\begin{proposition}[The corrected map on the quotient]\label{prop:corrected-map}
There is a well-defined map
\begin{equation}\label{eq:p-corrected}
 p_{\mathrm{norm}}^\vee:(T\extquot W)_\kappa\longrightarrow T/W,
 \qquad
 p_{\mathrm{norm}}^\vee=\Sc\circ(\mathfrak b^\vee_{\svee})^{-1}.
\end{equation}
For a label $[t,V]$, reconstruct its unipotent pair by Springer
theory and apply \eqref{eq:normalized-support}. This computes
$p_{\mathrm{norm}}^\vee$. In general this map depends on $V$
and does not factor through $\pr$.
\end{proposition}
\begin{proof}
The composite is well defined since the two maps in it are.
The reconstruction in the proof of Theorem~\ref{thm:main},
followed by \eqref{eq:normalized-support}, gives the specified
calculation. Example~\ref{ex:GL2} gives two labels with the same
value of $\pr$ and distinct values of $p_{\mathrm{norm}}^\vee$.
\end{proof}

\begin{remark}[Scope of the theorem]
Theorem~\ref{thm:main} is a parametrization of a set of enhanced
parameters. It does not assert an equivalence with a category of
smooth representations, an isomorphism of their Bernstein centers,
or equality with an analytic $R$-group algebra. In particular the
group $W_t$ is a Springer relative stabilizer and must not be
replaced by an analytic $R$-group without a separate comparison.
\end{remark}

\subsection{Changes of cover and central twists}\label{sec:functoriality}

We record two operations for which the algebras themselves can be
compared. They are useful when bringing a rigid seed construction
into the present conventions.

\begin{proposition}[Enlarging the finite central subgroup]\label{prop:enlarge-Z}
Let $Z\subset Z'\subset G^*$ be finite central subgroup schemes.
Inflate an enhancement along
\[
 \widehat{G^*/Z'}\longrightarrow\widehat{G^*/Z}.
\]
Both support maps commute with inflation. On a block obtained by
inflating its cuspidal datum, the induction algebras are isomorphic,
and their simple-module parametrizations can be made compatible.
This compares the inflated sector, not all enhancements at level $Z'$.
\end{proposition}
\begin{proof}
The homomorphism induces central surjections on the inverse-image
centralizers and on the quasi-Levis. Their unipotent varieties are
identified. Pullback identifies equivariant sheaves downstairs with
the full subcategory upstairs on which the central kernel acts
trivially. Indeed, in characteristic zero the equivariance data
descend precisely when this kernel acts trivially; this description
also proves full faithfulness.

The parabolic-induction correspondences are compatible with this
pullback. Thus the pullback of $\cF_t$ is the induction object
for the inflated cuspidal datum, and full faithfulness gives an
isomorphism of their endomorphism algebras. Every constituent of
the pullback remains kernel-trivial. The normalizer stabilizers
are the same, since the central kernel is already in the
quasi-Levi and pullback is fully faithful. Transport by conjugation
commutes with pullback, so the simple-module identifications are
coherent. Support compatibility is exactly the argument of
\cite[\S4.9]{DS}; the correcting cocharacter is compatible with
the central isogeny. Rigid duality identifies the inflated central
character with the same rigid class viewed at level $Z'$.
\end{proof}

\begin{proposition}[Twists central in the ambient dual group]\label{prop:ambient-twist}
Let $a:\WF\to Z(\hG)$ be a central cocycle for which twisting
preserves admissibility. The operation
$(\psi,u,\rho)\mapsto(a\psi,u,\rho)$ identifies the relevant
Springer algebras and support constructions, carrying a block to
the block of its twisted cuspidal datum. For fixed compatible
choices it transports the corresponding extended quotients.
\end{proposition}
\begin{proof}
For $g\in\hG$, the equality
$g(a(w)\psi(w))g^{-1}=a(w)\psi(w)$ is equivalent to
$g\psi(w)g^{-1}=\psi(w)$, since $a(w)$ is central. Hence
the full ordinary centralizers, their inverse images, and all
unipotent induction data are unchanged. Relative normalizers
and their action on local systems are unchanged as well.
The induction objects and their multiplication algebras may
therefore be identified directly. The same correcting
cocharacter works after multiplying both infinitesimal characters
by $a$, proving compatibility with normalized support too.
\end{proof}

The ambient-central hypothesis is essential. A twist central only
in the supporting Levi need not preserve the ambient centralizer;
see Remark~\ref{rem:jump}. Nor does changing the rigidification
on a fixed underlying inner form mean that enhancement labels
remain unchanged. Such a change is governed by the rigid character
transport, as in \cite[\S5.4]{DS}, and must be included in any
claimed natural comparison.

\subsection{Blocks with cuspidal support supplied by Yu data}\label{sec:yu}

The following application uses the cuspidal parameters constructed in
\cite{YuParameters}.

For each rigid Levi $\mathfrak m$ let $\mathscr Y(\mathfrak m)$
be a specified collection of Galois-side Yu data for which the
construction gives a relevant cuspidal rigid parameter
\[
 \operatorname{par}_{\mathfrak m}(D)\in\Phic^Z(\mathfrak m).
\]
The data include a common finite central subgroup, compatible
dual covers, full quotient characters, relative rigid realization
records, and the reference normalization prescribed in
\cite{YuParameters}. Let $\mathcal C_{\mathrm{Yu}}$ be the
unramified and normalizer saturation of these outputs, over all
rigid Levis of $\rig$.

\begin{definition}\label{def:Yu-blocks}
Let $\mathfrak B_{\mathrm{Yu}}^Z(\rig)$ be the set of ambient
inertial classes represented by $\mathcal C_{\mathrm{Yu}}$.
The Yu-supported parameter set is
\[
 \Phie^Z(\rig)_{\mathrm{Yu}}
 =\{x:[\Sc(x)]_{\mathrm{inert}}\in
                 \mathfrak B_{\mathrm{Yu}}^Z(\rig)\}.
\]
The saturation ensures that this is a union of whole blocks,
independently of the chosen presentation of an output.
\end{definition}

\begin{corollary}[The Galois quotient for Yu-supported blocks]
\label{cor:Yu-blocks}
For every constructed cuspidal output, its rigid Galois block has
the parametrization of Theorem~\ref{thm:main}. In particular
\begin{equation}\label{eq:Yu-blocks}
 \Phie^Z(\rig)_{\mathrm{Yu}}
 \simeq\bigsqcup_{\svee\in\mathfrak B_{\mathrm{Yu}}^Z(\rig)}
 (T_{\svee}\extquot W_{\svee})_{\kappa_{\svee}}.
\end{equation}
If the cuspidal construction is exhaustive for every relevant
rigid Levi, then \eqref{eq:Yu-blocks} is \eqref{eq:all-blocks}.
\end{corollary}
\begin{proof}
Apply Theorem~\ref{thm:main} to each output and its unramified
orbit, and discard repeated inertial classes. If every relevant
cuspidal parameter is an output, the support of any enhanced
parameter belongs to this union by existence of cuspidal support.
\end{proof}

The cuspidal exhaustion result of \cite[Corollary~11.4]{YuParameters}
applies to its strict quotient-character data in the tame prime-to-Weyl
range with fixed reference normalization. When applied to each relevant
Levi, it gives the last assertion of Corollary~\ref{cor:Yu-blocks}
in that range.

\begin{remark}[What ``tame'' means here]
Yu-supported refers to the origin of the cuspidal support in a
tame twisted-Levi construction. It does not mean that the
parameter is trivial on wild inertia; positive-depth parameters
usually are not. The prime restriction belongs to the cuspidal
construction or its exhaustion theorem, not to Dillery--Schwein's
block parametrization.
\end{remark}

\section{Wild reduction and cuspidal support}\label{sec:reduction}
\label{sq:scope}

We construct the parameter-side square of \eqref{eq:two-squares}.
For a fixed extension of the wild parameter and an intrinsic tame
cuspidal enhanced seed, its four maps are
\begin{equation}\label{sq:main-square}
\begin{CD}
 \Phie^Z(\rig)^{\svee}
 @>{\mathfrak b_{\svee}^{\vee}}>> \scrE_{\svee}^{\vee}\\
 @V{\red}V{\sim}V
 @VV{(f_{\calD}^{\vee},\iota_{\calD}^{\vee},a_{\calD})_*}V\\
 \Phie^Z(\rigH)^{\ssH}
 @>{\mathfrak b_{\ssH}^{\vee}}>> \scrE_{\ssH}^{\vee}.
\end{CD}
\end{equation}
Here $H$ plays the role of $G^0$; $M\subset G$ and $M_0\subset H$
will be the respective \emph{supporting parabolic Levis}. In particular,
$H$ must not be confused with the initial twisted Levi inside $M$.
The horizontal maps are the rigid generalized Springer parametrizations
of \cite[Theorem~4.8.7]{DS}, with compatible geometric conventions.
Their first coordinates record $\uSc$, the unnormalized support.

The left arrow will have the explicit formula
\begin{equation}\label{sq:reduction-formula}
 (\psi,u,\rho)\longmapsto
 \bigl(\varkappa^{-1}j^{-1}(\psi),u,\rho\bigr)
\end{equation}
after fixing the wild restriction. The enhancement in this formula is
transported through an equality of groups in a common cover. Thus the
left arrow is constructed before using either horizontal parametrization.
The extra symbol $a_{\calD}$ on the right records actual algebra
isomorphisms; it can be suppressed only after compatible homogeneous
coordinates have been fixed.

We first identify the groupoids of parameters and their supporting
Levis, then compute the maps on the enhanced orbits and the full Springer
algebras. The fixed-form statement is Theorem~\ref{sq:theorem}.
Theorem~\ref{sq:Sigma-matching} identifies the smaller group and its
rigid realization for the marked towers of Definition~\ref{sq:Sigma-input}.

\subsection{Wild reduction data and an explicit equivalence}
\label{sq:wild}

Let $G^*$ be tame and quasi-split over a nonarchimedean local field $F$,
with dual group $\hG$ and pinned $L$-group $\LG{G^*}$; write $G$ for a
specified inner form when necessary. We use the commuting-unipotent
model $(\psi,u)$, where $\psi:\WF\to\LG{G^*}$ is admissible and
$u\in Z_{\hG}(\psi)$ is unipotent. The homomorphism $\psi$ is the
restriction of an $\mathrm{SL}_2$-parameter to $\WF$, before inserting
any $|w|^{1/2}$ diagonal factor.

Fix an extendible finite-image homomorphism
$\lambda:\PF\to\hG$. Since $G^*$ is tame, this is the dual-group part
of its wild parameter. Suppose
\begin{equation}\label{sq:wild-Levi}
             \hH=Z_{\hG}(\lambda(\PF))
             \quad\text{is a connected Levi subgroup.}
\end{equation}
Choose one extension $\psi_*:\WF\to\LG{G^*}$ and put
$\cE_\lambda=\hH\psi_*(\WF)$. Normality of $\PF$ makes this a group;
its extension structure is independent of $\psi_*$. Indeed, two
extensions having the same actual wild restriction differ, at each
$w$, by an element of $\hH$.

\begin{definition}[Wild reduction datum]\label{sq:datum}
A datum $\calD$ consists of the following specified choices:
\begin{enumerate}[label=\textup{(D\arabic*)},leftmargin=*]
\item A tame quasi-split $H^*$ and an admissible tame $L$-embedding
$j:\LG{H^*}\to\LG{G^*}$ identifying $\widehat{H^*}=\hH$ and having image
$\cE_\lambda$. Tameness means $j(1,p)=(1,p)$ for $p\in\PF$.
\item A continuous full cocycle
$\varkappa\in Z^1(\WF,Z(\hH))$, for the action defined by $j$, with
$\varkappa|_{\PF}=\lambda$. Its full tame and unramified parts are retained.
\item A compatible common central cover
\[
 \pZ:\hG^+=\widehat{G^*/Z}\longrightarrow\hG,
 \qquad \hH^+=\pZ^{-1}(\hH)=\widehat{H^*/Z},
\]
with the indicated root-data marking. For fixed actual inner forms,
include the relative rigid realization described in
Section~\ref{sq:rigid}.
\end{enumerate}
The embeddings of dual groups will henceforth be suppressed.
\end{definition}

\begin{proposition}[Availability of the parameter-theoretic datum]
\label{sq:existence}
If $p\nmid |W_{\mathrm{abs}}(G)|$, condition~\eqref{sq:wild-Levi}
holds for every extendible wild parameter. A tame splitting as in
\textup{(D1)} exists. One can choose the cocycle in \textup{(D2)} with
finite $p$-power image in $Z(\hH)^\circ$.
\end{proposition}
\begin{proof}
The full-centralizer assertion is \cite[Lemma~2.2.9(iv)]{DF}; the tame
splitting is \cite[Lemma~2.2.4]{DF}. The central cocycle construction
also appears in \cite[Lemma~2.2.5 and \S2.2.6]{DF}. For completeness, the
prescribed-restriction argument for the last assertion is as follows.

The image of $\lambda$ lies in $Z(\hH)$: it centralizes its full
centralizer, and $Z_{\hG}(\hH)=Z(\hH)$ for a Levi subgroup.
Thus this follows directly
from~\eqref{sq:wild-Levi}. The order of $\pi_0Z(\hH)$ is prime to $p$
in the stated range. Hence the image lies in
$A[p^\infty]$, where $A=Z(\hH)^\circ$.
It is equivariant for the action given by $j$, since that embedding
has image $\cE_\lambda$.

Give $A[p^\infty]$ the discrete topology. The group $\IF/\PF$ is
pro-prime-to-$p$, so positive-degree continuous cohomology with these
coefficients vanishes, by averaging on finite quotients. The quotient
$\WF/\IF\simeq\Z$ has cohomological dimension one. Consequently
$H^2(\WF/\PF,A[p^\infty])=0$. Inflation--restriction extends
$\lambda$ to a cocycle $\varkappa$; coboundaries restrict trivially to
$\PF$, so its restriction is the specified homomorphism itself.
Its inertia values are finite by continuity. These values, its value
at Frobenius, and their translates under the finite action generate a
finite $p$-subgroup containing every value of $\varkappa$.
See also \cite[Lemma~9.2]{YuParameters} for this extension argument.
\end{proof}

The proposition concerns the dual data. A full common-cover marking and
an actual rational embedding in a specified rigid form must be retained
when the theorem is expressed for that form. Theorems~5.1 and~5.6 of \cite{YuParameters} give a rational reference
and relative realization for the marked towers considered below.

\begin{lemma}[Cancellation, including all conjugators]
\label{sq:groupoid}\label{lem:inertia-cut}
For a datum $\calD$, the formula
\begin{equation}\label{sq:forward}
       J_{\calD}(\psi_0)(w)=\varkappa(w)j(\psi_0(w))
\end{equation}
identifies the groupoid of tame parameters for $\LG{H^*}$ with the groupoid
of parameters for $\LG{G^*}$ having actual wild restriction $\lambda$.
Morphisms in the latter are conjugations that preserve this condition.
For every such parameter, with $\psi=J_{\calD}(\psi_0)$,
\begin{align}
 Z_{\hG}(\psi)&=Z_{\hH}(\psi_0),\label{sq:Weil-centralizer}\\
 Z_{\hG}(\psi,u)&=Z_{\hH}(\psi_0,u).
 \label{sq:unipotent-centralizer}
\end{align}
The assertions hold for full algebraic groups, and for their full
inverse images under $\pZ$.
\end{lemma}
\begin{proof}
The cocycle identity makes~\eqref{sq:forward} a homomorphism and its
wild restriction is $\lambda$. Conversely any extension of $\lambda$
has image in $\cE_\lambda=j(\LG{H^*})$. Pull it back and divide by
$\varkappa$; the resulting parameter has trivial wild restriction.
These operations are inverse. Central twisting preserves admissibility:
a power of a Frobenius image which kills the finite outer action is
multiplied by a commuting central semisimple element. Semisimplicity is
therefore preserved in both directions, as are the continuity and
inertia conditions. The fixed admissible embedding preserves these
conditions as well.

Any element of $\hG$ conjugating two parameters with actual wild
restriction $\lambda$ must centralize $\lambda(\PF)$, hence lies in
$\hH$. Within $\hH$, the factors $\varkappa(w)$ are central and cancel
from the conjugacy equations. This proves equality of all transporter
sets, in particular~\eqref{sq:Weil-centralizer}. A commuting unipotent
lies in this group; keep it unchanged to obtain
\eqref{sq:unipotent-centralizer}. Taking full inverse images proves the
cover assertions. Equivalently one may retain the same commuting
$\mathrm{SL}_2$-homomorphism throughout.
\end{proof}

Every parameter with wild restriction \emph{conjugate} to $\lambda$
can be put in this slice. Two such choices differ by an element of
$\hH$. Thus the lemma descends to conjugacy classes without a residual
quotient by a normalizer of the wild image. The actual homomorphism
$\lambda$, rather than just its image subgroup, is essential here.

\subsection{Cuspidal supports and their ambient Levi envelopes}
\label{sq:Levis}

To compare cuspidal support we must identify its ambient Levi. Write
$\cM_0\subset\LG{H^*}$ for a Levi $L$-subgroup, put
$\hM_0=\cM_0\cap\hH$, and set
\begin{equation}\label{sq:envelope}
 A_0=(Z(\hM_0)^{\WF})^\circ,
 \qquad \cM=Z_{\LG{G^*}}(A_0),\qquad \hM=Z_{\hG}(A_0).
\end{equation}
The action on the center is that defined by $\cM_0$.
Since $\cM_0=Z_{\LG{H^*}}(A_0)$, the expression for $\cM$ is a Levi
$L$-subgroup: its fibers over $\WF$ are nonempty, supplied by
$j(\cM_0)$. Notice also that $Z(\hH)\subset\hM_0$.

\begin{lemma}[Levi envelopes and invariant centers]\label{sq:centers}
The groups in~\eqref{sq:envelope} satisfy
\begin{equation}\label{sq:intersection}
 \hM\cap\hH=\hM_0,\qquad
 \cM\cap j(\LG{H^*})=j(\cM_0),\qquad
 (Z(\hM)^{\WF})^\circ=A_0.
\end{equation}
The map $J_{\calD}$ carries parameters in $\cM_0$ into $\cM$.
For these parameters it identifies their full centralizers inside
$\hM_0$ and $\hM$, respectively. It preserves discreteness between
these two supporting Levis.
\end{lemma}
\begin{proof}
The intersection identities follow by centralizing $A_0$. Both
$j(\cM_0)$ and $\varkappa$ centralize $A_0$, proving the factorization.
The torus $A_0$ is central in $\hM$ and fixed by a section of $\cM$,
so it lies in $(Z(\hM)^{\WF})^\circ$. Conversely $\lambda(\PF)$ lies
in $\hM_0\subset\hM$. Thus $Z(\hM)$ centralizes $\lambda$ and is
contained in $\hH$, hence in $Z(\hM_0)$. Its invariant connected torus
is therefore contained in $A_0$. This proves the last equality.
Choose the same commuting $\mathrm{SL}_2$ enrichment on both sides.
Intersect~\eqref{sq:Weil-centralizer} with its centralizer and with
$\hM$. This compares the full-parameter centralizers inside the two
supporting Levis. Discreteness is finiteness of this reductive
centralizer modulo the invariant center; the connected invariant
centers just proved equal give equivalence of this condition.
One must use the $\mathrm{SL}_2$ centralizer here: the centralizer of
the single unipotent can have a nontrivial unipotent radical.
\end{proof}

Let $Q^+=Z_{\hG}(\psi)^+=Z_{\hH}(\psi_0)^+$. Rigid enhancements
are local systems on unipotent orbits in $Q^+$, equivalently
irreducible representations of $\pi_0(Z_{Q^+}(u^+))$; the superscript
$+$ denotes full inverse image in the common cover. We use the
Dillery--Schwein quasi-cuspidal support throughout.

\begin{proposition}[Compatibility of cuspidal support]\label{sq:support}
The map $\inc$ induced by Lemma~\ref{sq:groupoid} identifies cuspidal
data in $\cM_0$ with their cuspidal data in the envelope $\cM$.
For arbitrary enhanced parameters in the wild slice it carries
unnormalized support to the corresponding envelope support.
Every cuspidal datum occurring in the wild slice is obtained in this
way, up to conjugacy.
\end{proposition}
\begin{proof}
For a cuspidal datum in $\cM_0$, the common centralizer and local
system give the same cuspidality condition on both sides. Discreteness
is Lemma~\ref{sq:centers}. This proves the first assertion.

For a general parameter, its quasi-cuspidal support is a triple
$(K^+,v^+,\epsilon)$ in the \emph{same} group $Q^+$ on both sides.
Put $S=\pZ(Z(K^+)^\circ)$. The two ambient support Levis are
$Z_{\LG{G^*}}(S)$ and $Z_{\LG{H^*}}(S)$, respectively, by
\cite[Construction~4.6.2 and Lemma~4.6.4]{DS}.
Their intersection with the common centralizer is the same $K^+$.
To see that these are the envelopes of~\eqref{sq:envelope}, denote
the second Levi by $\cM_0$. Choose a Jacobson--Morozov homomorphism
$\gamma:\mathrm{SL}_2\to(K^+)^\circ$ for $v^+$ and project it
through $\pZ$. It commutes with $S$, since $S$ comes from the center
of $K^+$. Its support parameter is discrete, hence
\[
 (Z(\hM_0)^{\WF})^\circ
    = Z_{\hM_0}(\psi_0,\pZ\gamma(\mathrm{SL}_2))^\circ.
\]
The torus $S$ is contained in this group. Conversely the invariant
center, lifted through the central isogeny on its identity component,
is central in $K^+$, so is contained in $S$. Thus $S=A_0$.
This gives precisely~\eqref{sq:envelope} and proves compatibility.

Conversely apply the inverse of Lemma~\ref{sq:groupoid} to a parameter
in the wild slice and compute its quasi-cuspidal support in $Q^+$.
The preceding construction recovers its ambient support. This proves
exhaustion of the cuspidal data which occur; it does not assert that
an arbitrary parabolic Levi of $G$ is the envelope of an arbitrarily
chosen Levi of $H$.
\end{proof}

\section{Orbits and Springer induction under wild reduction}
\label{sq:orbits}

Fix a tame cuspidal enhanced parameter
$t_{0,*}=[\psi_{0,*},v,\epsilon]_{\hM_0}$ and its image
$t_*=[\varkappa j\psi_{0,*},v,\epsilon]_{\hM}$.
Let $\ssH=[\cM_0,t_{0,*}]_{\hH}$ and
$\svee=[\cM,t_*]_{\hG}$. Put
\[
 U_M=\Xnr(\cM)=\bigl((Z(\hM)^{\IF})_{\Fr}\bigr)^\circ,\qquad
 U_0=\Xnr(\cM_0)=\bigl((Z(\hM_0)^{\IF})_{\Fr}\bigr)^\circ.
\]
Write
$T=U_Mt_*$ and $T_0=U_0t_{0,*}$ for the enhanced unramified orbits.
Their self-twist subgroups are denoted $U_{M,t_*}$ and $U_{0,t_{0,*}}$.
The inclusion of centers in Lemma~\ref{sq:centers} induces
$r:U_M\to U_0$.

\begin{proposition}[The orbit isomorphism, with self-twists]
\label{sq:orbit}
The map $r$ is an isogeny, and
\begin{equation}\label{sq:selftwists}
 U_{M,t_*}=r^{-1}(U_{0,t_{0,*}}).
\end{equation}
The self-twist groups therefore fit into the exact sequence
\begin{equation}\label{sq:selftwist-sequence}
 1\longrightarrow\ker r\longrightarrow U_{M,t_*}
   \xrightarrow{\,r\,} U_{0,t_{0,*}}\longrightarrow 1 .
\end{equation}
In particular there is an algebraic isomorphism
\begin{equation}\label{sq:f}
 f_{\calD}^{\vee}:T\xrightarrow{\sim}T_0,
 \qquad a\,t_*\longmapsto r(a)t_{0,*}.
\end{equation}
It is the restriction of the independently defined reduction map.
\end{proposition}
\begin{proof}
The actions on the two centers factor through finite groups. Over
$\C$, invariants and coinvariants for a finite action are identified
on Lie algebras by averaging. Thus the differential of $r$ is the
map on $\WF$-invariant Lie algebras of the centers, which is an
isomorphism by~\eqref{sq:intersection}. A map of complex tori with
isomorphic differential is an isogeny.

Twisting the formula~\eqref{sq:forward} by an unramified central
cocycle of $\cM$ commutes with reduction and gives $r(a)$ on the
source. Lemma~\ref{sq:groupoid}, intersected with $\hM$, says that
an isomorphism of the enhanced parameters on the large side exists
exactly when the corresponding isomorphism on the small side exists.
This proves~\eqref{sq:selftwists}. Since $r$ is surjective, restricting
it to this inverse image gives~\eqref{sq:selftwist-sequence}.
Quotienting the isogeny by these self-twist groups
gives an isomorphism
$U_M/U_{M,t_*}\simeq U_0/U_{0,t_{0,*}}$, proving the assertion about
the orbit torsors. The self-twist groups are finite; see
\cite[Lemma~4.5.7]{DS}.
\end{proof}

Write the subgroups preserving the enhanced unramified orbits as
\[
 N_G^T=\{n\in N_{\hG}(\cM):nT=T\},
 \qquad
 N_H^{T_0}=\{h\in N_{\hH}(\cM_0):hT_0=T_0\}.
\]
The full inertial Weyl groups are $W=N_G^T/\hM$ and
$W_0=N_H^{T_0}/\hM_0$. Their actions on the enhanced orbits can
have translation terms in coordinates centered at the respective seeds.

\begin{proposition}[Exhaustion of the inertial normalizers]
\label{sq:normalizer}
Put $\mathcal N=N_H^{T_0}$. Inside $\hG$ one has
\[
 N_G^T=\hM\,\mathcal N,\qquad \hM\cap\mathcal N=\hM_0.
\]
Thus every class in $W$ has a representative in $\mathcal N$, and
passing to the quotients gives an isomorphism
\begin{equation}\label{sq:iota}
 \iota_{\calD}^{\vee}:W\xrightarrow{\sim}W_0
\end{equation}
for which $f_{\calD}^{\vee}$ is equivariant, including its affine
translation terms. Consequently point stabilizers correspond.
\end{proposition}
\begin{proof}
An element $h\in\mathcal N$ preserves
$A_0=(Z(\hM_0)^{\WF})^\circ$ as a set, so normalizes its ambient
centralizer $\cM$. Centrality of $\varkappa$ in $\hH$ makes conjugation
by $h$ commute with $J_{\calD}$. Surjectivity of $r$ and
Proposition~\ref{sq:orbit} then show that $hT=T$.
The kernel of the resulting map $\mathcal N\to W$ is
$\mathcal N\cap\hM=\hM_0$.

For surjectivity, represent $w\in W$ by $n\in N_{\hG}(\cM)$.
Since $n$ preserves $T$, choose $a\in U_M$ and $m\in\hM$ such that
$mn$ carries a representative of $t_*$ to the representative of
$a t_*$, with its enhancement. Unramified twisting does not change
the actual wild restriction. Therefore
\[
              h=mn\in Z_{\hG}(\lambda(\PF))=\hH.
\]
It still normalizes $\cM$, hence normalizes
$\cM\cap j(\LG{H^*})=j(\cM_0)$, and reduction shows that it preserves
$T_0$. Thus $h\in\mathcal N$ and represents $w$.
Since $n=m^{-1}h$, this proves $N_G^T=\hM\,\mathcal N$ and
identifies both inertial Weyl groups with $\mathcal N/\hM_0$.

Equivariance follows from the same common representatives and the
cancellation formula. Explicitly, if $w t_*=a_wt_*$, then
$\iota(w)t_{0,*}=r(a_w)t_{0,*}$, with the elements understood modulo
self-twists. Hence the translation part is transported as well.
\end{proof}

The proof uses the entire enhanced orbit and adjusts an inertial
intertwiner inside its supporting Levi. A comparison of Weyl groups
of the connected centralizers at a single basepoint would not prove
this proposition.

\subsection{Springer induction objects and their coherent transports}
\label{sq:Springer}

For $t\in T$, choose a representative with fixed wild restriction
$\lambda$, and let $t_0=f_{\calD}^{\vee}(t)$. Put
\begin{align*}
 Q_t^+&=Z_{\hG}(\psi_t)^+
           =Z_{\hH}(\psi_{0,t_0})^+,\\
 K_t^+&=Q_t^+\cap\hM^+
           =Q_t^+\cap\hM_0^+.
\end{align*}
The cuspidal unipotent and its local system $\cE_t$ are the same on
both sides. In a nilpotent model choose a common quasi-parabolic
$P_t^+$ with quasi-Levi $K_t^+$ and set
\begin{equation}\label{sq:ind-object}
 \cF_t=\Ind_{K_t^+\subset P_t^+}^{Q_t^+}
                   \IC(\cO_t,\cE_t),\qquad
 \cA_t=\End(\cF_t).
\end{equation}
The conventions for induction, Fourier transform where used, and
Springer labels are fixed to be those of \cite{DS}. In particular no
identification with differently normalized AMS labels is being assumed.

We use left modules: a simple summand $\mathcal L$ corresponds to
$\Hom(\mathcal L,\cF_t)$, with endomorphisms acting by postcomposition.
Write $\cA_{G,t}$ and $\cA_{H,t_0}$ when keeping both sides visible.

\begin{proposition}[Comparison of actual stabilizer algebras]
\label{sq:algebras}
The common induction correspondence gives isomorphisms
\begin{equation}\label{sq:a}
 a_t:\cA_{G,t}\xrightarrow{\sim}\cA_{H,f_{\calD}^{\vee}(t)}.
\end{equation}
They respect the homogeneous stabilizer labels through
$\iota_{\calD}^{\vee}|_{W_t}$ and commute with orbit transports on
simple modules. They identify the summand multiplicity spaces.
\end{proposition}
\begin{proof}
The groups, quasi-parabolics, orbits and local systems in
\eqref{sq:ind-object} are identical in the common cover. Thus the
induction diagrams and their induced objects are identical. Their
endomorphisms, with composition, give~\eqref{sq:a}.
If the chosen presentation uses the Fourier-associated cuspidal local
system, apply the same Fourier construction in the common Lie algebra;
that associated system agrees on both sides as well. No assertion that
it equals the original local system is needed.
The geometric relative normalizer is
\[
 N_{Q_t^+}(K_t^+,\cO_t,\cE_t)/K_t^+.
\]
Its identification with the inertial point stabilizer
\cite[Lemma~4.7.2]{DS} uses an adjusted normalizer representative
fixing the Weil homomorphism. Proposition~\ref{sq:normalizer} uses
the same representatives on the two sides. Therefore the
identifications with $W_t$ and $(W_0)_{t_0}$ agree with $\iota$.

For transport from $t$ to $wt$, take a common representative in
$\mathcal N$, adjust it in $\hM_0$ to the chosen representative at
$wt$, and lift it to $\hH^+$. Denote the lift by $n$. Conjugation
carries $P_t^+$ to $nP_t^+n^{-1}$, which can differ from $P_{wt}^+$
even when $wt=t$.

In matched twisted group algebra coordinates, choose the coherent
transports of \cite[Lemma~4.8.6]{DS} once for the common Springer data
on each orbit. Part~(3) of that lemma identifies the Springer summand
at $t$ with its transported class at $wt$. The common representative
acts on the same unipotent pairs and local systems in $Q_t^+$ on both
sides. Thus the chosen transport gives the Springer transport for
each side, including when the chosen quasi-parabolic changes.
Part~(2) gives the composition law, and part~(1) makes a stabilizer
loop act by an inner automorphism, hence trivially on simple classes.
Writing $\tau$ for these transports, we obtain
\begin{equation}\label{sq:transport}
 a_{wt}\tau^G_{w,t}
       =\tau^H_{\iota(w),f(t)}a_t
       \quad\text{on simple spectra}.
\end{equation}
Finally, the common direct-sum decomposition of $\cF_t$ identifies
the multiplicity space $\Hom(\mathcal L,\cF_t)$ for each summand
$\mathcal L$.
\end{proof}

In homogeneous coordinates write
$\cA_{G,t}\simeq\C_{\kappa_{G,t}}[W_t]$ and
$\cA_{H,f(t)}\simeq\C_{\kappa_{H,f(t)}}[(W_0)_{f(t)}]$.
Here $\kappa$ denotes the Springer multiplier, reserving $\varkappa$ for
the full Weil cocycle. If
$a_t(e_w)=c_t(w)e_{\iota(w)}$, then multiplication gives
\begin{equation}\label{sq:cochain}
 \kappa_{G,t}(w,v)=
 \kappa_{H,f(t)}(\iota(w),\iota(v))
           \frac{c_t(w)c_t(v)}{c_t(wv)}.
\end{equation}
Matched geometric bases permit $c_t=1$. For independently prescribed
bases the cochains must be those obtained from $a_t$ and
\eqref{sq:transport}. Equality of the two cohomology classes alone
would leave the projective-label map undetermined.

\section{The Galois reduction square and its rigid sectors}
\label{sq:rigid}

For $J=G,H$ write
\[
 C_J^Z=\pi_0\bigl(\pZ^{-1}(Z(\widehat J)^{\WF})\bigr).
\]
A rigid form $\mathfrak j$ has a relevance character
$\zeta_{\mathfrak j}$ on this group. Relevance means equality of
characters after the natural central map to the refined component
group; that map need not be injective. These are the conventions of
\cite[\S5.4]{KalethaRigid} and \cite[\S5.2]{DS}.

The embedding $j$ and the common cover give a map
$C_G^Z\to C_H^Z$. At the cuspidal seed, the common enhancement
$\epsilon$ determines a character $\zeta_0$ of $C_H^Z$; it pulls back
to the ambient character $\zeta_{\rig}$. Retain the full character
$\zeta_0$. If fixed actual forms are desired, supply a rigid form
$\rigH$ with character $\zeta_0$ and a marked relative rigid
realization of $H\hookrightarrow G$ compatible with these dual data.
The supporting Levis and their rigidifications are the compatible
ones. This is a concrete realization requirement, not an additional
cocycle equality in the Springer calculation.

\begin{lemma}[The source rigid character is constant on a block]
\label{sq:relevance}
Every enhancement in the block $\svee$, reduced by
Lemma~\ref{sq:groupoid}, has central character $\zeta_0$ on $C_H^Z$.
Conversely every enhancement in $\ssH$ has ambient character
$\zeta_{\rig}$. Thus reduction restricts to the indicated matched
rigid sectors.
\end{lemma}
\begin{proof}
The image of $\pZ^{-1}(Z(\hH)^{\WF})$ is central in every common
$Q_t^+$. It acts by the same character on the cuspidal local system
and each simple summand of its induction object: the central group
acts trivially on the induction correspondence, and its action on
coefficients is the fixed scalar character. Equivalently use
\cite[Lemma~4.6.8]{DS}. Unramified twisting changes neither this
character nor the action on $Z(\hH)$, and conjugation by common
representatives in $\hH$ fixes its center. Thus the character is
constant on the whole block. Its pullback along $C_G^Z\to C_H^Z$
is the fixed ambient character, proving both assertions.
\end{proof}

For clarity, write the intrinsic extended quotients as
\[
 \scrE_{\svee}^{\vee}
  =\left(\coprod_{t\in T}\Irr(\cA_{G,t})\right)\!/W,
 \qquad
 \scrE_{\ssH}^{\vee}
  =\left(\coprod_{t_0\in T_0}\Irr(\cA_{H,t_0})\right)\!/W_0,
\]
where the quotients use the coherent orbit transports. In homogeneous
coordinates these are the usual twisted extended quotients.
Write $\Phie^Z(\LG{J^*})$ for enhanced parameter classes before
restricting to a specified rigid sector. In this reduction statement,
``depth zero'' means trivial wild restriction in the tame $L$-group.

\begin{theorem}[The rigid Galois depth-zero reduction square]
\label{sq:theorem}\label{thm:reduction}\label{thm:intro-square}
Fix a datum $\calD$, a tame cuspidal seed $t_{0,*}$ and its envelope
as above. The explicit reduction~\eqref{sq:reduction-formula} gives a
bijection of enhanced Galois Bernstein blocks
\begin{equation}\label{sq:R}
 \red:\Phie^Z(\LG{G^*})^{\svee}
                   \xrightarrow{\sim}\Phie^Z(\LG{H^*})^{\ssH}.
\end{equation}
The target block consists of tame parameters. The right vertical map
in~\eqref{sq:main-square} is
\begin{equation}\label{sq:EQ-map}
           [t,V]\longmapsto
           [f_{\calD}^{\vee}(t),V\circ a_t^{-1}].
\end{equation}
With matched Dillery--Schwein conventions, the square commutes and all
its arrows are bijections. With the compatible relative rigid
realization it restricts to the fixed-form square
\eqref{sq:main-square}. The refined groups are identified by
\[
 S_{\psi_0,u}^{+}\xrightarrow{\sim}S_{\psi,u}^{+},
 \qquad A_{\psi_0,u}^{+}\xrightarrow{\sim}A_{\psi,u}^{+}
\]
throughout; no assertion about equality of intrinsic AMS component
groups is required.
\end{theorem}
\begin{proof}
Every parameter in $\svee$ has wild restriction conjugate to
$\lambda$: its unnormalized support has the same Weil homomorphism,
and that support is conjugate to an unramified twist of $t_*$.
Choose a representative in the wild slice and apply
Lemma~\ref{sq:groupoid}. Its class and enhancement are independent
of this choice. Proposition~\ref{sq:support} and the orbit isomorphism
identify its support block with $\ssH$. The inverse construction
$\inc$ gives the reverse inclusion and proves~\eqref{sq:R}.
The Weil homomorphism in $\ssH$ is tame because its cuspidal support
is tame and unnormalized support retains that homomorphism.

Propositions~\ref{sq:orbit}, \ref{sq:normalizer}, and
\ref{sq:algebras} show that~\eqref{sq:EQ-map} is a well-defined
bijection, including orbit transports. Compute a horizontal map by
choosing the support $t$ and taking the multiplicity space of its
simple Springer summand in $\cF_t$. Reduction leaves that summand
and the entire induction object unchanged; it changes its algebra
coordinate by $a_t$. Consequently the resulting left module is
$V\circ a_t^{-1}$ and its support is $f_{\calD}^{\vee}(t)$.
This is commutativity, not a definition of the left arrow through
the other three arrows. Lemma~\ref{sq:relevance} gives the fixed-form
restriction. The refined-group assertion is
\eqref{sq:unipotent-centralizer} in the common cover.
\end{proof}

\begin{remark}[Fixed ambient relevance does not determine all source sectors]
For the union of blocks with wild class $[\lambda]$, a fixed character
on $C_G^Z$ may admit several characters on $C_H^Z$. Reduction then
has targets in the corresponding union of source rigid sectors.
For one enhanced Bernstein block the seed fixes the character, as
Lemma~\ref{sq:relevance} proves. This does not allow one to discard
its relative embedding record or identify different lifts of the
ambient character.
\end{remark}

\subsection{Normalized support and changes of the full cocycle}
\label{sq:choices}

\begin{proposition}[The support correction also transports]
\label{sq:normalized}
For the square of Theorem~\ref{sq:theorem}, both $\uSc$ and
$\Sc$ commute with the corresponding map on support data. The
ordinary projection of the extended quotient is still $\uSc$.
\end{proposition}
\begin{proof}
The assertion for $\uSc$ was proved above. The correcting
cocharacter for $\Sc$ is constructed in $Z(K_t^+)^\circ$ from
the unipotent and its support, as in \cite[Lemma~4.6.6]{DS}.
Choose the same construction in the common group $Q_t^+$.
Its image lies in the torus $A_0$ of the support envelope, hence in
the invariant centers of both supporting Levis. Evaluation at the
prescribed power of $q$ therefore gives correcting cocycles whose
classes correspond under $r$. Multiplication by $\varkappa$ commutes
with them. This proves normalized-support compatibility, including
the unchanged enhancement. The correction depends on the Springer
summand and need not be determined by the unnormalized orbit alone;
it does not change the meaning of the quotient projection.
\end{proof}

For $\lambda=1$, $j=1$ and $\varkappa=1$, the square is the
identity square. Example~\ref{ex:GL2} still distinguishes its normalized
and unnormalized support maps. Example~\ref{ex:GL4-wild} checks a
nontrivial wild reduction at points with jumping stabilizers.

\begin{proposition}[Changing the full central extension]
\label{sq:change}
Fix $j$ and replace $\varkappa$ by $\varkappa'=b\varkappa$, where
$b\in Z^1(\WF,Z(\hH))$ is tame. Then
\[
             R_{(j,\varkappa')}^+=b^{-1}\otimes R_{(j,\varkappa)}^+.
\]
The two squares are identified by this central twist on the $H$ side,
with the same induction algebras and their coherent transports.
A coboundary change gives conjugate parameters. A general tame change
need not preserve the same source inertial block.
\end{proposition}
\begin{proof}
The formula follows by division. Since $b$ is central in the ambient
$\hH$, it changes neither its Weil centralizers nor the local systems
or induction correspondences. Conjugation commutes with the twist,
so the same assertion holds for transports. A central coboundary is
implemented by conjugation by its defining element. If the image of
$[b]$ in $H^1(\WF,Z(\hM_0))$ belongs to $\Xnr(\cM_0)$, twisting
preserves the seed's enhanced unramified orbit. A class trivial on
inertia can lie in another component of $X_{\mathrm{w.nr}}(\cM_0)$
and change the block, as in Example~\ref{ex:weak-torus}.
\end{proof}

\subsection{Application to a marked ambient Yu tower}\label{sq:yu}

Suppose a parameter-side tower has been supplied \emph{in the ambient
$G$}, with compatible tame embeddings, full central cocycles, and the
same cover. Write its terminal composite as
\begin{equation}\label{sq:tower-composite}
        \psi_\Sigma=\varkappa_\Sigma j_\Sigma(\psi_0),
       \qquad \psi_0|_{\PF}=1,
       \qquad \varkappa_\Sigma\in Z^1(\WF,Z(\widehat{G^0})).
\end{equation}
The composite cocycle is the product of the transported stage
cocycles. Their values lie in $Z(\widehat{G^0})$, and the
homomorphism identity for the composite proves the cocycle identity.
All full tame and unramified factors, including a terminal central
factor, are included.

\begin{lemma}[A stagewise test for the full wild centralizer]
\label{sq:tower-test}
Assume $\lambda=\psi_\Sigma|_{\PF}$ centralizes
$\widehat{G^0}$. Suppose that, after each preceding cancellation,
the leading ramification image of the current stage has full
centralizer equal to the next smaller dual Levi. Then
\[
                 Z_{\hG}(\lambda)=\widehat{G^0}.
\]
If $j_\Sigma$ has the indicated compatible outer action, its image
is $\cE_\lambda$, so~\eqref{sq:tower-composite} supplies
\textup{(D1)--(D2)}.
\end{lemma}
\begin{proof}
An element centralizing $\lambda$ centralizes its highest leading
image and therefore lies in the first smaller dual Levi. Inside
that group the full central factor cancels from the wild
centralizer equations. Apply the next cut-down to the residual
wild homomorphism and continue. This gives containment in
$\widehat{G^0}$; the assumed centrality gives the reverse inclusion.
For the extension image, both $j_\Sigma(\LG{H^*})$ and $\cE_\lambda$
have connected dual group $\widehat{G^0}$. The factorization supplies
one common element over each $w$, since $\varkappa_\Sigma(w)$ and the
dual part of $\psi_0(w)$ belong to that group. Their fibers over
$\WF$ therefore agree.
\end{proof}

For a tower originally constructed only in the supporting $M$, the
lemma's ambient cut-down is an additional assertion. One can instead
start directly with the ambient wild centralizer of the cuspidal
parameter, apply Proposition~\ref{sq:existence}, and use the envelope
construction of Section~\ref{sq:Levis}. Theorem~\ref{sq:theorem} then
already gives an intrinsic parameter-side reduction square. To label
its lower group and choices by a prescribed marked tower $\Sigma$,
use the comparison below, which applies to its full quotient-character
data and intrinsic enhanced seed.

\section{Marked ambient quotient-character towers}
\label{sq:marked-Sigma}

We now specialize the Galois square to the full quotient-character
input of \cite{YuParameters}. Fix a rational
tame twisted Levi tower, its supporting parabolic Levis and an
intrinsic enhanced seed. Towers with the geometric configuration of
Kim--Yu/AFMO are included. The full quotient characters encode central
Weil cocycles through the torus correspondence; no supercuspidal
representation is part of this input. Ellipticity of the ambient
inclusions is not assumed. The source rigid form is the prescribed
actual subgroup when the seed is relevant to the induced rigidification
described below.

\subsection{The input and one simultaneous rigid presentation}

Write $H=G^0$, and fix an actual rational tame tower
\begin{equation}\label{sq:actual-tower}
        G^0\subsetneq G^1\subsetneq\cdots\subsetneq G^d=G.
\end{equation}
Let $M^0$ be a parabolic $F$-Levi of $G^0$, let $A$ be the maximal
$F$-split subtorus of $Z(M^0)^\circ$, and put
\begin{equation}\label{sq:actual-supports}
                 M^i=Z_{G^i}(A),\qquad M=M^d.
\end{equation}
Require $M^0$ to be elliptic in each $M^i$. Each $M^i$ is a
parabolic $F$-Levi of $G^i$ by construction. This is the supporting-Levi
configuration in \cite[\S4.1, before Definition~4.1.1]{AFMO-reduce}.
The sequence $M^0\subseteq\cdots\subseteq M^d$ can have repetitions.
For example, if $G^0=M^0=T$ is the split diagonal torus of
$G^1=\GL_2$, then $A=T$ and $M^1=T$. A positive-depth
$G^1$-generic character of $T(F)$ still contributes its full factor
to the supporting parameter at this step.
Assume $p\ne2$ and $p\nmid|W_{\mathrm{abs}}(G)|$.
The positive proper-stage depths satisfy
$0<r_0<\cdots<r_{d-1}$. A final character of $G(F)$ is retained
separately, including when its depth equals or exceeds $r_{d-1}$.
For $d=0$ there is only that central twist.

\begin{definition}[A marked quotient-character tower datum]
\label{sq:Sigma-input}
Fix the tower~\eqref{sq:actual-tower}, the supporting groups
\eqref{sq:actual-supports}, and the given ambient rigid inner twist
$(\xi,z_G):G^*\to G$, using a common finite central subgroup $Z$.
For each $0\leq i\leq d$, supply a full character $\bar\chi_i$ of
$D_{G^i}(F)=(G^i/G^i_{\mathrm{der}})(F)$ and put
\[
                   \chi_i=\bar\chi_i\circ q_{G^i}.
\]
For $i<d$, require $\chi_i$ to be $G^{i+1}$-generic of depth $r_i$.
Supply an intrinsic tame cuspidal enhanced seed in $M^0$, relevant to
its compatible induced rigidification, and keep the relative markings.
The resulting datum will be denoted $\Sigma^+$.

The full quotient characters are specified, including their tame and
unramified parts. If only their positive-filtration restrictions are
initially given, completing them is an additional choice. Its effect
on the central cocycle and enhanced seed is described in
Lemma~\ref{sq:Sigma-lift-change}; it is not discarded from $\Sigma^+$.
\end{definition}

\begin{proposition}[The prescribed tower in one rigid presentation]
\label{sq:Sigma-reference}
Choose a tame maximal $F$-torus $S\subset M^0$. There is an inner
modification $\Xi=\xi\circ\Ad(g)$ such that
\[
 K_i=\Xi^{-1}(G^i),\quad L_i=\Xi^{-1}(M^i),\quad
 S^*=\Xi^{-1}(S)
\]
are rational, all $K_i,L_i$ are quasi-split, and
$\Xi|_{S^*}$ is an $F$-isomorphism. The single rigid cocycle
\begin{equation}\label{sq:Sigma-z}
             z^g(w)=g^{-1}z_G(w)\sigma(g)
\end{equation}
takes values in $S^*$, has its original band in $Z$, and rigidifies all
these restrictions. The relative record for $K_0\hookrightarrow G^*$
is $(z^g,g)$ and realizes exactly the given $G^0\hookrightarrow G$.
The full quotient characters transport through the induced rational
isomorphisms $D_{K_i}\simeq D_{G^i}$.
\end{proposition}
\begin{proof}
Apply the simultaneous-reference theorem \cite[Theorem~5.6]{YuParameters} to the
nested chain $S\subset M^0\subset G^0\subset\cdots\subset G$,
omitting repetitions. Its underlying input is
\cite[Lemma~6.4]{KalethaCovers}, applied successively from the top of
the chain downwards. It gives the rational quasi-split references and
the rational restriction on $S^*$. Since
$\Ad(z^g)=\Xi^{-1}\sigma(\Xi)$ fixes $S^*$ pointwise, the identity
$Z_{G^*}(S^*)=S^*$ puts the cocycle in $S^*$ itself. Conjugation
fixes the common $Z$, including its gerbe-band map.

The preimage $A^*=\Xi^{-1}(A)$ is an $F$-split torus. Therefore
$L_i=Z_{K_i}(A^*)$ is a parabolic Levi in the quasi-split $K_i$ and
is quasi-split. It contains $S^*$, so the same cocycle restricts to
it. Equation~\eqref{sq:Sigma-z} is exactly the relative descent
identity for the given embedding and the original ambient twist.
In positive characteristic it is interpreted on a common flat
trivialization of the rigid gerbe; the scheme-theoretic centralizer
identity gives the same descent. Finally inner automorphisms act
trivially on abelian quotients, so the induced quotient isomorphisms
are rational and transport the supplied full characters.
\end{proof}

This proposition supplies one full marking of the dual root-data
rectangle and its common covers. It retains the actual rational
embedding, not only its stable conjugacy class. If a reference normalization $\mathfrak N$ has already been
fixed as in \cite{YuParameters}, transport
this presentation to $\mathfrak N$ with its enhancement and relative
record; reinitializing the enhancement while changing reference is not
one of the operations used here.

\subsection{The supporting-Levi rectangle and tame embeddings}

Use the dual marking from Proposition~\ref{sq:Sigma-reference}, writing
$\widehat{G^i}$ and $\widehat{M^i}$ for the duals of the reference
groups. The actual and reference abelian quotients have already been
identified. The following explains how the ambient and supporting
constructions can use the same embeddings.

\begin{lemma}[A common supporting torus]\label{sq:Sigma-envelope}
There are compatible realizations of the Levi $L$-subgroups
$\cM_i\subset{}^LG^i$ and a common torus
\[
             A^\vee=(Z(\widehat{M^0})^{\WF})^\circ
\]
such that
\begin{equation}\label{sq:Sigma-envelope-equalities}
 \cM_i=Z_{{}^LG^i}(A^\vee),\qquad
 (Z(\widehat{M^i})^{\WF})^\circ=A^\vee.
\end{equation}
In particular $\cM_d$ is the envelope of $\cM_0$ from
\eqref{sq:envelope}, and
$\widehat{M^i}=\widehat{G^i}\cap\widehat M$.
Tame $L$-embeddings
$j_i:{}^LG^i\to{}^LG^{i+1}$ can be chosen with this connected
marking; their restrictions give compatible tame embeddings
$\cM_i\to\cM_{i+1}$.
\end{lemma}
\begin{proof}
The common root-data marking identifies the dual inclusions for
$M^0\subset M^i$. The inclusions of invariant connected dual centers
are equalities because $M^0$ is elliptic in $M^i$: their dimensions
are the split ranks of the respective connected centers, and these
ranks agree by~\eqref{sq:actual-supports}. This is the fixed-center criterion of \cite[Lemma~4.6]{YuParameters}. Denote this common embedded torus by $A^\vee$.
For a parabolic Levi $L$-subgroup, its centralizer description using
its invariant connected center gives the first equality in
\eqref{sq:Sigma-envelope-equalities}. This also proves the intersection
assertion and identifies the ambient envelope. To make the compatibility explicit, in each quasi-split $K_i$ choose
a pinned maximal torus inside $L_i$ and a compatible Borel, with $L_i$
standard. Conjugate this torus to $S^*$ by an element of $L_i$ over
$F^s$. Such conjugation fixes $A^*$; its Weyl corrections belong to
the Weyl group of $L_i$. Dually those corrections act trivially on
$A^\vee$. Thus, in the common root datum on $S^*$, the based Weil
actions of every $K_i$ fix the same $A^\vee$ pointwise. This gives
the claimed compatible centralizer descriptions. It does not require
a simultaneously invariant Borel for all the twisted Levi inclusions.

Existence of compatible tame embeddings for the connected marked
inclusions follows from \cite[\S6.1, Corollary~6.7]{KalethaCovers}
with tame $\chi$-data. The root-field extensions are tame and
$p\ne2$, so the usual minimally ramified data are trivial on wild
inertia. Choose the pinned realizations so that the supporting Levis
are standard, as in the preceding centralizer description.
An element $x\in\cM_i$ centralizes $A^\vee$, and $j_i$ is the
identity on that marked torus. Consequently $j_i(x)$ centralizes
$A^\vee$ and lies in $\cM_{i+1}$. The restriction is an
$L$-embedding and is tame. This proves compatibility directly;
it does not identify an independently chosen supporting-Levi
embedding without its change-of-embedding correction.
\end{proof}

For each $i$ let
\[
  [\varkappa_i]=\mathcal L_{D_{G^i}}(\bar\chi_i)
           \in H^1(\WF,Z(\widehat{G^i})^\circ)
\]
and choose a representative. Write $J_i=j_{d-1}\circ\cdots\circ j_i$
for $i<d$, $J_d=1$, and let $\beta_i$ be its connected dual map.
In these coordinates set
\begin{equation}\label{sq:Sigma-product}
 j_\Sigma=J_0,\qquad
 \varkappa_\Sigma(w)=\prod_{i=0}^{d}\beta_i(\varkappa_i(w)).
\end{equation}
For $d=0$ use $j_\Sigma=1$. Every factor on the right lies in
$Z(\widehat{G^0})$: the center of a larger member of a nested Levi
tower is contained in every smaller member and centralizes it.

\begin{lemma}[The full composite and supporting restriction]
\label{sq:Sigma-product-lemma}
The product in~\eqref{sq:Sigma-product} is a central cocycle for the
$\WF$-action on $Z(\widehat{G^0})$ defined by $j_\Sigma$.
Given a tame parameter $\psi_0$ of $\cM_0$, define parameters
$\psi_i$ of $\cM_i$ successively by
\[
 \psi_{i+1}(w)=j_i|_{\cM_i}\bigl(\varkappa_i(w)\psi_i(w)\bigr)
 \qquad (0\leq i<d).
\]
The restricted composite, including the terminal twist, is
\begin{equation}\label{sq:Sigma-forward}
   \psi_\Sigma=\varkappa_d\psi_d=\varkappa_\Sigma j_\Sigma(\psi_0).
\end{equation}
The same formula works for arbitrary tame parameters of ${}^LG^0$.
Its wild restriction
$\lambda_\Sigma=\varkappa_\Sigma|_{\PF}$ is independent of $\psi_0$.
\end{lemma}
\begin{proof}
The action of $j_\Sigma(1,w)$ on
$\beta_i(Z(\widehat{G^i}))$ agrees with that of $J_i(1,w)$:
the difference is an element of $\beta_i(\widehat{G^i})$ and hence
acts trivially on its center. Thus each transported factor satisfies
the cocycle identity for the composite action. They commute, proving
the product identity. Expanding the successive construction gives
\eqref{sq:Sigma-forward}. Its restriction to the supporting tower
uses the embeddings of Lemma~\ref{sq:Sigma-envelope}. Moreover,
$\chi_i|_{M^i(F)}$ is the pullback of $\bar\chi_i$ along
$D_{M^i}\to D_{G^i}$, so functoriality of the torus correspondence
gives exactly the same factor through
$Z(\widehat{G^i})^\circ\to Z(\widehat{M^i})^\circ$.
If $M^i=M^{i+1}$, the restricted embedding is an isomorphism over
$\WF$ because it identifies the connected dual groups and each
fiber over $\WF$. The recursion still applies this isomorphism
to the full twist $\varkappa_i\psi_i$.
Finally all embeddings and $\psi_0$ are tame, giving the stated wild
restriction. This keeps the full terminal factor $\varkappa_d$.
\end{proof}

\subsection{The ambient cut-down and the exact extension image}

\begin{proposition}[Genericity identifies the prescribed $G^0$]
\label{sq:Sigma-cutdown}
For the marked quotient-character refinement,
\begin{equation}\label{sq:Sigma-cutdown-equality}
                  Z_{\hG}(\lambda_\Sigma)=\widehat{G^0}.
\end{equation}
For any parameter~\eqref{sq:Sigma-forward} one also has
\begin{equation}\label{sq:Sigma-extension}
 j_\Sigma({}^LG^0)
   =\widehat{G^0}\psi_\Sigma(\WF)
   =\cE_{\lambda_\Sigma}.
\end{equation}
These are full ambient assertions; no ellipticity of $G^0$ in $G$
is needed.
\end{proposition}
\begin{proof}
The prime-to-$p$ central-isogeny filtration calculation gives
$q_{G^i}(G^i(F)_{x,s})=D_{G^i}(F)_s$ for $s>0$. Consequently
the supplied quotient character has the same positive depth as
$\chi_i$; tame torus depth preservation gives depth $r_i$ for
$\varkappa_i$ when $i<d$. For a proper stage $i$, put
$B_i=\varkappa_i(I_F^{r_i})$ in
$Z(\widehat{G^i})^\circ$. Tame torus functoriality and the leading
filtration comparison identify the nontrivial root characters of
$B_i$ in $\widehat{G^{i+1}}$ with the normed-coroot characters
specified by $G^{i+1}$-genericity of $\chi_i$. Hence the connected
centralizer of $B_i$ is exactly $\widehat{G^i}$. The group $B_i$ is
finite of $p$-power order, so the prime-to-Weyl-order full-centralizer
result makes its full centralizer connected. Therefore
\[
          Z_{\widehat{G^{i+1}}}(B_i)=\widehat{G^i}.
\]
This is the root calculation of \cite[Lemma~4.4 and
Proposition~4.5]{YuParameters}. It works for the ambient tower here,
since ellipticity enters the separate discreteness assertion.

Remove $\varkappa_d$ first. Its values are central in $\hG$, so this
does not change the centralizer of the wild homomorphism. On
$I_F^{r_{d-1}}$ all earlier factors are trivial, because their depths
are strictly smaller; the embeddings are tame. The preceding
cut-down therefore puts every wild centralizer element in
$\widehat{G^{d-1}}$. Inside this group $\varkappa_{d-1}$ is central
and can be cancelled from the wild centralizer equations. Continue
downwards. This proves containment in $\widehat{G^0}$. The reverse
containment follows because every factor of $\lambda_\Sigma$ lies
in its center. For $d=0$ the whole wild image is ambient-central and
the same conclusion holds.

Formula~\eqref{sq:Sigma-forward} shows that
$\psi_\Sigma(w)$ and $j_\Sigma(1,w)$ differ by an element of
$\widehat{G^0}$ at every $w$. Their extension fibers therefore agree.
Combining this with~\eqref{sq:Sigma-cutdown-equality} and the
definition of $\cE_\lambda$ proves~\eqref{sq:Sigma-extension}.
In particular no second, independently normalized $L$-embedding is
needed to obtain the required extension image.
\end{proof}

\subsection{The prescribed rigid sector and the matching theorem}

Let $z_0$ and $z_{M^0}$ denote the one cocycle $z^g$ of
Proposition~\ref{sq:Sigma-reference} viewed in $K_0$ and $L_0$,
respectively. Their induced rigid forms are the actual $G^0$ and
$M^0$. Denote their relevance characters by $\zeta_0$ and
$\zeta_{M^0}$. Choose the intrinsic seed
$t_{0,*}=[\psi_{0,*},v,\epsilon]$ to be relevant to $z_{M^0}$.
This is the compatibility required in Definition~\ref{sq:Sigma-input}.

\begin{proposition}[Selection of the actual source form]
\label{sq:Sigma-rigid}
The source character selected by the common enhancement in
Section~\ref{sq:rigid} is exactly $\zeta_0$, and its relative
realization can be taken to be the actual marked inclusion
$G^0\hookrightarrow G$ with record $(z_0,g)$.
This conclusion retains the full character on $C_{G^0}^Z$.
\end{proposition}
\begin{proof}
There are compatible central maps
\[
               C_G^Z\longrightarrow C_{G^0}^Z
                         \longrightarrow C_{M^0}^Z.
\]
The pushforward of $z_{M^0}$ to $K_0$ is literally $z_0$;
its further pushforward is cohomologous to $z_G$ through $g$.
Functoriality of rigid duality therefore gives
\[
 \zeta_{M^0}|_{C_{G^0}^Z}=\zeta_0,
 \qquad \zeta_0|_{C_G^Z}=\zeta_{\rig},
\]
where restriction means pullback along these maps, not an assumed
inclusion of component groups. Relevance of $\epsilon$ to
$z_{M^0}$ gives the first character as its scalar action. Hence the
character extracted from the seed on $C_{G^0}^Z$ is $\zeta_0$.
Lemma~\ref{sq:relevance} propagates this equality across the block.
The relative cocycle identity~\eqref{sq:Sigma-z} already realizes the
prescribed subgroup, so no new rational embedding is selected.
The relative-realization theorem \cite[Theorem~5.1]{YuParameters} identifies its class
of marked reductions with this full source character, for the fixed
reference and ambient twist. Agreement only after pullback to $C_G^Z$
would not imply this conclusion.
\end{proof}

\begin{theorem}[Matching the marked ambient quotient-character datum]
\label{sq:Sigma-matching}
For $\Sigma^+$ of Definition~\ref{sq:Sigma-input}, compatible tame
embeddings and the retained full cocycles produce a wild reduction
datum $\calD_\Sigma$ for which all of the following hold:
\begin{enumerate}[label=\textup{(\roman*)},leftmargin=*]
\item Its smaller group, supporting Levi, common cover and relative
rigid realization are the prescribed $G^0$, $M^0$, $Z$ and marked
inclusions of the actual datum.
\item Its wild parameter and full central factor are
$\lambda_\Sigma$ and $\varkappa_\Sigma$ of
\eqref{sq:Sigma-product}. Equations
\eqref{sq:Sigma-cutdown-equality}--\eqref{sq:Sigma-extension} hold.
\item Its ambient supporting Levi is the prescribed $M$. For the
intrinsic seed $[\psi_{0,*},v,\epsilon]_{\widehat{M^0}}$, its cuspidal
enhanced parameter in $M$ is $[\psi_\Sigma,v,\epsilon]_{\widehat M}$,
where $\psi_\Sigma$ is the restricted composite
\eqref{sq:Sigma-forward} starting at $\psi_{0,*}$ and the enhancement
is transported through the common cover. The composite includes
every index of the given ambient tower.
\item The square~\eqref{sq:main-square} holds with
$\rigH$ equal to this actual rigid $G^0$ and with left arrow
$R_{\calD_\Sigma}^+$. Its right arrow uses the full normalizer and
Springer-algebra comparisons already proved in
Propositions~\ref{sq:normalizer} and~\ref{sq:algebras}.
\end{enumerate}
The seed in this theorem is the given enhanced parameter. Thus the
conclusion is commutativity of the Galois-side square with the prescribed
smaller group, supporting Levi and relative rigid realization. It
requires no depth-zero supercuspidal LLC.
\end{theorem}
\begin{proof}
Proposition~\ref{sq:Sigma-reference} supplies the actual marked
reference and common covers. Lemma~\ref{sq:Sigma-envelope} identifies
the supporting Levi envelope and compatible restrictions of the
embeddings. Lemma~\ref{sq:Sigma-product-lemma} and
Proposition~\ref{sq:Sigma-cutdown} verify \textup{(D1)--(D2)} of
Definition~\ref{sq:datum}, with the full factors of the given refinement.
Proposition~\ref{sq:Sigma-rigid} verifies the prescribed source sector
and relative part of \textup{(D3)}. Lemma~\ref{sq:Sigma-product-lemma}
identifies the restricted Weil parameter with the image of the seed
under $J_{\calD_\Sigma}$. Proposition~\ref{sq:support}, applied to
$\cM_0$ and its envelope $\cM_d$, transports $v$ and $\epsilon$ and
gives cuspidality in $M$. This proves \textup{(iii)} for the stated
composite, including equal successive supporting Levis.
When all the supporting inclusions are strict, this is also the
construction of \cite{YuParameters} with the restricted embeddings and quotient
characters. Indeed, for a common tame maximal torus in $M^0$, every
root of $M^{i+1}$ outside $M^i$ is a root of $G^{i+1}$ outside $G^i$,
since $M^{i+1}\cap G^i=M^i$. Its coroot has the same pairing with
the restricted leading functional. Thus the restricted character
has depth $r_i$ and satisfies the genericity criterion at each
proper step by \cite[Lemma~4.4]{YuParameters}.
Now apply Theorem~\ref{sq:theorem}. This proves all four assertions.
\end{proof}

\subsection{Changes of quotient characters and their lifts}

\begin{lemma}[Wild independence and the remaining tame factor]
\label{sq:Sigma-lift-change}
Keep the actual tower, connected marking and tame embeddings fixed.
Suppose two systems of full quotient characters agree on
$D_{G^i}(F)_{0+}$ for every $i$. Their wild homomorphisms
$\lambda_\Sigma$ are equal, while their full composite cocycles differ
by a tame central cocycle $b$ of $\widehat{G^0}$. Their reductions
satisfy $R'^+=b^{-1}\otimes R^+$.
If the pullbacks agree on $G^i(F)$ for every $i$, then $b$ is a
coboundary in $Z^1(\WF,Z(\widehat{G^0}))$. In that case the two
reduction maps agree on enhanced conjugacy classes.
\end{lemma}
\begin{proof}
The ratio of the quotient characters is trivial on the positive
filtration, so its tame-torus parameter is trivial on $\PF$
\cite[Lemma~2.4.1]{DF}. The product formula makes the ratio
of the composite cocycles a tame cocycle in $Z(\widehat{G^0})$.
Apply Proposition~\ref{sq:change}.

For the second assertion, fix $J=G^i$ and let $\eta$ be the ratio of
the two characters of $D_J(F)$. We show that $\eta\circ q_J=1$
forces its torus parameter to have trivial image in
$H^1(\WF,Z(\widehat J))$. Choose a $z$-extension
$1\to S\to\widetilde J\to J\to1$ with $S$ induced. The map
$\widetilde J(F)\to D_{\widetilde J}(F)$ is surjective since
$\widetilde J_{\mathrm{der}}$ is simply connected over the
nonarchimedean local field $F$ \cite[Lemma~2.4.1 and
Remark~2.4.2]{DF}. Hence the pullback of $\eta$ to
$D_{\widetilde J}(F)$ is trivial. By functoriality for tori its class
maps to zero in $H^1(\WF,\widehat{D_{\widetilde J}})$.
On dual centers we have an exact sequence
\[
 1\longrightarrow Z(\widehat J)\longrightarrow
 Z(\widehat{\widetilde J})=\widehat{D_{\widetilde J}}
 \longrightarrow\widehat S\longrightarrow1 .
\]
The map on $\WF$-fixed points at the right is surjective. Indeed,
$\widehat S^{\WF}$ is a connected torus because $S$ is induced,
and taking invariants for the finite Galois action on the complex
Lie algebras preserves surjectivity. The cohomology sequence thus
gives an injection
$H^1(\WF,Z(\widehat J))\hookrightarrow
  H^1(\WF,Z(\widehat{\widetilde J}))$. This proves the claim.

Each stage ratio is therefore a coboundary in the full center
$Z(\widehat{G^i})$. The maps $\beta_i$ on these centers are
$\WF$-equivariant for the composite action by
Lemma~\ref{sq:Sigma-product-lemma}. Their product in
\eqref{sq:Sigma-product} is a coboundary in $Z(\widehat{G^0})$.
Its defining central element implements the conjugacy of the
reduced parameters. A lift to the common cover centralizes
$\widehat{G^0}^{+}$, so this conjugacy also preserves the
enhancement.
\end{proof}

Consequently, if the full cocycle changes from $\varkappa_\Sigma$ to
$b\varkappa_\Sigma$, the reduced seed attached to a fixed ambient parameter
changes from $t_{0,*}$ to $b^{-1}t_{0,*}$. Membership of the image
of $[b]$ in $\Xnr(\cM_0)$ suffices for the new seed to lie in the
same enhanced unramified orbit. Agreement on the positive filtration
alone allows the block to change, even when $b|_{\IF}=1$
(Example~\ref{ex:weak-torus}). Lifts of the same characters on all
$G^i(F)$ give the same enhanced seed class.

\subsection{Recovery of the enhanced cuspidal seed}\label{sq:seed-comparison}
\begin{corollary}[Seed recovery]\label{sq:seed-matching}
Fix a datum $\calD$, a tame cuspidal enhanced seed $t_{0,*}$ and its
image
\begin{equation}\label{sq:seed-calibration}
 t_*=J^+_{\calD}(t_{0,*}).
\end{equation}
Then
\begin{equation}\label{sq:seed-recovery}
 R^+_{\calD}(t_*)=t_{0,*},\qquad
 R^+_{\calD}(a t_*)=r(a)t_{0,*}\quad(a\in U_M).
\end{equation}
The action of the full inertial normalizer on the enhanced orbit
corresponds under reduction:
\begin{equation}\label{sq:seed-equivariance}
 f^\vee_{\calD}(wt)=\iota^\vee_{\calD}(w)f^\vee_{\calD}(t).
\end{equation}
This includes the affine translation terms, the self-twists in
\eqref{sq:selftwist-sequence}, and the projective labels in the
specified rigid sectors. For Theorem~\ref{sq:Sigma-matching} the seed
recovered is the one supplied with the marked tower.
\end{corollary}
\begin{proof}
The Weil homomorphism of $t_*$ is $\varkappa j(\psi_{0,*})$. Cancelling
$\varkappa$ and pulling back along $j$ recovers $\psi_{0,*}$, and the
common cover identifies the unipotent element and its enhancement.
Proposition~\ref{sq:orbit} gives the twist formula and the exact sequence
of self-twists. Proposition~\ref{sq:normalizer} gives equivariance using
common normalizer representatives, and Proposition~\ref{sq:algebras}
identifies the projective labels with their transports. The rigid
character is the one in Lemma~\ref{sq:relevance}. For the tower apply
\eqref{sq:Sigma-forward}.
\end{proof}

\begin{remark}[Independently specified seeds]\label{sq:seed-independent}
An independently specified enhanced cuspidal parameter $t$ in the
ambient wild slice is the image of $t_{0,*}$ precisely when
$R^+_{\calD}(t)=t_{0,*}$ as enhanced classes in the specified rigid
sector. Keeping $j$ fixed and replacing $\varkappa$ by $b\varkappa$
changes the reduced seed attached to $t$ by $b^{-1}$, as in
Proposition~\ref{sq:change}. The full cocycle therefore enters the
comparison of any prescribed seeds.
\end{remark}

\clearpage
\part{Comparison of the two parametrizations}

\section{Transfer of block comparisons}
\label{sec:comparison}

Fix a block $\fs=[M,\sigma]_G$ of Theorem~\ref{thm:tame-global-extended-quotient}
and a rigid parameter block $\svee$ of Theorem~\ref{thm:main}.
Here $M$ denotes the supporting parabolic Levi written $L$ in
\S\ref{sec:tame}. We first describe the comparison at a single
level, then transport such a comparison along \eqref{eq:two-squares}.
Write the chosen representation-side parametrization as
\begin{equation}\label{eq:p-side}
 \mathfrak b_{\fs}:\Irr(G(F))^{\fs}
       \xrightarrow{\sim}(T^{\fs}\extquot W^{\fs})_{\kappa^{\mathrm p}}.
\end{equation}
Here $\kappa^{\mathrm p}$ is the multiplier of the actual
endpoint corner algebra, with the specified left-module
convention. By Theorem~\ref{thm:fiber-irr} and the tame transfer,
$\pr\mathfrak b_{\fs}$ records the endpoint central character.

\subsection{Ordinary orbit data}

\begin{proposition}[From an equivariant cuspidal bijection to orbit data]
\label{prop:orbit-LLC}
Assume a bijection between the cuspidal representation classes
and cuspidal rigid parameter classes in the supporting Levi,
on normalizer-stable sets containing the orbits under consideration.
Assume it is equivariant for the identified unramified character
tori and for the identified relative normalizer actions, with all
rigid and Whittaker transports included. Then it induces an
equivariant isomorphism
\begin{equation}\label{eq:orbit-comparison}
 f:T^{\fs}\xrightarrow{\sim}T^{\svee},\qquad
 \iota:W^{\fs}\xrightarrow{\sim}W^{\svee}.
\end{equation}
It identifies all point stabilizers and includes the affine
translation parts of the actions.
\end{proposition}
\begin{proof}
Twist equivariance and injectivity give equality of the
self-twist groups; surjectivity on the orbit gives the quotient
torus isomorphism. For an element of the relative normalizer,
the condition that it carry the base representation into its
unramified orbit is equivalent, by bijectivity and equivariance,
to the same condition on the parameter. Hence the inertial
normalizers agree. Applying the same argument to every point
gives the stabilizer assertion, and the action formula
\eqref{eq:affine-action} proves the statement about translations.
\end{proof}

The premise is stronger than a set bijection constructed from
Yu equivalence classes. In particular, an outer automorphism can
change the Whittaker datum or rigid presentation; the correction
to its enhancement action must be the one in the equivariance
premise. Equivariance with a different action does not prove
\eqref{eq:orbit-comparison}.

\subsection{Multipliers and prescribed projective labels}

At either level, choose $f,\iota$ as in \eqref{eq:orbit-comparison}. Write
\[
 A_t^{\mathrm p}=\C_{\kappa_t^{\mathrm p}}[W_t^{\fs}],
 \qquad
 A_{f(t)}^\vee=\C_{\kappa_{f(t)}^\vee}[W_{f(t)}^{\svee}].
\]
Denote their coherent algebra transports by
$\theta^{\mathrm p}_{w,t}$ and $\theta^\vee_{\iota(w),f(t)}$.
It is enough below for these transports and identities to hold
up to inner automorphisms, since the quotient uses simple classes.

\begin{hypothesis}[Prescribed endpoint comparison]\label{hyp:cocycle}
For each $t$, supply an algebra isomorphism
$a_t:A_t^{\mathrm p}\to A_{f(t)}^\vee$, inducing the group-label
identification $\iota_t$, such that
\begin{equation}\label{eq:transport-square}
 a_{wt}\theta^{\mathrm p}_{w,t}
 =\theta^\vee_{\iota(w),f(t)}a_t
 \quad\text{on simple spectra}.
\end{equation}
In homogeneous bases the multiplication condition says that there is a cochain
$c_t:W_t^{\fs}\to\C^\times$ with $c_t(1)=1$ such that
\begin{align}
 a_t(e_w^{\mathrm p})&=c_t(w)e_{\iota(w)}^\vee,
                                      \label{eq:cochain-map}\\
 \kappa_t^{\mathrm p}(w,w')
 &=\kappa_{f(t)}^\vee(\iota(w),\iota(w'))
       \frac{c_t(w)c_t(w')}{c_t(ww')}. \label{eq:cochain-identity}
\end{align}
The transport condition concerns the prescribed family, not only
the existence of some choices on orbit representatives.
\end{hypothesis}

\begin{theorem}[Conditional rigid block bijection]
\label{thm:block-comparison}
Assume \eqref{eq:orbit-comparison} and
Hypothesis~\ref{hyp:cocycle}. There is a bijection
\begin{equation}\label{eq:conditional-bijection}
 \mathcal L_{\fs}:\Irr(G(F))^{\fs}
              \xrightarrow{\sim}\Phie^Z(\rig)^{\svee}
\end{equation}
defined by
\begin{equation}\label{eq:composite}
 \mathcal L_{\fs}=(\mathfrak b^\vee_{\svee})^{-1}
                   \circ Q\circ\mathfrak b_{\fs},
 \qquad Q([t,V])=[f(t),V\circ a_t^{-1}].
\end{equation}
It is determined by the fixed comparison data. It respects the
endpoint orbit coordinate and the rigid relevance sector.
\end{theorem}
\begin{proof}
Equation~\eqref{eq:cochain-identity} is exactly the equality
obtained by applying $a_t$ to $e_w^{\mathrm p}e_{w'}^{\mathrm p}$
and comparing with $a_t(e_w^{\mathrm p})a_t(e_{w'}^{\mathrm p})$.
Thus \eqref{eq:cochain-map} defines the asserted isomorphism,
with no unstated inverse-cocycle convention. Transport of a left
module gives, explicitly,
\[
 e_{\iota(w)}^\vee\longmapsto
             c_t(w)^{-1}V(e_w^{\mathrm p}).
\]
Equation~\eqref{eq:transport-square} says these assignments
descend to orbit classes. Each is invertible, so $Q$ is a
bijection. Composing the two block parametrizations proves
\eqref{eq:conditional-bijection}. Its first coordinate is
$\bar f$, and the target block is relevant for $\rig$ by
Theorem~\ref{thm:main}.
\end{proof}

\begin{remark}[Why scalar choices matter]
Replacing $c_t$ by $c_t\chi_t$, for a character $\chi_t$ of
$W_t^{\fs}$, preserves \eqref{eq:cochain-identity} but changes
the transported label by $\chi_t^{-1}$. A family of cohomology
class equalities therefore does not select the particular map
\eqref{eq:composite}. Lemma~\ref{lem:orbit-transport} describes
what it does give after arbitrary choices. Theorem~\ref{thm:pinned-cross}
gives a finite test for choosing the prescribed map on simple classes at
each reduced stabilizer.
\end{remark}

\subsection{Transport of the comparison from depth zero}
\label{sec:transfer-proof}
We now combine the two squares in \eqref{eq:two-squares}. In the
notation of Theorem~\ref{thm:joint-transfer}, the maps on orbits and
groups underlying $d_{\rm p}$ and $d_\vee$ are
\[
 (f_\Sigma,\iota_\Sigma):(T,W)\longrightarrow(T_0,W_0),\qquad
 (f^\vee_{\calD},\iota^\vee_{\calD}):
       (T^\vee,W^\vee)\longrightarrow(T^\vee_0,W^\vee_0).
\]
The matched reduced group is $H=G^0$. The supporting Levi denoted
$M_0$ on the parameter side corresponds to $L^0$ on the representation
side. Theorem~\ref{sq:Sigma-matching} supplies this geometric
identification for the marked quotient-character tower. The cross-side
identification of the two reduced orbits and the comparison of their
algebras are the inputs of Theorem~\ref{thm:joint-transfer}.

\begin{proof}[Proof of Theorem~\ref{thm:joint-transfer}]
The transported maps on the original orbit torsors and inertial groups
are
\begin{equation}\label{eq:joint-orbits}
 f=(f_{\calD}^\vee)^{-1}f_0 f_\Sigma,\qquad
 \iota=(\iota_{\calD}^\vee)^{-1}\iota_0\iota_\Sigma .
\end{equation}
The three constituent orbit maps are equivariant for the corresponding
group isomorphisms, including the affine translations. Thus $f$ is
$\iota$-equivariant. For $t\in T$ put $t_0=f_\Sigma(t)$, $s_0=f_0(t_0)$
and $s=f(t)$, so $f_{\calD}^\vee(s)=s_0$. Write
\[
 u_t:A^{\rm p}_t\longrightarrow A^{\rm p}_{0,t_0},\qquad
 v_s:A^\vee_s\longrightarrow A^\vee_{0,s_0}
\]
for the two vertical algebra isomorphisms oriented towards depth zero.
For $u_t$ take the inverse of the map induced by
\eqref{eq:tame-endpoint-isom} on the centralizers of the neutral matrix
algebras in the matched corners. Lemma~\ref{lem:graded-corner-functoriality}
and Proposition~\ref{prop:tame-normalizer-isom} identify its grading
and its orbit transports. For $v_s$ use the common Springer induction
object of Proposition~\ref{sq:algebras}, with the same left-module
convention.

If $a^0_{t_0}:A^{\rm p}_{0,t_0}\to A^\vee_{0,s_0}$ is the prescribed
depth-zero comparison, put
\begin{equation}\label{eq:joint-algebra}
 a_t=v_s^{-1}\,a^0_{t_0}\,u_t:A^{\rm p}_t\xrightarrow{\sim}A^\vee_s .
\end{equation}
Each factor preserves the grading under its designated group map, so
$a_t$ carries degree $w$ to degree $\iota(w)$ for $w\in W_t$. All
three factors intertwine orbit transport on simple modules. Composing
these identities gives \eqref{eq:transport-square} for $a_t$, and
hence Theorem~\ref{thm:block-comparison} applies.

More explicitly, $V$ first goes to $V\circ u_t^{-1}$. Applying the
reduced comparison and then inverse Galois transport gives the module
$V\circ u_t^{-1}\circ (a^0_{t_0})^{-1}\circ v_s=V\circ a_t^{-1}$.
This proves $Q=d_\vee^{-1}Q_0d_{\rm p}$ on orbit classes. Invertibility
of $d_\vee$ proves uniqueness of this map. The two within-side
identities $d_{\rm p}b=b_0D_{\rm p}$ and
$d_\vee b^\vee=b_0^\vee D_\vee$ give \eqref{eq:joint-transfer}.
Finally, $v_s a_t u_t^{-1}=a^0_{t_0}$ proves the converse at the
level of the fixed algebra maps. Uniqueness on quotient sets allows
inner changes of the algebra coordinates, which act trivially on
simple classes.
\end{proof}

\begin{corollary}[The free locus]\label{cor:joint-free}
Suppose only the equivariant reduced orbit and group isomorphisms
$(f_0,\iota_0)$ are specified. Over the locus of $t\in T$ with $W_t=1$
there is a unique comparison of extended-quotient labels over
\eqref{eq:joint-orbits}. It commutes with the two reductions.
\end{corollary}
\begin{proof}
The maps in \eqref{eq:joint-orbits} and their factors preserve
stabilizers. At each of the four corresponding points the stabilizer
algebra is $\C$ and has one simple class. Thus the maps on orbits
already give all the maps on quotient labels over this locus.
\end{proof}

If a positive-depth construction independently prescribes the orbit
map or projective labels, \eqref{eq:joint-orbits} gives the test on
the base. For prescribed quotient maps $Q'$ over that base, the
condition is $d_\vee Q'=Q_0d_{\rm p}$. At the level of prescribed
homogeneous algebra maps it is the comparison with
\eqref{eq:joint-algebra}, with the allowed coordinate changes.

\subsection{Ordinary and normalized cuspidal support}

Define the representation-side corrected map, without assuming
any formula for it, by
\begin{equation}\label{eq:p-normalized-p}
 p_{\mathrm{norm}}^{\mathrm p}
 =\cusp\circ\mathfrak b_{\fs}^{-1}:
 (T^{\fs}\extquot W^{\fs})_{\kappa^{\mathrm p}}
                         \longrightarrow T^{\fs}/W^{\fs}.
\end{equation}
This is computed by reversing the parameter-removal
map and taking the original central character. In its connected
geometric charts the forward shift is
$(y,\sigma,\rho)\mapsto(y,\sigma-rh,\rho_0)$; the full
finite-sheet and residual-label data remain necessary.

\begin{proposition}[Exact support-compatibility criterion]
\label{prop:support-comparison}
The bijection $\mathcal L_{\fs}$ satisfies
$\Sc\circ\mathcal L_{\fs}=\bar f\circ\cusp$ if and only if
\begin{equation}\label{eq:correction-comparison}
 p_{\mathrm{norm}}^\vee\circ Q
       =\bar f\circ p_{\mathrm{norm}}^{\mathrm p}.
\end{equation}
Neither equality of cocycle classes nor equivariance of $f$
alone implies this identity.
\end{proposition}
\begin{proof}
Substitute \eqref{eq:composite} into the desired identity and
compose on the right with $\mathfrak b_{\fs}^{-1}$.
Equations~\eqref{eq:p-corrected} and
\eqref{eq:p-normalized-p} give exactly
\eqref{eq:correction-comparison}. The earlier hypotheses concern
the first-coordinate projection and twisted stabilizer algebras.
For example, in the scalar-Frobenius fiber of
Example~\ref{ex:GL2}, tensoring the two $S_2$-labels by the sign
character exchanges them without changing the base point or
cocycle. This operation changes normalized support. Taking this
operation on that orbit and the identity on other orbits gives
an isomorphism of quotient sets of the kind allowed by
Lemma~\ref{lem:orbit-transport}, but it does not preserve the
corrected map. Thus base and multiplier comparisons alone cannot
force the desired identity.
\end{proof}

\begin{corollary}[Transfer of normalized-support compatibility]
\label{cor:joint-support}
Assume Theorem~\ref{thm:joint-transfer}. Write $\cusp_0$ and
$\Sc_0$ for ordinary supercuspidal support and normalized parameter
support in the reduced blocks, and use a bar for the induced map on
orbit spaces. If
\begin{equation}\label{eq:joint-ordinary-support}
 \cusp_0D_{\rm p}=\bar f_\Sigma\,\cusp ,
\end{equation}
then $\Sc\,\mathcal L_\Sigma=\bar f\,\cusp$ holds if and only if
$\Sc_0\mathcal L_0=\bar f_0\,\cusp_0$ holds.
\end{corollary}
\begin{proof}
By Proposition~\ref{sq:normalized},
$\bar f_{\calD}^\vee\Sc=\Sc_0D_\vee$. On composing the desired
upper identity with $\bar f_{\calD}^\vee$, its left side becomes
$\Sc_0\mathcal L_0D_{\rm p}$ by \eqref{eq:joint-transfer}. Its
right side becomes $\bar f_0\cusp_0D_{\rm p}$ by
\eqref{eq:joint-orbits} and \eqref{eq:joint-ordinary-support}.
Now use surjectivity of $D_{\rm p}$ and injectivity of
$\bar f_{\calD}^\vee$.
\end{proof}
Equation~\eqref{eq:joint-ordinary-support} concerns the original
Bernstein center in its supercuspidal-support coordinates.
Theorem~\ref{thm:tame-global-extended-quotient} identifies the
endpoint central characters under $d_{\rm p}$. Proposition~\ref{prop:support-comparison}
also applies directly when one can compute the two corrected
maps on quotient labels.

Theorem~\ref{thm:main} characterizes the bounded
parameter subset by its bounded cuspidal labels. To prove
temperedness compatibility in \eqref{eq:conditional-bijection},
one must show that $Q\mathfrak b_{\fs}$ carries the tempered
representation subset onto precisely that subset. Compatibility
with standard modules, parabolic induction, or endoscopic
character identities is additional structure, not a consequence
of a bijection of extended quotients.

\subsection{Pinned comparison of the reduced blocks}\label{sec:pinned-mixed}

We now give sufficient conditions for the depth-zero input $Q_0$ in
Theorem~\ref{thm:joint-transfer} to carry prescribed Jordan labels to
prescribed Springer labels. Let the matched reduced group be $H=G^0$,
with supporting Levi $L^0$ and the specified rigid realization. Fix a
pinned correspondence $\ell_{\cP,L^0}$ on the relevant depth-zero
supercuspidal seeds of $L^0(F)$. Assume the twist and normalizer
equivariance of Proposition~\ref{prop:orbit-LLC}, with the prescribed
rigid and pinning transports. It gives $(f_0,\iota_0)$ on the actual
unramified orbits, including their self-twist quotients and affine
normalizer actions.

For $t_0\in T_0$ write $s_0=f_0(t_0)$ and $K=(W_0)_{t_0}$.
Use $\iota_{0,t_0}$ to grade the Springer algebra $A^\vee_{0,s_0}$ by
$K$. The two algebras to compare are now the endpoint commutant
$A^{\rm p}_{0,t_0}$ and this Springer algebra in the prescribed rigid
sector, both with the left convention. The required equality on $K$
is
\begin{equation}\label{eq:pinned-cross-class}
 [\kappa^{\rm p}_{0,t_0}]
        =\iota_{0,t_0}^*[\kappa^\vee_{0,s_0}]
             \qquad\text{in }H^2(K,\C^\times).
\end{equation}
It is an additional hypothesis even when the two effective multiplier
classes on $Q_{d,V}$ in \eqref{eq:pinned-effective} agree.

\begin{theorem}[Pinned comparison from depth-zero tests]\label{thm:pinned-cross}
Assume the reduced orbit data just described and
\eqref{eq:pinned-cross-class} at one representative of every orbit.
At each representative suppose that the finite Jordan and Springer
constructions independently prescribe pairs of simples of the two
stabilizer algebras. Suppose these pairs satisfy (a), (b) and
\eqref{eq:pinned-detecting} of Proposition~\ref{prop:pinned-trace}.
Assume that the pairs extend covariantly for orbit transports and
the stipulated changes of pinning, including returns, using the
actual graded algebra transports (up to inner change).
Then there is a unique map on quotient labels induced by graded
comparisons over $(f_0,\iota_0)$ which realizes the tests. It is
\begin{equation}\label{eq:pinned-Q}
 Q_0^{\rm pin}([t_0,E])=[s_0,(a^0_{t_0})_*E],
\end{equation}
where $a^0_{t_0}$ is any calibrated comparison. It is a bijection
covariant under the prescribed changes.

Suppose in addition that the representation-side map $b_0^{\rm pin}$
is supplied by Theorem~\ref{thm:pinned-removal}, including type and
Levi transports, and that the Galois parametrization $b_0^\vee$ is
natural for the corresponding changes of rigid data. The square
\begin{equation}\label{eq:pinned-reduced-square}
\begin{CD}
 P_0 @>{\mathcal L_0^{\rm pin}}>> Y_0\\
 @V{b_0^{\rm pin}}VV @VV{b_0^\vee}V\\
 \mathscr E_0 @>{Q_0^{\rm pin}}>> \mathscr E_0^\vee
\end{CD}
\end{equation}
defines a covariant bijection of reduced block classes. On a
supercuspidal block ($L^0=H$) it is the unramified extension of
$\ell_{\cP,H}$, provided both parametrizations have their specified
seed normalization.
\end{theorem}
\begin{proof}
Equation~\eqref{eq:pinned-cross-class} gives an initial graded
comparison at $t_0$. Apply Proposition~\ref{prop:pinned-trace} with
$C_1=A^{\rm p}_{0,t_0}$ and $C_0=A^\vee_{0,s_0}$. It proves existence
and uniqueness on simple classes. Lemma~\ref{lem:pinned-naturality}
then proves \eqref{eq:transport-square} on these classes and gives
covariance for the stipulated changes of pinning. Thus
\eqref{eq:pinned-Q} descends, and the inverse calibrated maps give
its inverse. A character annihilating $B(A^\vee_{0,s_0})$ fixes
every label, so all remaining cochain choices give the same map.

For a change of data write $U,V$ for the transports on $P_0,Y_0$
and $E^{\rm p},E^\vee$ for those on their quotients. Use primes for
the new data. The hypotheses and the first assertion give
\[
 (b_0^{\rm pin})' U=E^{\rm p}b_0^{\rm pin},\qquad
 (b_0^\vee)' V=E^\vee b_0^\vee,\qquad
 (Q_0^{\rm pin})' E^{\rm p}=E^\vee Q_0^{\rm pin}.
\]
Insertion in the composite $(b_0^\vee)^{-1}Q_0^{\rm pin}b_0^{\rm pin}$
proves $(\mathcal L_0^{\rm pin})'U=V\mathcal L_0^{\rm pin}$.
For $L^0=H$ the inertial Weyl group is trivial, so the orbit map
with the prescribed base point gives the stated extension of the seed.
\end{proof}

\begin{corollary}[Pinned transfer to the tame block]\label{cor:pinned-joint}
Under the hypotheses of Theorems~\ref{thm:joint-transfer}
and~\ref{thm:pinned-cross}, take $b^{\rm pin}$ from
Corollary~\ref{cor:pinned-tame}. The unique upper quotient comparison
and the resulting block bijection are
\[
 Q^{\rm pin}=d_\vee^{-1}Q_0^{\rm pin}d_{\rm p},\qquad
 \mathcal L_\Sigma^{\rm pin}
     =(b^\vee)^{-1}Q^{\rm pin}b^{\rm pin}
     =D_\vee^{-1}\mathcal L_0^{\rm pin}D_{\rm p}.
\]
They are covariant for changes of the pinned marked datum for which
the two reductions and the corner transports have the specified
covariance.
\end{corollary}
\begin{proof}
Apply \eqref{eq:joint-algebra} to any calibrated representatives of
\eqref{eq:pinned-Q}. Its action on simple classes is independent of
these representatives because both vertical maps transport the actual
corner modules. The two identities in \eqref{eq:two-squares} then
give the formulas, and the covariance identities compose.
\end{proof}

There is a corresponding freedom of coordinates in
\eqref{eq:pinned-reduced-square}. If $A$ relabels $\mathscr E_0$,
changing $b_0^{\rm pin}$ to $A b_0^{\rm pin}$ and $Q_0^{\rm pin}$ to
$Q_0^{\rm pin}A^{-1}$ leaves $\mathcal L_0^{\rm pin}$ unchanged.
A prescribed test pair must undergo this same relabelling. Thus a
choice can be assigned to either adjacent arrow once its action on
the tests has been specified.

For application to a general tame block, two finite comparisons
require independent input. On the representation side one must realize
the tests on $Q_{d,V}$ with the actual products
\eqref{eq:local-defect}, and carry the finite pinning-change law
\eqref{eq:pinned-jd-change} through sheets, corners and the full
normalizer. Across the two sides one must establish
\eqref{eq:pinned-cross-class} and realize the tests on $(W_0)_{t_0}$
using the source scalar extension and the rigid Springer sector.
The affine comparison and the supercuspidal seed alone do not supply
these tests. Once they are known at depth zero, the corollary
transports the resulting labels at every stabilizer. If an ambient
comparison is prescribed independently, agreement is tested by
\eqref{eq:joint-algebra} on simple classes. Compatibility with ordinary
cuspidal support is governed separately by
Proposition~\ref{prop:support-comparison}.

\section{Enhanced parameters for tame blocks with non-singular reduced support}
\label{sec:regular-blocks}

Let $\fs=[L,\sigma]_G$ be a tame Bernstein block with the fixed type of
\S\ref{sec:tame}. Write $H=G^0$ and
$\fs_{-1}=[L^0,\sigma_{-1}]_H$. The hypothesis concerns the actual
depth-zero supercuspidal representation $\sigma_{-1}$ of
\eqref{eq:tame-seed}, with the restriction of the ambient quadratic
twist used in that formula. We require that $\sigma_{-1}$ belong to a
Deligne--Lusztig packet for a maximally unramified elliptic torus
$S\subset L^0$ with an $F$-non-singular depth-zero character
\cite[\S2]{SXDepthZero}. Here $F$-non-singularity is the condition on
the norms of restricted coroots in \cite[Definition~3.1.1]{KalethaPackets}.
Testing only the roots of the connected finite parahoric quotient gives
the weaker $k_F$-non-singularity condition. Unramified twists are
trivial on the compact torus, so the same condition holds throughout
the inertial orbit of $\sigma_{-1}$.

For a tame $F$-non-singular elliptic pair in $L$, a Howe factorization
gives a depth-zero character $\phi_{-1}$ with the same restriction as
the original character to the torus in the simply connected cover of
$(L^0)_{\rm der}$ \cite[\S3.4]{KalethaPackets}. The coroot norms factor
through this torus. Thus a presentation in which $\phi_{-1}$ supplies
the actual reduced character satisfies the hypothesis. When changing
from an intrinsic Levi twist to the ambient twist, the test is on the
resulting character \cite[Remark~4.1.5]{AFMO-reduce}.

In particular, a regular Bernstein block admits the regular
presentation of \cite{AdlerMishraRegular,KalethaRegular}. In that
presentation the depth-zero torus character has trivial stabilizer on
$S(F)_0$ in $N_{L^0}(S)(F)/S(F)$. The ambient and intrinsic quadratic
characters on $S(F)_0$ are invariant under this normalizer: it fixes
the point of the reduced building attached to $S$, and the characters
extend to the corresponding reduced-point stabilizers
\cite[\S2.7 and Remark~4.1.5]{AFMO-reduce}. Their ratio therefore
preserves the stabilizer of the torus character. On the connected
finite quotient, tensoring by the ratio multiplies the inducing torus
character by its restriction, by the Deligne--Lusztig character formula
(the ratio is trivial on unipotent elements since $p$ is odd).
The actual reduced character is consequently regular and hence
$F$-non-singular \cite[Fact~3.1.5]{KalethaPackets}. This verifies
the hypothesis on reduced support. Applying the theorem to this
presentation also requires the quotient characters specified below.

We first compare the sets of enhanced parameters occurring in the
reduced block. Solleveld--Xu use the rigid version of AMS in
\cite[Lemma~8.1]{SXDepthZero}. The following observation identifies
the induction series needed here with those used in
\S\ref{sec:rigid}, including their component-group coefficients.

\begin{lemma}[Toral quasi-cuspidal support]\label{lem:regular-toral-series}
Let $Q$ be a possibly disconnected complex reductive group, let $T$ be
a maximal torus of $Q^\circ$, and put $K=Z_Q(T)$. For
$\epsilon\in\Irr(K/T)$, the unipotent pairs of AMS quasi-cuspidal
support $[K,1,\epsilon]_Q$ are precisely the simple summands of
geometric parabolic induction of the $K$-equivariant sheaf
$\epsilon$ at $1$. In particular the AMS and Dillery--Schwein
quasi-cuspidal series for this datum agree as sets of unipotent pairs.
\end{lemma}
\begin{proof}
Choose a Borel $B=TU$ of $Q^\circ$. Since $K$ fixes $T$ pointwise,
it preserves every root subgroup and normalizes $B$. Put $P=KB$.
Then $K^\circ=T$, $P/U=K$, and $P/B=K/T$. For a cocharacter of
$T$ in the chosen open chamber outside the fixed spaces of the
nonidentity actions of $N_Q(T)/K$ on $T$, the associated parabolic
is $P$ and its Levi is $K$. Consider the proper map
\[
  m:Q\times^P B\longrightarrow Q^\circ,\qquad [g,b]\longmapsto gbg^{-1}.
\]
The coefficient on its source is the local system given by $\epsilon$
through the finite cover $Q\times^B B\to Q\times^P B$. Each connected
component of $Q\times^P B$ is a copy of $Q^\circ\times^B B$,
and $m$ on it is the ordinary small Grothendieck--Springer map
followed by conjugation. Consequently its shifted direct image is
the intersection complex extending the direct image on the dense
regular semisimple locus. Over this locus the source identifies with
$Q\times^K T_{\rm reg}$: the unipotent radical $U$ uniquely conjugates
a regular semisimple element of $B$ into $T$. This description keeps
every component of $Q/P$. With the coefficient $\epsilon$
this is exactly the local system $q\pi_*(\widetilde{q\mathcal E})$
in \cite[\S5, (66)--(67)]{AMS}.

Over the unipotent variety the inverse image is $Q\times^P U$.
Proper base change identifies the restriction, with the perverse
shift for the unipotent variety, with the induction object of
\cite[\S2.4]{DS}. The criterion following (67) in \cite[\S5]{AMS}
tests membership by occurrence of the local system on the orbit in
the top Springer-fibre cohomology with these coefficients. For the
semismall map $Q\times^P U\to\mathcal U_{Q^\circ}$ this is exactly
the occurrence of the corresponding intersection complex as a direct
summand, as in \cite[Corollary~3.5.3]{DS}. All identifications are
$Q$-equivariant, so they retain $\epsilon$ on the full group $K/T$.
\end{proof}

\begin{theorem}[Enhanced parameters for a marked tame block with non-singular reduced support]
\label{thm:regular-block}
Let $\fs=[L,\sigma]_G$ and suppose that its reduced supercuspidal
representation $\sigma_{-1}$ is $F$-non-singular in the sense above.
Fix the AFMO block equivalence of \S\ref{sec:tame} and a marked
ambient quotient-character tower of Definition~\ref{sq:Sigma-input}
on the same rational tower, with the common cover and the prescribed
relative rigid realization. Writing $\vec{\phi}$ for the characters
of the fixed AFMO datum, require the full quotient characters of the
marked tower to satisfy
\[
  \bar\chi_i\circ q_{G^i}=\phi_i \quad\text{on }G^i(F)
  \qquad (0\leq i\leq d).
\]
This includes the terminal character $\phi_d$. Use the reduced seed
$\sigma_{-1}$ of \eqref{eq:tame-seed}, with the ambient quadratic
twist used in the AFMO reduction. Fix a realization of the Solleveld--Xu
correspondence for $H$, and use for $L^0$ the cuspidal correspondence
entering that construction. With these choices, taking cuspidal support
commutes with the correspondence \cite[Theorem~1(f)]{SXDepthZero}.
Take as the marked seed
\[
 t_{0,*}=\operatorname{rec}^{\rm SX}_{L^0}(\sigma_{-1}),\qquad
 t_* = J^+_{\calD_\Sigma}(t_{0,*})
\]
where the second map is taken in the supporting Levi. Write $\svee$
for the ambient inertial class of $t_*$. There is a unique bijection
completing the following square for these choices:
\begin{equation}\label{eq:regular-block-square}
\begin{CD}
 \Irr(G(F))^{\fs} @>{\mathcal L_\Sigma^{\rm ns}}>>
       \Phie^Z(\rig)^{\svee}\\
 @V{D_{\rm p}}VV @VV{R^+_{\calD_\Sigma}}V\\
 \Irr(H(F))^{\fs_{-1}} @>{\operatorname{rec}^{\rm SX}_H}>>
       \Phie^Z(\rigH)^{\ssH}.
\end{CD}
\end{equation}
Here $\ssH$ is the inertial class of $t_{0,*}$ in $H$.
Explicitly, for $\Pi\in\Irr(G(F))^\fs$, write
$\operatorname{rec}^{\rm SX}_H(D_{\rm p}\Pi)=[\psi_0,u,\rho]$.
Then
\begin{equation}\label{eq:regular-block-param}
 \mathcal L_\Sigma^{\rm ns}(\Pi)=
 [\,w\longmapsto\varkappa_\Sigma(w)j_\Sigma(\psi_0(w)),\,u,\,\rho\,].
\end{equation}
\end{theorem}
\begin{proof}
By $F$-non-singularity of $\sigma_{-1}$, Solleveld--Xu give the
bijection on the whole reduced block in
\cite[(10.12), Theorem~10.4]{SXDepthZero}. Its supporting cuspidal
correspondence agrees with Kaletha's by their Proposition~10.7, and
their Theorem~1(f) identifies its cuspidal support. Use the same
finite central subgroup $Z$ in their rigid component groups and in
\eqref{eq:rigid-groups}. For every unramified twist of $t_{0,*}$ the
supporting parameter is discrete in the supporting Levi and trivial
on $\SL_2(\C)$. Hence its connected Weil centralizer there is a torus:
discreteness makes the centralizer finite modulo the invariant
center. Its full inverse image in the common cover has toral identity
component as well. In the ambient Weil centralizer $Q$ in that cover,
the supporting group is a quasi-Levi. Its identity component is
therefore a maximal torus $T$ of $Q^\circ$, and the group itself is
$Z_Q(T)$. The inducing orbit is $\{1\}$, with the given enhancement
on $Z_Q(T)/T$. Lemma~\ref{lem:regular-toral-series} applies to every
irreducible representation of this finite group and identifies the
two unnormalized support series at every twist. The normalizing
cocharacters change support by a connected unramified twist, so the
resulting inertial blocks coincide. This identifies the lower right
set in \eqref{eq:regular-block-square}.

The left arrow is the AFMO bijection on irreducibles. By the character
matching hypothesis, the central Weil cocycle in
the lift is formed from quotient characters inducing the given AFMO
characters. For the fixed embeddings, Lemma~\ref{sq:Sigma-lift-change}
shows that any other quotient lifts of those same characters give the
same map on enhanced conjugacy classes.
Theorem~\ref{sq:Sigma-matching} identifies the right arrow and its
inverse $J^+_{\calD_\Sigma}$ on the two indicated rigid blocks.
Thus $\mathcal L_\Sigma^{\rm ns}=
J^+_{\calD_\Sigma}\operatorname{rec}^{\rm SX}_H D_{\rm p}$ is forced
and is bijective. Lemma~\ref{sq:groupoid} identifies the full
centralizers of the Weil part and $u$ in the common cover, which gives
the enhancement $\rho$ in \eqref{eq:regular-block-param}. Relevance
for the given rigid form follows from Lemma~\ref{sq:relevance}.
\end{proof}

In the extended-quotient coordinates of \eqref{eq:two-squares} the
result gives
\begin{equation}\label{eq:regular-quotient-square}
 Q_0^{\rm ns}=b_0^\vee\operatorname{rec}^{\rm SX}_H b_0^{-1},
 \qquad Q^{\rm ns}=d_\vee^{-1}Q_0^{\rm ns}d_{\rm p}
                  = b^\vee\mathcal L_\Sigma^{\rm ns}b^{-1}.
\end{equation}
These are maps on simple classes, including the classes at singular
orbits. To identify them with a specified graded stabilizer-algebra
comparison over the prescribed orbit map, one must compare the
endpoint and Springer normalizations in \S\ref{sec:pinned-mixed}.
Solleveld--Xu describe the choices affecting their enhancements in
\cite[Lemma~10.3]{SXDepthZero}. An independently normalized
positive-depth Kaletha seed $t_{\rm Kal}(\sigma)$ in the marked wild
slice gives the same
based target precisely when its reduction in the supporting Levi is
$t_{0,*}$ as an enhanced class. For equality of target blocks it
suffices that these seeds be inertially equivalent in the reduced
ambient group. These tests follow from Corollary~\ref{sq:seed-matching}
and account for the full central factors, the tame embeddings and the
ambient twist in \cite[Remark~4.1.5]{AFMO-reduce}.

\begin{example}[A regular seed with disconnected finite dual centralizer]
\label{ex:regular-disconnected}
Let $G=\SL_6$ over a local field with residue field $k=\mathbb F_{11}$.
In the dual of the hyperspecial reductive quotient, $\mathrm{PGL}_6$,
let $s$ be the projective class of the companion matrix of
$X^6+X^3+1=\Phi_9(X)$. Since $11$ has order $6$ modulo $9$, the
six primitive ninth roots are distinct and form one Frobenius orbit.
The diagonal torus in an eigenbasis therefore descends to an elliptic
maximal torus $T^*$, and $Z_{\mathrm{PGL}_6}(s)^\circ=T^*$.

A projective centralizer element lifts to a linear map that conjugates
the companion matrix to a scalar multiple of itself. The admissible
scalars are exactly $\mu_3$. Indeed, after choosing a primitive ninth
root, they are the translations preserving
$\{1,2,4,5,7,8\}\subset\mathbb Z/9$, namely $0,3,6$. Each is realized
by permuting the eigenlines, and the kernel is the diagonal torus.
Thus $\pi_0 Z_{\mathrm{PGL}_6}(s)\simeq\mu_3$ with Frobenius
$c\mapsto c^{11}=c^{-1}$. Its Frobenius-fixed subgroup is trivial.
The stabilizer of $s$ in the geometric Weyl group of $T^*$ is this
component group. By finite-torus duality the corresponding character
of the dual elliptic torus in $\SL_6(k)$ consequently has trivial
rational Weyl stabilizer. At the hyperspecial vertex this is a regular
depth-zero seed \cite[Fact~3.4.18 and Definition~3.4.19]{KalethaRegular}.
Here the rational points of the centralizer coincide with $T^*(k)$,
although the algebraic centralizer has three components.
\end{example}

\begin{example}[A positive-depth twist]\label{ex:regular-central-twist}
Take $G=\GL_2$, $Z=1$, and a positive-depth character $\xi$ of
$F^\times$, with Weil character $\widehat\xi$. Consider the regular
block $[T,\xi\boxtimes\xi]_G$ for the diagonal torus $T$. There is
only a terminal quotient character $\xi\circ\det$, so $H=G$,
$D_{\rm p}$ is tensoring by $(\xi\circ\det)^{-1}$, and $J^+_{\calD}$
multiplies the Weil part by $\widehat\xi$. The two constituents with
ordinary cuspidal support
$\{\xi|\cdot|^{1/2},\xi|\cdot|^{-1/2}\}$ are $\xi\circ\det$ and its
Steinberg twist. In $W_F\times\SL_2(\C)$ notation their parameters are,
respectively,
\[
 (\widehat\xi|\cdot|^{1/2})\oplus(\widehat\xi|\cdot|^{-1/2}),
 \qquad \widehat\xi\boxtimes\operatorname{Std}_2 .
\]
Both enhancements are trivial. Here the two parameters are
distinguished by their $\SL_2$ factors already in the reduced block.
The torus seed is regular and the scalar orbit has stabilizer $S_2$,
as in Example~\ref{ex:GL2}.

One also obtains a positive-depth block outside the regular case.
Let $S$ be the unramified elliptic norm-one torus in $\SL_2$, let $k$
have odd order $q$, and let $k_2/k$ be quadratic. Inflate the unique
quadratic character of $S(k)$ to $S(F)$. The coroot norm on residue
fields is $a\mapsto a/a^q$ from $k_2^\times$ onto $S(k)$, so the
character is $F$-non-singular. The rational Weyl involution acts by
inversion and fixes the character. A depth-zero supercuspidal
constituent $\tau$ of the associated Deligne--Lusztig packet is
therefore non-singular and nonregular. In $\SL_2\times\GL_2$ take
the block with supporting Levi $\SL_2\times T$ and supercuspidal
representation $\tau\boxtimes(\xi\boxtimes\xi)$. Its only
positive-depth factor is the same terminal determinant character on
$\GL_2$. Theorem~\ref{thm:regular-block} applies to the entire block,
including both products of $\tau$ with the constituents above.
\end{example}

\section{Tori, normalizers and singular orbits}\label{sec:checks}
We illustrate the role of the connected unramified character torus,
the full inertial normalizer and the choice of projective labels.

\subsection{Toral blocks and self-twists}
\begin{example}[A torus]\label{ex:weak-torus}
For an ambient torus the supporting Levi is the group itself
and $W=1$. Every nonempty relevant block is its unramified
orbit, the algebra $\cA_t$ is $\C$, and both support maps are
the identity. The action of $\Xnr$ on the underlying torus parameter
is free by translation in $H^1(\WF,\widehat T)$.

Let $T=\operatorname{Res}^{1}_{E/F}\mathbb G_m$ for a tamely ramified
quadratic extension $E/F$. Inertia acts on $\widehat T=\C^\times$
through inversion. Hence $\widehat T^{\IF}=\mu_2$,
$X_{\mathrm{w.nr}}({}^LT)=\mu_2$, and $\Xnr({}^LT)=1$. The cocycle
$b$ trivial on $\IF$ with $b(\Fr)=-1$ gives the nonidentity class.
Indeed, a coboundary trivial on $\IF$ comes from an element of
$\mu_2$, on which Frobenius acts trivially. Twisting the trivial
cocycle by $b$ therefore exchanges two distinct singleton blocks
with the same choice of enhancement.
\end{example}

\begin{example}[A supercuspidal block]
If the supporting Levi is the ambient group, then $W=1$ and
the block consists of unramified twists of its given cuspidal
enhanced parameter. Its induction series has one simple
constituent. Thus the general theorem introduces no extra
projective labels in a supercuspidal block.
\end{example}

\begin{example}[Different self-twist groups in a quadratic induction]
Take $Z=1$, $G=\GL_2$, and $H=\operatorname{Res}_{E/F}\mathbb G_m$
for the unramified quadratic extension $E/F$. Let $\sigma$ generate
$\operatorname{Gal}(E/F)$ and choose a character $\theta$ of $W_E$
whose restriction to $P_F=P_E$ differs from that of $\theta^\sigma$.
For $\psi=\Ind_{W_E}^{W_F}\theta$, the full wild centralizer is the
diagonal torus $\widehat H$ and Frobenius exchanges its two factors.
The invariant connected center is scalar, so $M_0=H$ and $M=G$.
Use the usual monomial $L$-embedding and the central cocycle given
by $\theta$ under the torus correspondence, so that $\psi$ reduces
to the trivial cocycle of $H$.

In the coordinate $(x,y)\mapsto xy$ on $(\widehat H)_{\Fr}$,
$r$ sends $z$ to $z^2$. The torus seed has trivial self-twist
group. An unramified twist fixing $\psi$ must fix $\theta$ upon
restriction to $W_E$, since the two wild characters are distinct.
Conversely, every twist trivial on $W_E$ fixes $\psi$ by induction.
The ambient self-twist group is therefore $\mu_2=\ker r$, as in
\eqref{sq:selftwist-sequence}. Both inertial Weyl groups are trivial.
\end{example}

\subsection{Compact inducing data and finite quotients}
\label{sec:finite-quotients}
Keep the depth-zero datum $(K_L,\tau_{\rb},H,\tau_H)$ and
$\Gamma_H=H/K_L$ of \S\ref{sec:notations}.

\begin{corollary}[The self-twist quotient torsor]\label{cor:extension-torsor}
Put $\widehat\Gamma_H=\Hom(\Gamma_H,\C^\times)$ and regard its characters
as characters of $H$. The set of isomorphism classes
\[
\mathcal E_H(\tau_H)=\{[\tau_H\otimes\chi]:\chi\in\widehat\Gamma_H\}
\]
is a torsor for
\[
\widehat\Gamma_H/X_H^{\nr}(\tau_H)\simeq T^\fs.
\]
Compact induction identifies it with the unramified orbit of $\pi$.
\end{corollary}
\begin{proof}
Twisting is transitive on the displayed set and has stabilizer exactly
$X_H^{\nr}(\tau_H)$. The quotient therefore acts simply transitively.
Lemma~\ref{lem:Ts-as-HomLambda} and induction in stages identify the result
with the unramified orbit of $\pi$.
\end{proof}

\begin{remark}\label{rem:genuine-extensions}
The equality $H=L(F)_{[x]}$ expresses invariance of $\tau_{\rb}$; it does
not imply that $\tau_H$ is an extension of $\tau_{\rb}$. It is a genuine
extension exactly when $\tau_H|_{K_L}\simeq\tau_{\rb}$, equivalently when
the Clifford multiplicity space has dimension one. In that case the
self-twist group on $H/K_L$ is trivial and the extension torsor is
$\widehat\Gamma_H$. In general the self-twist quotient must be retained.
\end{remark}

There is no general identification of $W^\fs$ with a normalizer quotient
inside the reductive quotient of an arbitrary point of $G$ containing the
chosen Levi compact datum. The following example explains the obstruction.

\begin{example}\label{ex:finite-quotient-obstruction}
Let $G=\mathrm{GL}_2(F)$, $L=T$ the diagonal torus and $\pi=1$. Then
$W^\fs\simeq S_2$. Choose the image $y$ of an enlarged torus point in the
standard apartment with root coordinate $\alpha(y)=1/3$. The enlarged
point stabilizer is the standard Iwahori subgroup $I$, and
\[
K_L=T(\mathcal O_F),\qquad I/I^+\simeq T(\ff).
\]
The compact datum is $\tau_{\rb}=1$, whose $T(F)$-orbit is a singleton.
Nevertheless,
\[
N_I(T(\mathcal O_F))/T(\mathcal O_F)=1,
\qquad N_{T(\ff)}(T(\ff))/T(\ff)=1.
\]
Indeed, $T(\mathcal O_F)$ is Zariski dense in $T$, so any normalizer is a
monomial matrix. An invertible anti-diagonal matrix over $\mathcal O_F$
has both off-diagonal entries units and cannot lie in $I$; the diagonal
normalizers in $I$ are exactly $T(\mathcal O_F)$. Thus neither displayed
quotient recovers $S_2$. This point is also zero-generic for the torus
embedding. A reflection cannot be corrected to fix it by an integral
valuation translation: the required root displacement is $2/3$.
\end{example}

This example distinguishes the inertial Weyl group used in the
parametrization from normalizer quotients in a finite reductive
quotient at an arbitrary building point.

\subsection{A split family with singular orbits}

\begin{example}[A split $\GL_4$ family with jumping stabilizers]\label{ex:GL4-wild}\label{ex:GL4-joint}
Take $G=\mathrm{GL}_4$ and two Weil characters $\chi_1,\chi_2$
whose wild restrictions $\lambda_1,\lambda_2$ are distinct. Put
\[
 \lambda=\operatorname{diag}(\lambda_1,\lambda_1,
                            \lambda_2,\lambda_2),\qquad
 \hH=\mathrm{GL}_2(\C)\times\mathrm{GL}_2(\C).
\]
Use the block-diagonal embedding and
$\varkappa=\operatorname{diag}(\chi_1 I_2,\chi_2 I_2)$.
For the diagonal supporting Levi $M_0=M$, start with the trivial
tame seed. Then $T=T_0=(\C^\times)^4$, $f$ is the identity in the
matched coordinates, and $W=W_0=S_2\times S_2$. A permutation
interchanging the two pairs cannot preserve the wild restriction,
even after an unramified twist, so it is absent from both inertial
normalizers. Centralizers can jump when the two Frobenius coordinates
within a pair coincide. At every such jump the two full centralizers
are still equal, and the corresponding Springer stabilizer algebras
are the same ordinary group algebras. This checks the family argument
at singular as well as generic points; a calculation only at a generic
point would miss the nontrivial stabilizer labels.

For a comparison with representations, take $p>3$, $Z=1$,
$H=G^0=\GL_2\times\GL_2$ and $\chi_2=1$, and let $\chi_1$ correspond
under local class field theory to a positive-depth character $\xi$ of
$F^\times$. Use the full quotient character
$(g_1,g_2)\mapsto\xi(\det g_1)$ on $H(F)$ and the trivial terminal
character. The roots outside $H$ join the two pairs of coordinates.
On their coroots the character restricts to $\xi$ or $\xi^{-1}$, so
the leading character is nontrivial at its depth for each such root.
This gives the split one-step generic tower. On both sides use the
diagonal torus as supporting Levi. On the representation side take
the depth-zero principal block of $H(F)$ with base character $1$,
a $0$-generic embedding and the corresponding fixed ambient-twisted
type. The upper base representation is the one determined by
\eqref{eq:tame-seed}. On the Galois side use the trivial torus seed
and the block-diagonal embedding.

Local class field theory on the split torus identifies the two
reduced orbits with $(\C^\times)^4$ using compatible Frobenius and
uniformizer conventions. The two transpositions in
$S_2\times S_2$ exchange the coordinates within each pair. For
$z=(z_1,z_2,z_3,z_4)$ write
\[
 I(z)=\{i\in\{1,2\}: z_{2i-1}=z_{2i}\}.
\]
The stabilizer is $\prod_{i\in I(z)}S_2$. One can calculate the
algebras themselves at depth zero. The torus compact type is trivial.
The quotient of its full monomial normalizer in $H(F)$ by
$T(\mathcal O_F)$ is $\Z^4\rtimes(S_2\times S_2)$, split by the
permutation matrices. The trivial type extends trivially to this
normalizer, so its multiplication lines have the ordinary group
products. The neutral torus algebra has one-dimensional fibres.
The endpoint comparison with the raw normalizer algebra consequently
gives $\C[\prod_{i\in I(z)}S_2]$ at $z$. On the parameter side each
equal pair gives ordinary toral Springer induction in $\GL_2(\C)$,
and each unequal pair gives induction in a torus. With a fixed
ordinary Springer convention the endomorphism algebra is the tensor
product of the same group algebras.

Choose in each $\GL_2$ factor the comparison sending its permutation
generator to its ordinary Springer generator. These comparisons
commute with transport: within each factor a transposition either
moves a point with trivial stabilizer or fixes the point and acts
trivially on simple classes. Taking products gives a reduced family
$a^0_z$ as in Theorem~\ref{thm:joint-transfer}. Formula
\eqref{eq:joint-algebra} transports it to the upper blocks. Thus the
numbers of endpoint labels when neither, one, or both of the equalities
$z_1=z_2$ and $z_3=z_4$ hold are respectively $1$, $2$ and $4$.
At a double collision the two transpositions still act on separate
factors and their product in either order is the same in both algebras.
Twisting the comparison in one factor by sign permutes the labels
involving that factor. This also illustrates Proposition~\ref{prop:pinned-trace}.
In the displayed homogeneous bases, if the trivial stabilizer module
on each side is prescribed as a pair, its trace is $1$ at every
group element. Hence $B_{\rm test}=B(\C[K])=K$ for
$K=\prod_{i\in I(z)}S_2$, and that one test fixes all labels over $z$.
This calculation uses the stated test pair as input.
Example~\ref{ex:GL2} explains why checking
normalized support distinguishes such choices.
\end{example}

\section*{Acknowledgments}
The author used ChatGPT to explore proof strategies, develop and revise
parts of the arguments, and edit the exposition. The author takes
responsibility for the statements and proofs in this paper.

\bibliographystyle{alpha}
\bibliography{depth-reduction-extended-quotients}

@incollection{ABPHecke,
 author={Aubert, Anne-Marie and Baum, Paul and Plymen, Roger},
 title={The {Hecke} algebra of a reductive {$p$}-adic group: a view from noncommutative geometry},
 booktitle={Noncommutative Geometry and Number Theory},
 series={Aspects of Mathematics}, volume={E37},
 publisher={Vieweg}, address={Wiesbaden}, year={2006}, pages={1--34},
 doi={10.1007/978-3-8348-0352-8_1},
 note={Preprint \href{https://arxiv.org/abs/math/0502234}{arXiv:math/0502234}}
}

@misc{ABPPrincipal,
 author={Aubert, Anne-Marie and Baum, Paul and Plymen, Roger},
 title={Extended quotients in the principal series of reductive {$p$}-adic groups},
 year={2011},
 howpublished={\href{https://arxiv.org/abs/1110.6596}{arXiv:1110.6596}}
}

@misc{ABPSGeometric,
 author={Aubert, Anne-Marie and Baum, Paul and Plymen, Roger and Solleveld, Maarten},
 title={Geometric structure and the local {Langlands} conjecture},
 year={2012},
 howpublished={\href{https://arxiv.org/abs/1211.0180}{arXiv:1211.0180}},
 note={Revised 2013}
}

@incollection{ABPSConjectures,
 author={Aubert, Anne-Marie and Baum, Paul and Plymen, Roger and Solleveld, Maarten},
 title={Conjectures about {$p$}-adic groups and their noncommutative geometry},
 booktitle={Around {Langlands} Correspondences},
 series={Contemporary Mathematics}, volume={691},
 publisher={American Mathematical Society}, year={2017}, pages={15--51},
 doi={10.1090/conm/691/13892},
 note={Preprint \href{https://arxiv.org/abs/1508.02837}{arXiv:1508.02837}}
}

@misc{AubertXuGtwo,
 author={Aubert, Anne-Marie and Xu, Yujie},
 title={The explicit local {Langlands} correspondence for {$G_2$}},
 year={2022},
 howpublished={\href{https://arxiv.org/abs/2208.12391}{arXiv:2208.12391}},
 note={Revised 2023}
}

@misc{DS,
 author={Dillery, Peter and Schwein, David},
 title={A stacky generalized {Springer} correspondence and rigid enhancements of {$L$}-parameters},
 year={2024},
 eprint={2308.11752v2},
 archivePrefix={arXiv},
 howpublished={\href{https://arxiv.org/abs/2308.11752v2}{arXiv:2308.11752v2}},
 note={Version of 10 June 2024; all numbered references use this version},
 url={https://arxiv.org/abs/2308.11752v2}
}

@article{AMS,
 author={Aubert, Anne-Marie and Moussaoui, Ahmed and Solleveld, Maarten},
 title={Generalizations of the {Springer} correspondence and cuspidal {Langlands} parameters},
 journal={Manuscripta Mathematica}, volume={157}, year={2018}, pages={121--192}
}

@article{KalethaRigid,
 author={Kaletha, Tasho},
 title={Rigid inner forms of real and {$p$}-adic groups},
 journal={Annals of Mathematics}, volume={184}, number={2}, year={2016}, pages={559--632},
 doi={10.4007/annals.2016.184.2.6}
}

@incollection{Haines,
 author={Haines, Thomas J.},
 title={The stable {Bernstein} center and test functions for {Shimura} varieties},
 booktitle={Automorphic Forms and {Galois} Representations, Vol.~2},
 series={London Mathematical Society Lecture Note Series},
 volume={415}, publisher={Cambridge University Press}, year={2014},
 pages={118--186}, doi={10.1017/CBO9781107297524.005}
}

@article{Dillery,
 author={Dillery, Peter}, title={Rigid inner forms over local function fields},
 journal={Advances in Mathematics}, volume={430}, year={2023}, pages={109204},
 doi={10.1016/j.aim.2023.109204}
}

@misc{YuParameters,
  author        = {Mishra, Manish},
  title         = {Construction and exhaustion of enhanced supercuspidal
                   {Langlands} parameters from {Galois}-side {Yu} data},
  year          = {2026},
  eprint        = {2611.37564},
  archivePrefix = {arXiv},
  howpublished  = {\href{https://arxiv.org/abs/2611.37564}{arXiv:2611.37564}}
}

@misc{DF,
 author={Dat, Jean-Fran\c{c}ois and Fintzen, Jessica},
 title={Parametrization and reduction to depth zero of {$\overline{\mathbb Z}[1/p]$}-blocks of tame {$p$}-adic groups},
 year={2026}, eprint={2602.12333}, archivePrefix={arXiv},
 howpublished={\href{https://arxiv.org/abs/2602.12333}{arXiv:2602.12333}},
 note={Numbered references use the 105-page authors' version, dated 24 March 2026, available at \url{https://webusers.imj-prg.fr/~jean-francois.dat/recherche/publis/tamegroups.pdf}}
}

@misc{KalethaCovers,
 author={Kaletha, Tasho},
 title={On {$L$}-embeddings and double covers of tori over local fields},
 year={2021}, eprint={1907.05173v2}, archivePrefix={arXiv},
 howpublished={\href{https://arxiv.org/abs/1907.05173v2}{arXiv:1907.05173v2}},
 note={Version of 10 February 2021; numbered references use this version}
}

@misc{AFMO-reduce,
 author={Adler, Jeffrey D. and Fintzen, Jessica and Mishra, Manish and Ohara, Kazuma},
 title={Reduction to depth zero for tame {$p$}-adic groups via {Hecke} algebra isomorphisms},
 year={2024}, eprint={2408.07805v1}, archivePrefix={arXiv},
 howpublished={\href{https://arxiv.org/abs/2408.07805v1}{arXiv:2408.07805v1}},
 note={Version of 14 August 2024; numbered references use this version}
}

@misc{AFMO-structure,
  author={Adler, Jeffrey D. and Fintzen, Jessica and Mishra, Manish and Ohara, Kazuma},
  title={Structure of {Hecke} algebras arising from types},
  year={2024}, howpublished={\href{https://arxiv.org/abs/2408.07801v1}{arXiv:2408.07801v1}}
}

@incollection{KY17,
  author={Kim, Ju-Lee and Yu, Jiu-Kang}, title={Construction of tame types},
  booktitle={Representation Theory, Number Theory, and Invariant Theory},
  publisher={Springer}, address={Cham}, year={2017}, pages={337--357}
}

@article{SolEndo,
  author={Solleveld, Maarten}, title={Endomorphism algebras and {Hecke} algebras for reductive $p$-adic groups},
  journal={Journal of Algebra}, volume={606}, year={2022}, pages={371--470},
  note={With the correction in \href{https://arxiv.org/abs/2005.07899v3}{arXiv:2005.07899v3}}
}

@article{Sol2025,
  author={Solleveld, Maarten}, title={Parameters of {Hecke} algebras for {Bernstein} components of $p$-adic groups},
  journal={Indagationes Mathematicae}, volume={36}, year={2025}, pages={124--170}
}

@incollection{AMS-GHA,
  author={Aubert, Anne-Marie and Moussaoui, Ahmed and Solleveld, Maarten},
  title={Graded {Hecke} algebras for disconnected reductive groups},
  booktitle={Geometric Aspects of the Trace Formula}, series={Simons Symposia},
  publisher={Springer}, year={2018}, pages={23--84},
  note={References include the erratum in \href{https://arxiv.org/abs/1607.02713v2}{arXiv:1607.02713v2}}
}

@misc{Ohara25,
  author={Ohara, Kazuma}, title={On parameters of {Hecke} algebras for $p$-adic groups},
  year={2025}, howpublished={\href{https://arxiv.org/abs/2505.16040v3}{arXiv:2505.16040v3}}
}

@article{Sol-red,
  author={Solleveld, Maarten}, title={On the classification of irreducible representations of affine {Hecke} algebras with unequal parameters},
  journal={Representation Theory}, volume={16}, year={2012}, pages={1--87},
  note={Numbering follows \href{https://arxiv.org/abs/1008.0177v3}{arXiv:1008.0177v3}}
}

@article{Lus-unip,
  author={Lusztig, George}, title={Classification of unipotent representations of simple $p$-adic groups, {II}},
  journal={Representation Theory}, volume={6}, year={2002}, pages={243--289},
  note={\href{https://arxiv.org/abs/math/0111248}{arXiv:math/0111248}}
}

@article{Sol-unip,
  author={Solleveld, Maarten}, title={A local {Langlands} correspondence for unipotent representations},
  journal={American Journal of Mathematics}, volume={145}, number={3}, year={2023}, pages={673--719},
  note={Numbering follows \href{https://arxiv.org/abs/1806.11357v2}{arXiv:1806.11357v2}}
}

@article{Sol-survey,
  author={Solleveld, Maarten}, title={Affine {Hecke} algebras and their representations},
  journal={Indagationes Mathematicae}, volume={32}, number={5}, year={2021}, pages={1005--1082},
  note={Numbering follows \href{https://arxiv.org/abs/2009.03007v2}{arXiv:2009.03007v2}}
}

@unpublished{JDII,
  author = {Mishra, Manish},
  title = {Enhanced {Jordan} decomposition for {Bruhat--Tits} quotients and
           affine {Hecke} algebras at depth zero},
  note = {Manuscript, 8 November 2026},
  year = {2026}
}

@article{AdlerMishraRegular,
  author={Adler, Jeffrey D. and Mishra, Manish},
  title={Regular {Bernstein} blocks},
  journal={Journal f{\"u}r die reine und angewandte Mathematik},
  volume={775}, year={2021}, pages={71--86}
}

@article{AubertXuBlocks,
  author={Aubert, Anne-Marie and Xu, Yujie},
  title={Hecke algebras for {$p$}-adic reductive groups and local {Langlands} correspondences for {Bernstein} blocks},
  journal={Advances in Mathematics}, volume={436}, year={2024}, pages={109384},
  doi={10.1016/j.aim.2023.109384}
}

@article{KalethaRegular,
  author={Kaletha, Tasho},
  title={Regular supercuspidal representations},
  journal={Journal of the American Mathematical Society}, volume={32},
  number={4}, year={2019}, pages={1071--1170}, doi={10.1090/jams/925}
}

@article{KalethaPackets,
  author={Kaletha, Tasho}, title={Supercuspidal {$L$}-packets},
  journal={Annals of Mathematics}, year={2026},
  note={To appear. Numbering follows
    \href{https://arxiv.org/abs/1912.03274v2}{arXiv:1912.03274v2}}
}

@article{SXDepthZero,
  author={Solleveld, Maarten and Xu, Yujie},
  title={Hecke algebras and local {Langlands} correspondence for non-singular depth-zero representations},
  journal={Forum of Mathematics, Sigma}, volume={14}, year={2026},
  pages={e51}, doi={10.1017/fms.2026.10198}
}
\end{document}